\documentclass[11pt]{amsart}
\usepackage{amsmath, amssymb, amsthm, amsfonts, mathrsfs}
\usepackage{enumitem}
\setbox0=\hbox{$+$}
\newdimen\plusheight
\plusheight=\ht0
\def\+{\;\lower\plusheight\hbox{$+$}\;}

\setbox0=\hbox{$-$}
\newdimen\minusheight

\minusheight=\ht0
\def\-{\;\lower\minusheight\hbox{$-$}\;}

\setbox0=\hbox{$\cdots$}
\newdimen\cdotsheight
\cdotsheight=\plusheight
\def\cds{\lower\cdotsheight\hbox{$\cdots$}}

\theoremstyle{definition}

\usepackage{graphicx, epsfig}
\theoremstyle{definition}

\renewcommand{\parallel}{\mathrel{/\mkern-5mu/}}
\makeatletter

\newcommand{\newparallel}{\mathrel{\mathpalette\new@parallel\relax}}
\newcommand{\new@parallel}[2]{%
  \begingroup
  \sbox\z@{$#1T$}
  \resizebox{!}{\ht\z@}{\raisebox{\depth}{$\m@th#1/\mkern-5mu/$}}%
  \endgroup}
\makeatletter
\newcommand{\notparallel}{%
  \mathrel{\mathpalette\not@parallel\relax}%
}
\newcommand{\not@parallel}[2]{%
  \ooalign{\reflectbox{$\m@th#1\smallsetminus$}\cr\hfil$\m@th#1\parallel$\cr}%
}
\makeatother
\renewcommand{\pmod}[1]{\,(\textup{mod}\,#1)}
\numberwithin{equation}{section}
\theoremstyle{plain}
\newtheorem{theorem}{Theorem}[section]
\newtheorem{corollary}[theorem]{Corollary}
\newtheorem{example}[theorem]{Example}
\newtheorem{convention}[theorem]{Convention}
\newtheorem{proposition}[theorem]{Proposition}
\newtheorem{definition}[theorem]{Definition}
\newtheorem{lemma}[theorem]{Lemma}
\newtheorem{claim}[theorem]{Claim}
\usepackage{titlesec}
\usepackage{geometry} 
\titleformat*{\section}{\LARGE\bfseries}
\makeatletter
\renewcommand\section{\@startsection{section}{1}{\z@}%
                                  {-3.5ex \@plus -1ex \@minus -.2ex}%
                                 {2.3ex \@plus.2ex}%
                                 {\normalfont\Large\bfseries}}
\renewcommand\subsection{\@startsection{subsection}{1}{\z@}%
                                  {-3.5ex \@plus -1ex \@minus -.2ex}%
                                 {2.3ex \@plus.2ex}%
                                 {\normalfont\large\bfseries}}
\renewcommand\subsubsection{\@startsection{subsubsection}{1}{\z@}%
                                  {-3.5ex \@plus -1ex \@minus -.2ex}%
                                 {2.3ex \@plus.2ex}%
                                 {\normalfont\normalsize\bfseries}}

\makeatother

\makeatletter
\renewcommand\contentsnamefont{\bfseries}
\def\@starttoc#1#2{\begingroup
\setTrue{#1}%
\par\removelastskip\vskip\z@skip
\@startsection{}\@M\z@{\linespacing\@plus\linespacing}%
{.5\linespacing}{
\contentsnamefont}{#2}%
\ifx\contentsname#2%
\else \addcontentsline{toc}{section}{#2}\fi
\makeatletter
\@input{\jobname.#1}%
\if@filesw
\@xp\newwrite\csname tf@#1\endcsname
\immediate\@xp\openout\csname tf@#1\endcsname \jobname.#1\relax
\fi
\global\@nobreakfalse \endgroup
\addvspace{32\p@\@plus14\p@}%
\let\tableofcontents\relax
}
\def\contentsname{Contents}
\def\l@section{\@tocline{2}{.5ex}{0mm}{5pc}{}}
\def\l@subsection{\@tocline{2}{0pt}{2em}{5pc}{}}
\def\l@subsubsection{\@tocline{2}{0pt}{3em}{5pc}{}}
\def\l@subsubsubsection{\@tocline{2}{0pt}{4em}{5pc}{}}
\usepackage{titlesec}
\usepackage{color}

\begin{document}
\title{Classification of Hyperbolic Dehn fillings II: Quadratic case}
\dedicatory{Dedicated to the memory of Walter Neumann}
\author{BoGwang Jeon}
\maketitle

\begin{abstract}
This paper is subsequent to \cite{CHDF}. In this paper, we extend the classification of hyperbolic Dehn fillings with sufficiently large coefficients by addressing the remaining case not covered in \cite{CHDF}. Specifically, by considering the case in which the two cusp shapes lie in the same quadratic field, we obtain the complete classification under a mild assumption satisfied by most manifolds.

The content of this paper is not limited to the classification of hyperbolic Dehn fillings. Along the way, we also classify key types of automorphisms of the holonomy variety of a two-cusped hyperbolic $3$-manifold and uncover an intriguing hidden structure in the complex volume of certain manifolds. Concrete examples illustrating these phenomena were discovered by S. Oh, and we elaborate on them in this paper. 

All the results presented here appear to be effective. In the third paper of this series, joint with S. Oh, we will provide examples confirming the optimality of our results.
\end{abstract}

\tableofcontents

\section{Introduction}\label{intro}
\subsection{Main result I}
In \cite{CHDF}, the author proposed a way to classify Dehn fillings of a two-cusped hyperbolic $3$-manifold by a single invariant. One of the main results of \cite{CHDF} is stated as follows:

\begin{theorem}\cite{CHDF}\label{23070303}
Let $\mathcal{M}$ be a two-cusped hyperbolic $3$-manifold.\footnote{Here and throughout the paper, all hyperbolic $3$-manifolds are assumed to be orientable and of finite volume.} There exists a finite set $\mathcal{S}\big(\subset \mathrm{GL}_2(\mathbb{Q})\times \mathrm{GL}_2(\mathbb{Q})\big)$ such that, for any two Dehn fillings $\mathcal{M}_{p'_1/q'_1, p'_2/q'_2}$ and $\mathcal{M}_{p_1/q_1, p_2/q_2}$ of $\mathcal{M}$ satisfying
\begin{equation}\label{25101403}
\mathrm{pvol}_{\mathbb{C}}\;\mathcal{M}_{p'_1/q'_1, p'_2/q'_2}=\mathrm{pvol}_{\mathbb{C}}\;\mathcal{M}_{p_1/q_1, p_2/q_2}, 
\end{equation}
$(p'_1/q'_1, p'_2/q'_2)$ is equal to either 
\begin{equation*}
(\sigma(p_1/q_1), \rho(p_2/q_2))\quad \text{or}\quad (\sigma(p_2/q_2), \rho(p_1/q_1))
\end{equation*} 
for some $(\sigma, \rho)\in \mathcal{S}$. Further, if $|p'_i|+|q'_i|$ and $|p_i|+|q_i|$ ($i=1,2$) are sufficiently large, then 
\begin{equation*}
\mathrm{vol}_{\mathbb{C}}\;\mathcal{M}_{p'_1/q'_1, p'_2/q'_2}=\mathrm{vol}_{\mathbb{C}}\;\mathcal{M}_{p_1/q_1, p_2/q_2}.
\end{equation*}
\end{theorem}

Here, $\mathrm{vol}_{\mathbb{C}}$ (resp. $\mathrm{pvol}_{\mathbb{C}}$) denotes the \textit{complex volume} (resp. \textit{pseudo complex volume}) of a manifold. The latter was introduced in \cite{CHDF}, and is defined as
\begin{equation}\label{2510140111}
\text{pvol}_{\mathbb{C}}\;\mathcal{M}_{p_1/q_1, p_2/q_2}:=\text{vol}_{\mathbb{C}}\;\mathcal{M}-\dfrac{\pi}{2}\sum_{i=1}^2 \lambda^i_{p_1/q_1,p_2/q_2}
\pmod{\sqrt{-1}\pi^2\mathbb{Z}}
\end{equation} 
where $\lambda^i_{p_1/q_1,p_2/q_2}$ ($i=1,2$) are the complex lengths of the core geodesics of $\mathcal{M}_{p_1/q_1, p_2/q_2}$.\footnote{See Section \ref{25111201} for the precise definition of $\lambda^i_{p_1/q_1,p_2/q_2}$ $(i=1,2$).}

Roughly, the pseudo complex volume captures the dominant part of the complex volume, so provides a good approximation to it. The definition is also very natural in the sense of \cite{nz, yoshida}, and it is expected that the two invariants always exhibit the same behavior. That is, we postulate the following equivalence for any $(p'_1/q'_1, p'_2/q'_2)$ and $(p_1/q_1, p_2/q_2)$:
\begin{center}
$\mathrm{pvol}_{\mathbb{C}}\mathcal{M}_{p'_1/q'_1, p'_2/q'_2}=\mathrm{pvol}_{\mathbb{C}}\mathcal{M}_{p_1/q_1, p_2/q_2}\;\;\Longleftrightarrow\;\;\mathrm{vol}_{\mathbb{C}}\mathcal{M}_{p'_1/q'_1, p'_2/q'_2}=\mathrm{vol}_{\mathbb{C}}\mathcal{M}_{p_1/q_1, p_2/q_2}$. 
\end{center}
The forward direction has already been established in Theorem \ref{23070303} for sufficiently large coefficients, and the analogue of the equivalence in the one-cusped case (also for sufficiently large coefficients) was proved in \cite{JO}. Moreover, extensive experimental evidence from the SnapPy census supports this equivalence, and no counterexamples are yet known.

From various perspectives, we are convinced that the pseudo complex volume provides an effective tool for analyzing the deeper structure of the complex volume. Since each $e^{\lambda^i_{p_1/q_1,p_2/q_2}}$ is an algebraic number, it can be readily handled via basic operations and, in particular in our case, used to study their multiplicative relations.\footnote{For instance, see the definition preceding Theorem \ref{22050203}, which gives rise to the dichotomy stated in the theorem.} Further, and perhaps more importantly, this framework allows us to employ a wide range of concepts and techniques from number theory and algebraic geometry, as initiated in \cite{CHDF}. Because the behavior and structure of the complex volume generally parallel those of the pseudo complex volume, the approach uncovers many finer, previously unknown properties of the former as well.\footnote{Setting aside our main results, the analysis in Example \ref{25110702}, for example, is a typical illustration that shows the nobility and power of our approach. In that example, the complex volumes of certain manifolds are shown to decompose in several distinct and natural ways, using the notion of pseudo complex volume and techniques adapted to that framework. We believe that the result in this example could hardly have been attainable otherwise.} In this paper, we continue to explore along these lines and present several interesting - and often somewhat surprising - examples and phenomena.

Returning to Theorem \ref{23070303}, when the two cusp shapes of $\mathcal{M}$ do not lie in the same quadratic field, one can obtain a more effective statement, provided the filling coefficients are sufficiently large. By abuse of notation, we still denote by $\mathcal{S}$ the subset from Theorem \ref{23070303}, under this assumption. Then a partial quantitative form of Theorem \ref{23070303} is given as follows:  

\begin{theorem}\cite{CHDF}\label{20103003}
Adopting the same notation and assumptions as in Theorem \ref{23070303}, suppose two cusp shapes $\tau_1, \tau_2$ of $\mathcal{M}$ do not belong to the same quadratic field, and $|p'_i|+|q'_i|,|p_i|+|q_i|$ are sufficiently large ($i=1,2$). 
\begin{enumerate} 
\item Assume that neither $\tau_1$ nor $\tau_2$ is contained in $\mathbb{Q}(\sqrt{-1})$ or $\mathbb{Q}(\sqrt{-3})$. If $\tau_1$ and $\tau_2$ are relatively independent (resp. dependent), then $|\mathcal{S}|\leq 1$ (resp. $2$). 

\item If $\tau_1\in \mathbb{Q}(\sqrt{-1})$ and $\tau_2\notin\mathbb{Q}(\sqrt{-3})$ (resp. $\tau_1\in \mathbb{Q}(\sqrt{-3})$ and $\tau_2\notin \mathbb{Q}(\sqrt{-1})$), then $|\mathcal{S}|\leq 2$ (resp. $3$). 

\item If $\tau_1\in \mathbb{Q}(\sqrt{-1})$ and $\tau_2\in \mathbb{Q}(\sqrt{-3})$, then $|\mathcal{S}|\leq 6$. 
\end{enumerate}
\end{theorem}

Here, the cusp shapes $\tau_1$ and $\tau_2$ are said to be {\it relatively dependent} if there exists $a,b,c,d\in \mathbb{Q}$ such that $\tau_1=\frac{a+b\tau_2}{c+d\tau_2}$ (otherwise, they are {\it relatively independent}). 

This theorem naturally generalizes its one-cusped analogue (see Corollary \ref{23071103}) and it appears to be effective. In \cite{CHDFIII}, the author, in a joint work with S. Oh, provides examples for each case of the theorem to confirm this fact.

In this paper, we focus on the remaining case, namely, when $\tau_1$ and $\tau_2$ are both elements of the same quadratic field. The following presents one of the main theorems of the paper. 

\begin{theorem}\label{23070505}
Adopting the same notation and assumptions as in Theorem \ref{23070303}, suppose two cusp shapes of $\mathcal{M}$ belong to the same quadratic field, and $|p'_i|+|q'_i|,|p_i|+|q_i|$ are sufficiently large ($i=1,2$). We further suppose one of the following is true: 
\begin{itemize}
\item the potential function of $\mathcal{M}$ has a non-trivial term of degree $4$; or
\item two cusps of $\mathcal{M}$ are SGI to each other.
\end{itemize}
If $\tau_1, \tau_2\in\mathbb{Q}(\sqrt{-3})$ (resp. $\mathbb{Q}(\sqrt{-2})$ or $\mathbb{Q}(\sqrt{-1})$), then $|\mathcal{S}|\leq 18$ (resp. $3$ or $8$). Otherwise, $|\mathcal{S}|\leq 2$. 
\end{theorem}

The conditions highlighted with dots in the above theorem, we believe, cover most two-cusped hyperbolic $3$-manifolds. In fact, among the many SnapPy census manifolds that we have tested so far, there is only one example, $v2788$ (in the SnapPy notation), which does not fall under the conditions of Theorem \ref{23070505}. 

The theorem also proves to be effective possibly except for the case where $\tau_1, \tau_2 \in \mathbb{Q}(\sqrt{-2})$. For instance, if $\mathcal{M}$ is the so-called Napoleon's $3$-manifold described in \cite{calegari}, and four of its cusps are left unfilled, then the remaining two cusps admit infinitely many families of $18$ distinct Dehn fillings, each with the same pseudo complex volume (and hence the same complex volume). Thus Theorem \ref{23070505} may be regarded as optimal in a broad sense when $\tau_1, \tau_2 \in \mathbb{Q}(\sqrt{-3})$.

If $\tau_1, \tau_2 \in \mathbb{Q}(\sqrt{-2})$, although we do not know of any example satisfying the required condition with $|\mathcal{S}| = 3$, we have found an example satisfying neither of the conditions but still with $|\mathcal{S}| = 3$.\footnote{Unlike the cases $\tau_1, \tau_2\in\mathbb{Q}(\sqrt{-1})$ or $\mathbb{Q}(\sqrt{-3})$, the situation where $\tau_1, \tau_2\in \mathbb{Q}(\sqrt{-2})$ is regarded as particularly exceptional and intriguing. First, this case does not occur in the one-cusped setting (i.e., it does not admit non-trivial Dehn filling pairs with the same pseudo complex volume). Moreover, while $\mathcal{S}$ in the other cases admits a group structure, it does not seem to possess one when $\tau_1, \tau_2 \in \mathbb{Q}(\sqrt{-2})$. As we illustrate throughout the paper, many exceptional phenomena are realized via this case. See, for instance, Examples \ref{25101509} and \ref{25110702}. } This again is achieved by $v2788$. See Example \ref{25110702}. (The manifold $v2788$ is particularly interesting and important, and we will return to it at several points later.) 

Further examples illustrating the effectiveness of each case in Theorem \ref{23070505} will also be discussed in \cite{CHDFIII}.

For manifolds whose potential functions have a vanishing term of degree $4$ including $v2788$, we can still get a bound for $|\mathcal{S}|$ but it seems to be slightly off the optimal value and is therefore excluded at this stage.

\subsection{Refinement of Theorem \ref{23070303}}

To prove Theorem \ref{23070303}, the author resolved a weak version of the so called Zilber-Pink conjecture for $\mathcal{X}\times \mathcal{X}$ where $\mathcal{X}$ is the holonomy variety of a two-cusped hyperbolic $3$-manifold.

Before stating an updated version of Theorem \ref{23070303}, let us introduce the following convention first. Note that, for each $\lambda^i_{p_1/q_1,p_2/q_2}$ ($i=1,2$) given as in \eqref{2510140111}, it is well defined modulo $2\pi \sqrt{-1}$. We say that $\lambda^1_{p_1/q_1,p_2/q_2}$ and $\lambda^2_{p_1/q_1,p_2/q_2}$ are {\it linearly independent} over $\mathbb{Q}$ if there exists $a_1, a_2\in \mathbb{Q}\backslash\{0\}$ such that 
\begin{equation}\label{25110607}
a_1\lambda^1_{p_1/q_1,p_2/q_2}+a_2\lambda^2_{p_1/q_1,p_2/q_2}\equiv 0 \pmod{\pi \sqrt{-1}\mathbb{Q}}.
\end{equation}
Otherwise, they are {\it linearly dependent} over $\mathbb{Q}$. 

Now Theorem \ref{23070303} is further refined as follows:

\begin{theorem}\label{22050203}(Theorem 6.17 and Corollary 6.19 in \cite{CHDF})
Let $\mathcal{M}$ be a two-cusped hyperbolic $3$-manifold and $\mathcal{X}$ be its holonomy variety. Let $\mathcal{M}_{p'_{1}/q'_{1},p'_{2}/q'_{2}}$ and $\mathcal{M}_{p_{1}/q_{1},p_{2}/q_{2}}$ be two Dehn fillings of $\mathcal{M}$ with $|p'_i|+|q'_i|$ and $|p_i|+|q_i|$ sufficiently large ($i=1,2$) having the same pseudo complex volume. Then there exists a finite set $\mathcal{H}$, depending only on $\mathcal{M}$, such that a Dehn filling point associated to the pair of the two Dehn fillings on $\mathcal{X}\times \mathcal{X}\subset \mathbb{C}^8(:=(M_1, L_1, \dots, M'_2, L'_2))$ is contained in some element of $\mathcal{H}$. More precisely, we have the following dichotomy. 
\begin{enumerate}
\item If $\lambda^1_{p_1/q_1, p_2/q_2}$ and $\lambda^2_{p_1/q_1, p_2/q_2}$ are linearly dependent over $\mathbb{Q}$, then so are $\lambda^1_{p'_1/q'_1, p'_2/q'_2}$ and $\lambda^2_{p'_1/q'_1, p'_2/q'_2}$. In particular, $H(\in \mathcal{H})$ containing the Dehn filling point in this case is defined by equations of the form   
\begin{equation}\label{220708011}
\begin{gathered}
M_1^{a_1}L_1^{b_1}M_2^{a_2}L_2^{b_2}=(M'_1)^{a'_1}(L'_1)^{b'_1}(M'_2)^{a'_2}(L'_2)^{b'_2}=M_1^{a_3}L_1^{b_3}(M'_1)^{a'_3}(L'_1)^{b'_3}=1, \\ 
M_1^{c_1}L_1^{d_1}M_2^{c_2}L_2^{d_2}=(M'_1)^{c'_1}(L'_1)^{d'_1}(M'_2)^{c'_2}(L'_2)^{d'_2}=M_1^{c_3}L_1^{d_3}(M'_1)^{c'_3}(L'_1)^{d'_3}=1,  
\end{gathered}
\end{equation} 
and $(\mathcal{X}\times \mathcal{X})\cap H$ is a $1$-dimensional anomalous subvariety of $\mathcal{X}\times \mathcal{X}$. Moreover, let
\begin{equation}\label{25101605}
A_j=\begin{pmatrix} a_j & b_j \\ c_j & d_j \end{pmatrix},
\qquad
A'_j=\begin{pmatrix} a'_j & b'_j \\ c'_j & d'_j \end{pmatrix}
\quad (1\le j\le 3).
\end{equation}
Then all $A_j$ and $A'_j$ are invertible, and 
\begin{equation}\label{25101501}
A_1(1/\tau_1)=A_2(1/\tau_2),\quad
A'_1(1/\tau_1)=A'_2(1/\tau_2),\quad
A_3(1/\tau_1)=A'_3(1/\tau_2),
\end{equation}
together with
\begin{equation}\label{25101504}
(A_1^T)^{-1}(p_1/q_1)=(A_2^T)^{-1}(p_2/q_2),\;
(A_1^{\prime\,T})^{-1}(p'_1/q'_1)=(A_2^{\prime\,T})^{-1}(p'_2/q'_2),\;
(A_3^T)^{-1}(p_1/q_1)=(A_3^{\prime\,T})^{-1}(p'_1/q'_1).
\end{equation}

\item If $\lambda^1_{p_1/q_1, p_2/q_2}$ and $\lambda^2_{p_1/q_1, p_2/q_2}$ are linearly independent over $\mathbb{Q}$, then so are $\lambda^1_{p'_1/q'_1, p'_2/q'_2}$ and $\lambda^2_{p'_1/q'_1, p'_2/q'_2}$. In particular, $H(\in \mathcal{H})$ containing the Dehn filling point is defined by equations of the following form 
\begin{equation}\label{20071501}
\begin{gathered}
(M'_1)^{a'_1}(L'_1)^{b'_1}=M_1^{a_1}L_1^{b_1}M_2^{a_2}L_2^{b_2},\;\;(M'_2)^{a'_2}(L'_2)^{b'_2}=M_1^{a_3}L_1^{b_3}M_2^{a_4}L_2^{b_4},\\
(M'_1)^{c'_1}(L'_1)^{d'_1}=M_1^{c_1}L_1^{d_1}M_2^{c_2}L_2^{d_2},\;\;(M'_2)^{c'_2}(L'_2)^{d'_2}=M_1^{c_3}L_1^{d_3}M_2^{c_4}L_2^{d_4}, 
\end{gathered}
\end{equation}
and $(\mathcal{X}\times \mathcal{X})\cap H$ is a $2$-dimensional anomalous subvariety of $\mathcal{X}\times \mathcal{X}$. Moreover, let $$A'_j:=\left(\begin{array}{cc}
a'_{j} & b'_{j}\\
c'_{j} & d'_{j}
\end{array}\right)\; (1\leq j\leq 2) \text{ and }  
A_j:=\left(\begin{array}{cc}
a_{j} & b_{j}\\
c_{j} & d_{j}
\end{array}\right)\; (1\leq j\leq 4)$$
and assume all $A_j$ and $A'_j$ are non-trivial (i.e., $A_j$ and $A_j$ are not the zero matrix). Then they are invertible and 
\begin{equation}\label{24080403}
\begin{gathered}
A'_1(1/\tau_1)=
A_1(1/\tau_1)=
A_2(1/\tau_2),\; 
A'_2(1/\tau_2)=
A_3(1/\tau_1)=
A_4(1/\tau_2)
\end{gathered}
\end{equation}
together with 
\begin{equation}\label{25101505}
(A_1^{\prime\,T})^{-1}(p'_1/q'_1)=(A_1^T)^{-1}(p_1/q_1)
=(A_2^T)^{-1}(p_2/q_2),\;
(A_2^{\prime\,T})^{-1}(p'_2/q'_2)=(A_3^T)^{-1}(p_1/q_1)
=(A_4^T)^{-1}(p_2/q_2).
\end{equation}
\end{enumerate}
\end{theorem}

In the statement of the theorem, $(\mathcal{X} \times \mathcal{X}) \cap H$ is termed \textit{anomalous}, meaning that $\mathcal{X} \times \mathcal{X}$ and $H$ do not intersect normally. Generically, since $\dim (\mathcal{X} \times \mathcal{X}) = \dim H = 4$ in the ambient space of dimension $8$, the expected dimension of $(\mathcal{X} \times \mathcal{X}) \cap H$ is $0$. However, as its dimension is $2$ in this case, which is rather unusual, it is called an {\it anomalous subvariety} of $\mathcal{X} \times \mathcal{X}$. The term was coined by E. Bombieri, D. Masser and U. Zannier in \cite{za}. For its precise definition and related theorems, we refer the reader to that paper (although they will not be needed further in this paper).

For $A:=\left(\begin{array}{cc}
a & b\\
c & d
\end{array}\right)\in\mathrm{GL}_2(\mathbb{Q})$, the action of $A(x/y)$ is defined as $\frac{a+bx}{c+dy}$ in \eqref{25101501}-\eqref{25101504} and \eqref{24080403}-\eqref{25101505}. Since $(p_i, q_i)$ and $(p'_i, q'_i)$ ($i=1,2$) are coprime pairs, once $A'_j, A_j$ and $(p_i,q_i)$ are given, the pairs $(p'_i, q'_i)$ are thus uniquely determined by either \eqref{25101504} or \eqref{25101505}. In other words, $\mathcal{S}$ is induced from $\mathcal{H}$, and Theorem \ref{23070303} from Theorem \ref{22050203}. 

As for the first case of Theorem \ref{22050203}, please also refer to Theorem \ref{22050501}, which is a preparatory result for it and can be seen as a special case of the Zilber-Pink conjecture over $\mathcal{X}$. Indeed, Theorem \ref{22050203} is obtained by repeatedly applying Theorem \ref{22050501}. 

Among the two cases in Theorem \ref{22050203}, the second one is more generic. That is, $\lambda^1_{p_1/q_1, p_2/q_2}$ and $\lambda^2_{p_1/q_1, p_2/q_2}$ are typically linearly independent over $\mathbb{Q}$.\footnote{This is also indicated by the statement of Theorem \ref{22050501}. See the remark after Example \ref{25102101}.} In fact, Theorem \ref{20103003} can likewise be derived from Theorem \ref{22050203}(2). To be precise, if $\tau_1$ and $\tau_2$ do not belong to the same quadratic field, then clearly there are no non-trivial $A'_j$ ($1\leq j\leq 2$) and $A_j$ ($1\leq j\leq 4$) satisfying \eqref{24080403}. That is, either $A'_1=A_1=0$ or $A_2=0$ (resp. $A'_2=A_4=0$ or $A_3=0$). However, if $A'_1=A_1=0$ (resp. $A'_2=A_4=0$), then the first (resp. last) two equations in \eqref{20071501} are 
\begin{equation*}
M_2^{a_2}L_2^{b_2}=M_2^{c_2}L_2^{d_2}=1\quad (\textnormal{resp. }M_1^{a_3}L_1^{b_3}=M_1^{c_3}L_1^{d_3}=1), 
\end{equation*}
contradicting the fact that each coordinate of a Dehn filling point is not a root of unity. Consequently, $A_2=A_3=0$ and \eqref{20071501} is simplified to 
\begin{equation*}
\begin{gathered}
(M'_1)^{a'_1}(L'_1)^{b'_1}=M_1^{a_1}L_1^{b_1},\;(M'_2)^{a'_4}(L'_2)^{b'_4}=M_2^{a_4}L_2^{b_4},\;(M'_1)^{c'_1}(L'_1)^{d'_1}=M_1^{c_1}L_1^{d_1},\;(M'_2)^{c'_4}(L'_2)^{d'_4}=M_2^{c_4}L_2^{d_4}.
\end{gathered}
\end{equation*}
In this case, the problem is treated analogously to the one-cusped case, which was covered in \cite{CHDF, JO} and will be reviewed in Section \ref{23070513}.

On the other hand, if $\tau_1$ and $\tau_2$ belong to the same quadratic field and, moreover, every $A'_j$ and $A_j$ is nonzero, the situation becomes quite different. There is no obvious way to apply the previous technique, and it turns out that a deeper analysis of the underlying structure of $\mathcal{X}$ is required. We return to this point after stating our second main result below.

\subsection{Main result II}\label{25101710}
\subsubsection{Refinement of Theorem \ref{22050203}}
Our second main result is a quantitative version of Theorem \ref{22050203} for the case where $\tau_1$ and $\tau_2$ belong to the same quadratic field, which, as remarked above, is the hardest case to deal with. Since the two cases in Theorem \ref{22050203} are treated separately, we present its effective version as two theorems below.

\begin{theorem}\label{250825011}
Let $\mathcal{M}$ be a two-cusped hyperbolic $3$-manifold whose two cusp shapes belong to the same quadratic field, and assume that either 
\begin{itemize}
\item the potential function of $\mathcal{M}$ has a non-trivial term of degree $4$; or
\item two cusps of $\mathcal{M}$ are SGI to each other. 
\end{itemize}
Let $\mathcal{M}_{p_1/q_1, p_2/q_2}$ be a Dehn filling of $\mathcal{M}$ with sufficiently large $|p_i|+|q_i|$ ($i=1,2$). If the complex lengths of its core geodesics of $\mathcal{M}_{p_1/q_1, p_2/q_2}$ are linearly dependent over $\mathbb{Q}$, then the number of Dehn fillings of $\mathcal{M}$ whose pseudo complex volume is equal to that of $\mathcal{M}_{p_1/q_1, p_2/q_2}$ is at most $9$ or $4$, depending on whether $\tau_1, \tau_2 \in \mathbb{Q}(\sqrt{-3})$ or $\mathbb{Q}(\sqrt{-1})$, respectively. Otherwise, it is $1$.
\end{theorem}

\begin{theorem}\label{25101401}
Let $\mathcal{M}$ be a two-cusped hyperbolic $3$-manifold whose two cusp shapes belong to the same quadratic field, and $\mathcal{M}_{p_1/q_1, p_2/q_2}$ be a Dehn filling of $\mathcal{M}$ with sufficiently large $|p_i|+|q_i|$ ($i=1,2$). If the complex lengths of its core geodesics of $\mathcal{M}_{p_1/q_1, p_2/q_2}$ are linearly independent over $\mathbb{Q}$, then the number of Dehn fillings of $\mathcal{M}$ whose pseudo complex volume is equal to that of $\mathcal{M}_{p_1/q_1, p_2/q_2}$ is at most $18, 3$ or $8$, depending on whether $\tau_1, \tau_2\in \mathbb{Q}(\sqrt{-3}), \mathbb{Q}(\sqrt{-2})$ or $ \mathbb{Q}(\sqrt{-1})$. Otherwise, it is $2$. 
\end{theorem}

More detailed versions of Theorems \ref{250825011} and \ref{25101401} can be found in Theorems \ref{25082501} and \ref{25101601}-\ref{25101610}, respectively. 

Note that the two conditions emphasized with dots in Theorem \ref{250825011} exactly coincide with those given in Theorem \ref{23070505}. Once the above two theorems are established, Theorem \ref{25101401} follows as simply as a corollary of them.
 
The proof of Theorem \ref{250825011} is relatively straightforward compared to that of Theorem \ref{25101401}. As noted earlier, the proof of Theorem \ref{22050203}(1) involves multiple applications of Theorem \ref{22050501}. Hence, by first making Theorem \ref{22050501} effective and then applying it repeatedly, Theorem \ref{250825011} follows naturally. We review the proof of Theorem \ref{22050203}(1) and outline the strategy for its quantification in detail in Section \ref{25102903}, preceding the proof of Theorem \ref{22050501} in Section \ref{25111402}. Conceptually and technically, however, this case does not pose any serious difficulty, nor does it differ substantially from the arguments in \cite{CHDF}. Therefore, we focus primarily on Theorem \ref{25101401} in the remainder of this section.

The proof of Theorem \ref{25101401} essentially relies on the classification of every $H$ defined in \eqref{20071501} with all $A_j$ nontrivial. Initially, as discussed in an earlier version of \cite{CHDF}, we believed that either $A_2 = A_3 = 0$ or $A_1 = A_4 = 0$ would always hold, and no $H$ with all $A_j\neq 0$ could arise. To simplify notation, let us, for the moment, call $H$ of \textit{Type I} (resp. \textit{Type II}) if $A_2 = A_3 = 0$ (resp. $A_1 = A_4 = 0$), and of \textit{Type III} otherwise. We will provide the precise definition in Section \ref{23070513}.

As will be laid out throughout the paper, particularly in Sections \ref{25012001}-\ref{25011802}, the existence of $H$ of Type III seems very unlikely, since the defining criteria are rather restrictive.\footnote{As explained earlier, the fact that $(\mathcal{X}\times \mathcal{X})\cap H$ is an \emph{anomalous} subvariety of $\mathcal{X}\times\mathcal{X}$ already suggests that this is not a typical phenomenon.} In the sections, it will be shown that these constraints translate into several equations in fewer variables, whose solutions generally do not exist. Nevertheless, thanks largely to the favorable properties of hyperbolic $3$-manifolds, we are able to identify the precise conditions under which these equations yield solutions, and so fully classify Type III elements. It turns out that Type III elements can exist only when the cusp shapes of $\mathcal{M}$ lie in $\mathbb{Q}(\sqrt{-1})$, $\mathbb{Q}(\sqrt{-2})$, or $\mathbb{Q}(\sqrt{-3})$, and even then only in highly rigid forms with at most finitely many possibilities.

Having this in mind, we still expected that the analysis, i.e., the classification of Type III elements, would remain purely theoretical rather than admitting any actual example realizing the phenomenon. However, interestingly and somewhat unexpectedly, a concrete example was found with the help of SnapPy by Sunul Oh. The manifold $v2788$, mentioned earlier, again illustrates the phenomenon. 
 
\begin{example}[S. Oh] \label{25101509}
{\normalfont Let $\mathcal{M}$ be $v2078$, $\mathcal{X}$ be its holonomy variety, and $H$ be defined by either 
\begin{equation}\label{25100902}
\begin{gathered}
(M'_1)^{2}=L_1L_2,\;L'_1=M_1^{-1}M_2^{-1},\;(M'_2)^{2}=L_1L_2^{-1},\;L'_2=M_1^{-1}M_2 
\end{gathered}
\end{equation}
or 
\begin{equation}\label{25100901}
(M'_1)^{2}=M_1L_1M_2^{-1},\;(L'_1)^{2}=M_1^{-2}L_1L_2^{-1},\;
(M'_2)^{2}=M_1M_2L_2^{-1},\;(L'_2)^{2}=L_1M_2^2L_2. 
\end{equation}
Then $(\mathcal{X}\times \mathcal{X})\cap H$ is a $2$-dimensional anomalous subvariety of $\mathcal{X}\times \mathcal{X}$.}
\end{example}

The description of $H$ given in either \eqref{25100902} or \eqref{25100901} corresponds precisely to that predicted theoretically in our analysis; see Theorem \ref{23070521} or \ref{24040801}, and Example \ref{25110401} or \ref{25110402} respectively.\footnote{Again, Example \ref{25101509} was not known when the first version of the paper was written. As stated there, we initially doubted the existence of such an example, and at least to us, the discovery of Example \ref{25101509} came as a big surprise.}

In truth, there are more $H$ fitting into the scenario of Example \ref{25101509}. The best way to handle them is to endow the set of all such $H$ with a group structure (which we shall do in Section \ref{25110409}), rather than listing all possible equation forms. See Example \ref{25111910} for a complete description of the group - and hence all the possible shapes of $H$ - for the above case.

Due to their length and the limited scope of the introduction, we do not explicitly state Theorems \ref{23070521} and \ref{24040801}, nor, more broadly, the classification of $H$ in Theorem \ref{22050203}(2) here. However this discovery is, in our view, among the most noteworthy results of the paper. They will be established in Sections \ref{24081804}-\ref{24081805}, and the final outcomes are all summarized in Theorems \ref{24090201}, \ref{25012901}, and \ref{24092702}. 

\subsubsection{How to classify $H$?}


To analyze the structure of $H$, we first associate it to a finite order element of $\mathrm{GL}_4(\mathbb{Q})$. In fact, if $A'_1$ and $A'_2$ in \eqref{20071501} are normalized to the identity, then the rest coefficients $A_j$ ($1\leq j\leq 4$) together determine a single element of $\mathrm{GL}_4(\mathbb{Q})$. Furthermore, as $(\mathcal{X}\times\mathcal{X})\cap H$ is a $2$-dimensional anomalous subvariety of $\mathcal{X}\times \mathcal{X}$, alternatively, $H$ is interpreted as an element of $\mathrm{Aut}\;\mathcal{X}$, the automorphism group of $\mathcal{X}$. 

At first glance, one might hope the possibility of appealing to a general theory applicable to our case. Indeed, according to Minkowski's theorem \cite{min, min2}, there exists a uniform bound depending only on $n$ such that the number of elements in every finite subgroup of $\mathrm{GL}_n(\mathbb{Q})$ less than this bound. For instance, when $n=2$ and $n=4$, the bounds are known to be $12$ and $1152$ respectively \cite{feit}; however, in light of Theorem \ref{20103003} and Corollary \ref{23071103}, these bounds are excessively large and far from optimal for our purposes. Moreover, classifying every finite subgroup of $\mathrm{GL}_4(\mathbb{Q})$ seems to be very difficult.

To address the problem in a more practical manner, we adopt a step-by-step approach based on the minimal polynomial of $H$. For any element of finite order in $\mathrm{GL}_4(\mathbb{Q})$, its minimal polynomial, having degree at most $4$, is one of those listed in \eqref{23062301}. For each of these in the list, we will examine and count all the possible forms of $H$ associated with it.\footnote{For example, when the minimal polynomial is $x^2+1$, there are $5$ distinct forms of $H$, as explicitly described in Theorem \ref{23070521}.} 
 
To elaborate further on this, as an initial step, we use the local analytic representation of $\mathcal{X}$ to deduce a necessary condition that $H$ must satisfy. This condition is formulated as a system of two homogeneous equations in two variables, involving a $2 \times 2$ matrix transformation. The resulting matrix, which we will term the \textit{primary matrix} associated with $H$,\footnote{The precise definition will be given in Definition \ref{25120501}.} is derived from the coefficient matrix of $H$ and encodes essential information about it. Utilizing the minimal polynomial constraint on $H$, as well as several key features of $\mathcal{X}$, we precisely determine all the possible entries of this $2 \times 2$ matrix, and, consequently, obtain all the possible concrete forms of $H$.

The formulation, analysis, and solution of this system of homogeneous equations will be the central job of the paper, and a substantial portion, particularly Sections \ref{25012001}-\ref{25011802}, is devoted to this task. Generically, equations of this type may have no solutions; however, due to the highly symmetric properties of $\mathcal{X}$, it is shown that only a few solutions exist in very rigid types.

\subsection{Remarks on the approach}

Although this paper builds upon the results of \cite{CHDF}, the proofs of many theorems have a very different flavor, with many technical details distinct from those in \cite{CHDF}. For example, in this paper, we do not use any heavy number-theoretic or algebraic-geometric machinery that was critically employed in \cite{CHDF}. Instead, the tools used here include more fundamental and structural results, such as the diagonalization theorem and the primary decomposition theorem in linear algebra.
 
From a computational standpoint, since we pursue an effective statement in every theorem, the reader may inevitably find the presentation of some parts of the paper is highly technical. We try to avoid repeating routine computations; however, whenever necessary, we include full details, even down to a single step. 

From a conceptual standpoint, we hope that the paper introduces a new perspective and reveals previously unknown structures in hyperbolic Dehn filling and the holonomy variety, which may ultimately help decipher the mystery of both the geometric and arithmetic invariants of hyperbolic 3-manifolds, along with their interrelation.

It would be preferable if the reader is already acquainted with the basic background, or at least with the statements of the main results, presented in \cite{CHDF}. Nevertheless, we will review the necessary material and highlight the key points in Section \ref{23070513}.

\subsection{Outline of the paper}

Section \ref{prem} is a preliminary overview. We summarize the background and several results from \cite{CHDF} that will be revisited in this paper, and then discuss concepts and theorems from linear algebra needed for the proofs presented here.

In Section \ref{25101507}, we prove Theorem \ref{250825011}. The proof of this theorem is quite rudimentary compared with that of Theorem \ref{25101401} and relies only on the material reviewed in the previous section.

Sections \ref{23071105}-\ref{24081805} are devoted to classifying Types I-III elements in $\mathrm{Aut}\;\mathcal{X}$. These sections constitute the main body of the paper and are preparatory steps for proving Theorem \ref{25101401}. 

In Section \ref{23071105}, as a warm-up, we remove the two polynomials, $x^4-x^3+x^2-x+1$ and $x^4+x^3+x^2+x+1$, from the list of minimal polynomials corresponding to Types I-III elements in $\mathrm{Aut}\;\mathcal{X}$. 

Sections \ref{25012001}-\ref{25011802} form the technical core and outline the key strategy for the classification. First, in Section \ref{25012001}, for each Type I, II or III element in $\mathrm{GL}_4(\mathbb{Q})$, we set up a system of homogeneous equations, whose solvability gives a necessary condition for the element to lie in $\mathrm{Aut}\,\mathcal{X}$. Then, in Section \ref{25011802}, we analyze these equations and determine when they admit solutions. 

In Sections \ref{24081804}-\ref{24081805}, for each polynomial of degree $\leq 4$ that is a product of cyclotomic polynomials, we explore all possible forms of Types I-III elements contained in both $\mathrm{GL}_4(\mathbb{Q})$ and $\mathrm{Aut}\,\mathcal{X}$ with that polynomial as their minimal polynomial.

In Section \ref{25012005}, building on the results of Sections \ref{24081804}-\ref{24081805}, we completely classify the finite subgroups of $\mathrm{Aut}\;\mathcal{X}$ consisting of elements of Types I-III. 

In Section \ref{25012003}, we verify Theorem \ref{25101401} by amalgamating the results of the preceding sections and, combining this with Theorem \ref{250825011}, finally establish Theorem \ref{23070505} in Section \ref{25102801}.

\subsection{Acknowledgement}
The first version of this paper dealt only with Theorem \ref{25101401}; the author mistakenly overlooked the first case of Theorem \ref{22050203}, and thus omitted Theorem \ref{250825011} in that version.

I would like to thank Sunul Oh for many helpful conversations and collaborations, and in particular for allowing me to include Example \ref{25101509} in this paper.

This work was supported by the National Research Foundation of Korea (NRF) grant funded by the Korea government (RS-2025-23323903).
\newpage
\section{Preliminaries}\label{prem}

\subsection{Essence from the precedent}\label{23070513}
 
In this section, we will briefly go over some essential elements from \cite{CHDF}. Rather than repeating all the technical details covered in \cite{CHDF}, we focus on the principal concepts and results that are frequently used and play crucial roles in this paper.

\subsubsection{Holonomy variety: algebraic \& analytic}\label{25111201}
We first review the notions of holonomy variety as well as Dehn filling point, which are key to understanding the statements of Theorems \ref{23070303} and \ref{22050203}. For a more detailed account, we refer to Sections 2.1-2.2 of \cite{CHDF}. 

For a given $n$-cusped hyperbolic $3$-manifold $\mathcal{M}$, let $\mathcal{T}$ be a geometric ideal triangulation and $G(\mathcal{T})$ be the \textit{gluing variety} induced from $\mathcal{T}$ in the sense of Thurston \cite{thu}. If $T_i$ denotes a torus cross-section of the $i^{\text{th}}$-cusp of $\mathcal{M}$ with the chosen meridian-longitude pair $m_i, l_i$ $(1\leq i\leq n)$, then each point of $G(\mathcal{T})$ gives rise to a Euclidean similarity structure on $T_i$, thus defining the so-called the \textit{holonomy} map 
\begin{equation*}
\pi_1(T_i)\longrightarrow \text{Aff}(\mathbb{C})=\{az+b\;:\;a\neq 0, b\in \mathbb{C}\}.
\end{equation*}
The derivatives of the holonomies of $l_i$ and $m_i$ then yield rational functions $M_i$ and $L_i$ respectively, on $G(\mathcal{T})$, and the \textit{holonomy variety} $\mathcal{X}$ of $\mathcal{M}$ is defined as the Zariski closure of the image of the map
\begin{equation*}
\xi\;:\;G(\mathcal{T})\longrightarrow (M_1, L_1, \dots, M_n, L_n).
\end{equation*}
It is known that the holonomy variety is independent of the triangulation $\mathcal{T}$ but depends only on the chosen meridian-longitude pair of $T_i$.

In general, the holonomy variety $\mathcal{X}$ may have several irreducible components, but we are only interested in so called the \textit{geometric component} of it, which is an $n$-dimensional variety in $\mathbb{C}^{2n}\big(:=(M_1, L_1, \dots, M_n, L_n)\big)$ containing $(1, \dots, 1)$ that corresponds the complete hyperbolic metric structure of $\mathcal{M}$. By abuse of notation, we still denote this component by $\mathcal{X}$ and call it the \textit{holonomy variety} of $\mathcal{M}$.   

The {\it analytic holonomy variety} $\mathcal{X}_{an}$ of $\mathcal{M}$ is obtained by taking the logarithm of each coordinate of $\mathcal{X}$, that is, 
$$u_i:=\log M_i\;\;\text{and}\;\; v_i:=\log L_i,\quad (1\leq i\leq n),$$
and is a complex manifold in $\mathbb{C}^{2n}\big(:=(u_1, v_1, \dots, u_n, v_n)\big)$, locally biholomorphic to a small neighborhood of a branch containing $(1, \dots, 1)$ in $\mathcal{X}$.\footnote{In \cite{CHDF}, $\mathcal{X}_{an}$ was denoted by $\log\mathcal{X}$.} In particular, it exhibits various symmetric properties, which are summarized below: 
\begin{theorem}\cite{jeon3, nz, PP} \label{potential}
\begin{enumerate}
\item Using $u_1, \dots, u_n$ as generators, each $v_i$ is represented of the form 	
\begin{equation*}
u_i\cdot\tau_i(u_1,\dots,u_n)
\end{equation*}
over $\mathcal{X}_{an}$ where $\tau _i(u_1,\dots,u_n)$ is a holomorphic function with $\tau_i(0,\dots,0)=\tau_i\in \mathbb{C}\setminus\mathbb{R}$ and $\mathrm{Im}\,\tau_i>0$ for $1\leq i\leq n$.

\item There is a holomorphic function $\Phi(u_1,\dots,u_n)$ such that 
\begin{equation*}
v_i=\frac{1}{2}\frac{\partial \Phi}{\partial u_i}\quad \text{and}\quad \Phi(0,\dots,0)=0, \quad (1\leq i\leq n). 
\end{equation*}
In particular, $\Phi(u_1,\dots,u_n)$ is even in each argument and so its Taylor expansion is of the following form:
\begin{equation*}
\Phi(u_1,\dots,u_n)=(\tau_1u_1^2+\cdots+\tau_n u_n^2)+(c_{4,\dots,0}u_1^4+\cdots+c_{0,\dots,4}u_n^4)+\text{(higher order)}.\\
\end{equation*}

\item Each $v_i$ is not linear. That is, 
\begin{equation*}
v_i(0, \dots, 0,  u_i, 0, \dots, 0)\neq \tau_i u_i, \quad (1\leq i\leq n).
\end{equation*}
\end{enumerate}
\end{theorem}
\begin{proof}
The first two are given in \cite{nz, PP}. For the last one, see Lemma 2.3 in \cite{jeon3}. 
\end{proof}

We call $\tau_i$ the \textit{cusp shape} of the $i^{\text{th}}$-cusp of $\mathcal{M}$ with respect to $m_i, l_i$, and $\Phi(u_1,\dots,u_n)$ as the \textit{Neumann-Zagier potential function}, or simply the \textit{potential function}, of $\mathcal{M}$ with respect to $m_i, l_i$ $(1\leq i\leq n)$.

\subsubsection{Dehn filling points}
For $\mathcal{M}$ and $\mathcal{X}$ as above, a \textit{Dehn filling point} corresponding to $\mathcal{M}_{p_1/q_1,\dots,p_n/q_n}$ is an intersection point between $\mathcal{X}$ and the variety defined by  
\begin{equation} \label{Dehn eq}
M_i^{p_i}L_i^{q_i}=1,\quad (1\leq i\leq n).  
\end{equation}
Obviously, the above equations are derived from the relations added to the fundamental group to $\mathcal{M}$, via the Seifert-van Kampen theorem, to obtain that of $\mathcal{M}_{p_1/q_1,...,p_n/q_n}$. Consequently, such a point gives rise to the discrete faithful representation 
\begin{equation*}
\phi\;:\;\pi_1(\mathcal{M}_{p_1/q_1, \dots, p_n/q_n})\longrightarrow \mathrm{PSL}_2\mathbb{C}
\end{equation*}
associated with $\mathcal{M}_{p_1/q_1,...,p_n/q_n}$. Following \cite{CHDF}, let 
\begin{equation}\label{25101701}
\mathcal{M}_{p_1/q_1,\dots, p_n/q_n}\quad \text{and}\quad  \mathcal{M}_{p'_1/q'_1, \dots, p'_n/q'_n}
\end{equation}
be two Dehn fillings of $\mathcal{M}$ with the same pseudo complex volume. By regarding each filling as a point in its copy of $\mathcal{X}$, the resulting point in $\mathcal{X} \times \mathcal{X}$ is called a \textit{Dehn filling point associated to the pair in \eqref{25101701}}.

In the context of $\mathcal{X}_{an}$, \eqref{Dehn eq} is specified as 
\begin{equation}\label{22063001}
p_iu_i+q_iv_i= \pm 2\pi \sqrt{-1},\quad  (1\leq i\leq n).  
\end{equation}
By abuse of notation, we still call an intersection point between \eqref{22063001} and $\mathcal{X}_{an}$ a \textit{Dehn filling point} associated to $\mathcal{M}_{p_1/q_1, \dots, p_n/q_n}$. If
\begin{equation*}
(u_1, v_1, \dots, u_n, v_n)=(\xi_{m_1}, \xi_{l_1}, \dots, \xi_{m_n}, \xi_{l_n})
\end{equation*}
is a Dehn filling point on $\mathcal{X}_{an}$ associated to $\mathcal{M}_{p_1/q_1, \dots, p_n/q_n}$, then the following quantity
$$\pm (r_i\xi_{m_i}+s_i\xi_{l_i}),$$ where $r_i, s_i\in \mathbb{Z}$ satisfy $p_is_i-q_ir_i=1$ and $|e^{r_i\xi_{m_i}+s_i\xi_{l_i}}|>1$, is the \textit{complex length} of the $i^{\text{th}}$-core geodesic of $\mathcal{M}_{p_1/q_1, \dots p_n/q_n}$. This appears as $\lambda^i_{p_1/q_1, \dots, p_n/q_n}$ in \eqref{2510140111} and the statements of Theorem \ref{22050203}, as well as in many others.

\subsubsection{Auxiliary theorems I}

In this sub-subsection and the next, we review several auxiliary theorems from \cite{CHDF}. In particular, those presented here are relevant to Theorem \ref{22050203}(1) and will be adapted to verify Theorem \ref{250825011}.

We also note that, from now until the end, all results will be presented in the setting of the analytic holonomy variety rather than its algebraic counterpart. This formulation provides a simpler and more efficient framework in which to state and verify them.

First, the following theorem combines Theorem 2.23 and Lemma 4.1 in \cite{CHDF}, and it can be regarded as a special case of the Zilber-Pink conjecture for $\mathcal{X}$ (or $\mathcal{X}_{an}$).
    
\begin{theorem}\label{22050501}\cite{CHDF}
Let $\mathcal{M}$ be a two-cusped hyperbolic $3$-manifold and $\mathcal{X}_{an}$ be its analytic holonomy set. Then there exists a finite set $\mathcal{A}(\subset \mathrm{GL}_2(\mathbb{Q}))$ of matrices such that, for any Dehn filling $\mathcal{M}_{p_1/q_1,p_2/q_2}$ of $\mathcal{M}$ with $|p_i|+|q_i|$ ($i=1,2$) sufficiently large, if the complex lengths $\lambda^i_{p_1/q_1,p_2/q_2}$ ($i=1,2$) of the core geodesics of $\mathcal{M}_{p_1/q_1, p_2/q_2}$ are linearly dependent over $\mathbb{Q}$, then a corresponding Dehn filling point associated with $\mathcal{M}_{p_1/q_1,p_2/q_2}$ is contained in a $1$-dimensional analytic subset of $\mathcal{X}_{an}$, defined by 
\begin{equation*}
\left(\begin{array}{c}
u_1\\
v_1
\end{array}\right)
=A\left(\begin{array}{c}
u_2\\
v_2
\end{array}\right)
\end{equation*}
for some $A\in \mathcal{A}$. Moreover, it further satisfies  
\begin{equation}\label{25110608}
\lambda^1_{p_1/q_1,p_2/q_2}\equiv(\det A)\lambda^2_{p_1/q_1,p_2/q_2} \pmod{\pi \sqrt{-1}\mathbb{Q}}
\end{equation}
and 
\begin{equation}\label{25120801}
\left(\begin{array}{cc}
p_{1} & q_{1}\\
\end{array}\right)=
\left(\begin{array}{cc}
p_{2} & q_{2}\\
\end{array}\right)A^{-1}.
\end{equation}
\end{theorem}

In short, if a Dehn filling of $\mathcal{M}$ satisfies the condition stated in the theorem, then a point associated with it is not randomly distributed over $\mathcal{X}_{an}$ but rather lies on one of finitely many analytic subsets of it. 

\begin{example}\label{25102101}
{\normalfont A typical example illustrating the above theorem is the Whitehead link complement $\mathcal{W}$. Since there exists an isometry exchanging the two cusps, the complex lengths of the core geodesics of $\mathcal{W}_{p/q,p/q}$ are equal for every co-prime pair $(p,q)$,\footnote{Here, we choose the canonical choice of meridian-longitude pairs for the two cusps} and a Dehn filling point associated to $\mathcal{W}_{p/q,p/q}$ lies in a variety defined by $$\left(\begin{array}{c}
u_1\\
v_1
\end{array}\right)
=\left(\begin{array}{c}
u_2\\
v_2
\end{array}\right).$$}
\end{example}

\noindent\textbf{Remark. }Another consequence of Theorem \ref{22050501} is that the two complex lengths of a Dehn filling are {\it generically} linearly independent. Indeed, if they were linearly dependent, then a Dehn filling point corresponding to that filling would lie over one of at most finitely many $1$-dimensional subvarieties of $\mathcal{X}$. Since the union of these $1$-dimensional subsets forms a proper Zariski-closed subset of $\mathcal{X}$, a randomly chosen Dehn filling point is not likely to lie over it. 

For instance, in the previous example, the complex lengths of the two core geodesics of every Dehn filling $\mathcal{W}_{p_1/q_1,p_2/q_2}$ are linearly independent over $\mathbb{Q}$ unless $p_1/q_1=p_2/q_2$.
\\

As a special case of Theorem \ref{22050501}, when $\mathcal{X}_{an}$ in the statement is given as the product of two analytic curves, the structure of $\mathcal{A}$ in Theorem \ref{22050501} is explicitly described as follows:
\begin{theorem}\label{25111301}
Let $\mathcal{M}$ be a two-cusped hyperbolic $3$-manifold and $\mathcal{X}_{an}$ be its analytic holonomy variety, given as the product of two analytic curves $\mathcal{C}_1\times\mathcal{C}_2$ in $\mathbb{C}^2\times \mathbb{C}^2(:=(u_1, v_1, u_2, v_2))$. Suppose $A_i$ ($i=1,2$) are elements of $\mathrm{GL}_2(\mathbb{Q})$ such that each 
\begin{equation}\label{25111303}
\left(\begin{array}{c}
u_1\\
v_1
\end{array}\right)
=A_i\left(\begin{array}{c}
u_2\\
v_2
\end{array}\right)
\end{equation}
is a $1$-dimensional analytic subset of $\mathcal{X}_{an}$ $(i=1,2)$. If the cusp shapes of $\mathcal{M}$ are contained in $\mathbb{Q}(\sqrt{-3})$ (resp. $\mathbb{Q}(\sqrt{-1})$, then $(A_2^{-1}A_1)^6=I$ (resp. $(A_2^{-1}A_1)^4=I$). Otherwise, $A_1=\pm A_2$. In particular, if $A_2^{-1}A_1\neq  \pm I$, then it is uniquely determined up to inversion and multiplication of $\pm 1$. 
\end{theorem}

The proof of the theorem directly follows from the next lemma, established in \cite{CHDF}. 
\begin{lemma}\label{20102501}(Lemma 4.2 in \cite{CHDF})
For $\mathcal{M}$ and $\mathcal{X}_{an}$ as in Theorem \ref{25111301}, we further assume that $\mathcal{C}_1$ and $\mathcal{C}_2$ coincide and let $\mathcal{C}:=\mathcal{C}_1(=\mathcal{C}_2)$. Suppose 
\begin{equation*}
\left(\begin{array}{c}
u_1\\
v_1
\end{array}\right)
=A\left(\begin{array}{c}
u_2\\
v_2
\end{array}\right)
\end{equation*}
is a $1$-dimensional subset of $\mathcal{C}\times \mathcal{C}$ for some $A\in\mathrm{GL}_2(\mathbb{Q})$. If the cusp shapes of $\mathcal{M}$ are contained in $\mathbb{Q}(\sqrt{-3})$ (resp. $\mathbb{Q}(\sqrt{-1})$, then $A^6=I$ (resp. $A^4=I$). Otherwise, $A=\pm I$. In particular, if $A \pm I$, then it is uniquely determined up to inversion and multiplication of $\pm 1$. 
\end{lemma}

\begin{proof}[Proof of Theorem \ref{25111301}]
Under the assumptions on the structure of $\mathcal{X}_{an}$ and the equation in \eqref{25111303}, each $A_i$ can be alternatively regarded as an isomorphism from $\mathcal{C}_1$ to $\mathcal{C}_2$. It follows that $A_2^{-1}A_1$ can be viewed as an automorphism of $\mathcal{C}_1$, and the conclusion then follows from Lemma \ref{20102501}.
\end{proof}

Note that the condition on $\mathcal{X}_{an}$ in Theorem \ref{25111301} corresponds precisely to the definition of the notion of `SGI' mentioned in the first main result, Theorem \ref{23070505}. Now we formally define it as follows:
\begin{definition}
\normalfont{Let $\mathcal{M}$ be a two-cusped hyperbolic $3$-manifold, and let $\mathcal{X}_{an}$ be its analytic holonomy variety. We say that the two cusps of $\mathcal{M}$ are \textit{strongly geometrically isolated (SGI)} if $\mathcal{X}_{an}$ is the product of two analytic curves $\mathcal{C}_1\times\mathcal{C}_2$ in $\mathbb{C}^2\times \mathbb{C}^2(:=(u_1, v_1, u_2, v_2))$; equivalently, each $v_i$ depends only on $u_i$ ($i=1,2$).}
\end{definition}

The above definition is due to Neumann-Reid in \cite{rigidity}, which was further studied and developed by D. Calegari in \cite{calegari1, calegari}. Basically, it means that the two cusps of $\mathcal{M}$ deform independently without affecting each other. 

Another immediate corollary of Lemma \ref{20102501} is the one-cusped analogue of our main results:\footnote{If the order of $A$ in Lemma \ref{20102501} is $6$ (resp. $4$), then $A^3=-I$ (resp. $A^2=-I$). Hence, $A^i(p/q)=A^{i+3}(p/q)$ (resp. $A^i(p/q)=A^{i+2}(p/q)$) for any coprime $p, q$ and for $0\leq i\leq 2$ (resp. $0\leq i\leq 1$). Consequently, the order of $\sigma$ is halved in Corollary \ref{23071103}.}

\begin{corollary}(Theorem 1.9 in \cite{CHDF})\label{23071103}
Let $\mathcal{M}$ be a one-cusped hyperbolic $3$-manifold with the cusp shape $\tau$. Suppose $\mathcal{M}_{p'/q'}$ and $\mathcal{M}_{p/q}$ be two Dehn fillings of $\mathcal{M}$ of the same pseudo complex volume with sufficiently large $|p'|+|q'|$ and $|p|+|q|$. If $\tau\in \mathbb{Q}(\sqrt{-3})$ (resp. $\mathbb{Q}(\sqrt{-1})$), then there exists $\sigma$ of order $3$ (resp. $2$) such that $p'/q'=\sigma^i(p/q)$ for some $0\leq i\leq 2$ (resp. $0\leq i\leq 1$). Otherwise, $p'/q'=p/q$. In particular, $\sigma$ is uniquely determined up to its inverse. 
\end{corollary}

In fact, what was proved in \cite{CHDF} (Theorem 4.3 therein) is more general than the above corollary, while the statement itself, appearing as Corollary \ref{23071103}, has been classically known (cf. \cite{JO}).\\
 
\subsubsection{Auxiliary theorems II}\label{25110409}
In this sub-subsection, we restate Theorem \ref{22050203}(2) in the context of the analytic holonomy variety $\mathcal{X}_{an}$ and present its several refinements of it. We also clarify the notions of Types I, II and III elements introduced in Section \ref{25101710}. The results and definitions discussed here are relevant to Theorem \ref{25101401} and are seen as preparatory steps toward it.

If we consider the analytic holonomy set corresponding to $\mathcal{X}_{an}\times \mathcal{X}_{an}$, the two identical copies of $\mathcal{X}_{an}$, embedded in $\mathbb{C}^{8}\big(:=(u_1, v_1, \dots, u'_2, v'_2)\big)$, then for $H$ as in \eqref{20071501}, the analytic counterpart of $(\mathcal{X}\times \mathcal{X})\cap H$ is the complex manifold defined by 
\begin{equation}\label{20091903}
\begin{gathered}
\left(\begin{array}{cccc}
a'_1 & b'_1 & 0 & 0\\
c'_1 & d'_1 & 0 & 0\\
0 & 0 & a'_4 & b'_4\\
0 & 0 & c'_4 & d'_4 
\end{array}\right)
\left(\begin{array}{c}
u'_1\\
v'_1\\
u'_2\\
v'_2
\end{array}\right)
=\left(\begin{array}{cccc}
a_1 & b_1 & a_2 & b_2\\
c_1 & d_1 & c_2 & d_2\\
a_3 & b_3 & a_4 & b_4\\
c_3 & d_3 & c_4 & d_4 
\end{array}\right)
\left(\begin{array}{c}
u_1\\
v_1\\
u_2\\
v_2
\end{array}\right).
\end{gathered}
\end{equation}
To simplify the problem, let us assume 
\begin{equation}\label{22022007}
\left(\begin{array}{cc}
a'_1 & b'_1 \\
c'_1 & d'_1 
\end{array}\right)=
\left(\begin{array}{cc}
a'_4 & b'_4 \\
c'_4 & d'_4 
\end{array}\right)=I
\end{equation} 
and 
\begin{equation}\label{22022005}
\left(\begin{array}{cccc}
a_1 & b_1 & a_2 & b_2\\
c_1 & d_1 & c_2 & d_2\\
a_3 & b_3 & a_4 & b_4\\
c_3 & d_3 & c_4 & d_4 
\end{array}\right)\in \text{Mat}_{4\times 4}(\mathbb{Q}).
\end{equation}
Denoting \eqref{22022005} by $M$ and following the same notation given in \cite{CHDF}, we let 
\begin{equation}\label{22022010}
A_i:=\left(\begin{array}{cc}
a_i & b_i \\
c_i & d_i 
\end{array}
\right)\quad (1\leq i\leq 4) 
\end{equation}
and represent $M$ as 
\begin{equation}\label{21072901}
\left(\begin{array}{cc}
A_1 & A_2\\
A_3 & A_4
\end{array}
\right).
\end{equation}

Using the above notation, a more quantitative version of Theorem \ref{22050203}(2) is stated below:

\begin{theorem}(Corollary 8.12 in \cite{CHDF})\label{23070517}
Let $\mathcal{M}$ be a two-cusped hyperbolic $3$-manifold and $\mathcal{X}_{an}$ be its analytic holonomy set. Let $\mathcal{M}_{p'_{1}/q'_{1}, p'_{2}/q'_{2}}$ and $\mathcal{M}_{p_{1}/q_{1}, p_{2}/q_{2}}$ be two Dehn fillings of $\mathcal{M}$ with sufficiently large $|p'_i|+|q'_i|, |p_i|+|q_i|$ ($i=1,2$) satisfying
\begin{equation*}
\mathrm{pvol}_{\mathbb{C}}\;\mathcal{M}_{p'_{1}/q'_{1}, p'_{2}/q'_{2}}=\mathrm{pvol}_{\mathbb{C}}\;\mathcal{M}_{p_{1}/q_{1}, p_{2}/q_{2}}. 
\end{equation*}
We further suppose the complex lengths of the two core geodesics of $\mathcal{M}_{p_{1}/q_{1}, p_{2}/q_{2}}$ are linearly independent over $\mathbb{Q}$. Then there exists a finite set $\mathcal{G}(\subset \mathrm{GL}_4(\mathbb{Q}))$ of matrices such that a Dehn filling point associated to the pair is contained in a $2$-dimensional analytic subset of $\mathcal{X}_{an}\times \mathcal{X}_{an}$, defined by 
\begin{equation}\label{22022101}
\left(\begin{array}{c}
u'_1\\
v'_1\\
u'_2\\
v'_2
\end{array}\right)
=M\left(\begin{array}{c}
u_1\\
v_1\\
u_2\\
v_2
\end{array}\right)
\end{equation}
where $M\in \mathcal{G}$. Further, when $M$ is represented as in \eqref{21072901}, either one of the following holds:
\begin{enumerate}
\item if $\det A_2=0$ or $\det A_3=0$, then 
\begin{equation}\label{23120801}
\det A_1=\det A_4=1\quad \text{and}\quad  A_2=A_3=0;
\end{equation}
\item if $\det A_1=0$ or $\det A_4=0$, then 
\begin{equation}\label{23120802}
\det A_3=\det A_2=1\quad \text{and}\quad A_1=A_4=0;
\end{equation}
\item if $\det A_j\neq 0\;(1\leq j\leq 4)$, then 
\begin{equation}\label{23071302}
\det A_1=\det A_4,\;\;  \det A_2=\det A_3,\;\; \det A_1+\det A_3=1,\;\; (\det A_1) A_1^{-1}A_2=-(\det A_3)A_3^{-1}A_4.
\end{equation}
\end{enumerate}
\end{theorem}

The following is an updated version of the relation appearing in \eqref{24080403}.

\begin{theorem}(Theorem 7.1 in \cite{CHDF})\label{22031403}
Let $\mathcal{M}, \mathcal{X}_{an}, \mathcal{M}_{p'_{1}/q'_{1}, p'_{2}/q'_{2}}$ and $\mathcal{M}_{p_{1}/q_{1}, p_{2}/q_{2}}$ be the same as above. Suppose a Dehn filling point associated to the pair is contained in a subset of $\mathcal{X}_{an}\times \mathcal{X}_{an}$ defined by \eqref{22022101} where $M$ is as in \eqref{21072901}. Then there exist $k_j\in \mathbb{Q}$ ($1\leq j\leq 4$) such that 
\begin{equation}\label{22031407}
\begin{gathered}
\left( \begin{array}{cc}
p'_1 & q'_1
\end{array}\right)
A_1=k_1\left( \begin{array}{cc}
p_1 & q_1\\
\end{array}\right), \quad \left( \begin{array}{cc}
p'_1 & q'_1
\end{array}\right)
A_2
=k_2\left( \begin{array}{cc}
p_2 & q_2
\end{array}\right), \\
\left( \begin{array}{cc}
p'_2  & q'_2\\
\end{array}\right)
A_3
=k_3\left( \begin{array}{cc}
p_1 & q_1
\end{array}\right), \quad 
\left( \begin{array}{cc}
p'_2  & q'_2\\
\end{array}\right)
A_4
=k_4 \left( \begin{array}{cc}
p_2 & q_2
\end{array}\right),
\end{gathered}
\end{equation}
and
\begin{equation}\label{25120101}
k_1+k_2=k_3+k_4=1.
\end{equation}
Moreover, if $k_j\neq 0$ for $1\leq j\leq 4$, it further satisfies
\begin{equation}\label{25120102}
\frac{\det A_1}{k_1}+\frac{\det A_3}{k_3}=\frac{\det A_2}{k_2}+\frac{\det A_4}{k_4}=1.
\end{equation} 
\end{theorem}


The following includes the precise definitions of the notions introduced heuristically in Section \ref{25101710}:
 
\begin{definition}\label{24060701}
\normalfont{Let $M\in \mathrm{GL}_4(\mathbb{Q})$ be a matrix represented as in \eqref{21072901}. 
\begin{enumerate}
\item We say $M$ is of block diagonal or anti-diagonal matrix if $A_2=A_3=0$ or $A_1=A_4=0$, respectively. 
\item We say $M$ is of Type I, II, or III matrix if it satisfies \eqref{23120801}, \eqref{23120802}, or \eqref{23071302} respectively.
\end{enumerate}}
\end{definition}

\noindent \textbf{Remark. }For $M\in \mathrm{Aut}\,\mathcal{X}$ given as in \eqref{21072901}, $M$ being of block diagonal matrix is equivalent to $M$ being of Type I. This is because if $M$ is of block diagonal matrix, then, as $M$ has a finite order, both $A_1$ and $A_4$ are also matrices of finite order. Consequently, each of $\det A_1$ and $\det A_4$ is either $1$ or $-1$. Using the fact that $\tau_1,\tau_2\in \mathbb{C}\setminus\mathbb{R}$ and \eqref{24080403}, it follows $\det A_1, \det A_4>0$ (see Corollary \ref{25010913}(3)), i.e., $M$ is of Type I. \\

Clearly the product of two Type I (or Type II) elements is of Type I, and the product of Type I and Type II elements is of Type II. In the following lemma, we extend this property to include elements of Type III. The lemma will be particularly useful for simplifying the proofs of various arguments in Section \ref{25012005} to analyze the structure of a subgroup of $\mathrm{Aut}\,\mathcal{X}$.

\begin{lemma}\label{25020503}
If $M$ (resp. $N$) is a matrix of Type III (resp. Type I or II), then $NM$ is of Type III. 
\end{lemma}
\begin{proof}
We prove only the second case, as the proof for the first case is analogous. Suppose $M$ and $N$ are given as 
\begin{equation*}
\left( \begin{array}{cc}
A_1 & A_2\\
A_3 & A_4
\end{array}\right)\;\; \text{and}\;\; 
\left( \begin{array}{cc}
0 & B_2\\
B_3 & 0
\end{array}\right)
\end{equation*}
respectively. Then 
\begin{equation*}
NM=\left( \begin{array}{cc}
B_2A_3 & B_2A_4\\
B_3A_1 & B_3A_2
\end{array}\right), 
\end{equation*}
and  
\begin{equation*}
\begin{gathered}
\det (B_2A_3)=\det A_3=\det A_2=\det (B_3A_2), \; \det (B_2A_4)=\det A_4=\det A_1=\det (B_3A_1),\\
\det (B_2A_3)+\det (B_3A_1)=\det A_3+\det A_1=1. 
\end{gathered}
\end{equation*} 
Since $\frac{\det A_1}{1-\det A_1}A_3A_1^{-1}A_2=A_4$, we obtain
\begin{equation*}
\frac{\det (B_2A_3)}{1-\det (B_2A_3)}B_3A_1(B_2A_3)^{-1}B_2A_4=\frac{1-\det A_1}{\det A_1}B_3A_1A_3^{-1}A_4=B_3A_2. 
\end{equation*}
Thus $NM$ is of Type III. 
\end{proof}

\begin{convention}
{\normalfont Throughout the remainder of the paper, we will no longer refer to the (algebraic) holonomy variety $\mathcal{X}$ and will instead focus on its analytic counterpart, $\mathcal{X}_{an}$. For brevity, and with a slight abuse of notation, we shall denote $\mathcal{X}_{an}$ simply by $\mathcal{X}$ from now on.}
\end{convention}

\subsection{Linear Algebra}

Let $\mathcal{M}$ be a two-cusped hyperbolic $3$-manifold and $ \mathcal{X}$ be its analytic holonomy set. If the manifold defined by \eqref{22022101} is a $2$-dimensional anomalous subset of $\mathcal{X}\times \mathcal{X}$, then, equivalently, the following linear map
\begin{align*}
\mathcal{X}&\longrightarrow \mathcal{X}\\
\bold{x}&\mapsto  M\bold{x}
\end{align*}
is an element of $\mathrm{Aut}\,\mathcal{X}$, the automorphism group of $\mathcal{X}$. Thus one may view   
\begin{equation}\label{23070511}
\Big\{M\in \mathrm{GL}_4(\mathbb{Q})\;:\;\left(\begin{array}{c}
u'_1\\
v'_1\\
u'_2\\
v'_2
\end{array}\right)
=M\left(\begin{array}{c}
u_1\\
v_1\\
u_2\\
v_2
\end{array}\right)\text{ is a }2\text{-dimensional anomalous subset of }\mathcal{X}\times \mathcal{X}\}
\end{equation}
as a subgroup of $\mathrm{Aut}\,\mathcal{X}$, which, moreover, is finite according to Theorem 1.4 in \cite{za}. 

For $M\in \mathrm{GL}_n(\mathbb{Q})$ in general, by the \textit{minimal polynomial} $\mathfrak{m}_M(x)\in \mathbb{Q}[x]$ of $M$, it means the monic polynomial of the smallest degree such that $\mathfrak{m}_M(M)=0$. The following is one of the basic theorems in linear algebra. 
\begin{theorem}[Cayley-Hamiton]
For $M\in \mathrm{GL}_n(\mathbb{Q})$, the minimal polynomial of $\mathfrak{m}_M(x)$ divides the characteristic polynomial $\mathfrak{\chi}_M(x)$ of $M$. 
\end{theorem}

For $M$ as in \eqref{23070511}, since it is of finite order and the degree of its minimal polynomial $\mathfrak{m}_M(x)$ is at most $4$, the candidates for $\mathfrak{m}_M(x)$ are determined by considering all the possible combinations of cyclotomic polynomials, as listed below:
\begin{equation}\label{23062301}
\begin{gathered}
x\pm 1, \;\; x^2\pm 1, \;\; x^2\pm x+1,  \;\; x^3\pm 1, \;\;  x^3\pm x^2+x\pm 1,\;\;x^3\pm 2x^2+2x\pm 1,\\
x^4\pm 1,\;\; x^4\pm x^2+1,\; \; x^4\pm x^3+x^2\pm x+1. 
\end{gathered}
\end{equation}

If $\mathfrak{m}_M(x)$ is either $x-1$ or $x+1$, then $M$ is obviously either $I$ or $-I$ respectively. 

However, not every polynomial in \eqref{23062301} arises as the minimal polynomial of some $M$ in \eqref{23070511}. In Section \ref{23071105}, we further refine \eqref{23062301} by dropping the last two polynomials, i.e. $x^4\pm x^3+x^2\pm x+1$, from the list.

After removing these two, we will explore the remaining cases and explicitly classify the shapes of $M$ for each minimal polynomial. This will be one of the central goals of the paper.\footnote{See Theorems \ref{24090201}, \ref{25012901}, and \ref{24092702} for the complete classification of $M$ of Types I, II, and III, respectively.} 

The concepts and theorem introduced below are fundamental in linear algebra.

The \textit{companion matrix} of a monic polynomial 
\begin{equation*}
f=x^n+a_{n-1}x^{n-1}+\cdots+a_1x+a_0
\end{equation*}
is the $n\times n$ matrix
\begin{equation*}
C(f):=\left(\begin{array}{cccccc}
0 & 0 & 0 & \cdots & 0 & -a_0\\
1 & 0 & 0 & \cdots & 0 & -a_1\\
0 & 1 & 0 & \cdots & 0 & -a_2\\
\vdots & \vdots & \vdots & \ddots & \vdots & \vdots\\
0 & 0 & 0 & \cdots & 0 & -a_{n-1}
\end{array}\right).
\end{equation*}
Given matrices $C_1$ and $C_2$, their \textit{direct sum} is the block diagonal matrix 
\begin{equation*}
C_1\oplus C_2
=\left(\begin{array}{cc}
C_1 & 0 \\
0 & C_2
\end{array}\right).
\end{equation*}
We denote the direct sum of $j$ copies of a $C$ by $C^{[j]}$. 

The following theorem is both well-known and pivotal in linear algebra. 
\begin{theorem}[Primary Decomposition Theorem]\label{23070515}
If the characteristic polynomial of $M$ is given by 
\begin{equation*}
\mathfrak{\chi}_M(x)=p_1(x)^{e_1}\cdots p_r(x)^{e_r}
\end{equation*}
where $p_i$ are distinct monic irreducible polynomials, then $M$ is similar to 
\begin{equation*}
C(p_1)^{[e_1]}\oplus \cdots \oplus C(p_r)^{[e_r]}.
\end{equation*}
\end{theorem}

\begin{convention}
{\normalfont For later use, we introduce 
$$C_1\widetilde{\oplus} C_2$$
to denote the block anti-diagonal matrix 
\begin{equation*}
\left(\begin{array}{cc}
0 & C_1 \\
C_2 & 0
\end{array}\right).
\end{equation*}}
\end{convention}

\subsection{Change of Variables}\label{25011001}
One useful technique for analyzing the analytic holonomy set is a change of variables. For instance, consider the following transformation 
\begin{equation*}
\begin{gathered}
\left(\begin{array}{c}
\tilde{u}'_i\\
\tilde{v}'_i
\end{array}\right)
=\tilde{A}_i\left(\begin{array}{c}
u'_i\\
v'_i
\end{array}\right), \quad 
\left(\begin{array}{c}
\tilde{u}_i\\
\tilde{v}_i
\end{array}\right)
=\tilde{A}_i\left(\begin{array}{c}
u_i\\
v_i
\end{array}\right)  
\end{gathered}
\end{equation*} 
where $\tilde{A}_i\in \mathrm{GL}_2(\mathbb{Q})$ ($i=1,2$). Under this change of variables, the equation   
\begin{equation*}
\left(\begin{array}{c}
u'_1\\
v'_1\\
u'_2\\
v'_2
\end{array}\right)
=M\left(\begin{array}{c}
u_1\\
v_1\\
u_2\\
v_2
\end{array}\right)
\end{equation*}
is transformed into 
\begin{equation*}
\left(\begin{array}{c}
\tilde{u}'_1\\
\tilde{v}'_1\\
\tilde{u}'_2\\
\tilde{v}'_2
\end{array}\right)
=\left(\begin{array}{cc}
 \tilde{A}_1 & 0\\
 0 & \tilde{A}_2
\end{array}\right)
M\left(\begin{array}{cc}
\tilde{A}_1 & 0\\
 0 & \tilde{A}_2
\end{array}\right)^{-1}
\left(\begin{array}{c}
\tilde{u}_1\\
\tilde{u}_1\\
\tilde{u}_2\\
\tilde{u}_2
\end{array}\right),
\end{equation*}
and the symmetries
\begin{equation*}
\frac{\partial {v}'_2}{\partial {u}'_1}=\frac{\partial {v}'_1}{\partial {u}'_2}\quad\text{and}\quad \frac{\partial {v}_2}{\partial {u}_1}=\frac{\partial {v}_1}{\partial {u}_2}
\end{equation*}
are turned into
\begin{equation}\label{23062307}
\frac{\det \tilde{A}_1}{\det \tilde{A}_2}\frac{\partial \tilde{v}'_2}{\partial \tilde{u}'_1}=\frac{\partial \tilde{v}'_1}{\partial \tilde{u}'_2}\quad\text{and}\quad \frac{\det \tilde{A}_1}{\det \tilde{A}_2}\frac{\partial \tilde{v}_2}{\partial \tilde{u}_1}=\frac{\partial \tilde{v}_1}{\partial \tilde{u}_2}
\end{equation}
by the chain rule. 

The technique will be frequently invoked to simplify many proofs throughout the paper. In particular, for $M$ given as in \eqref{21072901}, we often conjugate $M$ to
\begin{equation}\label{24100301}
\left(\begin{array}{cc}
 I & 0\\
 0 & A_2
\end{array}\right)
M\left(\begin{array}{cc}
 I & 0\\
 0 & A_2
\end{array}\right)^{-1}
=\left(\begin{array}{cc}
A_1 & I\\
A_2A_3 & A_2A_4A_2^{-1}
\end{array}\right)
\end{equation}
and work with \eqref{24100301}, rather than the original $M$, as it offers significant computational advantages in various aspects. (See the discussion following Proposition \ref{24082301}.) 

Typically, when $\mathcal{X}$ is the analytic holonomy variety of a cusped hyperbolic $3$-manifold, a change of variables on $\mathcal{X}$ may result in a manifold that no longer belongs to the category of analytic holonomy sets. For instance, this occurs when $\frac{\det \tilde{A}_1}{\det \tilde{A}_2} \neq 1$, as established in \eqref{23062307}. To accommodate such manifolds, we further enlarge the category of analytic holonomy sets. By extracting the essential properties described in Theorem \ref{potential}, we generalize the definition as follows:
\begin{definition}\label{23120604}
\normalfont{We call an analytic set $\mathcal{X}$ defined by holomorphic functions 
\begin{equation}\label{23120805}
v_i(u_1,u_2)=\tau_i u_i+\sum_{\substack{l=3\\ \text{odd}}}^{\infty}\Big(\sum_{j+k=l} c^i_{jk}u_1^ju_2^k\Big)
\end{equation}
($i=1,2$) in $\mathbb{C}^4(:=(u_1, v_1, u_2, v_2))$ a \textit{generalized analytic holonomy set} if it satisfies
\begin{itemize}
\item $\tau_i\in \mathbb{C}\setminus \mathbb{R}$ and $\mathrm{Im}\,\tau_i>0$ ($i=1,2$);
\item $v_1(u_1,0)\neq \tau_1 u_1$ and $v_2(0,u_2)\neq \tau_2 u_2$;
\item $v_1(u_1, u_2)=v_1(u_1, -u_2), \; v_1(u_1, u_2)=-v_1(-u_1, u_2),\;  v_2(u_1, u_2)=v_2(-u_1, u_2), \\ 
v_2(u_1, u_2)=-v_2(u_1, -u_2)$;
\item $\frac{\partial v_1}{\partial u_2}=a\frac{\partial v_2}{\partial u_1}$ for some $a\in \mathbb{C}$;
\item every subgroup of $\mathrm{Aut}\,\mathcal{X}$ contained in $\mathrm{GL}_4(\mathbb{Q})$ is a finite set. 
\end{itemize}
We denote the set of all generalized analytic holonomy sets by $\mathfrak{Ghol}$.}
\end{definition}

Consequently, the analytic holonomy set $\mathcal{X}$ of any cusped hyperbolic $3$-manifold is contained in $\mathfrak{Ghol}$, and any manifold, obtained via a change of variables from $\mathcal{X}$, is also an element of $\mathfrak{Ghol}$. 

\begin{convention}
{\normalfont For $\mathcal{X} \in \mathfrak{Ghol}$ as in \eqref{23120805}, by abuse of notation and for the sake of brevity, we still refer to $\tau_1$ and $\tau_2$ as the \textit{cusp shapes} of $\mathcal{X}$, and assume that they belong to the same quadratic field throughout the paper, unless otherwise stated.}
\end{convention}

\newpage
\section{Proof of Theorem \ref{250825011}}\label{25101507}
In this section, we prove Theorem \ref{250825011}. More precisely, we establish the following detailed version of it.

\begin{theorem}\label{25082501}
Let $\mathcal{M}$ be a two-cusped hyperbolic $3$-manifold whose cusp shapes belong to the same quadratic field, and let $\mathcal{M}_{p'_1/q'_1, p'_2/q'_2}$ and $\mathcal{M}_{p_1/q_1, p_2/q_2}$ be two Dehn fillings of $\mathcal{M}$ having the same pseudo complex volume with sufficiently large $|p'_i|+|q'_i|$ and $|p_i|+|q_i|$ ($i=1,2$). Assume that the complex lengths of the core geodesics of $\mathcal{M}_{p_1/q_1, p_2/q_2}$ are linearly dependent over $\mathbb{Q}$.
\begin{enumerate}
\item Suppose the potential function of $\mathcal{M}$ contains a non-trivial term of degree $4$. Then the following statements hold. 
\begin{enumerate}
\item If $\tau_1, \tau_2\in \mathbb{Q}(\sqrt{-3})$, then there exist unique $\sigma_1, \sigma_2\in \mathrm{GL}_2(\mathbb{Q})$, each of order $3$, such that 
\begin{equation}\label{25102607}
p'_1/q'_1=\sigma_1^i(p_1/q_1),\quad p'_2/q'_2=\sigma_2^i(p_2/q_2)
\end{equation}
for some $0\leq i\leq 2$. 
\item Suppose $\tau_1, \tau_2\in \mathbb{Q}(\sqrt{-1})$.  
\begin{enumerate}
\item There exist unique $\sigma_1, \sigma_2\in \mathrm{GL}_2(\mathbb{Q})$, each of order $2$, such that $(p'_1/q'_1, p'_2/q'_2)$ is one of  
\begin{equation}\label{25102608}
(\sigma_1^i(p_1/q_1), \sigma_2^i(p_2/q_2)), \;\;(p_1/q_1, \sigma_2^i(p_2/q_2)) 
 \;\;\text{or}\;\;(\sigma_1^i(p_1/q_1), p_2/q_2)
\end{equation}
for some $0\leq i\leq 1$. In each case, there are only two possibilities corresponding to $i=0,1$. 

\item There exist $\sigma_1, \sigma_2\in \mathrm{GL}_2(\mathbb{Q})$, both non-cyclic, such that $(p'_1/q'_1, p'_2/q'_2)$ is   
\begin{equation}\label{25102609}
(\sigma_1(p_1/q_1),\sigma_2(p_2/q_2)). 
\end{equation}
In particular, there are at most two possible choices for $(\sigma_1, \sigma_2)$. 

\end{enumerate}
\item Otherwise, $(p'_1/q'_1, p'_2/q'_2)=(p_1/q_1,p_2/q_2)$. 
\end{enumerate}  

\item Suppose two cusps of $\mathcal{M}$ are SGI to each other. If $\tau_1, \tau_2\in \mathbb{Q}(\sqrt{-3})$ (resp. $\mathbb{Q}(\sqrt{-1})$), then there exist unique $\sigma_1, \sigma_2\in \mathrm{GL}_2(\mathbb{Q})$, both of order $3$ (resp. $2$), such that 
\begin{equation}\label{25112403}
p'_1/q'_1=\sigma_1^i(p_1/q_1),\quad p'_2/q'_2=\sigma_2^j(p_2/q_2)
\end{equation}
for some $0\leq i, j\leq 2$ (resp. $0\leq i, j\leq 1$). Otherwise, $(p'_1/q'_1, p'_2/q'_2)=(p_1/q_1,p_2/q_2)$. 
\end{enumerate}
\end{theorem}

\subsection{Preparations}\label{25102903}

Before presenting the proof, we revisit Theorem \ref{22050203}(1), rephrase it in terms of the analytic holonomy set. The proof of Theorem \ref{25082501} will then follow naturally by quantifying this reformulation, relying primarily on either Theorem \ref{25111301} or Lemma \ref{20102501}. 

First, in the analytic setting, the equations in \eqref{220708011} are equivalent to 
\begin{equation}\label{25102201}
\begin{gathered}
A_1\left(\begin{array}{c}
u_1 \\
v_1  
\end{array}\right)+
A_2\left(\begin{array}{c}
u_2 \\
v_2  
\end{array}\right)
=A'_1\left(\begin{array}{c}
u'_1 \\
v'_1  
\end{array}\right)+
A'_2\left(\begin{array}{c}
u'_2 \\
v'_2  
\end{array}\right)=
A_3\left(\begin{array}{c}
u_1 \\
v_1  
\end{array}\right)
+A'_3\left(\begin{array}{c}
u'_1 \\
v'_1  
\end{array}\right)=0
\end{gathered}
\end{equation} 
over $\mathcal{X}\times \mathcal{X}$ in $\mathbb{C}^8(:=(u_1, v_1, \dots, u'_2, v'_2))$. 

By Theorem \ref{22050203}(1), since the complex lengths of the core geodesics of $\mathcal{M}_{p_1/q_1,p_2/q_2}$ are linearly dependent over $\mathbb{Q}$, the same holds for those of $\mathcal{M}_{p'_1/q'_1,p'_2/q'_2}$. By Theorem \ref{22050501}, this implies 
\begin{equation}\label{25102203}
\begin{gathered}
A_1\left(\begin{array}{c}
u_1 \\
v_1  
\end{array}\right)+
A_2\left(\begin{array}{c}
u_2 \\
v_2  
\end{array}\right)=0\quad (\text{resp. }A'_1\left(\begin{array}{c}
u'_1 \\
v'_1  
\end{array}\right)+A'_2\left(\begin{array}{c}
u'_2 \\
v'_2  
\end{array}\right)=0)
\end{gathered}
\end{equation}
defines a $1$-dimensional analytic subset $\mathcal{C}$ (resp. $\mathcal{C}'$) of $\mathcal{X}$ in $\mathbb{C}^4$ with coordinates $(u_1, v_1, u_2, v_2)$ (resp. $(u'_1, v'_1, u'_2, v'_2)$). 

Projecting each $\mathcal{C}$ (resp. $\mathcal{C}'$) to $\mathbb{C}^2$ with coordinates $(u'_1, u'_2)$ (resp. $(u_1, u_2)$), the last equation 
\begin{equation}\label{25102204}
A_3\left(\begin{array}{c}
u_1 \\
v_1  
\end{array}\right)+A'_3\left(\begin{array}{c}
u'_1 \\
v'_1  
\end{array}\right)=0
\end{equation}
in \eqref{25102201} can be regarded as a $1$-dimensional analytic subset of $\mathcal{C}\times \mathcal{C}'$ in $\mathbb{C}^4$ (with coordinates $(u_1, v_1, u'_1, v'_1)$). 

As a special case, if two cusps of $\mathcal{M}$ are SGI to each other, then $\mathcal{C}$ and $\mathcal{C}'$ above are exactly the same. In this case, the order of $-A_3^{-1}A'_3$ is $4$ or $6$ for $\tau_1, \tau_2\in \mathbb{Q}(\sqrt{-1})$ or $\mathbb{Q}(\sqrt{-3})$ respectively, and $2$ otherwise, by Lemma \ref{20102501}. Hence, the desired $\sigma_1$ in \eqref{25112403} is obtained from \eqref{25101504} in Theorem \ref{22050203}. Analogously, one finds $\sigma_2$ in \eqref{25112403}.  

Consequently, Theorem \ref{25082501}(2) is readily derived except for the uniqueness part, which will be verified in Lemma \ref{24090301} below. Hence, in the remainder of this section, we solely focus on the first case of the theorem.

In the following lemma, as a preliminary step toward the proof of Theorem \ref{25082501}(1), we derive two equations relating coefficients of the potential function of $\mathcal{M}$ to the eigenvalues of $A_i$ and $A'_i$. 

\begin{lemma}\label{25102901}
Having the same notation and assumptions as above, to simplify the proof, we may, by changing basis if necessary, further assume that 
\begin{equation}\label{25101607}
A'_1=I\quad \text{and}\quad A'_2=-I. 
\end{equation}
Set  
\begin{equation}\label{25102602}
\left(\begin{array}{cc}
a & b\\
c & d
\end{array}\right):=
-A_2^{-1}A_1, \;\;\left(\begin{array}{cc}
\alpha & \beta\\
\gamma & \delta
\end{array}\right):=
-(A'_3)^{-1}
A_3,
\end{equation}
and let
\begin{equation*}
c_{4,0}u_1^4+c_{2,2}u_1^2u_2^2+c_{0,4}u_2^4
\end{equation*}
be the terms of the potential function of $\mathcal{M}$ of homogeneous degree $4$. Then
\begin{equation}\label{25100705}
\begin{aligned}
2c_{4,0}(a+b\tau)^3+c_{2,2}(a+b\tau)&=(d-b\tau)(2c_{4,0}+c_{2,2}(a+b\tau)^2),\\
(\delta -\beta \tau)(2c_{4,0}+c_{2,2}(a+b\tau)^2)&=(\alpha +\beta \tau)^3(2c_{4,0}+c_{2,2})
\end{aligned}
\end{equation}
where $\tau:=\tau_1(=\tau_2)$.\footnote{The equality $\tau_1=\tau_2$ follows from the assumption \eqref{25101607}.}
\end{lemma}
\begin{proof}
By the assumption in \eqref{25101607}, since 
\begin{equation*}
\left(\begin{array}{c}
u'_1 \\
v'_1  
\end{array}\right)=\left(\begin{array}{c}
u'_2 \\
v'_2  
\end{array}\right)
\end{equation*}
defines a $1$-dimensional anomalous subvariety of $\mathcal{X}$, we have $c_{4, 0}=c_{0, 4}$. By the notation introduced, the first equation in \eqref{25102203} is restated as 
\begin{equation}\label{25082502}
\left(\begin{array}{c}
u_2\\
v_2
\end{array}\right)
=\left(\begin{array}{cc}
a & b\\
c & d
\end{array}\right)
\left(\begin{array}{c}
u_1\\
v_1
\end{array}\right).
\end{equation}
As discussed above, since \eqref{25082502} defines a $1$-dimensional anomalous subvariety of $\mathcal{X}$, by substituting $u_2=au_1+bv_1$ into $v_2=cu_1+dv_1$, we get 
\begin{equation}\label{25102205}
\begin{aligned}
&\tau (au_1+bv_1)+2c_{4,0}(au_1+bv_1)^3+c_{2,2}(au_1+bv_1)u_1^2+\cdots\\
=&cu_1+d(\tau u_1+2c_{4,0}u_1^3+c_{2,2}u_1(au_1+bv_1)^2+\cdots).
\end{aligned}
\end{equation}
Extracting the terms of degrees $3$ in $u_1$ from \eqref{25102205} and comparing their coefficients, we attain
\begin{equation}\label{25100101}
\begin{aligned}
2c_{4,0}(a+b\tau)^3+c_{2,2}(a+b\tau)&=(d-b\tau)(2c_{4,0}+c_{2,2}(a+b\tau)^2), 
\end{aligned}
\end{equation}
which gives the first equality in \eqref{25100705}. 

The second equality in \eqref{25100705} can be derived in a similar manner from the fact that \eqref{25102204} defines a $1$-dimensional subvariety of $\mathcal{C}\times \mathcal{C}'$, and hence the proof is omitted here.
\end{proof}

In the above proof, note that 
\begin{equation*}
b\tau^2+(a-d)\tau-c=0\Longrightarrow \tau=\frac{d-a\pm \sqrt{(a-d)^2+4bc}}{2b}
\end{equation*}
and, as $\mathrm{Im}\,\tau>0$, it follows that 
\begin{equation}\label{25111001}
a+b\tau=\frac{a+d+\sqrt{(a+d)^2+4(ad-bc)}}{2} \quad (\text{resp. }\frac{a+d-\sqrt{(a+d)^2-4(ad-bc)}}{2})
\end{equation}
for $b>0$ (resp. $b<0$). Similarly, 
\begin{equation}\label{25111002}
d-b\tau=\frac{a+d-\sqrt{(a+d)^2+4(ad-bc)}}{2} \quad (\text{resp. }\frac{a+d+\sqrt{(a+d)^2-4(ad-bc)}}{2})
\end{equation}
for $b>0$ (resp. $b<0$). That is, $a+b\tau$ and $d-b\tau$ are the eigenvalues of $\left(\begin{array}{cc}
a & b\\
c & d
\end{array}\right)$. 

Analogously, $\alpha+\beta\tau$ and $\delta-\beta\tau$ are the eigenvalues of $\left(\begin{array}{cc}
\alpha & \beta\\
\gamma & \delta
\end{array}\right)$. 

\begin{convention}
\normalfont{To simplify notation, for \(-A_2^{-1}A_1\) as given in \eqref{25102602}, we denote
\[
(a+d)^2-4(ad-bd)
\]
by \(\mathrm{Disc}\,(-A_2^{-1}A_1)\) and call it the \textit{discriminant} of \(-A_2^{-1}A_1\).}
\end{convention}

The following lemma concerns the rigidity of a matrix of the above type. It shows that the matrix is uniquely determined, up to multiplication by its determinant and inversion, by one of its eigenvalues.

\begin{lemma}\label{24090301}
Let 
\begin{equation}\label{25110601}
A:=\left(\begin{array}{cc}
a & b\\
c & d
\end{array}\right)\quad \text{and}\quad  B:=\left(\begin{array}{cc}
\alpha & \beta\\
\gamma & \delta
\end{array}\right) 
\end{equation}
be two elements in $\mathrm{GL}_2(\mathbb{Q})$ such that 
\begin{equation}\label{24123101}
1/\tau=A(1/\tau)=B(1/\tau)
\end{equation}
for some $\tau\in \mathbb{Q}(\sqrt{D})\setminus \mathbb{Q}$ where $D$ is a negative integer. If $A$ and $B$ have a common eigenvalue, then $B$ is either $A$ or $(\det A)A^{-1}$.
\end{lemma}
\begin{proof}
By \eqref{24123101},  
\begin{equation}\label{25102601}
\begin{gathered}
\frac{a-d}{b}=\frac{\alpha-\delta}{\beta}, \quad \frac{\sqrt{\mathrm{Disc}\,A}}{b}=\frac{\pm\sqrt{\mathrm{Disc}\,B}}{\beta}.
\end{gathered}
\end{equation}
By the assumption, since $A$ and $B$ have a common eigenvalue, it follows that $\mathrm{tr}\,A=\mathrm{tr}\,B$ and $\det A=\det B$. Combining these with \eqref{25102601}, we attain $b=\pm \beta$ and so $a-d=\pm (\alpha-\delta)$. A straightforward computation shows that if $b=\beta$ and $a-d=\alpha-\delta$ (resp. $b=-\beta$ and $a-d=-\alpha+\delta$), then $B=A$ (resp. $B=(\det A)A^{-1}$).
\end{proof}

In the proof of Theorem \ref{25082501}, for given $c_{4,0}$ and $c_{2,2}$, we first determine the possible eigenvalues of \eqref{25100101} using \eqref{25100705}. Once these eigenvalues are obtained, the possible forms of the matrices in \eqref{25100101} are then determined as well, thanks to Lemma \ref{24090301}.

The lemma will be broadly applicable, not only to the proof of Theorem \ref{25082501} below but also throughout the paper, particularly in Section \ref{25012005}.

\subsection{Proof}\label{25111402}
Now we complete the proof of Theorem \ref{25082501}. As observed earlier, it is enough to prove the first case of the theorem.   

Before the proof, we note that, in the setting of Theorem \ref{25082501}(1), 
\begin{equation}\label{25121701}
p'_1/q'_1=\big(A_3^{\prime\,T}(A_3^T)^{-1}\big)(p_1/q_1)
\end{equation}
and
\begin{equation*}
p'_2/q'_2=\big(A_2^{\prime\,T}(A_1^{\prime\,T})^{-1}A_3^{\prime\,T}(A_3^T)^{-1}A_1^T(A_2^T)^{-1}\big)(p_2/q_2)
\end{equation*}
hold by \eqref{25101504}. If $A'_1$ and $A'_2$ are normalized as in \eqref{25101607}, the latter one further reduces to 
\begin{equation}\label{25121702}
p'_2/q'_2=\big(A_3^{\prime\,T}(A_3^T)^{-1}A_1^T(A_2^T)^{-1}\big)(p_2/q_2).
\end{equation}

In the proof, we will also make use of the following relation, established in Lemma 9.1 of \cite{CHDF}:
\begin{equation*}
1+\frac{\det A_1
}{\det A_2} 
=\frac{\det A_3}
{\det A'_3}\Big(1+\frac{\det A'_1}{\det A'_2}\Big).
\end{equation*}
In our case, by the assumptions in \eqref{25101607}-\eqref{25102602}, it is equivalent to 
\begin{equation}\label{22051402}
1+(ad-bc)=2(\alpha\delta-\beta\gamma).
\end{equation}

\begin{proof}[Proof of Theorem \ref{25082501}(1)]
To simplify notation, we denote $-A_2^{-1}
A_1$ and $-(A'_3)^{-1}A_3$ by $A$ and $B$ respectively as given in \eqref{25110601}. By letting 
\begin{equation*}
z:=a+b\tau\quad \text{and}\quad w:=\alpha+\beta \tau, 
\end{equation*}
we rewrite \eqref{25100705} as 
\begin{equation}\label{25090105}
\begin{aligned}
\left(\begin{array}{cc}
z^3-\overline{z} & z-\overline{z}z^2\\  
w^3-\overline{w} & w^3-\overline{w}z^2
\end{array}\right)
\left(\begin{array}{cc}
2c_{4,0}\\
c_{2,2}
\end{array}\right)
=\left(\begin{array}{c}
0\\  
0
\end{array}\right).
\end{aligned}
\end{equation}
If the determinant of the coefficient matrix in \eqref{25090105} is nonzero, then $c_{4,0}=c_{2,2}=0$, contradicting the assumption. Thus its determinant is $0$ and so
\begin{equation}\label{25093001}
\begin{aligned}
(z^3-\overline{z})(w^3-\overline{w}z^2)&=(z-\overline{z}z^2)(w^3-\overline{w})\\
\Longrightarrow 
w^3(z+\overline{z})(z^2-1)&=\overline{w}(z^3+z)(z^2-1). 
\end{aligned} 
\end{equation}

\begin{enumerate}
\item If $z^2=1$, then $A=\pm I$. Without loss of generality, if we set $A=I$, then $\mathcal{C}$ and $\mathcal{C}'$ are the same algebraic curve and so $B$ defines an automorphism of this algebraic curve. By Lemma \ref{20102501}, the order of $B$ is either $6$ or $4$, and accordingly, $\tau_1, \tau_2$ are contained in either $\mathbb{Q}(\sqrt{-3})$ or $\mathbb{Q}(\sqrt{-1})$ respectively. Consequently, by \eqref{25121701}-\eqref{25121702}, it falls into either \eqref{25102607} (when $\tau_1, \tau_2 \in \mathbb{Q}(\sqrt{-3})$) or \eqref{25102608} (when $\tau_1, \tau_2 \in \mathbb{Q}(\sqrt{-1})$).

\item If $z^2\neq 1$ and $z+\overline{z}=0$, then $z^3+z=0\Longrightarrow z=\pm \sqrt{-1}$, implying $A$ is of order $4$. By \eqref{22051402},  
\begin{equation}\label{25102604}
1+|z|^2=2|w|^2\Longrightarrow |w|=1,
\end{equation}
and so $\overline{w}=\frac{1}{w}$, and \eqref{25090105} is reduced to 
\begin{equation*}
\begin{aligned}
\Big(w^3-\frac{1}{w}\Big)2c_{4,0}+ \Big(w^3+\frac{1}{w}\Big)c_{2,2}=0.
\end{aligned}
\end{equation*}
\begin{enumerate}
\item If $c_{4,0}=0$, then $w^3+\frac{1}{w}=0$, contradicting the fact that $w(=\alpha+\beta\tau$) lies in $\mathbb{Q}(\sqrt{-1})$. 

\item If $c_{2,2}=0$, then $w^4=1\Longrightarrow w=\pm 1, \pm \sqrt{-1}$, implying $B=\pm I,\pm A$ by Lemma \ref{24090301}. Thus it falls into a case in \eqref{25102608} by \eqref{25121701}-\eqref{25121702}. 

\item If $c_{4,0}, c_{2,2}\neq 0$, then 
\begin{equation}\label{25110603}
\frac{w^4-1}{w^4+1}=-\frac{c_{2,2}}{2c_{4,0}}.
\end{equation} 
For given $c_{2,2}$ and $2c_{4,0}$, it follows that if $w$ is a root satisfying \eqref{25110603}, then the other roots are $-w$ and $\pm \sqrt{-1}w$. By Lemma \ref{24090301}, for each fixed $w$, there are at most two possibilities for $B$. We claim that only one of the two satisfies the required condition unless $B^4=I$ (equivalently, $w=\pm 1$ or $\pm \sqrt{-1}$). Adapting the notation introduced in Section \ref{25102903}, suppose both 
\begin{equation*}
\left(\begin{array}{c}
u'_1 \\
v'_1  
\end{array}\right)=B\left(\begin{array}{c}
u_1 \\
v_1  
\end{array}\right)\;\;\text{and}\;\;
\left(\begin{array}{c}
u'_1 \\
v'_1  
\end{array}\right)=(\det B)B^{-1}\left(\begin{array}{c}
u_1 \\
v_1  
\end{array}\right)
\end{equation*}
are $1$-dimensional subsets of $\mathcal{C}\times \mathcal{C}'$. By the remark preceding Lemma \ref{25102901}, if both $B$ and $(\det B)B^{-1}$ are viewed as isomorphisms $\mathcal{C}\rightarrow \mathcal{C}'$, then $\frac{1}{\det B}B^2$ is an automorphism of $\mathcal{C}$. By Lemma \ref{20102501}, since the cusp shapes are contained in $\mathbb{Q}(\sqrt{-1})$, we have $B^4=I$. (Hence, the situation corresponds to the one just discussed above.)

Otherwise, if neither $w=\pm 1$ nor $w=\pm \sqrt{-1}$ is a root of \eqref{25110603}, then $B$ is uniquely determined for each $w$ satisfying \eqref{25110603}. As noted above, since there are two roots satisfying \eqref{25110603} up to sign, it falls into the case described in \eqref{25102609} by \eqref{25121701}-\eqref{25121702}.

\end{enumerate}

\item Lastly, suppose $z^2\neq 1$ and $z+\overline{z}\neq 0$, and so 
\begin{equation}\label{25100801}
\begin{aligned}
w^3(z+\overline{z})=\overline{w}(z^3+z)
\Longrightarrow 
\frac{z(1+z^2)}{z+\overline{z}}=\frac{w^3}{\overline{w}}.
\end{aligned} 
\end{equation}
By \eqref{22051402} (or \eqref{25102604}), we have 
\begin{equation}\label{25102605}
1+|z|^2=2|w|^2. 
\end{equation}
Letting $z:=re^{i\theta}$, \eqref{25100801} implies 
\begin{equation*}
\begin{aligned}
\Big|\frac{z(1+z^2)}{z+\overline{z}}\Big|=|w|^2, 
\end{aligned}
\end{equation*}
and, combining this with \eqref{25102605}, it follows that 
\begin{equation*}
\begin{aligned}
\Big|\frac{z(1+z^2)}{z+\overline{z}}\Big|&=\frac{1+r^2}{2}\\
\Longrightarrow
\Big(\frac{\cos \theta+r^2\cos 3\theta}{2\cos \theta}\Big)^2+\Big(\frac{\sin \theta+r^2\sin 3\theta}{2\cos \theta}\Big)^2&=\frac{(1+r^2)^2}{4}\\
\Longrightarrow (r^2-1)^2(1-\cos^2\theta)&=0. 
\end{aligned}
\end{equation*}
Thus either $r=1$ or $\cos \theta=\pm 1$. 

\begin{enumerate}
\item If $\cos \theta=\pm 1$, then $z$ is real, and $w\in \mathbb{R}$ is either real or purely imaginary. This, in turn, implies that $A=kI$, and $B=kI$ or $B^2=kI$ for some $k\in\mathbb{R}$. Consequently, by \eqref{25121701}-\eqref{25121702}, it falls into one of the cases described in \eqref{25102608}. 

\item If $r=1$, then $\overline{z}=\frac{1}{z}$ and $\overline{w}=\frac{1}{w}$, so, by \eqref{25093001},   
\begin{equation*}
\begin{aligned}
\frac{z(1+z^2)}{z+\overline{z}}=\frac{w^3}{\overline{w}}
\Longrightarrow 
z=\pm w^2.
\end{aligned}
\end{equation*}
In this case, we will derive a contradiction by showing that $\alpha+\delta\notin\mathbb{Q}$. We first claim 
\begin{claim}\label{25082901}
If $z=w^2$ (resp. $z=-w^2$), then 
\begin{equation*}
A=(\mathrm{tr}\,B)B-2I\quad 
(\text{resp. }-(\mathrm{tr}\,B)B+2I).
\end{equation*}
\end{claim}
\begin{proof}[Proof of Claim \ref{25082901}]
Without loss of generality, assume $b, \beta>0$ and $z=w^2$; we will only consider this case. Since $z=w^2$, 
\begin{equation}\label{25082902}
\begin{aligned}
\frac{\mathrm{tr}\,A+\sqrt{\mathrm{Disc}\,A}}{2}&=\frac{(\mathrm{tr}\,B)^2-2+(\mathrm{tr}\,B)\sqrt{\mathrm{Disc}\,B}}{2}\\
\Longrightarrow 
\mathrm{tr}\,B>0, \quad \mathrm{tr}\,A&=(\mathrm{tr}\,B)^2-2.
\end{aligned}
\end{equation}
Recall that $\tau$ is a root of $bx^2+(a-d)x-c=\beta x^2+(\alpha-\delta)x-\gamma=0$, and so
\begin{equation}\label{25090102}
\frac{b}{\beta}=\frac{a-d}{\alpha-\beta}=\frac{c}{\gamma}, \quad \frac{d-a+\sqrt{\mathrm{Disc}\,A}}{2b}=\frac{\delta-\alpha+\sqrt{\mathrm{Disc}\,B}}{2\beta}(=\tau).
\end{equation}
Combining the above with \eqref{25082902}, we further obtain 
\begin{equation*}
\frac{b}{\beta}=\frac{a-d}{\alpha-\delta}=\frac{c}{\gamma}=\alpha+\delta(=\mathrm{tr}\,B). 
\end{equation*}
Together with the last equation in \eqref{25082902}, we conclude  
\begin{equation*}
a=\alpha^2+\alpha\delta-2=\alpha (\mathrm{tr}\,B)-2,\quad d=\delta (\mathrm{tr}\,B)-2, 
\end{equation*} 
and hence
\begin{equation*}
A=(\mathrm{tr}\,B)B-2I. 
\end{equation*}
\end{proof}
By the claim, 
\begin{equation*}
\begin{aligned}
\det A(=1)=\det\big((\mathrm{tr}\,B)B-2I\big)\Longrightarrow
3=(\mathrm{tr}\,B)^2.
\end{aligned}
\end{equation*} 
But this contradicts the fact that $\mathrm{tr}\,B\in \mathbb{Q}$ (as $B\in \mathrm{GL}_2(\mathbb{Q})$). 
\end{enumerate}
\end{enumerate}
\end{proof}

As illustrated in Example \ref{25102101}, for a given two-cusped hyperbolic $3$-manifold, if the complex lengths of the two core geodesics of its Dehn filling are linearly dependent, then they are typically equal, i.e., $\frac{a_1}{a_2}$ in \eqref{25110607} is usually $-1$. In the proof of Theorem \ref{25082501}, this condition is equivalent to $\det A = 1$, thanks to \eqref{25110608} in Theorem \ref{22050501}. In fact, most cases in the above proof satisfy this condition, and the only potential exception is the one classified in \eqref{25102609}. However, we doubt the existence of a concrete example falling under \eqref{25102609} with $\det A\neq 1$.

Concerning this, one noteworthy phenomenon again comes from the manifold $v2788$. 
\begin{example}\label{25110702}
{\normalfont Let $\mathcal{M}$ be $v2788$, and $\lambda^i_{p_1/q_1,p_2/q_2}$ ($i=1,2$) be the complex lengths of the core geodesics of $\mathcal{M}_{p_1/q_1,p_2/q_2}$. Then one can observe and verify   
\begin{equation}\label{25110609}
\textnormal{pvol}_{\mathbb{C}}\,\mathcal{M}_{p/q,p/q}=\textnormal{pvol}_{\mathbb{C}}\,\mathcal{M}_{-2q/p, (2p-2q)/(p+2q)}(=\textnormal{pvol}_{\mathbb{C}}\,\mathcal{M}_{(2p-2q)/(p+2q),-2q/p})
\end{equation}
for every coprime pair $(p,q)$ such that $(\frac{2p-2q}{3},\frac{p+2q}{3})$ is also a coprime pair. In particular, the manifold admits a symmetry between two cusps and thus 
\begin{equation}\label{25111005}
\begin{aligned}
\lambda^1_{p/q,p/q}&\equiv\lambda^2_{p/q,p/q}\pmod{\pi \sqrt{-1}\mathbb{Z}}.
\end{aligned}
\end{equation}
We further confirm\footnote{Consequently, it follows that  
$$\lambda^1_{p/q,p/q}+\lambda^2_{p/q,p/q}\equiv
\lambda^1_{-2q/p, (2p-2q)/(p+2q)}+\lambda^2_{-2q/p, (2p-2q)/(p+2q)}\pmod{\pi \sqrt{-1}\mathbb{Z}}$$
which is equivalent to \eqref{25110609}.}
\begin{equation}\label{25110701}
\begin{aligned}
\lambda^1_{-2q/p, (2p-2q)/(p+2q)}&\equiv 3\lambda^2_{-2q/p, (2p-2q)/(p+2q)}
\pmod{\pi \sqrt{-1}\mathbb{Z}},\\
2\lambda^1_{p/q,p/q}&\equiv \lambda^1_{-2q/p, (2p-2q)/(p+2q)}\pmod{\pi \sqrt{-1}\mathbb{Z}}. 
\end{aligned}
\end{equation}

As remarked in Section \ref{intro}, the behavior of the complex volume is governed by the pseudo complex volume. Here, using the volume formula in \cite{nz}, one can similarly show that the complex volumes of $\mathcal{M}_{p/q,p/q}$ and $\mathcal{M}_{-2q/p, (2p-2q)/(p+2q)}$ are also decomposed into linear combinations of two numbers, where each number corresponds to the complex volume change at each cusp, and the ratios between these numbers are exactly matching those appearing in \eqref{25111005}-\eqref{25110701}.}
\end{example}
The above example is interesting in the sense that, to the best of our knowledge, it is the only example admitting Dehn fillings whose core geodesics have linearly dependent complex lengths with a nontrivial ratio; that is, $a_1/a_2$ in \eqref{25110607} is not equal to $-1$. In some respects, this phenomenon is related to the fact that, as mentioned earlier, $v2788$ is the unique manifold we have found so far for which the terms of homogeneous degree $4$ in the potential function vanish. (Consequently, it does not satisfy the assumptions of Theorem \ref{23070505} or \ref{250825011}.)

One of the central mysteries in the study of hyperbolic $3$-manifolds concerns their complex volumes (resp. volumes), particularly the questions of whether these numbers are linearly dependent or independent \cite{thu1}, and why distinct hyperbolic $3$-manifolds can share the same complex volume (resp. volume) \cite{walter}. We hope that Example \ref{25110702}, and more broadly the analysis carried out in this paper, offers new and deeper insight into the structure of the complex volume (resp. volume), especially with regard to these questions. 

A more detailed study, including rigorous proofs of the results in Examples \ref{25101509} and \ref{25110702} as well as several additional perspectives, will appear in \cite{CHDFIII}.\\

\noindent\textbf{Warning. }In Theorem \ref{23070517}, we show that two Dehn fillings of $\mathcal{M}$ satisfying the required conditions are connected via an element of $\mathrm{Aut}\,\mathcal{X}$ of Type I, II, or III. The converse, however, does not necessarily hold: even if two Dehn fillings of $\mathcal{M}$ correspond via such an element, the required assumption need not be satisfied; that is, their core geodesics may all have linearly dependent complex lengths. The manifold $v2788$ is the one that precisely falls into this category.\\

\newpage
\section{Removing $x^4\pm x^3+x^2\pm x+1$}\label{23071105}
From now until Section \ref{24081805}, this part is devoted to classifying the Types I-III elements in $\mathrm{Aut}\,\mathcal{X}$ where $\mathcal{X} \in \mathfrak{Ghol}$ with respect to the minimal polynomial listed in \eqref{23062301}. 

In this section, as an initial step, we prove that not every polynomial in the list appears as the minimal polynomial of an element of $\mathrm{Aut}\,\mathcal{X}$. In particular, we remove the following two from the list:
\begin{equation}\label{25030501}
x^4 +x^3 + x^2 +x + 1, \quad x^4 - x^3 + x^2 - x + 1.
\end{equation}

Recall that the minimal polynomial of $M\in\mathrm{Aut}\,\mathcal{X}$ is one of those in \eqref{25030501}, then it is similar to either
\begin{equation}\label{23070601}
\left(\begin{array}{cccc}
0 & 0 & 0 & -1\\
1 & 0 & 0 & -1\\
0 & 1 & 0 & -1\\
0 & 0 & 1 & -1 
\end{array}\right)\quad \text{or}\quad 
\left(\begin{array}{cccc}
0 & 0 & 0 & -1\\
1 & 0 & 0 & 1\\
0 & 1 & 0 & -1\\
0 & 0 & 1 & 1 
\end{array}\right),
\end{equation}
thanks to the primary decomposition theorem, Theorem \ref{23070515}. 

We now state and prove 
\begin{proposition}\label{23062802}
Let $M$ be an element in $ \mathrm{Aut}\,\mathcal{X}$ where $\mathcal{X}\in \mathfrak{Ghol}$. Suppose $M$ is a matrix of Type I, II or III described in Definition \ref{24060701}. Then $\mathfrak{m}_M(x)$ is neither $x^4+x^3+x^2+x+1$ nor $x^4-x^3+x^2-x+1$. 
\end{proposition}

For the proof, the material introduced in Section \ref{prem} would be enough to establish the result. Nonetheless, as a warm-up before tackling the main argument, it would be instructive to see how the basic ideas from linear algebra can be applied directly.

\begin{proof}
Without loss of generality, suppose $\mathfrak{m}_M(x)=x^4+x^3+x^2+x+1$,\footnote{Since $\mathfrak{m}_{-M}(x)=x^4+x^3+x^2+x+1$, it is enough to consider this case.} hence $M$ is similar to the first matrix, which we denote by $C$, in \eqref{23070601}. That is, there exists an invertible matrix $S$ such that $M=S^{-1}CS$ and  \eqref{22022101} is equivalent to 
\begin{equation}\label{23062901}
S\left(\begin{array}{c}
u'_1\\
v'_1\\
u'_2\\
v'_2
\end{array}\right)
-CS\left(\begin{array}{c}
u_1\\
v_1\\
u_2\\
v_2
\end{array}\right)
=\left(\begin{array}{c}
0\\
0\\
0\\
0
\end{array}\right). 
\end{equation}
Let
\begin{equation}\label{23070301}
S:=\left(\begin{array}{cccc}
a_1 & b_1 & c_1 & d_1\\
a_2 & b_2 & c_2 & d_2\\
a_3 & b_3 & c_3 & d_3\\
a_4 & b_4 & c_4 & d_4 
\end{array}\right),
\end{equation}
and rewrite \eqref{23062901} as 
\begin{equation}\label{23062903}
\left(\begin{array}{cccc}
a_1 & b_1 & c_1 & d_1\\
a_2 & b_2 & c_2 & d_2\\
a_3 & b_3 & c_3 & d_3\\
a_4 & b_4 & c_4 & d_4 
\end{array}\right)
\left(\begin{array}{c}
u'_1\\
v'_1\\
u'_2\\
v'_2
\end{array}\right)
-\left(\begin{array}{cccc}
-a_4 & -b_4 & -c_4 & -d_4 \\
a_1-a_4 & b_1-b_4 & c_1-c_4 & d_1-d_4\\
a_2-a_4 & b_2-b_4 & c_2-c_4 & d_2-d_4\\
a_3-a_4 & b_3-b_4 & c_3-c_4 & d_3-d_4
\end{array}\right)
\left(\begin{array}{c}
u_1\\
v_1\\
u_2\\
v_2
\end{array}\right)=
\left(\begin{array}{c}
0\\
0\\
0\\
0
\end{array}\right). 
\end{equation}
Note that the Jacobian of \eqref{23062903} at $(u'_1, u'_2, u_1, u_2)=(0,0,0,0)$ is 
\begin{equation}\label{23062905}
\left(\begin{array}{cccc}
a_1+b_1\tau_1 & c_1+d_1\tau_2 & a_4+b_4\tau_1 & c_4+d_4\tau_2 \\
a_2+b_2\tau_1 & c_2+d_2\tau_2 & (-a_1+a_4)+(-b_1+b_4)\tau_1 & (-c_1+c_4)+(-d_1+d_4)\tau_2\\
a_3+b_3\tau_1 & c_3+d_3\tau_2 & (-a_2+a_4)+(-b_2+b_4)\tau_1 & (-c_2+c_4)+(-d_2+d_4)\tau_2\\
a_4+b_4\tau_1 & c_4+d_4\tau_2 & (-a_3+a_4)+(-b_3+b_4)\tau_1 & (-c_3+c_4)+(-d_3+d_4)\tau_2
\end{array}\right).
\end{equation}

Let 
\begin{equation*}
\bold{v_j}:=(a_j+b_j\tau_1, c_j+d_j\tau_2)\quad \text{for}\quad 1\leq j\leq 4.
\end{equation*}
 
\begin{enumerate}
\item Suppose $\bold{v_1}$ and $\bold{v_2}$ are linearly independent, and so let
\begin{equation*}
\bold{v_3}=\alpha\bold{v_1}+\beta\bold{v_2},\quad 
\bold{v_4}=\gamma\bold{v_1}+\delta\bold{v_2}
\end{equation*}
for some $\alpha,\beta,\gamma,\delta\in \mathbb{Q}(\tau_1)$. Then \eqref{23062905} is represented as
\begin{equation}\label{23062909}
\left(\begin{array}{cccc}
1 & 0 & \gamma & \delta \\
0 & 1 & \gamma-1 & \delta \\
\alpha & \beta & \gamma & \delta-1 \\
\gamma & \delta & -\alpha+\gamma & -\beta+\delta
\end{array}\right)
\left(\begin{array}{cccc}
a_1+b_1\tau_1 & c_1+d_1\tau_2 & 0 & 0\\
a_2+b_2\tau_1 & c_2+d_2\tau_2 & 0 & 0\\
0 & 0 & a_1+b_1\tau_1 & c_1+d_1\tau_2 \\
0 & 0 & a_2+b_2\tau_1 & c_2+d_2\tau_2 
\end{array}\right).
\end{equation}
Since \eqref{23062903} is defined a complex manifold of dimension $2$, the rank of \eqref{23062909}, and thus the rank of the first matrix in \eqref{23062909} is $2$. Consequently, it follows that 
\begin{equation}\label{25011507}
\alpha(\gamma, \delta)+\beta(\gamma-1, \delta)=(\gamma, \delta-1)\;\;\text{and}\;\; \gamma(\gamma, \delta)+\delta(\gamma-1, \delta)=(-\alpha+\gamma, -\beta+\delta).
\end{equation}
We analyze \eqref{25011507} case by case. 
\begin{enumerate}
\item If $\beta=0$, then the first equation in \eqref{25011507} reduces to $\alpha(\gamma, \delta)=(\gamma, \delta-1)$, implying either $\gamma=0$ or $\alpha=1$. However, if $\gamma\neq 0$ and $\alpha=1$, then $\delta=\delta-1$, which is a contradiction. So $\gamma=0$, and two equations in \eqref{25011507} are $\alpha(0, \delta)=(0, \delta-1)$ and $\delta(-1, \delta)=(-\alpha, \delta)$. However, in this case, one can check there are no $\alpha$ and $\delta$ satisfying the equations.  

\item Otherwise, if $\beta\neq 0$, then $\gamma\neq 0$ and, since $\delta\neq 0$, \eqref{25011507} induces
\begin{equation}\label{25031203}
\begin{gathered}
\alpha+\beta=\frac{\beta+\gamma}{\gamma}=\frac{\delta-1}{\delta}, \quad \gamma+\delta=\frac{-\alpha+\gamma+\delta}{\gamma}=\frac{\delta-\beta}{\delta}\\
\Longrightarrow
\frac{\beta}{\gamma}=\frac{-1}{\delta}, \quad \frac{-\alpha+\delta}{\gamma}=\frac{-\beta}{\delta}
\Longrightarrow
\gamma=-\beta\delta, \quad \alpha=\delta-\beta^2.
\end{gathered}
\end{equation}
Since $\delta(\gamma+\delta)=\delta-\beta$, we get
\begin{equation*}
\delta^2-\beta\delta^2=\delta-\beta\Longrightarrow \beta(\delta^2-1)=\delta^2-\delta, 
\end{equation*}
which implies $\delta=1$ or $\beta=\frac{\delta}{\delta+1}$. If $\delta=1$, then 
\begin{equation*}
\delta(\alpha+\beta)=\delta-1\Longrightarrow \alpha+\beta=0\Longrightarrow \beta^2-\beta-1=0\;(\text{by the last equation in }\eqref{25031203}).
\end{equation*}
But this contradicts the fact that $\beta$ is an element of an imaginary quadratic field $\mathbb{Q}(\tau_1)$. If $\beta=\frac{\delta}{\delta+1}$, then 
\begin{equation*}
\begin{aligned}
\delta(\alpha+\beta)=\delta-1
&\Longrightarrow \delta(\delta-\beta^2+\beta)=\delta(\delta-\frac{\delta^2}{(\delta+1)^2}+\frac{\delta}{\delta+1})=\delta-1\\
&\Longrightarrow \delta^4+\delta^3+\delta^2+\delta+1=0, 
\end{aligned}
\end{equation*}
and we again get the same contradiction. 
\end{enumerate}

\item  Now we suppose $\bold{v_1}$ and $\bold{v_2}$ are linearly dependent. Without loss of generality, if we assume $\bold{v_1}$ and $\bold{v_3}$ are linearly independent. Then 
\begin{equation*}
\bold{v_2}=\alpha\bold{v_1}, \quad \bold{v_4}=\beta\bold{v_1}+\gamma\bold{v_3}
\end{equation*}
for some $\alpha,\beta,\gamma\in \mathbb{C}$, and hence 
\eqref{23062905} is represented as 
\begin{equation}\label{25111805}
\left(\begin{array}{cccc}
1 & 0 & \beta & \gamma \\
\alpha & 0 & \beta-1 & \gamma \\
0 & 1 & -\alpha+\beta & \gamma \\
\beta & \gamma & \beta & \gamma-1 
\end{array}\right)
\left(\begin{array}{cccc}
a_1+b_1\tau_1 & c_1+d_1\tau_2 & 0 & 0\\
a_3+b_3\tau_1 & c_3+d_3\tau_2 & 0 & 0\\
0 & 0 & a_1+b_1\tau_1 & c_1+d_1\tau_2 \\
0 & 0 & a_3+b_3\tau_1 & c_3+d_3\tau_2 
\end{array}\right).
\end{equation}
Seen as above, if \eqref{23062901} defines a complex manifold of dimension $2$, then the rank of the first matrix in \eqref{25111805} is $2$, which implies
\begin{equation*}
\begin{gathered}
\alpha(\beta, \gamma)=(\beta-1, \gamma)\quad\text{and}\quad \beta(\beta, \gamma)+\gamma(-\alpha+\beta, \gamma)=(\beta, \gamma-1).
\end{gathered}
\end{equation*}
But one can easily check that there are no $\alpha,\beta$ and $\gamma$ satisfying the above equations. 
\end{enumerate}
\end{proof}

\newpage
\section{Strategy for classifying $M$: I}\label{25012001}
In this section and the next, we outline the strategy for classifying elements of Types I-III in $\mathrm{Aut}\,\mathcal{X}$ where $\mathcal{X} \in \mathfrak{Ghol}$. As an initial step, we derive a necessary condition for $M\in \mathrm{GL}_4(\mathbb{Q})$ to be an element of $\mathrm{Aut}\,\mathcal{X}$. More precisely, by calculating the second-lowest degree terms in the equations appearing in \eqref{22022101}, we deduce a system of two homogeneous equations in $u_1$ and $u_2$, involving a $2$ by $2$ matrix transformation directly arising from $M$. Then, in Section \ref{25011004}, we diagonalize the matrix to transform the equation into a much simpler form, easier to handle. Finally, in Section \ref{25011802}, several methods for resolving the equation are proposed. In particular, we establish criteria for the matrix under which a solution to the system of equations exists. 

\subsection{Formulation}\label{25011801}

The following proposition is the key starting point of the strategy. The first claim of the proposition is simply a restatement of the equalities in \eqref{24080403}, which is attained by considering the linear terms in \eqref{22022101}. The second claim, in which we set up the central equation, is derived by comparing the terms of the second-lowest degree in \eqref{22022101}.

\begin{proposition}\label{24020501}
Let $\mathcal{X}$ be an element in $\mathfrak{Ghol}$ parametrized by $v_i$ ($i=1,2$) as in \eqref{23120805}. Let $M\in \mathrm{Aut}\,\mathcal{X}$ be given by 
\begin{equation}\label{25010811}
\left(\begin{array}{cc}
A_1 & A_2\\
A_3 & A_4
\end{array}\right)
\Bigg(=\left(\begin{array}{cccc}
a_1 & b_1 & a_2 & b_2\\
c_1 & d_1 & c_2 & d_2\\
a_3 & b_3 & a_4 & b_4\\
c_3 & d_3 & c_4 & d_4 
\end{array}\right)\Bigg)\in \mathrm{GL}_4(\mathbb{Q}).
\end{equation}
Then the following statements hold. 
\begin{enumerate}
\item $\tau_1(a_j+b_j\tau_1)=c_j+d_j\tau_1$ for $1\leq j\leq 2$ and $\tau_2(a_j+b_j\tau_2)=c_j+d_j\tau_2$ for $3\leq j\leq 4$. 
\item If $\Theta_i(u_1, u_2)$ be the homogeneous polynomial of the second smallest degree of $v_i(u_1, u_2)$ ($i=1,2$), then  
\begin{equation}\label{24032305}
\left(\begin{array}{c}
\Theta_1\Big(\left(\begin{array}{cc}
\omega_1 & \omega_2\\
\omega_3 & \omega_4
\end{array}
\right)\bold{U}\Big) \\
\Theta_2\Big(\left(\begin{array}{cc}
\omega_1 & \omega_2\\
\omega_3 & \omega_4
\end{array}
\right)\bold{U}\Big)
\end{array}
\right)=\left(\begin{array}{cc}
\overline{\omega_1} & \overline{\omega_2}\\
\overline{\omega_3} & \overline{\omega_4}
\end{array}\right)\left(\begin{array}{c}
\Theta_1(\bold{U}) \\
\Theta_2(\bold{U})
\end{array}
\right)
\end{equation}
where $\bold{U}:=\left(\begin{array}{c}
u_1\\
u_2
\end{array}\right)$, $\omega_j$ ($1\leq j\leq 4$) is an eigenvalue of $A_j$ such that $\mathrm{Im}\,\omega_j\geq 0$ (resp. $\mathrm{Im}\,\omega_j< 0$) for $b_j\geq0$ (resp. $b_j<0$), and $\overline{\omega_j}$ is the conjugate of $\omega_j$. 
\end{enumerate}
\end{proposition}
\begin{proof}
By the assumption, since \eqref{22022101} is a $2$-dimensional analytic subset of $\mathcal{X}\times \mathcal{X}$, the following two identities 
\begin{equation}\label{23100401}
\begin{aligned}
&v'_1\big(a_1u_1+b_1v_1(u_1, u_2)+a_2u_1+b_2v_1(u_1, u_2), a_3u_1+b_3v_1(u_1, u_2)+a_4u_2+b_4v_2(u_1, u_2)\big)\\
=&c_1u_1+d_1v_1(u_1, u_2) +c_2u_2+d_2v_2(u_2, u_2)
\end{aligned}
\end{equation}
and 
\begin{equation}\label{25010911}
\begin{aligned}
&v'_2\big(a_1u_1+b_1v_1(u_1, u_2)+a_2u_1+b_2v_1(u_1, u_2), a_3u_1+b_3v_1(u_1, u_2)+a_4u_2+b_4v_2(u_1, u_2)\big)\\
=&c_3u_1+d_3v_1(u_1, u_2) +c_4u_2+d_4v_2(u_2, u_2)
\end{aligned}
\end{equation}
hold for any $u_1$ and $u_2$.

The first claim follows immediately from the coefficients of the terms of degree $1$ in \eqref{23100401}-\eqref{25010911}. 

To prove the second claim, we compare the coefficients of the terms of the second-smallest degree, and get 
\begin{equation}\label{24032303}
\begin{aligned} 
&\Theta_1\big((a_1+b_1\tau_1) u_1+(a_2+b_2\tau_2)u_2, (a_3+d_3\tau_1)u_1+(a_4+b_4\tau_2)u_2\big)\\
=&(d_1-b_1\tau_1)\Theta_1(u_1, u_2)+(d_2-b_2\tau_2)\Theta_2(u_1, u_2),  
\end{aligned}
\end{equation}
and
\begin{equation}\label{24032304}
\begin{aligned} 
&\Theta_2\big((a_1+b_1\tau_1) u_1+(a_2+b_2\tau_2)u_2, (a_3+d_3\tau_1)u_1+(a_4+d_4\tau_2)u_2\big)\\
=&(d_3-b_3\tau_1)\Theta_1(u_1, u_2)+(d_4-b_4\tau_2)\Theta_2(u_1, u_2). 
\end{aligned}
\end{equation}
By letting $\omega_j:=a_j+b_j\tau_1$ ($j=1,2$) and $\omega_j:=a_j+b_j\tau_2$ ($j=3,4$), the equations \eqref{24032303}-\eqref{24032304} are combined and summarized as in \eqref{24032305}. The fact that $\omega_j$ is an eigenvalue of $A_j$ with the required conditions follows from the first claim and the observations in \eqref{25111001}-\eqref{25111002}. (If $b_j=0$, it is clear from the first claim that $c_j=0, \omega_j=a_j=d_j$,  and $A_j=a_j I$.)
\end{proof}

Summing up, for a given $M$, we obtain two matrices 
\begin{equation}\label{25010915}
\left(\begin{array}{cc}
\omega_1 & \omega_2\\
\omega_3 & \omega_4
\end{array}
\right)\quad \text{and}\quad  
\left(\begin{array}{cc}
\overline{\omega_1} & \overline{\omega_2}\\
\overline{\omega_3} & \overline{\omega_4}
\end{array}\right)
\end{equation}
associated with it, and, if $M$ is realized as an element of $\mathrm{Aut}\,\mathcal{X}$ for some $\mathcal{X}\in \mathfrak{Ghol}$, then there exist homogeneous polynomials $\Theta_i$ ($i=1,2$) satisfying \eqref{24032305}. Said differently, the existence of $\Theta_i$ meeting \eqref{24032305} is a necessary condition for the existence of $\mathcal{X}\in \mathfrak{Ghol}$ such that $M\in \mathrm{Aut}\,\mathcal{X}$. Hence, the problem is reduced to solving the equation for $\Theta_i$ in \eqref{24032305}, given the data \eqref{25010915}. 


The following term arises from the above proposition and will be frequently used throughout the paper.
\begin{definition}\label{25120501}
\normalfont{We call the first matrix in \eqref{25010915} the \textit{primary matrix associated} to $M$ and denote it by $P_M$. Similarly, we call the second one in \eqref{25010915} the \textit{conjugate matrix} of $P_M$ and denote it by $\overline{P_M}$. More generally, for 
\begin{equation}\label{24090601}
\left(\begin{array}{cc}
a_1+b_1\sqrt{D} & a_2+b_2\sqrt{D}\\
a_3+b_3\sqrt{D} & a_4+b_4\sqrt{D}
\end{array}\right)\in \mathrm{GL}_2(\mathbb{Q}(\sqrt{D}) )
\end{equation}
where $D$ is a negative integer, we call 
\begin{equation*}
\left(\begin{array}{cc}
a_1-b_1\sqrt{D} & a_2-b_2\sqrt{D}\\
a_3-b_3\sqrt{D} & a_4-b_4\sqrt{D}
\end{array}\right) 
\end{equation*}
the \textit{conjugate matrix} of \eqref{24090601}. }
\end{definition}

As mentioned previously, we will often encounter situations where, when $M$ is given as in \eqref{25010811}, a minimal polynomial condition on $M$ forces each $A_j$ ($2\leq j\leq 4$) to be represented in terms of $A_1$. In such cases, the following lemma provides an intimate relationship between the arithmetic of $M$ and its primary matrix. 

\begin{lemma}\label{250111011}
Let $\mathcal{X}$ be an element in $\mathfrak{Ghol}$ and $M\in \mathrm{GL}_4(\mathbb{Q})$ be an element in $\mathrm{Aut}\,\mathcal{X}$ given as in \eqref{25010811}. Suppose each $A_j$ ($2\leq j\leq 4$) is expressed as a linear function of $A_1$. More precisely, 
\begin{equation*}
M=\left(\begin{array}{cc}
A_1 & a_2 I+b_2A_1 \\
a_3 I+b_3A_1 & a_4 I+b_4A_1
\end{array}\right)
\end{equation*}
for some $a_j, b_j\in \mathbb{Q}$ for $2\leq j\leq 4$. Then the primary matrix associated with $M$ has a finite order. 
\end{lemma}
\begin{proof}
By Proposition \ref{24020501}, the primary matrix associated with $M$ is given as
\begin{equation}\label{25011201}
\left(\begin{array}{cc}
\omega_1 & a_2+b_2\omega_1 \\
a_3+b_3\omega_1 & a_4+b_4\omega_1
\end{array}\right)
\end{equation}
where $\omega_1$ is an eigenvalue of $A_1$. Note that, for any polynomial $f$ and $g$, $f(A_1)=g(A_1)$ if and only if $f(\omega_1)=g(\omega_1)$. Thus it follows that $n$-th power of \eqref{25011201} is the primary matrix associated with $M^n$. As the order of $M$ is finite, the conclusion follows.  
\end{proof}

If \eqref{25011201} is an element of finite order, then the two eigenvalues of \eqref{25011201} are roots of unity. Since the degree of $\omega_1$ is at most $2$,  there are only finitely many possible choices for these eigenvalues. By computing their exact values, one can explicitly determine the possible forms of $M$.\footnote{Note that once we find the exact value of $\omega_i$, the explicit form of $A_i$ is determined by Proposition \ref{24020501}(1) and Lemma \ref{24090301}.} We will further generalize the above lemma in Lemma \ref{25011101} and often invoke both in Sections \ref{24081804}-\ref{24081805}. 

The following corollary of Proposition \ref{24020501} appears to share a similar flavor with a statement from Theorem \ref{23070517}. But, the corollary removes the assumption reparding the existence of a Dehn filling point associated with a pair of Dehn fillings having the same pseudo-complex volume. 

\begin{corollary}\label{25010913}
Let $\mathcal{X}$ be an element in $\mathfrak{Ghol}$ and $M\in \mathrm{GL}_4(\mathbb{Q})$ be an element in $\mathrm{Aut}\,\mathcal{X}$ given as in \eqref{25010811}. Then one of the following is true:
\begin{enumerate}
\item $A_1=A_4=0$;

\item $A_2=A_3=0$;

\item $\det A_j> 0$ for all $1\leq j\leq 4$.
\end{enumerate}
\end{corollary}
\begin{proof}
If $\det A_1=0$, then $A_1=0$ by the first claim in Proposition \ref{24020501}. By Claim 7.5 in \cite{CHDF}, if $A_1=0$, then $\det A_4=0$ and thus $A_4=0$. 

By symmetry, if either $\det A_2=0$ or $\det A_3=0$, we similarly obtain $A_2=A_3=0$. 

If $\det A_j\neq 0$, then $\det A_j> 0$ follows from the fact that the two cusp shapes of $\mathcal{X}$ are non-real, contained in $\mathbb{Q}(\sqrt{\mathrm{Disc}\,A_j})$.
\end{proof}

\subsection{Transformation}\label{25011004}

To address \eqref{24032305} in general, we first transform the equation into an equivalent form, which is much easier to handle with, by adopting some elementary theory of linear algebra. 

Let $\zeta_i$ ($i=1,2$) be the eigenvalues of $\overline{P_M}$, that is, the roots of 
\begin{equation*}
x^2-(\overline{\omega_1}+\overline{\omega_4})x+\big(\overline{\omega_1}\overline{\omega_4}-\overline{\omega_2}\overline{\omega_3}\big)=0. 
\end{equation*}
If $\zeta_1\neq\zeta_2$ and $\omega_2\neq 0$, then $\overline{P_M}$ is factored as 
\begin{equation*}
\begin{gathered}
\left(\begin{array}{cc}
\overline{\omega_2} & \overline{\omega_2}\\
\zeta_1-\overline{\omega_1} & \zeta_2-\overline{\omega_1}
\end{array}\right)
\left(\begin{array}{cc}
\zeta_1 & 0\\
0 & \zeta_2
\end{array}\right)
\left(\begin{array}{cc}
\overline{\omega_2} & \overline{\omega_2}\\
\zeta_1-\overline{\omega_1} & \zeta_2-\overline{\omega_1}
\end{array}\right)^{-1}.
\end{gathered}
\end{equation*}
Thus, if we set 
\begin{equation*}
\left(\begin{array}{c}
g_1(\bold{U})\\
g_2(\bold{U})
\end{array}\right)
:=\left(\begin{array}{cc}
\overline{\omega_2} & \overline{\omega_2}\\
\zeta_1-\overline{\omega_1} & \zeta_2-\overline{\omega_1}
\end{array}\right)^{-1}
\left(\begin{array}{c}
\Theta_1(\bold{U})\\
\Theta_2(\bold{U})
\end{array}\right)
\end{equation*} 
where $\bold{U}:=\left(\begin{array}{c}
u_1\\
u_2
\end{array}\right)$, then \eqref{24032305} is turned into
\begin{equation}\label{24032501}
\begin{gathered}
\left(\begin{array}{c}
g_1(P_M\bold{U})\\
g_2(P_M\bold{U})
\end{array}\right)
=\left(\begin{array}{cc}
\zeta_1 & 0\\
0 & \zeta_2
\end{array}\right)
\left(\begin{array}{c}
g_1(\bold{U})\\
g_2(\bold{U})
\end{array}\right).
\end{gathered}
\end{equation}
To solve the equation, we factor $P_M$ as 
\begin{equation*}
\begin{gathered}
\left(\begin{array}{cc}
\omega_2 & \omega_2\\
\lambda_1-\omega_1 & \lambda_2-\omega_1
\end{array}\right)
\left(\begin{array}{cc}
\lambda_1 & 0\\
0 & \lambda_2
\end{array}\right)
\left(\begin{array}{cc}
\omega_2 & \omega_2\\
\lambda_1-\omega_1 & \lambda_2-\omega_1
\end{array}\right)^{-1}
\end{gathered}
\end{equation*}
and set a new basis as
\begin{equation}\label{24050602}
\begin{gathered}
\left(\begin{array}{c}
\tilde{u}_1\\
\tilde{u}_2
\end{array}\right)
:=\left(\begin{array}{cc}
\omega_2 & \omega_2\\
\lambda_1-\omega_1 & \lambda_2-\omega_1
\end{array}\right)^{-1}
\left(\begin{array}{c}
u_1\\
u_2
\end{array}\right).
\end{gathered}
\end{equation}
If $\tilde{g}_i (\tilde{u}_1, \tilde{u}_2)$ ($i=1,2$) are defined as
\begin{equation*}
g_i\big(\omega_2\tilde{u}_1+\omega_2\tilde{u}_2, (\lambda_1-\omega_1)\tilde{u}_1+(\lambda_2-\omega_1)\tilde{u}_2\big), 
\end{equation*}
then \eqref{24032501} (and so \eqref{24032305} as well) is now equivalent to 
\begin{equation}\label{24050601}
\begin{gathered}
\left(\begin{array}{c}
\tilde{g_1}(\lambda_1\tilde{u}_1, \lambda_2\tilde{u}_2)\\
\tilde{g_2}(\lambda_1\tilde{u}_1, \lambda_2\tilde{u}_2)
\end{array}\right)
=\left(\begin{array}{cc}
\zeta_1 & 0\\
0 & \zeta_2
\end{array}\right)
\left(\begin{array}{c}
\tilde{g_1}(\tilde{u}_1, \tilde{u}_2)\\
\tilde{g_2}(\tilde{u}_1, \tilde{u}_2)
\end{array}\right).
\end{gathered}
\end{equation}
In the next section, particularly in Lemmas \ref{23112902} and \ref{24010902}, we will explore the general solutions satisfying \eqref{24050601}. 

Recall that one of the important symmetric properties that $\Theta_i$ ($i=1,2$) possesses is the following (by Definition \ref{23120604}):
\begin{equation}\label{23120811}
\Theta_1(u_1, u_2)=\Theta_1(u_1, -u_2)\quad\text{and}\quad \Theta_2(u_1, u_2)=-\Theta_2(u_1, -u_2).
\end{equation}
Thus if $\iota:=\left(\begin{array}{cc}
1 & 0\\
0 & -1
\end{array}\right)$, then 
\begin{equation}\label{23120807}
\left(\begin{array}{c}
\Theta_1(P_M\iota\bold{U})\\
\Theta_2(P_M\iota\bold{U})
\end{array}\right)
=\overline{P_M}
\left(\begin{array}{c}
\Theta_1(\iota\bold{U})\\
\Theta_2(\iota\bold{U})
\end{array}\right)
=\overline{P_M}\iota
\left(\begin{array}{c}
\Theta_1(\bold{U})\\
\Theta_2(\bold{U})
\end{array}\right).
\end{equation}
Note that $\overline{P_M}\iota$ is the conjugate matrix of $P_M\iota$. The fact that $\Theta_i$ ($i=1,2$) satisfy both \eqref{24032305} and \eqref{23120807} places a strong restriction on the structure of $P_M$. For instance, if $\frac{\lambda_1}{\lambda_2}$ is not a root of unity, it will be shown that either $P_M$ commutes with $P_M\iota$ or $\mathrm{tr}\,(P_M \iota)=0$ (see Proposition \ref{24072003}). This observation, together with its further consequences (e.g., Proposition \ref{24082301}), will simplify many proofs by significantly reducing the amount of computation required in Sections \ref{24081804}-\ref{25012005}.

More generally, let $\mathcal{G}$ be defined as the group generated by $P_M$ and $\iota$. Then
\begin{equation*}
\left(\begin{array}{c}
\Theta_1(P\bold{U})\\
\Theta_2(P\bold{U})
\end{array}\right)
=\overline{P}
\left(\begin{array}{c}
\Theta_1(\bold{U})\\
\Theta_2(\bold{U})
\end{array}\right)
\end{equation*}
holds for any $P\in \mathcal{G}$ and its conjugate matrix $\overline{P}$. Sometimes, we prefer to work with elements in $\mathcal{G}$ other than $P_M$, especially when their eigenvalues and eigenvectors are easier to handle compared to those of $P_M$. 

In this respect, Lemma \ref{250111011} is generalized as follows:

\begin{lemma}\label{25011101}
Let $\mathcal{X}$ be an element in $\mathfrak{Ghol}$, $M\in \mathrm{GL}_4(\mathbb{Q})$ be an element in $\mathrm{Aut}\,\mathcal{X}$ given as in \eqref{25010811}, and $P_M$ be a primary matrix associated with $M$. Suppose each $A_j$ ($2\leq j\leq 4$) is expressed as a function of $A_1$, then any element in $\langle P_M, \iota \rangle $ has a finite order. 
\end{lemma}

This lemma will be particularly useful in completing the proofs in Section \ref{24081804}. 

We also remark that, by \eqref{24050602}, if we set 
 \begin{equation*}
 \tilde{\Theta}_i(\tilde{u}_1, \tilde{u}_2):=\Theta_i\big(\omega_2\tilde{u}_1+\omega_2\tilde{u}_2, (\lambda_1-\omega_1)\tilde{u}_1+(\lambda_2-\omega_1)\tilde{u}_2\big)\quad (i=1,2),
 \end{equation*}
then $a\frac{\partial \Theta_1(u_1, u_2)}{\partial u_2}=\frac{\partial \Theta_2(u_1, u_2)}{\partial u_1}$ is equivalent to   
\begin{equation*}
\begin{aligned}
a\bigg(\frac{\partial \tilde{u}_1}{\partial u_2}\frac{\partial \tilde{\Theta}_1}{\partial \tilde{u}_1}+\frac{\partial \tilde{u}_2}{\partial u_2}\frac{\partial \tilde{\Theta}_1}{\partial \tilde{u}_2}\bigg)=&\frac{\partial \tilde{u}_1}{\partial u_1}\frac{\partial \tilde{\Theta}_2}{\partial \tilde{u}_1}+\frac{\partial \tilde{u}_2}{\partial u_1}\frac{\partial \tilde{\Theta}_2}{\partial \tilde{u}_2}\\
\Longrightarrow 
a\bigg(-\omega_2\frac{\partial \tilde{\Theta}_1}{\partial \tilde{u}_1}+\omega_2\frac{\partial \tilde{\Theta}_1}{\partial \tilde{u}_2}\bigg)=&(\lambda_2-\omega_1)\frac{\partial \tilde{\Theta}_2}{\partial \tilde{u}_1}-(\lambda_1-\omega_1)\frac{\partial \tilde{\Theta}_2}{\partial \tilde{u}_2}
\end{aligned}
\end{equation*}
by \eqref{24050602}. Since $\left(\begin{array}{cc}
\overline{\omega_2} & \overline{\omega_2}\\
\zeta_1-\overline{\omega_1} & \zeta_2-\overline{\omega_1}
\end{array}\right)\left(\begin{array}{cc}
\tilde{g}_1\\
\tilde{g}_2
\end{array}\right)=
\left(\begin{array}{cc}
\tilde{\Theta}_1\\
\tilde{\Theta}_2
\end{array}\right)$, it is further reduced to 
\begin{equation*}
\begin{aligned}
&a\bigg(-\omega_2\overline{\omega_2}\frac{\partial \tilde{g}_1}{\partial \tilde{u}_1}-\omega_2\overline{\omega_2}\frac{\partial \tilde{g}_2}{\partial \tilde{u}_1}+\omega_2\overline{\omega_2}\frac{\partial \tilde{g}_1}{\partial \tilde{u}_2}+\omega_2\overline{\omega_2}\frac{\partial \tilde{g}_2}{\partial \tilde{u}_2}\bigg)\\
=&(\lambda_2-\omega_1)\Big((\zeta_1-\overline{\omega_1})\frac{\partial \tilde{g}_1}{\partial \tilde{u}_1}+(\zeta_2-\overline{\omega_1})\frac{\partial \tilde{g}_2}{\partial \tilde{u}_1}\Big)-(\lambda_1-\omega_1)\Big((\zeta_1-\overline{\omega_1})\frac{\partial \tilde{g}_1}{\partial \tilde{u}_2}+(\zeta_2-\overline{\omega_1})\frac{\partial \tilde{g}_2}{\partial \tilde{u}_2}\Big), 
\end{aligned}
\end{equation*}
and thus 
\begin{equation}\label{24050605}
\begin{aligned}
&\big((\lambda_1-\omega_1)(\zeta_1-\overline{\omega_1})+a\omega_2\overline{\omega_2}\big)\frac{\partial \tilde{g}_1}{\partial \tilde{u}_2}+\big((\lambda_1-\omega_1)(\zeta_2-\overline{\omega_1})+a\omega_2\overline{\omega_2}\big)\frac{\partial \tilde{g}_2}{\partial \tilde{u}_2}\\
=&\big((\lambda_2-\omega_1)(\zeta_1-\overline{\omega_1})+a\omega_2\overline{\omega_2}\big)\frac{\partial \tilde{g}_1}{\partial \tilde{u}_1}+\big((\lambda_2-\omega_1)(\zeta_2-\overline{\omega_1})+a\omega_2\overline{\omega_2}\big)\frac{\partial \tilde{g}_2}{\partial \tilde{u}_1}.
\end{aligned}
\end{equation}
We summarize what have been discussed so far in the following:
\begin{proposition}\label{24051906}
Let $\Theta_i$ ($i=1,2$) be homogeneous polynomials of the same degree $\geq 3$ satisfying 
\begin{equation}\label{24051905}
\left(\begin{array}{c}
\Theta_1(P\bold{U}) \\
\Theta_2(P\bold{U})
\end{array}
\right)=\overline{P}\left(\begin{array}{c}
\Theta_1(\bold{U}) \\
\Theta_2(\bold{U})
\end{array}
\right)
\end{equation}
where $\bold{U}:=\left(\begin{array}{c}
u_1\\
u_2
\end{array}\right)$ and $P:=\left(\begin{array}{cc}
\omega_1 & \omega_2(\neq 0)\\
\omega_3 & \omega_4
\end{array}\right)\in \mathrm{GL}_2(\mathbb{Q}(\sqrt{D}))$ for some negative integer $D$. Let $\lambda_i$ (resp. $\zeta_i$) ($i=1,2$) be the eigenvalues of $P$ (resp. $\overline{P}$) such that $\lambda_1\neq \lambda_2$. Let 
\begin{equation}\label{24082302}
\bold{\tilde{U}}:=\left(\begin{array}{c}
\tilde{u}_1\\
\tilde{u}_2
\end{array}\right):=S_{\lambda}^{-1}\bold{U}, \;\; 
\left(\begin{array}{c}
g_1(\bold{U})\\
g_2(\bold{U})
\end{array}\right)
:=S_{\zeta}^{-1}
\left(\begin{array}{c}
\Theta_1(\bold{U})\\
\Theta_2(\bold{U})
\end{array}\right) \;\;\text{and}\;\; \left(\begin{array}{c}
\tilde{g_1}(\tilde{\bold{U}})\\
\tilde{g_2}(\tilde{\bold{U}})
\end{array}\right)
:=\left(\begin{array}{c}
g_1(S_{\lambda}\tilde{\bold{U}})\\
g_2(S_{\lambda}\tilde{\bold{U}})
\end{array}\right)
\end{equation} 
where 
\begin{equation*}
S_{\lambda}:=\left(\begin{array}{cc}
\omega_2 & \omega_2\\
\lambda_1-\omega_1 & \lambda_2-\omega_1
\end{array}\right)\;\;\text{and}\;\;  S_{\zeta}:=\left(\begin{array}{cc}
\overline{\omega_2} & \overline{\omega_2}\\
\zeta_1-\overline{\omega_1} & \zeta_2-\overline{\omega_1}
\end{array}\right).
\end{equation*}
Then \eqref{24051905} is equivalent to 
\begin{equation}\label{24091301}
\begin{gathered}
\left(\begin{array}{c}
\tilde{g_1}(\lambda_1\tilde{u}_1, \lambda_2\tilde{u}_2)\\
\tilde{g_2}(\lambda_1\tilde{u}_1, \lambda_2\tilde{u}_2)
\end{array}\right)
=\left(\begin{array}{cc}
\zeta_1 & 0\\
0 & \zeta_2
\end{array}\right)
\left(\begin{array}{c}
\tilde{g_1}(\tilde{u}_1, \tilde{u}_2)\\
\tilde{g_2}(\tilde{u}_1, \tilde{u}_2)
\end{array}\right).
\end{gathered}
\end{equation}
Further, $a\frac{\partial \Theta_1(u_1, u_2)}{\partial u_2}=\frac{\partial \Theta_2(u_1, u_2)}{\partial u_1}$ corresponds to \eqref{24050605}. 
\end{proposition}
We also note that, under the transformations given in \eqref{24082302}, \eqref{23120807} is turned into 
\begin{equation}\label{24091601}
\left(\begin{array}{c}
\tilde{g}_1(S_{\lambda}^{-1}P\iota S_{\lambda}\bold{\tilde{U}}) \\
\tilde{g}_2(S_{\lambda}^{-1}P\iota S_{\lambda}\bold{\tilde{U}})
\end{array}
\right)=S_{\zeta}^{-1}\overline{P}\iota S_{\zeta}
\left(\begin{array}{c}
\tilde{g}_1(\bold{\tilde{U}}) \\
\tilde{g}_2(\bold{\tilde{U}})
\end{array}\right). 
\end{equation}

In conclusion, combining the above with Proposition \ref{24020501}, for a given $M\in \mathrm{GL}_2(\mathbb{Q})$ and $P:=P_M$, the existence of $\tilde{g}_i$ ($i=1,2$) satisfying both \eqref{24091301} and \eqref{24091601} is a necessary condition for $M$ to be an element of $\mathrm{Aut}\,\mathcal{X}$. 

In the next section, we closely examine the equation in \eqref{24091301} and determine the explicit forms of $\tilde{g}_i$ ($i=1,2$) that solve it. 


\begin{convention}
{\normalfont For simplicity, let us denote both 
\begin{equation*}
\left(\begin{array}{cc}
1 & 0 \\
0 & -1 
\end{array}\right)\quad \text{and}\quad  
\left(\begin{array}{cc}
I & 0 \\
0 & -I 
\end{array}\right)
\Bigg(=\left(\begin{array}{cccc}
1 & 0 & 0 & 0 \\
0 & 1 & 0 & 0 \\
0 & 0 & -1 & 0 \\
0 & 0 & 0 & -1 
\end{array}\right)\Bigg)
\end{equation*}
by $\iota$. }
\end{convention}
Although the notation above may seem ambiguous, its meaning will be clear from context.

\section{Strategy for classifying $M$: II}\label{25011802}
In this section, we solve the equation formulated in Proposition \ref{24051906} and provide useful criteria that will be applied subsequently in Sections \ref{24081804}-\ref{24081805}. We also introduce a proposition that ensures the compatibility of two distinct elements of $\mathrm{Aut}\,\mathcal{X}$. This proposition will serve as a key tool, much like a Rosetta stone, in proving the main argument.

\subsection{Solution}
In the next two lemmas, we describe $\tilde{g_i}$ satisfying \eqref{24091301}, and find the relations that $\lambda_i$ and $\zeta_i$ ($i=1,2$) must meet. 

Depending on whether $\frac{\lambda_1}{\lambda_2}$ is a root of unity or not, the general solutions will exhibit different behaviors. For instance, if $\frac{\lambda_1}{\lambda_2}$ is not a root of unity, then $\tilde{g_i}$ in \eqref{24091301} is relatively simple to describe as will be shown in the following lemma:
\begin{lemma}\label{23112902}
Let $g(u_1, u_2)$ be a homogeneous polynomial satisfying
\begin{equation}\label{23120607}
 g(\lambda_1u_1, \lambda_2u_2)=\alpha g(u_1, u_2)
\end{equation} 
for some $\alpha, \lambda_i \in \mathbb{C}\setminus \{0\}$ ($i=1,2$). If $\frac{\lambda_2}{\lambda_1}$ is not a root of unity, then $g$ is a constant multiple of $u_1^ku_2^l$ for some $k,\l\in \mathbb{N}\cup \{0\}$ and $\alpha=\lambda_1^k\lambda_2^l$. 
\end{lemma}
\begin{proof}
Without loss of generality, we assume $g$ is factored as
\begin{equation*}
g(u_1, u_2)=\beta u_1^{k}u_2^{l}\prod_{j=1}^n(u_1+\beta_j u_2)
\end{equation*}
for some $k,\l\in \mathbb{N}\cup \{0\}$ and $\beta,\beta_j (\neq 0)\in \mathbb{C}$ ($1\leq j\leq n$). Then 
\begin{equation*}
g(\lambda_1 u_1, \lambda_2 u_2)=\beta \lambda_1^{k+n}\lambda_2^{l} u_1^{k}u_2^{l}\prod_{j=1}^n\Big(u_1+\frac{\lambda_2}{\lambda_1} \beta_ju_2\Big), 
\end{equation*}
and thus
\begin{equation*}
\prod_{j=1}^n(u_1+\beta_ju_2)=\prod_{j=1}^n\Big(u_1+\frac{\lambda_2}{\lambda_1}\beta_ju_2\Big), 
\end{equation*}
thanks to \eqref{23120607}. That is, 
\begin{equation}\label{23120608}
\{\beta_1, \dots, \beta_n\}=\Big\{\frac{\lambda_2}{\lambda_1}\beta_1, \dots, \frac{\lambda_2}{\lambda_1}\beta_n\Big\}. 
\end{equation}
Without loss of generality, let us suppose $\beta_j= \Big(\frac{\lambda_2}{\lambda_1}\Big)^{j-1}\beta_1$. As $\frac{\lambda_2}{\lambda_1}$ is not a root of unity by the assumption, there is no $j$ ($2\leq j\leq n$) such that $\beta_n=\beta_1$, contradicting the equality \eqref{23120608}. Hence the conclusion follows. 
\end{proof}

If $\frac{\lambda_1}{\lambda_2}$ is a root of unity, then $\tilde{g}_i$ in \eqref{24091301} can be described by the following lemma:

\begin{lemma}\label{24010902}
Let $g(u_1, u_2)$ be a homogeneous polynomial satisfying
\begin{equation*}
 g(\lambda_1u_1, \lambda_2u_2)=\alpha g(u_1, u_2)
\end{equation*} 
for some $\alpha, \lambda_i \in \mathbb{C}\setminus \{0\}$ ($i=1,2$). If $\frac{\lambda_1}{\lambda_2}$ is a root of unity of order $d$, then $g$ is a constant multiple of a polynomial of the following form:
\begin{equation}\label{24091501}
u_1^ku_2^l\prod_{i=1}^{n}(u_1^d+\alpha_iu_2^d)
\end{equation}
for some $k,\l\in \mathbb{N}\cup \{0\}$ and $\alpha_i\in \mathbb{C}$ ($1\leq i\leq n$), and $\alpha=\lambda_1^{k+dn}\lambda_2^l$.
\end{lemma}
\begin{proof}
Let $h$ be a polynomial such that 
\begin{itemize}
\item $h$ is a factor of $g$;
\item $h$ is neither divisible by $u_1$ nor $u_2$;
\item $h(\lambda_1 u_1, \lambda_2 u_2)=\beta h(u_1, u_2)$ for some constant $\beta\in \mathbb{C}$;
\item $h$ is a polynomial of the smallest degree satisfying the above three.
\end{itemize} 
If $h$ is defined by $\prod_{i=1}^{m}(u_1+\beta_iu_2)$, then
\begin{equation*}
\begin{gathered}
h(\lambda_1u_1, \lambda_2u_2)=\prod_{i=1}^{m}(\lambda_1u_1+\beta_i\lambda_2u_2)=\lambda_1^m\prod_{i=1}^{m}\Big(u_1+\beta_i\frac{\lambda_2}{\lambda_1}u_2\Big)
\end{gathered}
\end{equation*}
and, as it is equal to $\beta h(u_1, u_2)$, it follows that $\beta=\lambda_1^m$ and 
\begin{equation*}
\{\beta_i\;|\; 1\leq i\leq m\}=\Big\{\beta_i\frac{\lambda_2}{\lambda_1}\;|\; 1\leq i\leq m\Big\}. 
\end{equation*}
Since the order of $\frac{\lambda_2}{\lambda_1}$ is $d$, without loss of generality, if we let $\beta_{i+1}=\Big(\frac{\lambda_2}{\lambda_1}\Big)^{i}\beta_1$ for $1\leq i\leq d$, then 
\begin{equation*}
\prod_{i=1}^{d}(\lambda_1u_1+\beta_i\lambda_2u_2)=\lambda_1^d\prod_{i=1}^{d}\Big(u_1+\beta_i\frac{\lambda_2}{\lambda_1}u_2\Big)=\lambda_1^d\prod_{i=1}^{d}(u_1+\beta_iu_2).
\end{equation*}
Setting $h_1(u_1, u_2):=\prod_{i=1}^{d}(u_1+\beta_i u_2)$, it follows that $h_1$ is a polynomial satisfying the first three assumptions imposed on $h$ and thus $h=h_1$ by the fourth one. Moreover, as $\beta_i=(\frac{\lambda_2}{\lambda_1})^{i-1}\beta_1$ (and the order of $\frac{\lambda_2}{\lambda_1}$ is $d$), 
\begin{equation*}
\prod_{i=1}^{d}(u_1+\beta_iu_2)=u_1^d+\alpha_1 u_2^d
\end{equation*}
for some $\alpha_1\in \mathbb{C}$. Now, considering $\frac{g}{h}$, repeating the same process, we conclude that $g$ is a constant multiple of a function of the form given in \eqref{24091501}. 

The claim for $\alpha$ is then immediate from \eqref{24091501}.   
\end{proof}

Thanks to Lemmas \ref{23112902} and \ref{24010902}, we analyze the equation in \eqref{24050605} in more detail as follows. 

\begin{lemma}\label{24090112}
Let $g_i(u_1, u_2)$ $(i=1,2)$ be homogeneous functions given as
\begin{equation*}
u_1^{k_i}u_2^{l_i}\prod_{j=1}^{n_i} (u_1^d+\beta^i_j u_2^d), \quad (i=1,2)
\end{equation*}
satisfying $(k_1, l_1, n_1)\neq (k_2, l_2, n_2)$, $d\geq 2$, $\deg g_1=\deg g_2$ is odd, and 
\begin{equation}\label{24083101}
a_1\frac{\partial g_1}{\partial u_1}+a_2\frac{\partial g_1}{\partial u_2}+a_3\frac{\partial g_2}{\partial u_1}+a_4\frac{\partial g_2}{\partial u_2}=0
\end{equation}
for some $a_i\in \mathbb{C}$. Then the equation in \eqref{24083101} splits into
\begin{equation}\label{24090109}
a_1\frac{\partial g_1}{\partial u_1}+a_4\frac{\partial g_2}{\partial u_2}=a_2\frac{\partial g_1}{\partial u_2}+a_3\frac{\partial g_2}{\partial u_1}=0.
\end{equation}  
Furthermore, if $d\geq 3$, then either $a_1=a_4=0$ or $a_2=a_3=0$ is true. 
\end{lemma}
\begin{proof}
\begin{enumerate}
\item Suppose $a_i\neq 0$ for all $1\leq i\leq 4$. We first claim
\begin{claim}\label{24120205}
If \eqref{24083101} holds, then either $l_1=0$ or $l_2=0$, and similarly, either $k_1=0$ or $k_2=0$. 
\end{claim}
\begin{proof}[Proof of Claim \ref{24120205}]
Suppose $l_1, l_2\neq 0$. By expanding and rearranging the terms of $\frac{\partial g_i(u_1, u_2)}{\partial u_1}$ and $\frac{\partial g_i(u_1, u_2)}{\partial u_2}$ ($i=1,2$) in descending order of degree with respect to $u_1$, the leading-term exponents are determined as follows: 
\begin{itemize}
\item $(k_1 + n_1 d - 1, l_1)$ for $\frac{\partial g_1}{\partial u_1}$, and $(k_1 + n_1 d, l_1 - 1)$ for $\frac{\partial g_1}{\partial u_2}$;
\item $(k_2 + n_2 d - 1, l_2)$ for $\frac{\partial g_2}{\partial u_1}$, and $(k_2 + n_2 d, l_2 - 1)$ for $\frac{\partial g_2}{\partial u_2}$.
\end{itemize}
By the equation in \eqref{24083101}, since the leading-terms cancel each other, their exponents must be paired, and so 
\begin{equation}\label{24120203}
\begin{gathered}
(k_1 + n_1 d, l_1 - 1)=(k_2+n_2d, l_2-1)\Longrightarrow k_1+n_1d=k_2+n_2d, \;l_1=l_2.
\end{gathered}
\end{equation}

Now we rearrange the terms of $\frac{\partial g_i(u_1, u_2)}{\partial u_1}$ and $\frac{\partial g_i(u_1, u_2)}{\partial u_2}$ ($i = 1, 2$) in descending order of degree with respect to $u_2$, and repeat the same process. 

If $k_1, k_2 \neq 0$, as in the previous case, $k_1 = k_2$, which further implies $n_1=n_2$ by \eqref{24120203}. But this contradicts the assumption that $(k_1, l_1, n_1)\neq (k_2, l_2, n_2)$.  

Consequently, either $k_1$ or $k_2$ must be $0$, and, by symmetry, either $l_1$ or $l_2$ is $0$. This completes the proof of the claim. 
\end{proof}

If $k_1, l_1\neq 0$ and $k_2=l_2=0$, then (again by arranging the terms in \eqref{24083101} in decreasing order with respect to $u_1$) we obtain the following from the exponents of the first and last terms:
\begin{equation}\label{24090111}
\begin{gathered}
(k_1+n_1d, l_1-1)=(n_2d-1,0)\;\;\text{and}\;\;(k_1-1,n_1d+ l_1)=(0, n_2d-1),
\end{gathered}
\end{equation}
which imply $k_1=l_1=1,\; n_1d=n_2d-2 \Longrightarrow d=2, \; n_1-n_2=1$. But this contradicts the fact that $\deg g_1=2+2n_1=\deg g_2=2n_2$ is odd. 

Consequently, either $k_1,l_2\neq 0$ and $k_2=l_1=0$, or $k_1, l_2\neq 0$ and $k_2=l_1=0$. Without loss of generality, if we assume the first case, then 
\begin{equation*}
(k_1-1+n_1d,0)=(n_2d, l_2-1),\; (k_1-1, n_2d)=(0,n_2d+l_2-1)\Longrightarrow k_1=l_2=1, \;n_1=n_2. 
\end{equation*}
Letting $n:=n_1=n_2$, we finally get  
\begin{equation*}
\begin{aligned}
g_1(u_1, u_2)=u_1\sum^{n}_{i=0}\alpha_{1,i} u_1^{(n-i)d}u_2^{id}\;\; \text{and}\;\; g_2(u_1, u_2)=u_2\sum^{n}_{i=0}\alpha_{2,i} u_1^{(n-i)d}u_2^{id}\;\; \textnormal{for some }\alpha_{1,i}, \alpha_{2,i}\in \mathbb{C}, 
\end{aligned}
\end{equation*}
where their partial derivatives are given as   
\begin{equation*}
\begin{aligned}
\frac{\partial g_1}{\partial u_1}=\sum^{n}_{i=0}((n-i)d+1)\alpha_{1,i} u_1^{(n-i)d}u_2^{id},\quad  \frac{\partial g_1}{\partial u_2}=\sum^{n}_{i=1}id\alpha_{1,i} u_1^{(n-i)d+1}u_2^{id-1}
\end{aligned}
\end{equation*}
and
\begin{equation*}
\begin{aligned}
\frac{\partial g_2}{\partial u_1}=\sum^{n-1}_{i=0}(n-i)d\alpha_{2,i} u_1^{(n-i)d-1}u_2^{id+1},\quad  \frac{\partial g_2}{\partial u_2}=\sum^{n}_{i=0}(id+1)\alpha_{2,i} u_1^{(n-i)d}u_2^{id}. 
\end{aligned}
\end{equation*}

If $d\geq 3$, as every term of $\frac{\partial g_1}{\partial u_2}$ (resp. $\frac{\partial g_2}{\partial u_1}$) differs from all others, it follows that $\alpha_{1,i}=0$ for $1\leq i\leq n$ (resp. $\alpha_{2,j}=0$ for $0\leq j\leq n-1$). By \eqref{24083101}, this further implies each $g_i$ ($i=1, 2$) is the zero function, which leads to a contradiction.  

For $d=2$, clearly every term of $\frac{\partial g_1}{\partial u_1}$ (resp. $\frac{\partial g_1}{\partial u_2}$) corresponds to a term of $\frac{\partial g_2}{\partial u_2}$ (resp. $\frac{\partial g_2}{\partial u_1}$). Moreover, since there is no common term between $\frac{\partial g_1}{\partial u_1}$ and $\frac{\partial g_1}{\partial u_2}$, \eqref{24083101} splits into the equations in \eqref{24090109}. 

\item If one of the $a_i$ ($1\leq i\leq 4$) is zero, following similar steps as above, it can be shown exactly two of the $a_i$ ($1\leq i\leq 4$) must be zero. Further, as $\frac{\partial g_1}{\partial u_1}$ and $\frac{\partial g_1}{\partial u_2}$ (resp. $\frac{\partial g_2}{\partial u_1}$ and $\frac{\partial g_2}{\partial u_2}$) do not have a common term, we conclude that either $a_1=a_4=0$ or $a_2=a_3=0$. This completes the proof of the lemma. 
\end{enumerate}
\end{proof}

\noindent \textbf{Remark. } In the above lemma, if $n_i=0$, meaning each $g_i(u_1, u_2)$ $(i=1,2)$ is simply given as $u_1^{k_i}u_2^{l_i}$ with $\det g_1=\det g_2$, then it is straightforward to check that \eqref{24083101} implies either $a_1=a_4=0$ or $a_2=a_3=0$.  \\

Combining Lemmas \ref{23112902}-\ref{24090112} with the above remark, the complete solution of the equation proposed in Proposition \ref{24051906} is as follows:
\begin{proposition}\label{24010906}
Having the same notation as in Proposition \ref{24051906}, if we further suppose $\lambda_i\neq 0$ ($i=1,2$), then the following dichotomy holds for $\tilde{g}_i$ ($i=1,2$) satisfying \eqref{24091301}:
\begin{enumerate}
\item if $\frac{\lambda_1}{\lambda_2}$ is not a root of unity, then $\tilde{g}_i(\tilde{u}_1, \tilde{u}_2)$ ($i=1,2$) is a constant multiple of a polynomial of the following form:
\begin{equation*}
\tilde{u}_1^{k_i}\tilde{u}_2^{l_i} 
\end{equation*}
where $\zeta_i=\lambda_1^{k_i}\lambda_2^{l_i}$;
\item if $\frac{\lambda_1}{\lambda_2}$ is a root of unity of order $d$, then $\tilde{g}_i(\tilde{u}_1, \tilde{u}_2)$ ($i=1,2$) is a constant multiple of a polynomial of the following form:
\begin{equation*}
\tilde{u}_1^{k_i}\tilde{u}_2^{l_i} \prod_{j=1}^{n_i} (\tilde{u}_1^{d}+\beta_j\tilde{u}_2^{d})
\end{equation*}
where $\zeta_i=\lambda_1^{k_i}\lambda_2^{l_i}\lambda_1^{n_id}$. 
\item Further, if $\frac{\lambda_1}{\lambda_2}$ is either not a root of unity or a root of unity of order $d\geq 3$, then either 
\begin{equation*}
(\lambda_1-\omega_1)(\zeta_1-\overline{\omega_1})+(\det A_2)\omega_2\overline{\omega_2}=(\lambda_2-\omega_1)(\zeta_2-\overline{\omega_1})+(\det A_2)\omega_2\overline{\omega_2}=0
\end{equation*}
or 
\begin{equation*}
(\lambda_1-\omega_1)(\zeta_2-\overline{\omega_1})+(\det A_2)\omega_2\overline{\omega_2}=(\lambda_2-\omega_1)(\zeta_1-\overline{\omega_1})+(\det A_2)\omega_2\overline{\omega_2}=0
\end{equation*}
is true. 
\end{enumerate}
\end{proposition}

\subsection{Eigenvalues of $P_M$}

As discussed in the previous subsections, computing the eigenvalues of $P_M$ and $\overline{P_M}$ is critical for solving \eqref{24091301}. In this subsection, we present several useful criteria relevant to this task.  

Note that since $P_M$ and $\overline{P_M}$ are conjugate to each other, their characteristic polynomials are also conjugate. In this case, we show that the eigenvalues of the two matrices are intimately connected, as described below: 
\begin{lemma}\label{23122103}
Consider the following two equations:
\begin{equation}\label{23122101}
x^2+(a+b\sqrt{D})x+(c+d\sqrt{D})=0, \quad x^2+(a-b\sqrt{D})x+(c-d\sqrt{D})=0
\end{equation}
where $a, b, c, d\in \mathbb{Q}$ and $D$ is a negative integer. Let $\lambda_{i}$ and $\zeta_{i}$ ($i=1,2$) be the roots of the first and second equations respectively. 
\begin{enumerate}

\item If \eqref{23122101} are reducible over $\mathbb{Q}(\sqrt{D})$, then we have either $\lambda_1=\overline{\zeta_1}$ and $\lambda_2=\overline{\zeta_2}$, or $\lambda_1=\overline{\zeta_2}$ and $\lambda_2=\overline{\zeta_1}$. 

\item Suppose \eqref{23122101} are irreducible over $\mathbb{Q}(\sqrt{D})$ and, further, $\lambda_1=\zeta_1^{\alpha}\zeta_2^{\beta}$ for some $\alpha, \beta\in \mathbb{Z}$ satisfying $|\alpha+\beta|>1$. Then $\lambda_1\lambda_2$ is a root of unity. 
\end{enumerate}
\end{lemma}
\begin{proof}
\begin{enumerate}
\item Let 
\begin{equation*}
\lambda_{i}=\frac{a+b\sqrt{D}\pm \sqrt{(a+b\sqrt{D})^2-4(c+d\sqrt{D})}}{2}
\end{equation*}
and 
\begin{equation*}
\zeta_{i}=\frac{a-b\sqrt{D}\pm \sqrt{(a-b\sqrt{D})^2-4(c-d\sqrt{D})}}{2}
\end{equation*}
where $i=1,2$. Since $(a+b\sqrt{D})^2-4(c+d\sqrt{D})$ is the conjugate of $(a-b\sqrt{D})^2-4(c-d\sqrt{D})$, if \eqref{23122101} are reducible over $\mathbb{Q}(\sqrt{D})$, then $\lambda_{i}=\frac{a+b\sqrt{D}\pm (e+f\sqrt{D})}{2}$ for some $e+f\sqrt{D}\in \mathbb{Q}(\sqrt{D})$, and clearly $\zeta_{i}=\frac{a-b\sqrt{D}\pm (e-f\sqrt{D})}{2}$ follows.

\item We first claim $\lambda_2=\zeta_1^{\beta}\zeta_2^{\alpha}$. Since $\zeta_1+\zeta_2$ and $\zeta_1\zeta_2$ are elements of $\mathbb{Q}(\sqrt{D})$, by elementary means, one can check $\zeta_1^{\alpha}\zeta_2^{\beta}+\zeta_1^{\beta}\zeta_2^{\alpha}$ is also an element of $\mathbb{Q}(\sqrt{D})$ for any $\alpha, \beta\in \mathbb{Z}_{\geq 0}$. As $(\zeta_1^{\alpha}\zeta_2^{\beta})(\zeta_1^{\beta}\zeta_2^{\alpha})\in \mathbb{Q}(\sqrt{D})$, this implies $\zeta_1^{\alpha}\zeta_2^{\beta}$ and $\zeta_1^{\beta}\zeta_2^{\alpha}$ are the roots of some quadratic polynomial $f(x)$ defined over $\mathbb{Q}(\sqrt{D})$. Since the equations in \eqref{23122101} are irreducible over $\mathbb{Q}(\sqrt{D})$ by the assumption, $\zeta_1^{\alpha}\zeta_2^{\beta}(=\lambda_1)\notin \mathbb{Q}(\sqrt{D})$, which implies $f(x)$ is the minimal polynomial of $\lambda_1$ over $\mathbb{Q}(\sqrt{D})$. By the uniqueness of the minimal polynomial, $f(x)=0$ is essentially equal to the first equation in \eqref{23122101}, and so $\lambda_2=\zeta_1^{\beta}\zeta_2^{\alpha}$.  

Finally, as $\lambda_1\lambda_2=(\zeta_1\zeta_2)^{\alpha+\beta}$, $\lambda_1\lambda_2=\overline{\zeta_1\zeta_2}$ (by \eqref{23122101}) and $|\alpha+\beta|>1$, it is concluded that $\zeta_1\zeta_2$ is a root of unity. 
\end{enumerate}
\end{proof}

The following lemma will be applied repeatedly in Sections \ref{24081804}-\ref{24081805}. 

\begin{lemma}\label{24072601}
Let $D$ be a negative integer. Let $\lambda$ be a root of unity satisfying $\deg \lambda=4$ and $\lambda^2-a\lambda-b=0$ for some $a\in \mathbb{Q}(\sqrt{D})$ and $b=\pm 1$. Then the following statements hold. 
\begin{enumerate} 
\item If $b=1$, then $a$ is either $\pm \sqrt{-1}$ or $\pm \sqrt{-2}$. 
\item If $b=-1$, then $a$ is one of $\frac{-1\pm \sqrt{5}}{2}, \frac{1\pm \sqrt{5}}{2}, \pm \sqrt{3}$ or $\pm \sqrt{2}$. 
\end{enumerate}
\end{lemma}
\begin{proof}
Since $\lambda$ is a root of unity of degree $4$, it is a root of one of the following: 
\begin{equation*}
x^4+x^3+x^2+x+1=0,\quad  x^4-x^3+x^2-x+1=0, \quad x^4-x^2+1=0\quad\text{or}\quad x^4+1=0, 
\end{equation*}
The roots of $x^4+x^3+x^2+x+1=0$ (resp. $x^4-x^2+1=0$) are 
\begin{equation*}
\frac{\sqrt{5}-1\pm\sqrt{-2\sqrt{5}-10}}{4}, \; \frac{-\sqrt{5}-1\pm\sqrt{2\sqrt{5}-10}}{4}\;\; (\text{resp. } \frac{\sqrt{3}\pm \sqrt{-1}}{2},\;\frac{-\sqrt{3}\pm \sqrt{-1}}{2}). 
\end{equation*}
Similarly, the roots of $x^4-x^3+x^2-x+1=0$ (resp. $x^4+1=0$) are
\begin{equation*}
\frac{-\sqrt{5}+1\pm\sqrt{-2\sqrt{5}-10}}{4}, \; \frac{\sqrt{5}+1\pm\sqrt{2\sqrt{5}-10}}{4}\;\; (\text{resp. }\frac{\sqrt{2}\pm \sqrt{-2}}{2},\;  \frac{-\sqrt{2}\pm \sqrt{-2}}{2}).
\end{equation*}
Now the result follows by solving $\lambda^2-a\lambda-b=0$ for each of the possible values of $\lambda$ listed above.
\end{proof}

Along the way, we also state the following lemma, whose proof is analogous to that of the previous lemma, thus omitted here.

\begin{lemma}\label{24122301}
Let $\lambda$ be a root of unity satisfying $\deg \lambda=4$ and $\lambda^2-a\lambda -b=0$ where $a, b$ are both contained in either $\mathbb{Q}(\sqrt{-1})$ or $\mathbb{Q}(\sqrt{-3})$. Then the following statements hold. 
\begin{enumerate}
\item If $b$ is either $\pm \sqrt{-1}$ or $\frac{1\pm\sqrt{-3}}{2}$, then $a=0$. 

\item If $b$ is $\frac{-1\pm \sqrt{-3}}{2}$, then there is no $a\in \mathbb{Q}(\sqrt{-3})$ satisfying the required condition. 
\end{enumerate}
\end{lemma}

\subsection{Compatibility}\label{25011407}
In this subsection, we discuss the compatibility of two elements in $\mathrm{Aut}\,\mathcal{X}$. As an initial step, we examine $\tilde{g}_i$ ($i=1,2$) in \eqref{24091301}, which further satisfies \eqref{24091601}. For instance, in the following lemma, we show that, assuming $\frac{\lambda_1}{\lambda_2}$ is not a root of unity, the condition in \eqref{24091601} imposes very strong constraints on both $S_{\lambda}^{-1}P\iota S_{\lambda}$ and $S_{\zeta}^{-1}\overline{P}\iota S_{\zeta}$.  

\begin{lemma}\label{24052409}
Let $g_i(u_1, u_2)=u_1^{k_i}u_2^{l_i}$ where $k_i, l_i$ ($i=1,2$) are some non-negative integers satisfying 
\begin{itemize}
\item $k_1+l_1=k_2+l_2\geq 3$ and it is odd, 
\item $(k_1, l_1)\neq (k_2, l_2)$.  
\end{itemize}
If there exist $P, Q\in \mathrm{GL}_2(\overline{\mathbb{Q}})$ such that 
\begin{equation}\label{24052401}
\left(\begin{array}{c}
g_1(P\bold{U})\\
g_2(P\bold{U})
\end{array}\right)=
Q\left(\begin{array}{c}
g_1(\bold{U})\\
g_2(\bold{U})
\end{array}\right)
\end{equation}
where $\bold{U}:=\left(\begin{array}{c}
u_1\\
u_2
\end{array}\right)$, then both $P$ and $Q$ are either diagonal or anti-diagonal matrices.
\end{lemma}
\begin{proof}
Let $P:=\left(\begin{array}{cc}
\alpha_1 & \beta_1\\
\gamma_1 & \delta_1
\end{array}\right)$ and $Q:=\left(\begin{array}{cc}
\alpha_2 & \beta_2\\
\gamma_2 & \delta_2
\end{array}\right)$. Without loss of generality, if we assume $k_2\geq k_1$ and $l_1\geq l_2$, then \eqref{24052401} is expanded as
\begin{equation}\label{24082305}
\begin{aligned}
\left(\begin{array}{c}
(\alpha_1u_1+\beta_1u_2)^{k_1} (\gamma_1u_1+\delta_1u_2)^{l_1}\\
(\alpha_1u_1+\beta_1u_2)^{k_2} (\gamma_1u_1+\delta_1u_2)^{l_2}
\end{array}\right)
=\left(\begin{array}{c}
u_1^{k_1}u_2^{l_2}(\alpha_2u_2^{l_1-l_2}+\beta_2u_1^{k_2-k_1}) \\
u_1^{k_1}u_2^{l_2}(\gamma_2u_2^{l_1-l_2}+\delta_2u_1^{k_2-k_1}) 
\end{array}\right).
\end{aligned}
\end{equation}
By the assumption on $k_i$ and $l_i$ ($i=1,2$), one gets either $\beta_1=\gamma_1=0$ or $\alpha_1=\delta_1=0$. 
\begin{enumerate}
\item If $\beta_1=\gamma_1=0$, then \eqref{24082305} reduces to 
\begin{equation*}
\begin{aligned}
(\alpha_1u_1)^{k_1}(\delta_1u_2)^{l_1}=u_1^{k_1}u_2^{l_2}(\alpha_2u_2^{l_1-l_2}+\beta_2u_1^{k_2-k_1}) ,\quad (\alpha_1u_1)^{k_2}(\delta_1u_2)^{l_2}=u_1^{k_1}u_2^{l_2}(\gamma_2u_2^{l_1-l_2}+\delta_2u_1^{k_2-k_1}),  
\end{aligned}
\end{equation*}
thus $\beta_2=\gamma_2=0$. 

\item If $\alpha_1=\delta_1=0$, then \eqref{24082305} reduces to 
\begin{equation}\label{24082306}
\begin{aligned}
(\beta_1u_2)^{k_1} (\gamma_1u_1)^{l_1}=u_1^{k_1}u_2^{l_2}(\alpha_2u_2^{l_1-l_2}+\beta_2u_1^{k_2-k_1}), \quad (\beta_1u_2)^{k_2} (\gamma_1u_1)^{l_2}=u_1^{k_1}u_2^{l_2}(\gamma_2u_2^{l_1-l_2}+\delta_2u_1^{k_2-k_1}),
\end{aligned}
\end{equation}
so either $\alpha_2=0$ or $\beta_2=0$, and either $\gamma_2=0$ or $\delta_2=0$. If $\alpha_2=\gamma_2=0$ or $\beta_2=\delta_2=0$, then $(k_1, l_1)=(k_2, l_2)$, contradicting our assumption. Thus either $\alpha_2=\delta_2=0$ or $\beta_2=\gamma_2=0$ must hold. But, for the second case, \eqref{24082306} is further reduced to 
\begin{equation*}
\begin{aligned}
(\beta_1u_2)^{k_1} (\gamma_1u_1)^{l_1}=\alpha_2u_1^{k_1}u_2^{l_1}, \quad (\beta_2u_2)^{k_2} (\gamma_1u_1)^{l_2}=\delta_2u_1^{k_2}u_2^{l_2}, 
\end{aligned}
\end{equation*}
implying $k_i=l_i$ ($i=1,2$), which contradicts the fact that $k_i+l_i$ is odd.  
\end{enumerate}
In conclusion, either $\beta_i=\gamma_i=0$ or $\alpha_i=\delta_i=0$ $(i=1,2)$, and hence both $P$ and $Q$ are either diagonal or anti-diagonal matrices. 
\end{proof}

As a corollary, if $\frac{\lambda_1}{\lambda_2}$ (where $\lambda_1$ and $\lambda_2$ are defined in \eqref{24091301}) is not a root of unity, both $S_{\lambda}^{-1}P\iota S_{\lambda}$ and $S_{\zeta}^{-1}\overline{P}\iota S_{\zeta}$ in \eqref{24091601} are either diagonal or anti-diagonal.  

Another immediate and useful result obtained from the above lemma is as follows:
\begin{proposition}\label{24072003}
Let $P:=\left(\begin{array}{cc}
\omega_1 & \omega_2(\neq 0) \\
\omega_3 & \omega_4 
\end{array}\right)\in\mathrm{GL}_2(\mathbb{Q}(\sqrt{D}))$ where $D$ is a negative integer, and $\Theta_i(u_1, u_2)$ ($i=1,2$) be homogeneous polynomials of the same degree $\geq 3$. Suppose 
\begin{equation}\label{25011003}
\left(\begin{array}{c}
\Theta_1(\iota\bold{U}) \\
\Theta_2(\iota\bold{U})
\end{array}
\right)=\iota\left(\begin{array}{c}
\Theta_1(\bold{U}) \\
\Theta_2(\bold{U})
\end{array}\right)
\quad \text{and}\quad 
\left(\begin{array}{c}
\Theta_1(P\bold{U}) \\
\Theta_2(P\bold{U})
\end{array}
\right)=\overline{P}\left(\begin{array}{c}
\Theta_1(\bold{U}) \\
\Theta_2(\bold{U})
\end{array}\right)
\end{equation}
where $\bold{U}:=\left(\begin{array}{c}
u_1 \\
u_2
\end{array}\right)$. Let $P_1$ and $P_2$ be two elements in $\langle P, \iota \rangle $. If $\lambda_i$ ($i=1,2$) are two nonzero eigenvalues of $P_1$ such that $\frac{\lambda_1}{\lambda_2}$ is not a root of unity, then we have either $P_1P_2=P_2P_1$ or $\mathrm{tr}\,P_2=0$. 
\end{proposition}
\begin{proof}
To simplify the proof, we assume $P_1:=P$, and let $\tilde{g}_i$, $S_{\lambda}$ and $S_{\zeta}$ be the same as in Proposition \ref{24051906}. Since $\frac{\lambda_1}{\lambda_2}$ is not a root of unity, $\tilde{g}_i$ ($i=1,2$) in \eqref{24091301} is of the form $\tilde{u}_1^{k_i}\tilde{u}_2^{l_i}$ ($i=1,2$) by Proposition \ref{24010906}. Since $\lambda_1\neq \lambda_2$, it implies $\zeta_1\neq \zeta_2$, and thus $(k_1, l_1)\neq (k_2, l_2)$ (again by Proposition \ref{24010906}). Then 
\begin{equation*}
\left(\begin{array}{c}
\Theta_1(P_2\bold{U}) \\
\Theta_2(P_2\bold{U})
\end{array}
\right)=\overline{P_2}\left(\begin{array}{c}
\Theta_1(\bold{U}) \\
\Theta_2(\bold{U})
\end{array}\right)
\end{equation*}
is equivalent to 
\begin{equation*}
\left(\begin{array}{c}
\tilde{g}_1(S_{\lambda}^{-1}P_2 S_{\lambda}\bold{\tilde{U}})\\
\tilde{g}_2(S_{\lambda}^{-1}P_2 S_{\lambda}\bold{\tilde{U}})
\end{array}\right)=
S_{\zeta}^{-1}\overline{P_2} S_{\zeta}
\left(\begin{array}{c}
\tilde{g}_1(\bold{\tilde{U}})\\
\tilde{g}_2(\bold{\tilde{U}})
\end{array}\right), 
\end{equation*}
and, applying Lemma \ref{24052409}, both $S_{\lambda}^{-1}P_2 S_{\lambda}$ and $S_{\zeta}^{-1}\overline{P_2}S_{\zeta}$ in \eqref{24091601} are either diagonal or anti-diagonal matrices. For the first case, it means both $P_1$ and $P_2$ are simultaneously diagonalizable by $S_{\lambda}$, thus they commute. The second case obviously yields $\mathrm{tr}\, P_2=\mathrm{tr}\, \overline{P_2}=0$.
\end{proof}

By amalgamating the results obtained thus far in this section as well as Section \ref{25012001}, we finally arrive at the following proposition. This proposition provides a simple yet highly useful criterion. It will be particularly crucial for analyzing the structure of elements of $\mathrm{Aut}\,\mathcal{X}$ and for checking whether any two elements $M$ and $N\in \mathrm{GL}_4(\mathbb{Q})$ belong to the same subgroup of $\mathrm{Aut}\,\mathcal{X}$ or not. 


\begin{proposition}\label{24082301}
Let $\mathcal{X}$ be an element in $\mathfrak{Ghol}$, and $M, N\in \mathrm{GL}_4(\mathbb{Q})$ be elements in $\mathrm{Aut}\,\mathcal{X}$. Let $P_M$ and $P_N$ be the primary matrices associated with $M$ and $N$ respectively. If $P\in\langle P_M, P_N, \iota \rangle $ has a nonzero $(1,2)$-entry and distinct nonzero eigenvalues, then either $\mathrm{tr}\,(P\iota)=0$ or both eigenvalues of $P$ are roots of unity.
\end{proposition}
\begin{proof}
If $\Theta_i$ ($i=1,2$) are given as in Proposition \ref{24020501}, then clearly the equations in \eqref{25011003} holds for $P$. Let 
\begin{equation}\label{24082201}
\mathfrak{\chi}_P(x):=x^2+(a+b\sqrt{D})x+(c+d\sqrt{D})
\end{equation}
where $D$ is a negative integer, and $\lambda_i$ $(i=1, 2)$ be the eigenvalues of $P$. 
\begin{enumerate}
\item If $\frac{\lambda_1}{\lambda_2}$ is not a root of unity, by Proposition \ref{24072003}, we have either $\mathrm{tr}\,(P\iota)=0$ or $P(P\iota)=(P\iota) P$. If $P$ is given as $\left(\begin{array}{cc}
\omega_1 & \omega_2(\neq 0)\\
\omega_3 & \omega_4
\end{array}\right)$, then the latter condition implies 
\begin{equation*}
\begin{aligned}
\left(\begin{array}{cc}
\omega_1 & \omega_2\\
\omega_3 & \omega_4
\end{array}\right)
\left(\begin{array}{cc}
\omega_1 & -\omega_2\\
\omega_3 & -\omega_4
\end{array}\right)=&\left(\begin{array}{cc}
\omega_1 & -\omega_2\\
\omega_3 & -\omega_4
\end{array}\right)\left(\begin{array}{cc}
\omega_1 & \omega_2\\
\omega_3 & \omega_4
\end{array}\right)\\
\Longrightarrow 
\left(\begin{array}{cc}
\omega_1^2+\omega_2\omega_3 & -\omega_1\omega_2-\omega_2\omega_4 \\
\omega_3\omega_1+\omega_4\omega_3 & -\omega_3\omega_2-\omega_4^2
\end{array}\right)
 =&\left(\begin{array}{cc}
\omega_1^2-\omega_2\omega_3 & \omega_1\omega_2-\omega_2\omega_4\\
\omega_3\omega_1-\omega_4\omega_3 & \omega_3\omega_2-\omega_4^2
\end{array}\right).
\end{aligned}
\end{equation*}
Thus $\omega_2\omega_3= \omega_1\omega_2=\omega_4\omega_3=\omega_3\omega_2=0\Longrightarrow \omega_1=\omega_3=0$. But this contradicts the fact $\lambda_i$ are nonzero. Consequently, $\mathrm{tr}\,(P\iota)=0$.

\item Suppose $\frac{\lambda_1}{\lambda_2}$ is a root of unity. 
\begin{enumerate}
\item If \eqref{24082201} is reducible over $\mathbb{Q}(\sqrt{D})$, by Lemma \ref{23122103}, we have either $\lambda_1=\overline{\zeta_1}$ and $\lambda_2=\overline{\zeta_2}$, or $\lambda_1=\overline{\zeta_2}$ and $\lambda_2=\overline{\zeta_1}$. Since each $\zeta_i$ is of the form $\lambda_1^{\alpha_i}\lambda_2^{\beta_i}$ with $\alpha_i+\beta_i\geq 3$, as stated in Proposition \ref{24010906}, it follows that all $\lambda_i$ ($i=1,2$) are roots of unity. 

\item If \eqref{24082201} is irreducible over $\mathbb{Q}(\sqrt{D})$, then $\det P(=\lambda_1\lambda_2)$ is a root of unity by Lemma \ref{23122103}, concluding $\lambda_i$ ($i=1,2$) are all roots of unity. 
\end{enumerate}
\end{enumerate}
\end{proof}


As emphasized previously, Propositions \ref{24051906} and \ref{24082301}, along with Lemma \ref{25011101}, form our primary toolkit for tackling the main arguments in later sections.

When $M$ is given as in \eqref{25010811}, suppose that, after normalization, $M$ is given as 
\begin{equation*}
\left(\begin{array}{cc}
A_1 & I\\
A_2A_3 & A_2A_4A_2^{-1}
\end{array}\right). 
\end{equation*}
Then, using the minimal polynomial condition on $M$, one can compute the traces and determinants of $A_2A_3$ and $A_2A_4A_2^{-1}$ in terms of the trace and determinant of $A_1$. As a result, each entry of the primary matrix associated with $P_M$ is likewise expressed in terms of the trace and determinant of $A_1$. By applying either Lemma \ref{25011101} or Proposition \ref{24082301}, the exact determinant and trace of $A_1$ are determined, and hence the explicit form of $M$ as well, thanks to Lemma \ref{24090301}. For simple cases, where the degree of the minimal polynomial of $M$ is at most $3$, Lemma \ref{25011101} suffices; otherwise, Proposition \ref{24082301} may be required.

After analyzing all possible structures of $M$, we proceed to examine the structure of a subgroup of $\mathrm{Aut}\,\mathcal{X}$. As mentioned earlier, the crucial point here is determining whether two elements are compatible with each other, that is, whether they belong to the same subgroup. In this case, as expected, Proposition \ref{24082301} will be further adapted and serve as a cornerstone to support the main arguments.

\newpage

\section{Classification of $M$: I}\label{24081804}
Using the strategies outlined and the results obtained in Sections \ref{25012001}-\ref{25011802}, in this section and the next, we analyze the structure of $M$, determining its exact form when it appears as an element of $\mathrm{Aut}\,\mathcal{X}$ for some $\mathcal{X}\in \mathfrak{Ghol}$. The end results are all summarized in Theorems \ref{24090201}, \ref{25012901}, and \ref{24092702}, each of which concerns Type I, II, and III, respectively. 

These two sections constitute the most technical part of the paper. As mentioned earlier, because we pursue effective forms of the statements, the proofs here inevitably involve laborious computations and repeated citations of various lemmas proved in Section~\ref{25011802}. We try to simplify computations as much as possible and, in particular, to avoid repetitive ones; however, when necessary, we do not hesitate to include them. 

We encourage the reader to skip the proofs in these sections on a first reading and focus only on the statements of the theorems. Nonetheless, for readers who wish to examine a proof on a first reading, it suffices to skim the primitive case in Section~\ref{25111903}. In broad terms, every proof follows a strategy similar to that described in Sections~\ref{25012001}-\ref{25011802}, although they differ slightly in many details. 

Thanks to Proposition \ref{23062802}, the list in \eqref{23062301} is now reduced to 
\begin{equation*}
x^2\pm1, \; x^2\pm x+1, \; x^3\pm1, 
\; x^3\pm x^2+x\pm 1, \;x^3\pm2x^2+2x\pm1,\; x^4\pm1, \; x^4\pm x^2+1.
\end{equation*}
We first focus only on the following in this section (while the remaining cases will be addressed in the next section):
\begin{equation*}
\begin{gathered}
x^2\pm1,\; x^2\pm x+1,\; x^3\pm1,\; x^3\pm x^2+x\pm1,\;x^3\pm2x^2+2x\pm1.
\end{gathered}
\end{equation*}
In addition, rather than assuming that $M$ is of Type I, II or III, we adopt a more general approach for the following:
\begin{equation}\label{25010901}
x^2\pm1,\;\; x^2\pm x+1.
\end{equation}
That is, when $\mathfrak{m}_M(x)$ is one of \eqref{25010901}, we do not necessarily restrict $M$ to be of Type I, II or III, but instead, we allow $M$ to be an arbitrary matrix. This generality is essential because a broader framework is required to address the following remaining four later in Section \ref{24081805}:\footnote{For instance, when considering $M$ such that $\mathfrak{m}_M(x)=x^4-1$, we will reduce it to the case of $M^2$ to use the fact that $ \mathfrak{m}_{M^2}(x)=x^2-1$. However, even if $M$ is of Type III and satisfies $\mathfrak{m}_M(x)=x^4-1$, it does not necessarily mean that $M^2$ is also of Type III. Therefore, we need to analyze all the possible forms of $M$ in a broader setting.}
\begin{equation*}
x^4\pm1,\;\; x^4\pm x^2+1. 
\end{equation*}

Before beginning the proof, we recall that if $A_j$ ($1\leq j\leq 4$) and $M$ are as in \eqref{22022010}-\eqref{21072901}, and moreover, if both $M$ and $A_1$ are invertible, then a well-known formula for $M^{-1}$ is given as 
\begin{equation}\label{23070519}
\left(\begin{array}{cc}
A_1^{-1}+A_1^{-1}A_2(M/A_1)^{-1}A_3A_1^{-1} & -A_1^{-1}A_2(M/A_1)^{-1}\\
-(M/A_1)^{-1}A_3A_1^{-1} & (M/A_1)^{-1}
\end{array}
\right)
\end{equation}
where $M/A_1:=A_4-A_3A_1^{-1}A_2$.

\subsection{$x^2-1$}\label{25111903}
In this subsection, we consider $M$ such that $\mathfrak{m}_M(x)=x^2-1$ or, equivalently, $M=M^{-1}$. First, by normalizing $M$ and equating it with the inverse formula in \eqref{23070519}, we represent $M$ in terms of $A_1$. We then find the primary matrix associated with $M$, compute its eigenvalues, and explicitly determine all possible forms of $M$. This is the simplest type, but the proof of the theorem contains all the necessary information and strategies that will be subsequently employed in the proofs of many other cases. 

\begin{theorem}\label{23063002}
Let $\mathcal{X}$ be an element in $\mathfrak{Ghol}$ and $M$ be an element in $\mathrm{Aut}\,\mathcal{X}$ such that $\mathfrak{m}_M(x)=x^2-1$. Then the following statements hold:
\begin{enumerate}
\item if $M$ is of block diagonal (resp. anti-diagonal) matrix, then it is 
\begin{equation*}
\pm (I\oplus -I) \quad (\text{resp. }
A_2\widetilde{\oplus} A_2^{-1} \text{ for some }A_2\in \mathrm{GL}_2(\mathbb{Q}));
\end{equation*}
\item otherwise, it is either  
\begin{equation}\label{24120209}
\left(\begin{array}{cc}
\frac{I}{2} & A_2\\
\frac{3}{4}A_2^{-1} & -\frac{I}{2}
\end{array}\right)\;\; \text{or}\;\; 
\left(\begin{array}{cc}
-\frac{I}{2} & A_2\\
\frac{3}{4}A_2^{-1} & \frac{I}{2}
\end{array}\right)
\end{equation}
for some $A_2\in \mathrm{GL}_2(\mathbb{Q})$. 
\end{enumerate}
\end{theorem}
\begin{proof}
The first claim is straightforward, so we only consider the second one. That is, we assume $M$ is neither of block diagonal nor anti-diagonal matrix, and $\det A_j\neq 0$ for $1\leq j\leq 4$ by Corollary \ref{25010913}.

To simplify the problem, by change variables, let us normalize $M$ as 
\begin{equation}\label{25022101}
\left(\begin{array}{cc}
A_1 & I\\
A_3 & A_4
\end{array}
\right). 
\end{equation}
Since $M=M^{-1}$, we get 
\begin{equation*}
A_1=A_1^{-1}+A_1^{-1}(M/A_1)^{-1}A_3A_1^{-1}, \; I=-A_1^{-1}(M/A_1)^{-1}, \; A_3=-(M/A_1)^{-1}A_3A_1^{-1}, \; A_4=(M/A_1)^{-1}
\end{equation*}
by \eqref{23070519}, which imply 
\begin{equation*}
A_4=-A_1, \; A_1=A_1^{-1}-A_3A_1^{-1}\Longrightarrow A_3=I-A_1^2=-(\mathrm{tr}\, A_1)A_1+(1+\det A_1)I. 
\end{equation*}
Thus the primary matrix $P_M$ associated to $M$ is 
\begin{equation}\label{25011102}
\left(\begin{array}{cc}
\omega & 1\\
1-\omega^2 & -\omega
\end{array}\right)
\end{equation}
where $\omega=\frac{\mathrm{tr}\,A_1\pm\sqrt{\mathrm{Disc}\,A_1}}{2}$, and the eigenvalues $\lambda_{i}$ $(i=1,2)$ of $P_M\iota$ are the solutions of  
\begin{equation*}
\det \left(\begin{array}{cc}
\omega-x & -1\\
1-\omega^2 & \omega-x 
\end{array}\right)=x^2-2\omega x+1.
\end{equation*}
By Lemma \ref{25011101}, since the order of \eqref{25011102} is finite, $\lambda_{i}$ $(i=1,2)$ are roots of unity. 

\begin{enumerate}
\item If $\deg \lambda_{i}\leq 2$ ($i=1,2$), then they are the roots of either $x^2+1=0$ or $x^2\pm x+1=0$. If $\omega=0$, then it contradicts the fact that $\det A_1\neq 0$. Hence $\omega=\pm \frac{1}{2}$ and $M$ is of the form given in \eqref{24120209}. 

\item If $\deg \lambda_{i}=4$ ($i=1,2$), by Lemma \ref{24072601}, $2\omega$ is one of the following: $\frac{-1\pm \sqrt{5}}{2}, \frac{1\pm \sqrt{5}}{2}, \pm \sqrt{3}$ or $\pm \sqrt{2}$. But this contradicts the assumption that $\omega\in \mathbb{Q}(\sqrt{D})$ for some negative integer $D$.  
\end{enumerate}
This completes the proof of the theorem.
\end{proof}

\subsection{$x^2+1$}
In this subsection, we consider $M\in \mathrm{Aut}\,\mathcal{X}$ such that $\mathfrak{m}_M(x)=x^2+1$. The strategy of the proof closely parallels that of the previous subsection. 
\begin{theorem}\label{23070521}
Let $\mathcal{X}$ be an element in $\mathfrak{Ghol}$ and $M$ be an element in $\mathrm{Aut}\,\mathcal{X}$ such that $\mathfrak{m}_M(x)=x^2+1$. Then the following statements hold. 
\begin{enumerate}
\item If $M$ is of  block diagonal (resp. anti-diagonal) matrix, it is of the form
\begin{equation*}
A_1\oplus A_4 \quad (\text{resp. }
A_2\widetilde{\oplus} (-A_2^{-1}) )
\end{equation*}
for some $A_1, A_4\in \mathrm{GL}_2(\mathbb{Q})$ (resp. $A_2\in \mathrm{GL}_2(\mathbb{Q})$) satisfying $\mathfrak{m}_{A_1}(x)=\mathfrak{m}_{A_4}(x)=x^2+1$. 

\item Otherwise, it is of the following form:
\begin{equation}\label{24120208}
\left(\begin{array}{cc}
A_1 & A_2\\
-(\det A_2)A_2^{-1} & -A_2^{-1}A_1A_2
\end{array}\right) 
\end{equation}
for some $A_1, A_2\in \mathrm{GL}_2(\mathbb{Q})$ such that $\det A_1+\det A_2=1$. In particular, $\mathfrak{m}_{A_1}(x)$ is one of $x^2+\frac{1}{4}, x^2+\frac{1}{2}$ or $x^2+\frac{3}{4}$.
\end{enumerate}
\end{theorem}
\begin{proof}
Similar to the previous theorem, the first claim is straightforward, so we only focus on the second one. 

To simplify the problem, we normalize $M$ as in \eqref{25022101}. Since $M = -M^{-1}$, using \eqref{23070519} and applying the same procedure as in the previous theorem, we get $A_4=-A_1, A_3=-I-A_1^2$, and thus the primary matrix $P_M$ associated to $M$ is  
\begin{equation*}
\left(\begin{array}{cc}
\omega & 1\\
-\omega^2-1 & -\omega
\end{array}
\right)
\end{equation*}
where $\omega=\frac{\mathrm{tr}\, A_1\pm\sqrt{\mathrm{Disc}\, A_1}}{2}$. By Lemma \ref{25011101}, $P_M\iota$ has a finite order, hence the eigenvalues $\lambda_{i}$ ($i=1,2$) of it, which are the solutions of $x^2-2\omega x-1=0$, are the roots of unity. 
\begin{enumerate}
\item If $\deg \lambda_{i}\leq 2$ ($i=1,2$), then they are roots of $x^2-1=0$, $x^2+1=0$ or $x^2\pm x+1=0$. For the first case, we have $\omega=0$, and so $\mathrm{tr}\,A_1=\det A_1=0$. But it contradicts the fact that $\det A_1\neq 0$. For the second case, we get $\omega=\pm\sqrt{-1}$, which implies $\mathrm{tr}\,A_3=\det A_3=0$, again leading to the same contradiction. For the last case, it follows that $\omega=\pm \frac{\sqrt{-3}}{2}$, implying $\mathrm{tr}\,A_1=0$ and $\det A_1=\frac{3}{4}$. 

\item If $\deg \lambda_{i}=4$ ($i=1,2$), by Lemma \ref{24122301}, $2\omega$ is either $\pm \sqrt{-1}$ or $\pm \sqrt{-2}$. Thus $\mathrm{tr}\,A_1=0$ and $\det A_1$ is either $\frac{1}{4}$ or $\frac{1}{2}$. For each case, since the order of $\frac{\lambda_1}{\lambda_2}$ is at least $3$, $\det A_2$ is equal to either $\det A_1-1$ or $-\det A_1+1$ by Proposition \ref{24010906}. As $\det A_2>0$, the latter holds, that is, $\det A_1+\det A_2=1$. 

In conclusion, as $\det A_2\neq 0$, $\mathrm{tr}\,A_1=0$ and $\det A_1$ is one of $\frac{1}{4}, \frac{1}{2}$ or $\frac{3}{4}$. The original $M$, prior to normalization, is of the type given in \eqref{24120208}.  
\end{enumerate}
\end{proof}

\begin{example}\label{25110401}
{\normalfont For $\mathcal{M}$ as given in Example \ref{25101509} and $H$ as defined in \eqref{25100902}, if we normalize the coefficient matrix of $H$ by \eqref{22022007}, then, as an element of $\mathrm{GL}_{4}(\mathbb{Q})$, it is given by  
\begin{equation}\label{25112505}
\left(\begin{array}{cccc}
0 & \frac{1}{2} & 0  & \frac{1}{2}\\
-1 & 0 & -1 & 0\\
0 & \frac{1}{2} & 0 & -\frac{1}{2}\\
-1 & 0 & 1 & 0
\end{array}\right).
\end{equation}
Here, $A_1=A_2=A_3=-A_4=\left(\begin{array}{cc}
0 & \frac{1}{2} \\
-1 & 0 
\end{array}\right)$, and hence it satisfies the second criterion in the second case of the preceding theorem.}
\end{example}

\subsection{$x^2\pm x+1$}
In this subsection, we consider $M$ such that $\mathfrak{m}_M(x)$ is either $x^2-x+1$ or $x^2+x+1$. Note that if $\mathfrak{m}_M(x)=x^2-x+1$, then $\mathfrak{m}_{-M}(x)=x^2+x+1$, and vice versa. Hence we prove only one case, as the other follows immediately. Similar to the previous subsections, we follow the same outline. 

\begin{theorem}\label{24040801}
Let $\mathcal{X}$ be an element in $\mathfrak{Ghol}$ and $M$ be an element in $\mathrm{Aut}\,\mathcal{X}$ such that $\mathfrak{m}_M(x)=x^2-x+1$. Then $M$ is not of block anti-diagonal matrix, and we have the following dichotomy.  
\begin{enumerate}
\item If $M$ is of block diagonal matrix, then it is one of the following forms:
\begin{equation}\label{24102401}
A_1 \oplus A_4,\quad A_1 \oplus I \quad \text{or}\quad I \oplus A_4
\end{equation}
where $A_1, A_4\in \mathrm{GL}_2(\mathbb{Q})$ and $\mathfrak{m}_{A_1}(x)=\mathfrak{m}_{A_4}(x)=x^2-x+1$. 

\item Otherwise, $M$ is one of the following:  
\begin{equation}\label{25011011}
\left(\begin{array}{cc}
\frac{I}{2} & A_2\\
-\frac{3}{4}A_2^{-1} & \frac{I}{2}
\end{array}
\right), \quad  
\left(\begin{array}{cc}
A_1 & A_2\\
-\frac{A_2^{-1}}{2} & \frac{A_2^{-1}A_1^{-1}A_2}{2}
\end{array}\right) 
\quad \text{or}\quad
\left(\begin{array}{cc}
A_1 & A_2\\
-\frac{A_2^{-1}}{4} & \frac{3A_2^{-1}A_1^{-1}A_2}{4}
\end{array}\right)
\end{equation}
for some $A_1, A_2\in \mathrm{GL}_2(\mathbb{Q})$ satisfying $\det A_2=\frac{3}{4}, \mathfrak{m}_{A_1}(x)=x^2-x+\frac{1}{2}$ or $\mathfrak{m}_{A_1}(x)=x^2-x+\frac{3}{4}$ respectively.
\end{enumerate}
\end{theorem}
\begin{proof}
For a block diagonal or anti-diagonal matrix $M$, the conclusion is easily verified. Thus we consider only the case where $M$ is neither block diagonal nor anti-diagonal matrix, and $\det A_j\neq 0$ for $1\leq j\leq 4$ (by Corollary \ref{25010913}).

To simplify the problem, we normalize $M$ as in \eqref{25022101}. Since $M^{-1}=I-M$, using \eqref{23070519} and applying the same method in Theorem \ref{23063002}, it follows that $A_4=I-A_1, A_3=-A_1^2+A_1-I$, and thus the primary matrix $P_M$ associated to $M$ is 
\begin{equation*}
\left(\begin{array}{cc}
\omega & 1\\
-\omega^2+\omega-1 & 1-\omega
\end{array}\right)
\end{equation*}
where $\omega:=\frac{\mathrm{tr}\,A_1\pm\sqrt{\mathrm{Dics}\,A_1}}{2}$. If we denote two eigenvalues of $P_M\iota$ by $\lambda_{i}$ ($i=1,2$), then $\lambda_{i}$ are the solutions of $x^2-(2\omega-1)x-1=0$ and roots of unity by Lemma \ref{25011101}. 

We analyze each case separately, based on the degree of $\lambda_i$. 
\begin{enumerate}
\item First, if $\deg \lambda_{i}\leq 2$, then they are $\pm 1$, $\pm \sqrt{-1}$, $\frac{\pm 1+\sqrt{-3}}{2}$ or $\frac{\pm 1- \sqrt{-3}}{2}$.   
\begin{enumerate}
\item If $\lambda_{i}=\pm 1$, then $\omega=\frac{1}{2}$, and, consequently, the original $M$, before normalization, is of the first type given in \eqref{25011011}.

\item If $\lambda_{i}=\pm \sqrt{-1}$, then $2\omega-1=\pm 2 \sqrt{-1}$. Without loss of generality, let us assume $\omega=\sqrt{-1}+\frac{1}{2}$ and set $P_M=\left(\begin{array}{cc}
\sqrt{-1}+\frac{1}{2} & 1\\
\frac{1}{4} & -\sqrt{-1}+\frac{1}{2}
\end{array}\right)$. Then $P:=(P_M\iota) P_M(P_M\iota)=\left(\begin{array}{cc}
-\frac{3}{2} & 2\sqrt{-1}-1\\
-\frac{1}{4}-\frac{\sqrt{-1}}{2} & -\frac{3}{2}
\end{array}\right)$ and the eigenvalues of $P$ are $\frac{-3\pm\sqrt{5}}{2}$. But this contradicts the fact that the eigenvalues of $P$ must be roots of unity (by Lemma \ref{25011101}). 

\item If $\lambda_{i}=\frac{\pm 1+ \sqrt{-3}}{2}$, then $2\omega-1=\pm \sqrt{-3}$. In this case, $P_M=\left(\begin{array}{cc}
\frac{1+\sqrt{-3}}{2} & 1\\
0 & \frac{1-\sqrt{-3}}{2}
\end{array}\right)$, implying $\mathrm{tr}\,A_3=\det A_3=0$. But it contradicts the fact that $\det A_3\neq 0$.

Similarly, we get the same contradiction for $\lambda_i=\frac{\pm 1- \sqrt{-3}}{2}$. 
\end{enumerate}

\item Next if $\deg \lambda_{i}=4$ ($i=1,2$), by Lemma \ref{24072601}, $2\omega-1$ is either $\pm \sqrt{-1}$ or $\pm \sqrt{-2}$. If $2\omega-1=\pm \sqrt{-1}$ (resp. $\pm \sqrt{-2}$), then $\mathrm{tr}\,A_1=1, \det A_1=\frac{1}{2}$ (resp. $\det A_1=\frac{3}{4}$)  and $M$ is of the form 
\begin{equation*}
\left(\begin{array}{cc}
A_1 & I\\
-\frac{I}{2} & \frac{A_1^{-1}}{2}
\end{array}\right)\quad
(\text{resp. }\left(\begin{array}{cc}
A_1 &I\\
-\frac{I}{4} & \frac{3A_1^{-1}}{4}
\end{array}\right))
\end{equation*}
where $\mathfrak{m}_{A_1}(x)=x^2-x+\frac{1}{2}$ (resp. $x^2-x+\frac{3}{4}$). Consequently, the original $M$, prior to normalization, is the second (resp. third) type given in \eqref{25011011}. 
\end{enumerate}
\end{proof}

\noindent\textbf{Remark. } Note that the first matrix in \eqref{25011011} is obtained by multiplying $\iota$ with the first one in \eqref{24120209}. As summarized in Theorem \ref{24092702}, for each element of Type III, there exists a corresponding element (of the same type) paired in this manner. \\

The following corollary is immediate from the previous theorem. 
\begin{corollary}\label{24121701}
Let $\mathcal{X}$ be an element in $\mathfrak{Ghol}$ and $M$ be an element in $\mathrm{Aut}\,\mathcal{X}$ such that $\mathfrak{m}_M(x)=x^2+x+1$. Then $M$ is not of block anti-diagonal matrix, and we have the following dichotomy.  
\begin{enumerate}
\item If $M$ is of block diagonal matrix, then it is one of the following forms:
\begin{equation}\label{25112201}
A_1 \oplus A_4,\quad A_1 \oplus (-I)
\quad \text{or}\quad
(-I)\oplus A_4
\end{equation}
where $A_1, A_4\in \mathrm{GL}_2(\mathbb{Q})$ and $\mathfrak{m}_{A_1}(x)=\mathfrak{m}_{A_4}(x)=x^2+x+1$. 

\item Otherwise, $M$ is one of the following:  
\begin{equation*}
\left(\begin{array}{cc}
-\frac{I}{2} & A_2\\
-\frac{3}{4}A_2^{-1} & -\frac{I}{2}
\end{array}
\right),\quad  
\left(\begin{array}{cc}
A_1 & A_2\\
-\frac{A_2^{-1}}{2} & \frac{A_2^{-1}A_1^{-1}A_2}{2}
\end{array}\right) \quad \text{or}\quad
\left(\begin{array}{cc}
A_1 & A_2\\
-\frac{A_2^{-1}}{4} & \frac{3A_2^{-1}A_1^{-1}A_2}{4}
\end{array}\right)
\end{equation*}
for some $A_1, A_2\in \mathrm{GL}_2(\mathbb{Q})$ satisfying $\det A_2=\frac{3}{4}$, $\mathfrak{m}_{A_1}(x)=x^2+x+\frac{1}{2}$ or $\mathfrak{m}_{A_1}(x)=x^2+x+\frac{3}{4}$ respectively.
\end{enumerate}
\end{corollary}

\begin{example}\label{25110402}
{\normalfont Returning to Example \ref{25101509}, if $H$ is now as in \eqref{25100901}, then, by \eqref{22022007}, its coefficient matrix is
\begin{equation}\label{25110405}
\left(\begin{array}{cccc}
\frac{1}{2} & \frac{1}{2} & -\frac{1}{2}  & 0\\
-1 & \frac{1}{2} & 0 & -\frac{1}{2}\\
\frac{1}{2} & 0 & \frac{1}{2} & -\frac{1}{2}\\
0 & \frac{1}{2} & 1 & \frac{1}{2}
\end{array}\right). 
\end{equation}
In this case, $A_1$ and $A_2$ are   
\begin{equation*}
\left(\begin{array}{cc}
\frac{1}{2} & \frac{1}{2} \\
-1 & \frac{1}{2} 
\end{array}\right)\;\text{and}\;
\left(\begin{array}{cc}
-\frac{1}{2}  & 0\\
0 & -\frac{1}{2}
\end{array}\right), 
\end{equation*}
and so the matrix in \eqref{25110405} is of the third form given in \eqref{25011011}. }
\end{example}

\subsection{$x^3\pm 1$}\label{23122001}

In this subsection, we consider $M$ such that $\mathfrak{m}_M(x)$ is either $x^3-1$ or $x^3+1$. For the remainder of this section, we restrict $M$ to Type I, II or III, rather than allowing it to be arbitrary. The basic strategy of the proof is similar to those used in the previous cases. Once the result is established for $\mathfrak{m}_M(x)=x^3-1$, the other case follows immediately from the fact that $\mathfrak{m}_{-M}(x)=x^3+1$. 

\begin{theorem}\label{24010907}
Let $\mathcal{X}$ be an element in $\mathfrak{Ghol}$ and $M$ be an element in $\mathrm{Aut}\,\mathcal{X}$ such that $\mathfrak{m}_M(x)=x^3-1$ (resp. $x^3+1$). Then $M$ is not of Type II, and we have the following dichotomy.  
\begin{enumerate}
\item If $M$ is of Type I, then it is either 
\begin{equation*}
A_1 \oplus I 
\quad \text{or}\quad  
I \oplus A_4\quad (\text{resp. }
A_1 \oplus (-I)\quad \text{or}\quad  
(-I) \oplus A_4)
\end{equation*} 
for some $A_1, A_4\in \mathrm{GL}_2(\mathbb{Q})$ satisfying $\mathfrak{m}_{A_1}(x)=\mathfrak{m}_{A_4}(x)=x^2+x+1$ (resp. $x^2-x+1$). 

\item If $M$ is of Type III, then it is of the form 
\begin{equation}\label{25031101}
\left(\begin{array}{cc}
A_1 & A_2\\
-3A_2^{-1}A_1^2 & A_2^{-1}A_1A_2
\end{array}\right)
\end{equation}
for some $A_1, A_2\in \mathrm{GL}_2(\mathbb{Q})$ satisfying $\mathfrak{m}_{A_1}(x)=x^2-\frac{1}{2}x+\frac{1}{4}$ (resp. $x^2+\frac{1}{2}x+\frac{1}{4}$).
\end{enumerate}
\end{theorem}

\begin{proof}
First, given the assumption $\mathfrak{m}_{M}=x^3-1$, it is clear that $M$ is not of Type II. If $M$ of Type I and so given as $A_1\oplus A_4$, we get $A_1^3=A_4^3=I$. Since $A_1, A_4\in \mathrm{GL}_2(\mathbb{Q})$, it follows that each of $\mathfrak{m}_{A_1}(x)$ and $\mathfrak{m}_{A_4}(x)$ is either $x^2+x+1$ or $x-1$. However, if $\mathfrak{m}_{A_1}(x)=\mathfrak{m}_{A_4}(x)$, it contradicts the fact that the minimal polynomial of $M$ is $x^3-1$. 

Now we assume $M$ is of Type III. Since $M^2=M^{-1}$, we get 
\begin{equation*}
\left(\begin{array}{cc}
A_1^2+A_2A_3 & A_1A_2+A_2A_4 \\
A_3A_1+A_4A_3 & A_3A_2+A_4^2 
\end{array}\right)=
\left(\begin{array}{cc}
(\det A_1) A_1^{-1} & (\det A_2)A_3^{-1} \\
(\det A_2)A_2^{-1} & -(\det A_2)A_2^{-1}A_1A_3^{-1}
\end{array}\right)
\end{equation*}
by \eqref{23070519}, and thus 
\begin{equation*}
\begin{aligned}
A_1A_2+A_2A_4=(\det A_2)A_3^{-1} \Longrightarrow 
A_1A_2A_3-\frac{\det A_1}{\det A_2}A_2A_3A_1^{-1}A_2A_3=(\det A_2)I. 
\end{aligned}
\end{equation*}
Combining it with $A_2A_3=-A_1^2+(\det A_1) A_1^{-1}$, we further get 
\begin{equation*}
A_1^6+\big(1-3(\det A_1)\big) A_1^3+(\det A_1)^3I =0.
\end{equation*}
Thus $\mathfrak{m}_{A_1}(x)$ divides 
\begin{equation}\label{24112301}
x^6+\big(1-3(\det A_1)\big) x^3+(\det A_1)^3
\end{equation}
and, by elementary means, one can check \eqref{24112301} is factored as\footnote{For $x^2+(t-1)x+\frac{t^2-t+1}{3}$, by letting $m:=1-t$, we obtain $x^2-mx+\frac{m^2-m+1}{3}$. Thus the second polynomial is equivalent to the third one.} 
\begin{equation*}
\begin{aligned}
(x^2+x+\det A_1)(x^2-tx+\det A_1)(x^2+(t-1)x+\det A_1)
\end{aligned}
\end{equation*}
for some $t\in \mathbb{Q}$ satisfying $\det A_1=\frac{t^2-t+1}{3}$. To simplify the problem, by changing variables if necessary, let us normalize $M$ as  
\begin{equation}\label{24021902}
\left(\begin{array}{cc}
A_1 & I \\
A_2A_3 & A_2A_4A_2^{-1}
 \end{array}\right)=\left(\begin{array}{cc}
A_1 & I \\
A_2A_3 & -\frac{\det A_1}{\det A_2}A_2A_3A_1^{-1}
 \end{array}\right).  
\end{equation}

First, if $\mathfrak{m}_{A_1}(x)\;|\;x^2+x+\det A_1$, then $A_2A_3=(\det A_1-1)I$ and one finds $M^2+M+I=0$. But this contradicts the fact that $\mathfrak{m}_M(x)=x^3-1$. Hence $\mathfrak{m}_{A_1}(x)\;|\;x^2-tx+\frac{t^2-t+1}{3}$ for some $t\in \mathbb{Q}$. If $\mathfrak{m}_{A_1}(x)\neq x^2-tx+\frac{t^2-t+1}{3}$ (equivalently if $A_1\parallel I$), then $\frac{t^2-t+1}{3}=\frac{t^2}{4}$ and so $t=2$. However, this contradicts the assumption that $\det A_1\neq 1$. Thus $\mathfrak{m}_{A_1}(x)=x^2-tx+\frac{t^2-t+1}{3}$.  

Since 
\begin{equation*}
\begin{aligned}
A_2A_3=-A_1^2+(\det A_1)A_1^{-1}=(-t-1)A_1+\frac{t^2+2t+1}{3}I, 
\end{aligned}
\end{equation*}
\eqref{24021902} is equal to 
\begin{equation*}
\begin{gathered}
\left(\begin{array}{cc}
A_1 & I \\
(-t-1)A_1+\frac{t^2+2t+1}{3}I & \frac{t^2-t+1}{t^2-t-2}\big((-t-1)I+\frac{t^2+2t+1}{3}A_1^{-1}\big)\end{array}\right). 
\end{gathered}
\end{equation*} 
Thus the primary matrix $P_M\iota$ associated to $M$ is
\begin{equation*} 
\left(\begin{array}{cc}
\frac{t+\sqrt{\frac{-1}{3}}(2-t)}{2} & -1\\
\frac{t^2-t-2}{3}\Big(\frac{-1+\sqrt{-3}}{2}\Big) & \frac{t-1}{2}-\frac{(t+1)\sqrt{-3}}{6}
\end{array}\right)
\end{equation*}
where
\begin{equation}\label{24120201}
\mathfrak{\chi}_{P_M\iota}(x)=x^2-\big(t-\frac{1}{2}+\frac{1-2t}{6}\sqrt{-3}\big)x+\frac{1-\sqrt{-3}}{2}. 
\end{equation}
By Lemma \ref{25011101}, the solutions of \eqref{24120201} are roots of unity. If their algebraic degrees are at most $2$, then $t$ is either $2$ or $-1$; otherwise, $\frac{1}{2}$ by Lemma \ref{24122301}. However, $t=2$ and $-1$ contradict the fact that $\det A_1\neq 1$. Consequently, $t=\frac{1}{2}, \mathfrak{m}_{A_1}(x)=x^2-\frac{1}{2}x+\frac{1}{4}$ and $M$, before normalization, is of the form given in \eqref{25031101}. 
\end{proof}

\subsection{$x^3\pm x^2+x\pm 1$}

Now we consider $M$ such that $\mathfrak{m}_M(x)$ is either $x^3-x^2+x-1$ or $x^3+x^2+x+1$. 

\begin{theorem}\label{25010813}
Let $\mathcal{X}$ be an element in $\mathfrak{Ghol}$ and $M$ be an element in $\mathrm{Aut}\,\mathcal{X}$ such that $\mathfrak{m}_M(x)=x^3-x^2+x-1$ (resp.  $x^3+x^2+x+1$). Then $M$ is not of Type II, and we have the following dichotomy.  
\begin{enumerate}
\item If $M$ is of Type I, then $M$ is either 
\begin{equation*}
A_1 \oplus I\;\; \text{or}\;\;
I\oplus A_4\;\; 
(\text{resp. }
A_1 \oplus (-I)\;\; \text{or}\;\; 
(-I) \oplus A_4)
\end{equation*}
for some $A_1, A_4\in \mathrm{GL}_2(\mathbb{Q})$ satisfying $\mathfrak{m}_{A_1}(x)=\mathfrak{m}_{A_4}(x)=x^2+1$ (res. $x^2-1$).
\item If $M$ is of Type III, then $M$ is of the following form 
\begin{equation}\label{24120207}
\left(\begin{array}{cc}
A_1 & A_2\\
-A_2^{-1}A_1^2 & A_2^{-1}A_1A_2
\end{array}\right)
\end{equation}
where $\mathfrak{m}_{A_1}(x)=x^2-x+\frac{1}{2}$ (resp. $x^2+x+\frac{1}{2}$).
\end{enumerate}
\end{theorem}

Similar to the previous case, we only consider $M$ with $\mathfrak{m}_{M}(x)=x^3-x^2+x-1$, as the other case follows straightforwardly. 

\begin{proof}
For $M$ of Type I, the claim is immediate. For $M$ of Type II, since  $M^3+M$ (resp. $M^2+I$) is of Type II (resp. Type I), it follows that $M^3+M=M^2+I=0$, contradicting the assumption that $\mathfrak{m}_{M}(x)=x^3+x^2+x+1$. 

Now we assume $M$ is of Type III. Since $M^2+I=M+M^{-1}$, we get 
\begin{equation}\label{25120301}
\left(\begin{array}{cc}
A_1^2+A_2A_3+I & A_1A_2+A_2A_4 \\
A_3A_1+A_4A_3 & A_3A_2+A_4^2+I
\end{array}\right)=
\left(\begin{array}{cc}
A_1+(\det A_1) A_1^{-1} & A_2+(\det A_2)A_3^{-1} \\
A_3+(\det A_2)A_2^{-1} & A_4-(\det A_2)A_2^{-1}A_1A_3^{-1}
\end{array}\right),
\end{equation}
which, combinined with \eqref{23070519}, implies
\begin{equation*}
\begin{aligned}
A_1A_2+A_2A_4&=A_2+(\det A_2)A_3^{-1} \\
\Longrightarrow
A_1A_2A_3-\frac{\det A_1}{\det A_2}A_2A_3A_1^{-1}A_2A_3&=A_2A_3+(\det A_2)I.
\end{aligned}
\end{equation*}
If $t:=\mathrm{tr}\,A_1$ and $a:=\det A_1$, since $A_2A_3=-A_1^2-I+aA_1^{-1}$ (by \eqref{25120301}), the above equation is simplified as
\begin{equation}\label{25011405}
\begin{aligned}
A_1^6-2A_1^5+(2+a)A_1^4-4aA_1^3+(a^2+2a)A_1^2-2a^2 A_1+a^3 I&=0 \\
\Longrightarrow (t^2-2t)A^2_1+(2t-2ta)A_1+(-2+2a)aI&=0. 
\end{aligned}
\end{equation}
If $A_1\parallel I$, then $A_1$ is either $I$ or $0$, contradicting the assumption. If $A_1\notparallel I$ and $\mathfrak{m}_{A_1}(x)=x^2-tx+a$, the last equation in \eqref{25011405} is equivalent to $A_1^2-tA_1+aI$, implying
\begin{equation*}
\frac{2t-2ta}{t^2-2t}=-t \Longrightarrow a=\frac{t^2-2t+2}{2}.
\end{equation*}
To simplify the problem, by changing variables if necessary, let us normalize $M$ as  
\begin{equation}\label{24112203}
\left(\begin{array}{cc}
A_1 & I \\
A_2A_3 & A_2A_4A_2^{-1}
 \end{array}\right)=\left(\begin{array}{cc}
A_1 & I \\
A_2A_3 & -\frac{\det A_1}{\det A_2}A_2A_3A_1^{-1}
 \end{array}\right).  
\end{equation}
Since $\mathfrak{m}_{A_1}(x)=x^2-tx+\frac{t^2-2t+2}{2}$, 
\begin{equation*}
\begin{aligned}
A_2A_3=-A_1^2+A_1-I+(\det A_1) A_1^{-1}=-tA_1+\frac{t^2}{2}I, 
\end{aligned}
\end{equation*}
and thus $M$ and $P_M$ are 
\begin{equation}\label{24112201}
\begin{aligned}
\left(\begin{array}{cc}
A_1 & I \\
-tA_1+\frac{t^2}{2}I &  -\frac{t}{t-2}A_1+\frac{2t-2}{t-2}I
\end{array}\right)\quad \text{and}\quad 
\left(\begin{array}{cc}
\frac{t+(2-t)\sqrt{-1}}{2} & 1\\
\frac{t(t-2)\sqrt{-1}}{2} &  \frac{2-t}{2}+\frac{t\sqrt{-1}}{2}
\end{array}\right)
\end{aligned}
\end{equation} 
respectively. Note that, for $P_M$ given as in \eqref{24112201}, 
\begin{equation*}
\mathfrak{\chi}_{P_M\iota}(x)=x^2-\big((t-1)+(1-t)\sqrt{-1}\big)x-\sqrt{-1} 
\end{equation*}
and two eigenvalues of $P_M\iota $ are roots of unity (by Lemma \ref{25011101}). If the degree of the eigenvalues is $2$, then $t$ is either $2$ or $0$, contradicting the fact that $0<\det A_1< 1$. Otherwise, if the degree of the eigenvalues is $4$, then $t=1$ by Lemma \ref{24122301}, and so $M$ is of the form given in \eqref{24120207}. 

\end{proof}

\subsection{$x^3\pm 2x^2+2x\pm 1$}

The proof of the following theorem is very analogous to that of Theorem \ref{25010813}, so it will be omitted here.   
\begin{theorem}
Let $\mathcal{X}$ be an element in $\mathfrak{Ghol}$ and $M$ be an element in $\mathrm{Aut}\,\mathcal{X}$ such that $\mathfrak{m}_M(x)=x^3-2x^2+2x-1$ (resp.  $x^3+2x^2+2x+1$). Then $M$ is not of Type II, and we have the following dichotomy.  
\begin{enumerate}
\item If $M$ is of Type I, then $M$ is either 
\begin{equation*}
A_1 \oplus I\;\; \text{or}\;\; 
I\oplus A_4\;\; (\text{resp. }A_1 \oplus (-I)
\;\; \text{or}\;\; 
(-I) \oplus A_4)
\end{equation*}
for some $A_1, A_4\in \mathrm{GL}_2(\mathbb{Q})$ satisfying $\mathfrak{m}_{A_1}(x)=\mathfrak{m}_{A_1}(x)=x^2-x+1$ (resp. $x^2+x+1$).

\item If $M$ is of Type III, then $M$ is of the following form 
\begin{equation*}
\left(\begin{array}{cc}
A_1 & A_2\\
-\frac{A_2^{-1}A_1^2}{3} & A_2^{-1}A_1A_2
\end{array}\right)
\end{equation*}
where $\mathfrak{m}_{A_1}(x)=x^2-\frac{3}{2}x+\frac{3}{4}$ (resp. $x^2+\frac{3}{2}x+\frac{3}{4}$).
\end{enumerate}
\end{theorem}

\newpage
\section{Classification of $M$: II}\label{24081805}
In this section, we analyze the remaining cases from the previous section. Recall that what remains is the case where $M\in \mathrm{Aut}\,\mathcal{X}$ has minimal polynomial one of the following:
\begin{equation*}
x^4-1, \; x^4+1, \;  x^4+x^2+1, \; x^4-x^2+1. 
\end{equation*}
For such $M$, the minimal polynomial of $M^2$ is, respectively, 
\begin{equation*}
x^2-1, \; x^2+1, \;  x^2+x+1, \; x^2-x+1. 
\end{equation*}
Thus, our approach builds on the results from the previous section, and the problem is split into the following three cases based on the type of $M^2$:
\begin{itemize}
\item $M^2$ is of block diagonal matrix;
\item $M^2$ is of block anti-diagonal matrix;
\item $M^2$ is neither of block diagonal nor anti-diagonal matrix. 
\end{itemize}

The first two cases are relatively easy to handle. For instance, if \(M\) is of Type III and \(M^2\) is block diagonal, then the shape of \(M\) is generally of very rigid type, as shown in the following:

\begin{lemma}\label{24082308}
Let $M$ be a matrix in $\mathrm{Aut}\,\mathcal{X}$ where $\mathcal{X}\in \mathfrak{Ghol}$. If $M$ is of Type III and $M^2$ is of block diagonal matrix, then $M$ and $M^2$ are given as 
\begin{equation}\label{25021601}
\left(\begin{array}{cc}
A_1 & A_2\\
\frac{\det A_2}{\det A_1}A_2^{-1}A_1^2 & -A_2^{-1}A_1A_2
\end{array}\right)\quad \text{and}\quad 
\left(\begin{array}{cc}
\frac{1}{\det A_1}A^2_1 & 0\\
0 & \frac{1}{\det A_1}A_2^{-1}A_1^2A_2
\end{array}\right)
\end{equation}
for some $A_1, A_2\in \mathrm{GL}_2(\mathbb{Q})$ respectively. 
\end{lemma}
\begin{proof}  
By the assumption,  
\begin{equation*}
\left(\begin{array}{cc}
A_1 & A_2\\
A_3 & A_4
\end{array}\right)^2
=\left(\begin{array}{cc}
A^2_1+A_2A_3 & A_1A_2+A_2A_4\\
A_3A_1+A_4A_3 & A_3A_2+A_4^2
\end{array}\right)
=\left(\begin{array}{cc}
\tilde{A_1} & 0\\
0 & \tilde{A_4}
\end{array}\right) 
\end{equation*}
for some $\tilde{A_1}$ and $\tilde{A_4}$ in $\mathrm{GL}_2(\mathbb{Q})$. Since 
\begin{equation*}
A_4=-\frac{\det A_1}{\det A_2}A_3A_1^{-1}A_2\quad \text{and}\quad A_1A_2+A_2A_4=A_3A_1+A_4A_3=0, 
\end{equation*}
it follows that 
\begin{equation*}
\begin{gathered}
A_1A_2-\frac{\det A_1}{\det A_2}A_2A_3A_1^{-1}A_2=0 \Longrightarrow A_1^2=\frac{\det A_1}{\det A_2}A_2A_3. 
\end{gathered}
\end{equation*}
In conclusion, $M$ and $M^2$ are given as in \eqref{25021601}.  
\end{proof}

For $M^2$ as a block anti-diagonal matrix, the analysis is not difficult but we treat it separately, as the shape of $M$ varies slightly depending on its minimal polynomial. 

Hence, among the three cases listed above, what really matters is the last one. In this case, using the results from the previous section, we derive several equations involving $A_j$ ($1\leq j\leq 4$), and, by manipulating them, express the traces and determinants of $A_2A_3$ and $A_2A_4A_2^{-1}$ in terms of those of $A_1$. The exact values of the trace and determinant of $A_1$ are then determined by applying the results of Section \ref{25011407}, and these will ultimately allow us to describe $M$ explicitly in terms of $A_1$ and $A_2$.

In a broad sense, the general scheme of each proof is analogous to that of the theorems in the previous section, but the computations involved here are substantially more intricate and subtle. We present detailed computations as fully as possible; however, when a computation has already appeared and its repetition is routine, we omit it in order to maintain the overall clarity of exposition. 

As in the previous section, an impatient reader may skip this section on a first reading; for those who wish to see a proof, it suffices to go through the first one, Theorem \ref{24070302}.

Before proceeding further, we record the following lemma. This lemma will be repeatedly applied throughout the section, but its proof is quite elementary and therefore skipped here.

\begin{lemma}\label{24112505}
Let $f(x):=x^4+ax^2+b$ be a polynomial in $\mathbb{Q}(x)$. If $f(x)$ is reducible and has a factor of degree $2$, then it is factored into either 
\begin{equation*}
(x^2+mx+n)(x^2-mx+n)\quad \text{or}\quad (x^2+m)(x^2+n)
\end{equation*}
for some $m, n\in \mathbb{Q}$. 
\end{lemma}

\subsection{$x^4-1$}
In this subsection, we analyze $M$ such that $\mathfrak{m}_{M}(x)=x^4-1$. As $\mathfrak{m}_{M^2}(x)=x^2-1$, the proof of the theorem depends on Theorem \ref{23063002}. 
 
\begin{theorem}\label{24070302}
Let $M$ be a matrix in $\mathrm{Aut}\;\mathcal{X}$ (where $\mathcal{X}\in \mathfrak{Ghol}$) whose minimal polynomial is $x^4-1$. Then $M$ is of Type III and has the form 
\begin{equation}\label{25010701}
\left(\begin{array}{cc}
A_1 & A_2\\
-A_2^{-1}A_1^2 & A_2^{-1}A_1A_2
\end{array}\right)
\end{equation}
for some $A_1, A_2\in\mathrm{GL}_2(\mathbb{Q})$ such that $\mathfrak{m}_{A_1}(x)=x^2\pm x+\frac{1}{2}$. 
\end{theorem}
\begin{proof}

Suppose $M$ is of Type I, given by $A_1\oplus A_4$. If $\mathfrak{m}_{M}(x)=x^4-1$, then both $\mathfrak{m}_{A_1}(x)$ and $\mathfrak{m}_{A_4}(x)$ divide $x^4-1$, and thus the only possible consideration is either $\mathfrak{m}_{A_1}(x)=x^2-1$ and $\mathfrak{m}_{A_4}(x)=x^2+1$, or $\mathfrak{m}_{A_1}(x)=x^2+1$ and $\mathfrak{m}_{A_4}(x)=x^2-1$. However either case contradicts the fact that $\det A_1=\det A_4=1$. Similarly, one proves $M$ is not of Type II. 
 
Now we assume $M$ is of Type III. 
\begin{enumerate}
\item First, suppose $M^2$ is a block anti-diagonal form, that is, 
\begin{equation*}
\left(\begin{array}{cc}
A_1 & A_2\\
A_3 & A_4
\end{array}\right)^2
=\left(\begin{array}{cc}
A^2_1+A_2A_3 & A_1A_2+A_2A_4\\
A_3A_1+A_4A_3 & A_3A_2+A_4^2
\end{array}\right)
=\left(\begin{array}{cc}
0 & \tilde{A_2}\\
\tilde{A_2}^{-1} & 0
\end{array}\right)
\end{equation*}
for some $\tilde{A_2}$. Since 
\begin{equation}\label{24103001}
A_4^2=-A_3A_2, \quad A^2_1=-A_2A_3\quad\text{and}\quad  A_4=-\frac{\det A_1}{\det A_3}A_3A_1^{-1}A_2,
\end{equation}
it follows that $\det A_2=\det A_3$. As $(A_1A_2+A_2A_4)(A_3A_1+A_4A_3)=I$, we get
\begin{equation*}
\begin{gathered}
A_1A_2A_3A_1+A_2A_4A_3A_1+A_1A_2A_4A_3+A_2A_4A_4A_3=-4A_1^4=I. 
\end{gathered}
\end{equation*}
Thus $\mathfrak{m}_{A_1}(x)$ divides $x^4+\frac{1}{4}$, which factors as $(x^2+x+\frac{1}{2})(x^2-x+\frac{1}{2})$. By \eqref{24103001}, it is clear that $M$ is of the form given in \eqref{25010701}. 

\item Now we assume $M^2$ is a block diagonal form. By Lemma \ref{24082308}, $M$ and $M^2$ are of the forms
\begin{equation*}
\left(\begin{array}{cc}
A_1 & A_2\\
\frac{\det A_2}{\det A_1}A_2^{-1}A_1^2 & -A_2^{-1}A_1A_2
\end{array}\right)
\quad \text{and}\quad 
\left(\begin{array}{cc}
\frac{1}{\det A_1}A^2_1 & 0\\
0 & \frac{1}{\det A_1}A_2^{-1}A_1^2A_2
\end{array}\right)
\end{equation*}
respectively. Since $A^4_1-(\det A_1)^2I=0$, 
$\mathfrak{m}_{A_1}(x)$ divides either $x^2+\det A_1$ or $x^2-\det A_1$, implying either $M^2=-I$ or $M^2=I$ respectively. But this contradicts the assumption that $\mathfrak{m}_M(x)=x^4-1$.   

\item Lastly, if $M^2$ is neither diagonal nor anti-diagonal matrix, then it is either 
\begin{equation}\label{24120401}
\left(\begin{array}{cc}
\frac{I}{2} & \tilde{A_2}\\
\frac{3}{4}\tilde{A_2}^{-1} & -\frac{I}{2}
\end{array}\right)\quad \text{or}\quad 
\left(\begin{array}{cc}
-\frac{I}{2} & \tilde{A_2}\\
\frac{3}{4}\tilde{A_2}^{-1} & \frac{I}{2}
\end{array}\right)
\end{equation}
for some $\tilde{A_2}\in \mathrm{GL}_2(\mathbb{Q})$ satisfying $\det \tilde{A_2}=\frac{3}{4}$ by Theorem \ref{23063002}. Without loss of generality, we consider only the first case, as the second one can be treated similarly. 

Changing variables if necessary, we normalize $M^2$ and denote $M$ as\footnote{That is, $M^2=\left(\begin{array}{cc}
A_1 & A_2\\
A_3 & A_4
\end{array}\right)^2
=\left(\begin{array}{cc}
\frac{I}{2} & I\\
\frac{3}{4}I & -\frac{I}{2}
\end{array}\right)$.}
\begin{equation*}
\left(\begin{array}{cc}
\frac{I}{2} & I\\
\frac{3}{4}I & -\frac{I}{2}
\end{array}\right)\quad \text{and}\quad \left(\begin{array}{cc}
A_1 & A_2\\
A_3 & A_4
\end{array}\right)
\end{equation*}
respectively. Then\footnote{The original $M$, prior to normalization, is $
\left(\begin{array}{cc}
A_1 & A_2\tilde{A_2}\\
\tilde{A_2}^{-1}A_3 & A_4
\end{array}\right)$ and, since it is of Type III, we get the equalities in \eqref{24051602}.} 
\begin{equation}\label{24051602}
\frac{3}{4}\det A_2=1-\det A_1, \quad \det A_3=\frac{3}{4}(1-\det A_1), \quad A_4=-\frac{\det A_1}{1-\det A_1}A_3A_1^{-1}A_2
\end{equation}
and 
\begin{equation}\label{23072401}
A_1^2+A_2A_3=\frac{1}{2}I, \quad  A_1A_2+A_2A_4=I, \quad A_3A_1+A_4A_3=\frac{3}{4}I, \quad A_3A_2+A_4^2=-\frac{1}{2}I.
\end{equation}
We first claim
\begin{claim}\label{24082401}
\begin{enumerate}
\item $A_2^{-1}=\frac{3}{4}A_3^{-1}$.
\item $\det A_1=\frac{5}{8}$ and $\mathrm{tr}\,A_1=\pm\frac{3}{2}$. 
\item If $\mathrm{tr}\,A_1=-\frac{3}{2}$, then $\det A_2=\frac{1}{2}$ and $\mathrm{tr}\,A_2=-1$.
\end{enumerate}
\end{claim}
\begin{proof}[Proof of Claim \ref{24082401}]
\begin{enumerate}

\item By the first equation in  \eqref{23072401}, 
\begin{equation}\label{25120401}
A_2=\frac{1}{2}A_3^{-1}-A_1^2A_3^{-1}.
\end{equation}
Applying this together with the last equation in \eqref{24051602} to the second and third equations in \eqref{23072401}, we get
\begin{equation}\label{24051606}
\begin{aligned}
A_1A_2+(\frac{1}{2}A_3^{-1}-A_1^2A_3^{-1})(-\frac{\det A_1}{1-\det A_1}A_3A_1^{-1}A_2)&=I\\
\Longrightarrow 
\frac{-\frac{1}{2}(\det A_1)A_1^{-1}+A_1}{1-\det A_1}&=A_2^{-1}
\end{aligned}
\end{equation}
and 
\begin{equation}\label{25120403}
\begin{aligned}
A_3A_1-\frac{\det A_1}{1-\det A_1}A_3A_1^{-1}(\frac{1}{2}A_3^{-1}-A_1^2A_3^{-1})A_3&=\frac{3}{4}I\\
\Longrightarrow 
\frac{-\frac{1}{2}(\det A_1)A_1^{-1}+A_1}{1-\det A_1}&=\frac{3}{4}A_3^{-1}
\end{aligned}
\end{equation}
respectively. The conclusion of the first claim is immediate from \eqref{24051606} and \eqref{25120403}. 

\item Applying \eqref{25120401} to the last one in \eqref{23072401}, we get
\begin{equation}\label{24051601}
\begin{aligned}
\frac{1}{2}I-A_3A_1^2A_3^{-1}+\frac{(\det A_1)^2}{(1-\det A_1)^2}A_3(\frac{1}{2}A_1^{-2}-I)(\frac{1}{2}I-A_1^2)A_3^{-1}&=-\frac{1}{2}I \\
\Longrightarrow A_1^4-A_1^2-\frac{(\det A_1)^2}{4(1-2\det A_1)}I&=0.
\end{aligned}
\end{equation}
If $A_1\parallel I$, then one can check there is no $a\in \mathbb{Q}$ satisfying both $A_1=aI$ and \eqref{24051601}. If $A_1\notparallel I$ and $\mathfrak{\chi}_{A_1}(x)=\mathfrak{m}_{A_1}(x)$, by Lemma \ref{24112505}, the equation \eqref{24051601} is factored either as 
\begin{equation*}
\big(A_1^2+(\det A_1)I\big)\big(A_1^2-\frac{\det A_1}{4(1-2\det A_1)}I\big)=0
\end{equation*}
or 
\begin{equation}\label{24112509}
\big(A_1^2-\mathrm{tr}\,A_1+(\det A_1)I\big)\big(A_1^2+\mathrm{tr}\,A_1+(\det A_1)I\big)=0.
\end{equation}

For the first case, by \eqref{24051602}, \eqref{24051606} and \eqref{25120403}, one gets all $A_i$ $(2\leq i\leq 4)$ are scalar multiples of $A_1$. If $A_i=a_iA_1$ for some $a_i\in \mathbb{Q}$ ($2\leq i\leq 4$), then
\begin{equation}\label{24112507}
\begin{aligned}
\frac{3}{4}a_2^2\det A_1=1-\det A_1,\; a_3^2\det A_1=\frac{3}{4}(1-\det A_1)\Longrightarrow a_2a_3=\frac{1-\det A_1}{\det A_1}
\end{aligned}
\end{equation}  
by the first two equalities in \eqref{24051602}, and 
\begin{equation}\label{25021602}
A_1^2+A_2A_3=A_1^2+a_2a_3A_1^2=-(\det A_1)(1+a_2a_3)I =\frac{1}{2}I
\end{equation}
by the first equality in \eqref{23072401}. By \eqref{24112507}, it follows that 
\begin{equation*}
-(\det A_1)(1+a_2a_3)=-(\det A_1)(1+\frac{1-\det A_1}{\det A_1})=-1,  
\end{equation*}
which contradicts \eqref{25021602}. Consequently, the last equation in \eqref{24051601} is factored as in \eqref{24112509}, and solving 
$$-\frac{(\det A_1)^2}{4(1-2\det A_1)}I=(\det A_1)^2I$$ yields the conclusion of the second claim. 

\item Now we assume $\mathrm{tr}\,A_1=-\frac{3}{2}$, and find the trace and determinant of $A_2$. 

Combining \eqref{24051606} with $A_1+\frac{3}{2}I+\frac{5}{8}A_1^{-1}=0$, it follows that  
\begin{equation}\label{24051609}
\begin{gathered}
4(A_1+\frac{1}{2}I)=A_2^{-1}, 
\end{gathered}
\end{equation}
so 
\begin{equation*}
\mathrm{tr}\,A_2^{-1}=\mathrm{tr}\,\big(4(A_1+\frac{1}{2}I)\big)=-2,\;\;\frac{1}{\det A_2}=16\det(A_1+\frac{1}{2}I)=2\;\;\text{and}\;\; \mathrm{tr}\,A_2=-1.  
\end{equation*} 
\end{enumerate}
This completes the proof of the claim. 
\end{proof}

To simplify the computation, we further change basis of $M$ again and normalize it as 
\begin{equation*}
\left(\begin{array}{cc}
A_1 & I \\
A_2A_3 & A_2A_4A_2^{-1}
\end{array}\right),
\end{equation*}
which, by \eqref{24051602} and the above claim, equals 
\begin{equation*}
\left(\begin{array}{cc}
A_1 & I \\
\frac{3}{4}A_2^2 & -\frac{5}{4}A_2^2A_1^{-1}
\end{array}\right). 
\end{equation*}
Without loss of generality, we assume $\mathrm{tr}\,A_1=-\frac{3}{2}$ and find the primary matrix $P_M$ associated to $M$. Since 
\begin{equation*}
\mathrm{tr}\,A_2^2=\mathrm{tr}\,(-A_2-\frac{I}{2})=0,\; \; \det A_2^2=\frac{1}{4}\;\; \text{and} \;\; 4(A_1^2+\frac{1}{2}A_1)=A_1A_2^{-1} 
\end{equation*}
by Claim \ref{24082401} and \eqref{24051609},   
\begin{equation*}
\begin{aligned}
\mathrm{tr}\,(A_1A_2^{-1})=&\mathrm{tr}\,\big(4(A_1^2+\frac{1}{2}A_1)\big)=4\mathrm{tr}\,(-A_1-\frac{5}{8}I)=1
\end{aligned}
\end{equation*} 
and so
\begin{equation*}
\begin{aligned}
\mathrm{tr}\,(A_2^2A_1^{-1})=\mathrm{tr}\,\big((-A_2-\frac{I}{2})A_1^{-1}\big)=-\frac{\mathrm{tr}\,(A_1A_2^{-1})}{\det (A_1A_2^{-1})}+\frac{6}{5}=\frac{2}{5}. 
\end{aligned}
\end{equation*}
Finally, the eigenvalues of $A_2$ and $A_2^2A_1^{-1}$ are $\frac{\pm\sqrt{-1}}{2}$ and $\frac{1\pm3\sqrt{-1}}{5}$ respectively, hence the primary matrix $P_M$ associated to $M$ is\footnote{The signs in the formula are independent. That is, each $\pm $ can be chosen separately. The same convention will apply to analogous formulas in the subsequent sections.} 

\begin{equation*}
\left(\begin{array}{cc}
\frac{-3\pm\sqrt{-1}}{4} & 1 \\
\frac{\pm 3\sqrt{-1}}{8} & \frac{-1\pm3\sqrt{-1}}{4}
\end{array}\right).
\end{equation*}
If $P_M$ is, for instance, given as 
\begin{equation*} 
\left(\begin{array}{cc}
\frac{-3+\sqrt{-1}}{4} & 1 \\
\frac{-3\sqrt{-1}}{8} & \frac{-1+3\sqrt{-1}}{4}
\end{array}\right),  
\end{equation*}
then, since two nonzero eigenvalues of $P_M$ are distinct each other and $\mathrm{tr}\,P_M\iota\neq 0$, it follows the eigenvalues are all roots of unity by Proposition \ref{24082301}. However, this is a contradiction as $\det P_M=-\frac{\sqrt{-1}}{4}$ is not a root of unity. Similarly, one derives a contradiction for each of the remaining cases. 

In conclusion, there is no $M$ of Type III such that $M^2$ is of the form given in \eqref{24120401}. 
\end{enumerate}
\end{proof}

\subsection{$x^4+1$}
In this subsection, we consider $M$ satisfying $\mathfrak{m}_{M}(x)=x^4+1$. We follow the same scheme as in the previous subsection. 

\begin{theorem}\label{24092502}
Let $M$ be an element in $\mathrm{Aut}\;\mathcal{X}$ (where $\mathcal{X}\in \mathfrak{Ghol}$) such that $\mathfrak{m}_M(x)=x^4+1$. Then $M$ is not of Type I, and the following statements hold.  
\begin{enumerate}
\item If $M$ is of Type II, then it is of the following form 
\begin{equation*}
A_2\widetilde{\oplus}A_3 
\end{equation*}
for some $A_2, A_3\in \mathrm{GL}_2(\mathbb{Q})$ satisfying $\mathfrak{m}_{A_2A_3}(x)=x^2+1$.

\item If $M$ is of Type III, then $M$ is one of the following forms:
\begin{equation}\label{24010303}
\left(\begin{array}{cc}
A_1 & A_2\\
A_2^{-1}A_1^2 & -A_2^{-1}A_1A_2
\end{array}\right), \;\; 
\left(\begin{array}{cc}
A_1 & A_2\\
\frac{A_2^{-1}}{2} & A_2^{-1}A_1A_2
\end{array}\right) \;\; \text{or}\;\;
\left(\begin{array}{cc}
A_1 & A_2\\
\frac{A_2^{-1}}{4} & -\frac{3A_2^{-1}A_1^{-1}A_2}{4}
\end{array}\right) 
\end{equation}
for some $A_1, A_2\in \mathrm{GL}_2(\mathbb{Q})$ satisfying $\mathfrak{m}_{A_1}(x)$ is $x^2\pm x+\frac{1}{2}$, $x^2+\frac{1}{2}$ or $x^2\pm x+\frac{3}{4}$ respectively. In particular, $M^2$ is of Type I, Type II, or Type III respectively, and $\tau_i\in \mathbb{Q}(\sqrt{-1})$ in the first case, while $\mathbb{Q}(\sqrt{-2})$ in the second and third cases, for $i=1,2$.
\end{enumerate}
\end{theorem}
\begin{proof}
If $M$ is of Type I, given as $A_1 \oplus A_4$, then both $\mathfrak{m}_{A_1}(x)$ and $\mathfrak{m}_{A_4}(x)$, which are polynomials of degree at most $2$, divide $x^4+1$, leading to a contradiction. When $M$ is of Type II, it is clear that $M$ must be of the form presented above. 

Now we assume $M$ is of Type III. If $M^2$ is a block anti-diagonal form, then, following a similar procedure to that in the proof of Theorem \ref{24070302}, one finds $4A_1^4-I=0$ and thus $M$ is of the second form among given in \eqref{24010303}. 
\begin{enumerate}
\item Suppose $M^2$ is a block diagonal matrix. By Lemma \ref{24082308}, $M$ and $M^2$ are of the forms 
\begin{equation*}
\left(\begin{array}{cc}
A_1 & A_2\\
\frac{\det A_2}{\det A_1}A_2^{-1}A_1^2 & -A_2^{-1}A_1A_2
\end{array}\right)
\quad \text{and}\quad \left(\begin{array}{cc}
\frac{1}{\det A_1}A^2_1 & 0\\
0 & \frac{1}{\det A_1}A_2^{-1}A_1^2A_2
\end{array}\right)
\end{equation*}
respectively. Since $A^4_1+(\det A_1)^2I=0$, it follows that $\mathfrak{m}_{A_1}(x)$ is either $x^2+tx+\frac{t^2}{2}$ or $x^2-tx+\frac{t^2}{2}$ for some $t\in \mathbb{Q}$. Without loss of generality, let us assume $t>0$ and $\mathfrak{m}_{A_1}(x)=x^2-tx+\frac{t^2}{2}$. Changing variables if necessary, if $M$ is normalized as 
\begin{equation*}
\left(\begin{array}{cc}
A_1 & I\\
\frac{2-t^2}{t^2}A_1^2 & -A_1
\end{array}\right),  
\end{equation*}
then the primary matrix $P_M$ associated with $M$ is 
\begin{equation*}
\left(\begin{array}{cc}
t\omega & 1\\
(2-t^2)\omega^2 & -t\omega
\end{array}\right) 
\end{equation*}
where $\omega=\frac{1\pm\sqrt{-1}}{2}$. By Lemma \ref{25011101}, the order of $P_M\iota$ is finite, and thus the eigenvalues of $P_M\iota$, which are the solutions of 
\begin{equation*}
x^2-2t\omega x+2\omega^2=0, 
\end{equation*}
are roots of unity. As $2\omega^2=\pm\sqrt{-1}$, if the degree of the eigenvalues is $4$, then $t=0$ by Lemma \ref{24122301}. But this contradicts the fact that $\det A_1\neq 0$. Otherwise, if the degrees of the eigenvalues are at most $2$, then one can check $t=1$ is the desired one, which induces the original $M$, before normalization, being as the first type given in \eqref{24010303}. 

\item  Assume $M^2$ is neither diagonal nor anti-diagonal. Since $M^2\in \mathrm{Aut}\,\mathcal{X}$ and $\mathfrak{m}_{M^2}(x)=x^2+1$, by Theorem \ref{23070521}, $M^2$ is a matrix of the following form: 
\begin{equation}\label{241204011}
\left(\begin{array}{cc}
\tilde{A_1} & \tilde{A_2}\\
-(\det \tilde{A_2})\tilde{A_2}^{-1} & -\tilde{A_2}^{-1}\tilde{A_1}\tilde{A_2}
\end{array}\right)
\end{equation}
where $\det \tilde{A_1}+\det \tilde{A_2}=1$ and $\mathfrak{m}_{\tilde{A_1}}(x)$ is one of the following: $x^2+\frac{1}{4}, x^2+\frac{1}{2}$ or $x^2+\frac{3}{4}$.\footnote{Note that, if $\mathfrak{m}_{\tilde{A_1}}(x)$ is $x^2+\frac{1}{4}, x^2+\frac{1}{2}$ or $x^2+\frac{3}{4}$, then the cusp shapes of $\mathcal{X}$ are contained in $\mathbb{Q}(\sqrt{-1}), \mathbb{Q}(\sqrt{-2})$ or $\mathbb{Q}(\sqrt{-3})$ respectively. We will use this fact later in the proof. }

Changing variables if necessary, we normalize $M^2$ and denote $M$ by 
\begin{equation*}
\left(\begin{array}{cc}
\tilde{A_1} & I\\
-(\det \tilde{A_2})I & -\tilde{A_1}
\end{array}
\right)\quad \text{and}\quad \left(\begin{array}{cc}
A_1 & A_2\\
A_3 & A_4\\
\end{array}\right)
\end{equation*}
respectively. Then  
\begin{equation}\label{23072703}
A_1^2+A_2A_3=\tilde{A_1}, \;\; A_1A_2+A_2A_4=I, \;\; A_3A_1+A_4A_3=-(\det \tilde{A_2})I, \;\;  A_3A_2+A_4^2=-\tilde{A_1},
\end{equation}
and\footnote{As in \eqref{24051602}, the original $M$ is of the form $\left(\begin{array}{cc}
A_1 & A_2 \tilde{A_2}\\
 \tilde{A_2}^{-1}A_3 & A_4\\
\end{array}\right)$ (before normalization) and, since it is of Type III, we get the equalities in \eqref{24020605}.}  
\begin{equation}\label{24020605}
\det A_1=\det A_4, \;\; \det A_1+(\det A_2)(\det \tilde{A_2})=\det A_1+\frac{\det A_3}{\det \tilde{A_2}}=1, \;\; A_4=-\frac{\det A_1}{1-\det A_1}A_3A_1^{-1}A_2.
\end{equation}
We address the problem by dividing it into two cases depending on whether $A_2\parallel I$ or $A_2\notparallel I$. 

\begin{enumerate}
\item Suppose $A_2\notparallel I$. We derive a couple of equalities using \eqref{23072703}-\eqref{24020605}.  
\begin{claim}\label{24020607}
\begin{enumerate}
\item $\mathrm{tr}\,A_2=0, \det A_2=-\frac{1}{\mathrm{tr}\;A_2^{-1}A_1}=\frac{1-\det A_1}{\det \tilde{A_2}}$ and $\mathrm{tr}\;A_1=-\mathrm{tr}\;A_4$. 
\item $(\mathrm{tr}\;A_1)^2=2\det A_1-\det \tilde{A_2}$ and $\mathrm{tr}\, (A_2A_3)= \det \tilde{A_2}$.
\end{enumerate} 
\end{claim}
\begin{proof}[Proof of Claim \ref{24020607}]
\begin{enumerate}
\item Combining the second equation in \eqref{23072703} with the last one in \eqref{24020605}, we get
\begin{equation}\label{24020609}
\begin{aligned}
A_1A_2+A_2(-\frac{\det A_1}{1-\det A_1}A_3A_1^{-1}A_2)&=I\\
\Longrightarrow 
A_3=-\frac{1-\det A_1}{\det A_1}A_2^{-1}(I-A_1A_2)A_2^{-1}A_1&=-\frac{1-\det A_1}{\det A_1}(A_2^{-2}A_1-A_2^{-1}A_1^2),
\end{aligned}
\end{equation}
thus
\begin{equation}\label{24020610}
A_4=-\frac{\det A_1}{1-\det A_1}A_3A_1^{-1}A_2=A_2^{-1}-A_2^{-1}A_1A_2.
\end{equation}
Plugging \eqref{24020609}-\eqref{24020610} into the third equation in \eqref{23072703}, 
\begin{equation}\label{25112601}
A_2^{-3}A_1-(A_2^{-1}A_1)^2=\frac{\det \tilde{A_2}\det A_1}{1-\det A_1}I.
\end{equation}
Since 
\begin{equation*}
(A_2^{-1}A_1)^2-(\mathrm{tr}\,A_2^{-1}A_1)A_2^{-1}A_1+(\det A_2^{-1}A_1)I=0
\end{equation*}
and 
\begin{equation*}
\det A_2^{-1}A_1=\frac{\det \tilde{A_2}\det A_1}{1-\det A_1}\;(\text{by }\eqref{24020605}), 
\end{equation*}
combining them with \eqref{25112601}, it follows that 
\begin{equation}\label{24021504}
\begin{gathered}
A_2^{-3}A_1=(\mathrm{tr}\,A_2^{-1}A_1)A_2^{-1}A_1\Longrightarrow 
A_2^{-2}=(\mathrm{tr}\,A_2^{-1}A_1)I. 
\end{gathered}
\end{equation}
As $A_2\notparallel I$ by the assumption, we conclude 
\begin{equation*}
\begin{gathered}
\mathrm{tr}\,A_2=0, \;\; \det A_2=-\frac{1}{\mathrm{tr}\,A_2^{-1}A_1}=\frac{1-\det A_1}{\det \tilde{A_2}}, 
\end{gathered}
\end{equation*}
which are the first two equalities in the claim. The third one follows from $\mathrm{tr}\,A_2=0$ and \eqref{24020610}. 

\item To simplify notation, let $\tilde{a}:= \det \tilde{A_2}$ and $a:=\det A_1$. By \eqref{23072703},  
\begin{equation}\label{25120405}
A_1^2+A_2A_3=-A_3A_2-A_4^2, 
\end{equation}
and, applying \eqref{24020609}-\eqref{24020610}, \eqref{25120405} is transformed into
\begin{equation}\label{25111801}
\begin{aligned}
A_2^{-2}=\frac{1}{a}(A_2^{-1}A_1-A_1^2+A_2^{-2}A_1A_2-A_2^{-1}A_1^2A_2). 
\end{aligned}
\end{equation}
Since $A_2^{-2}=-\frac{\tilde{a}}{1-a}I$ and $A_1^2=(\mathrm{tr}\,A_1)A_1-(\det A_1)I$, \eqref{25111801} is equivalent to 
\begin{equation}\label{25120407}
-\frac{\tilde{a}}{1-a}I=\frac{1}{a}\big(A_2^{-1}A_1-(\mathrm{tr}\,A_1)A_1+2aI+A_1A_2^{-1}+(\mathrm{tr}\,A_1)\frac{1-a}{\tilde{a}}A_2^{-1}A_1A_2^{-1}\big), 
\end{equation}
and, using
\begin{equation*}
\begin{aligned}
A_1A_2^{-1}=(\mathrm{tr}\,A_1A_2^{-1})I-(\det A_1A_2^{-1}) (A_1A_2^{-1})^{-1}=-\frac{\tilde{a}}{1-a}(I+aA_2A_1^{-1}), 
\end{aligned}
\end{equation*}
\eqref{25120407} is further simplified as
\begin{equation}\label{25112603}
\begin{aligned}
-(\mathrm{tr}\,A_1)A_1+(2a-\tilde{a})I-(\mathrm{tr}\,A_1)aA_1^{-1}=0.
\end{aligned}
\end{equation}
Now $(\mathrm{tr}\,A_1)^2=2a-\tilde{a}$ follows from \eqref{25112603} and $A_1-(\mathrm{tr}\,A_1)I+aA_1^{-1}=0$. Since 
\begin{equation*}
A_2A_3=-\frac{1-a}{a}(A_2^{-1}A_1-A_1^2)
\end{equation*}
by \eqref{24020609}, we attain 
\begin{equation*}
\mathrm{tr}\, (A_2A_3)=-\frac{1-a}{a}\big(-\frac{\tilde{a}}{1-a}-\mathrm{tr}\, ((\mathrm{tr}\,A_1)A_1-aI)\big)=\tilde{a}, 
\end{equation*}
which completes the proof of the claim. 
\end{enumerate}
\end{proof}

To simplify the computation, we further change basis of $M$ and normalize it as  
\begin{equation}\label{24122901}
\left(\begin{array}{cc}
A_1 & I \\
A_2A_3 & A_2A_4A_2^{-1}
\end{array}\right). 
\end{equation}
Then the primary matrix $P_M$ associated to $M$ is\footnote{Note that the signs in \eqref{24120603} and \eqref{23112205} are independent.} 
\begin{equation}\label{24120603}
\left(\begin{array}{cc}
\frac{\mathrm{tr}\,A_1\pm \sqrt{-\tilde{a}-2a}}{2} & 1\\
\frac{\tilde{a}\pm \sqrt{\tilde{a}^2-4(1-a)^2}}{2} & \frac{-\mathrm{tr}\,A_1\pm \sqrt{-\tilde{a}-2a}}{2}
\end{array}\right), 
\end{equation}
and $\det P_M$ is either 
\begin{equation}\label{23112205}
-a-\frac{\tilde{a}\pm\sqrt{\tilde{a}^2-4(1-a)^2}}{2}\quad \text{or}\quad \pm \frac{\sqrt{\tilde{a}^2-4a^2}}{2}\pm\frac{\sqrt{\tilde{a}^2-4(1-a)^2}}{2}.
\end{equation}

\begin{enumerate} 
\item We first assume $\det P_M$ is a root of unity. 
\begin{enumerate}
\item Suppose $\det P_M$ is equal to the first type in \eqref{23112205}. Note that, since $\det P_M$ is quadratic, it is either $\pm 1, \pm \sqrt{-1}$ or $\frac{\pm 1\pm \sqrt{-3}}{2}$.   
\begin{itemize}
\item[-] If $\det P_M=\pm 1$, then $\tilde{a}=\pm 2-2a$ and $\tilde{a}^2=4(1-a)^2$. If $\tilde{a}=-2-2a$ and $\tilde{a}^2=4(1-a)^2$, then $a=0$, contradicting the fact that $a>0$. On the other hand, if $\tilde{a}=2-2a$, then \eqref{24120603} and Lemma \ref{24020501} imply $A_2A_3\parallel I$; more precisely, $A_2A_3=(1-a)I$. So, by the last equality in \eqref{24020605}, $A_2A_4A_2^{-1}=-\frac{\det A_1}{1-\det A_1}A_2A_3A_1^{-1}=-aA_1^{-1}.$ Consequently, the primary matrix $P_M$ associated to $M$ is 
\begin{equation*}
\left(\begin{array}{cc}
\frac{\sqrt{4a-2}+\sqrt{-2}}{2} & 1\\
1-a & \frac{-\sqrt{4a-2}+\sqrt{-2}}{2}
\end{array}\right), 
\end{equation*}
and the eigenvalues of $P_M\iota $ are $\frac{\sqrt{4a-2}\pm \sqrt{4a-6}}{2}$. Now, since two nonzero eigenvalues are distinct each other\footnote{If they are the same, then $a=\frac{3}{2}$, which is a contradiction.} and $\mathrm{tr}\,P_M\neq 0$, it follows from Proposition \ref{24082301} that the eigenvalues are roots of unity. One can check $a=\frac{3}{4}$ is the only value satisfying this condition. In conclusion, we derive $\tilde{a}=\frac{1}{2}$, $A_2A_3=\frac{I}{4}$ and $A_2A_4A_2^{-1}=-\frac{3}{4}A_1^{-1}$, which implies that \eqref{24122901} is $\left(\begin{array}{cc}
A_1 & I \\
\frac{I}{4} & -\frac{3}{4}A_1^{-1}
\end{array}\right)$, and the original $M$, prior to normalization, corresponds to the third form stated in \eqref{24010303}.

\item[-] If $\det P_M=\pm \sqrt{-1}$, then $\tilde{a}+2a=0$ and $\tilde{a}^2-4(1-a)^2=\pm 4$, implying $a= 1$ or $0$. But this contradicts the $a$ is one of $\frac{1}{4}, \frac{1}{2}$ or $\frac{3}{4}$. Similarly, for $\det P_M=\frac{\pm 1 \pm \sqrt{-3}}{2}$, we get the same contradiction. 
\end{itemize}
\item Suppose $\det P_M$ is of the second type in \eqref{23112205}. If $\det P_M=\pm 1$, then one can check $a$ is either $1$ or $0$, contradicting the fact that $a$ is one of $\frac{1}{4}, \frac{1}{2}$ or $\frac{3}{4}$. For $\det P_M=\pm\frac{1}{2}\pm \frac{\sqrt{-3}}{2}$ or $\pm \frac{\sqrt{2}}{2}\pm \frac{\sqrt{-2}}{2}$, one analogously obtains the same contradiction. 
\end{enumerate}

\item Next we assume $\det P_M$ is not a root of unity. We only consider the following two cases, as the remaining cases can be handled in a similar manner to one of these. 
\begin{enumerate}
\item First, we consider $P_M=
\left(\begin{array}{cc}
\frac{\mathrm{tr}\,A_1+\sqrt{-\tilde{a}-2a}}{2} & 1\\
\frac{\tilde{a}+\sqrt{\tilde{a}^2-4(1-a)^2}}{2} & \frac{-\mathrm{tr}\,A_1+ \sqrt{-\tilde{a}-2a}}{2}
\end{array}\right)$. If the eigenvalues of $P_M\iota$ are the same, then $\tilde{a}=2(1-a)$, contradicting the fact that $\det P_M$ is not a root of unity. As $\mathrm{tr}\,P_M\neq 0$ and $\det (P_M\iota)\neq 0$, it follows from Proposition \ref{24082301} that the eigenvalues are roots of unity, which contradicts the fact that $\det P_M$ is not a root of unity. 

\item Second, we set $P_M=\left(\begin{array}{cc}
\frac{\mathrm{tr}\,A_1+\sqrt{-\tilde{a}-2a}}{2} & 1\\
\frac{\tilde{a}+\sqrt{\tilde{a}^2-4(1-a)^2}}{2} & \frac{-\mathrm{tr}\,A_1- \sqrt{-\tilde{a}-2a}}{2}
\end{array}\right)$. If $\det P_M=0$, then one finds $a=\frac{1}{2}$. As $\mathbb{Q}(\sqrt{-\tilde{a}-2a})=\mathbb{Q}(\sqrt{\tilde{a}^2-4(1-a)^2})$ and $\tilde{a}$ is one of $\frac{1}{4}, \frac{1}{2}$ or $\frac{3}{4}$, it follows that $\tilde{a}=\frac{3}{4}$ and hence $\mathbb{Q}(\sqrt{-\tilde{a}-2a})=\mathbb{Q}(\sqrt{-7})$. However, this contradicts the fact that two cusp shapes of $\mathcal{X}$ are contained in $\mathbb{Q}(\sqrt{-1})$ (see the footnote after \eqref{241204011}), and hence $\det P_M\neq 0$. If $\lambda_i$ (resp. $\zeta_i$) for $i=1,2$ are the eigenvalues of $P_M$ (resp. $\overline{P_M}$), since $\det P_M$ is not a root of unity, we have either 
\begin{equation}\label{24100205}
\lambda_1=\overline{\zeta_1}, \lambda_2=\overline{\zeta_2}\quad \text{or}\quad  \lambda_2=\overline{\zeta_1}, \lambda_1=\overline{\zeta_2} 
\end{equation}
by Lemma \ref{23122103}. Further, as $\lambda_1+\lambda_2=0$, $\zeta_1=\lambda_1^k(-1)^l$ for some $k(\geq 3)$ and $l\in \mathbb{N}$ by Proposition \ref{24010906}. However, either condition in \eqref{24100205} induces that $\lambda_1$ (so $\lambda_2$ as well) is a root of unity, which is a contradiction. 
\end{enumerate}
\end{enumerate}

\item If $A_2\parallel I$, since $a+\tilde{a}(\det A_2)=1$ by \eqref{24020605}, $A_2$ is either $\sqrt{\frac{1-a}{\tilde{a}}}I$ or $-\sqrt{\frac{1-a}{\tilde{a}}}I$. Without loss of generality, let us assume the first one. 
\begin{claim}\label{24021501}
If $A_2=\sqrt{\frac{1-a}{\tilde{a}}}I$, then $A_3=-\sqrt{(1-a)\tilde{a}}I$ and $A_4=aA_1^{-1}$.
\end{claim}
\begin{proof}[Proof of Claim \ref{24021501}]
Note that the following two equations 
\begin{equation*}
A_3=-\frac{1-a}{a}(A_2^{-2}A_1-A_2^{-1}A_1^2) \quad \text{and}\quad A_2^{-2}=\mathrm{tr}\,(A_2^{-1}A_1)I
\end{equation*}
obtained in the proof of Claim \ref{24020607} (see \eqref{24020609} and \eqref{24021504}) hold without any assumption on $A_2$. Thus $\mathrm{tr}\,A_1=\sqrt{\frac{\tilde{a}}{1-a}}$ and
\begin{equation*}
\begin{aligned}
A_3=-\frac{1-a}{a}\big(\frac{\tilde{a}}{1-a}A_1-\sqrt{\frac{\tilde{a}}{1-a}}A_1^2\big)=-\sqrt{(1-a)\tilde{a}}I.
\end{aligned}
\end{equation*} 
From $A_4=-\frac{a}{1-a}A_3A_1^{-1}A_2$, $A_4=aA_1^{-1}$ follows, which completes the proof of the claim. 
\end{proof}
Since $A_1^2+A_2A_3=-A_3A_2-A_4^2$ by \eqref{23072703}, we get $A_1^2+(a-1)I=(1-a)I-a^2A_1^{-2}$. As $A_1+aA_1^{-1}=\sqrt{\frac{\tilde{a}}{1-a}}I$, we conclude $\frac{\tilde{a}}{1-a}=2$. But this contradicts the fact that $\mathrm{tr}\,A_1= \sqrt{\frac{\tilde{a}}{1-a}}\in \mathbb{Q}$. 
\end{enumerate}
In conclusion, if $M$ is of Type III such that $M^2$ is of the form given in \eqref{241204011}, then $M$ must be of the third form given in \eqref{24010303}. 
\end{enumerate}
\end{proof}

\subsection{$x^4+x^2+1$}

This subsection deals with $M$ satisfying $\mathfrak{m}_M(x)=x^4+x^2+1$. The overall strategy of the subsection is similar to that of the previous two subsections. In this case, since $\mathfrak{m}_{M^2}(x)$ is  $x^2+x+1$, we appeal to Corollary \ref{24121701} to prove the main theorem. 

\begin{theorem}\label{23070525}
Let $M$ be a matrix in $\mathrm{Aut}\,\mathcal{X}$ (where $\mathcal{X}\in \mathfrak{Ghol}$) whose minimal polynomial is $x^4+x^2+1$. Then the following statements holds. 
\begin{enumerate}
\item If $M$ is of Type I, then it is of the following form 
\begin{equation*}
A_1 \oplus A_4
\end{equation*}
for some $A_1, A_4\in \mathrm{GL}_2(\mathbb{Q})$ satisfying either $\mathfrak{m}_{A_1}(x)=x^2+x+1$ and $\mathfrak{m}_{A_4}(x)=x^2-x+1$, or $\mathfrak{m}_{A_1}(x)=x^2-x+1$ and $\mathfrak{m}_{A_4}(x)=x^2+x+1$.

\item If $M$ is of Type II, then it is of the following form 
\begin{equation*}
A_2\widetilde{\oplus}A_3
\end{equation*}
for some $A_2, A_3\in \mathrm{GL}_2(\mathbb{Q}$) satisfying $\mathfrak{m}_{A_2A_3}(x)=x^2+x+1$.

\item If $M$ is of Type III, then it is either 
\begin{equation}\label{24120412}
\left(\begin{array}{cc}
A_1 & A_2\\
3A_2^{-1}A_1^2 & -A_2^{-1}A_1A_2
\end{array}\right)\quad \text{or}\quad 
\left(\begin{array}{cc}
A_1 & A_2\\
\frac{A_2^{-1}}{4} & A_2^{-1}A_1A_2
\end{array}
\right)
\end{equation}
for some $A_1, A_2\in\mathrm{GL}_2(\mathbb{Q})$ satisfying $\mathfrak{m}_{A_1}(x)=x^2\pm\frac{1}{2}x+\frac{1}{4}$ or $\mathfrak{m}_{A_1}(x)=x^2+\frac{3}{4}$ respectively. In particular, if $M$ is of the first (resp. second) form given in \eqref{24120412}, then $M^2$ is of Type I (resp. Type III). Moreover, in both cases, we have $\tau_i\in \mathbb{Q}(\sqrt{-3})$ for $i=1,2$. 
\end{enumerate}
\end{theorem}
\begin{proof}
If $M$ is of Type I, expressed as $A_1\oplus A_4$, then both $\mathfrak{m}_{A_1}(x)$ and $\mathfrak{m}_{A_4}(x)$ divide $x^4+x^2+1$. Since $x^4+x^2+1$ is factored as $(x^2+x+1)(x^2-x+1)$, we get the conclusion as stated above. 

Suppose $M$ is of Type II, given as $A_2\widetilde{\oplus}A_3$. Since the minimal polynomial of $M^2=A_2A_3\oplus  A_3A_2$ is $x^2+x+1$, $M^2$ is one of the forms given in \eqref{25112201} by Corollary \ref{24121701}. However, the last two in \eqref{25112201} do not occur, as $\mathfrak{m}_{A_2A_3}(x)=\mathfrak{m}_{A_3A_2}(x)$. 

If $M$ is of Type III, then first note that $M^2$ is not of block anti-diagonal matrix by Corollary \ref{24121701}. 
\begin{enumerate}
\item If $M^2$ is of block diagonal matrix, by Lemma \ref{24082308}, $M$ and $M^2$ are given as 
\begin{equation*}
\left(\begin{array}{cc}
A_1 & A_2\\
\frac{\det A_2}{\det A_1}A_2^{-1}A_1^2 & -A_2^{-1}A_1A_2
\end{array}\right)
\quad \text{and}\quad \left(\begin{array}{cc}
\frac{1}{\det A_1}A^2_1 & 0\\
0 & \frac{1}{\det A_1}A_2^{-1}A_1^2A_2
\end{array}\right)
\end{equation*}
respectively. Since $A_1^4+(\det A_1)A_1^2+(\det A_1)^2I=0$, it follows that $\mathfrak{m}_{A_1}(x)$ is either $x^2+tx+t^2$ or $x^2-tx+t^2$ for some $t\in \mathbb{Q}$. Without loss of generality, assume $t>0$ and, changing variables if necessary, normalize $M$ as 
\begin{equation*}
\left(\begin{array}{cc}
A_1 & I\\
\frac{1-t^2}{t^2}A_1^2 & -A_1
\end{array}\right), 
\end{equation*}
in which case the primary matrix $P_M$ associated with $M$ is 
\begin{equation}\label{25011503}
\left(\begin{array}{cc}
t\omega & 1\\
(1-t^2)\omega^2 & -t\omega
\end{array}\right) 
\end{equation}
where $\omega=\frac{1\pm\sqrt{-3}}{2}$. By Lemma \ref{25011101}, the order of $P_M\iota$ is finite, so the eigenvalues of \eqref{25011503}, which are the solutions of 
\begin{equation*}
x^2-2t\omega x+\omega^2=0, 
\end{equation*}
are roots of unity. Since $t>0$, the degrees of the eigenvalues are at most $2$ by Lemma \ref{24122301}, and one can verify the only possible value for $t$ is $\frac{1}{2}$. Therefore the original $M$, before normalization, is of the first form given in \eqref{24120412}. 

\item Now we assume $M^2$ is neither of block diagonal nor anti-diagonal matrix. By Corollary \ref{24121701}, $M^2$ is one of the following forms: 
\begin{equation}\label{24120403}
\left(\begin{array}{cc}
-\frac{I}{2} & \tilde{A_2}\\
-\frac{3}{4}\tilde{A_2}^{-1} & -\frac{I}{2}
\end{array}
\right),\;\; 
\left(\begin{array}{cc}
\tilde{A_1} & \tilde{A_2}\\
-\frac{\tilde{A_2}^{-1}}{2} &-\tilde{A_2}^{-1}(I+\tilde{A_1})\tilde{A_2}
\end{array}\right)\;\; \text{or}\;\; 
\left(\begin{array}{cc}
\tilde{A_1} & \tilde{A_2}\\
-\frac{\tilde{A_2}^{-1}}{4} & -\tilde{A_2}^{-1}(I+\tilde{A_1})\tilde{A_2}
\end{array}\right)
\end{equation}
where $\det \tilde{A_2}=\frac{3}{4}$, $\mathfrak{m}_{\tilde{A_1}}(x)=x^2+x+\frac{1}{2}$ or $\mathfrak{m}_{\tilde{A_1}}(x)=x^2+x+\frac{3}{4}$ respectively. 
\begin{enumerate}
\item We assume $M^2$ is of the first form given in \eqref{24120403}. By changing variables if necessary, we normalize $M^2$ and denote $M$ as 
\begin{equation*}
\left(\begin{array}{cc}
-\frac{I}{2} & I\\
-\frac{3}{4}I & -\frac{I}{2}
\end{array}\right)\quad \text{and}\quad \left(\begin{array}{cc}
A_1 & A_2\\
A_3 & A_4
\end{array}\right)
\end{equation*}
respectively.\footnote{That is, $M^2=\left(\begin{array}{cc}
A_1 & A_2\\
A_3 & A_4
\end{array}\right)^2
=\left(\begin{array}{cc}
-\frac{I}{2} & I\\
-\frac{3}{4}I & -\frac{I}{2}
\end{array}\right)$.} 
Similar to the proof of Claim \ref{24082401}, one obtains
$\det A_1=\frac{3}{4}$ and $A_3=\frac{1}{4}A_2^{-1}$. Consequently, $M$ is of the following form
\begin{equation*}
\left(\begin{array}{cc}
A_1 & A_2\\
\frac{A_2^{-1}}{4} & A_2^{-1}A_1A_2
\end{array}\right), 
\end{equation*}
which corresponds to the second type given in \eqref{24120412}. 

\item Now we consider the second case in \eqref{24120403}. Similar to the above cases, by changing variables if necessary, we assume $M$ and $M^2$ are given as
\begin{equation}\label{25112103}
\left(\begin{array}{cc}
A_1 & A_2\\
A_3 & A_4
\end{array}
\right)\quad \text{and}\quad \left(\begin{array}{cc}
\tilde{A_1} & I\\
-\frac{1}{2}I & -(I+\tilde{A_1})
\end{array}\right)
\end{equation}
respectively. 

\begin{enumerate}
\item First, if we assume $A_2\notparallel I$, then, following the analogous steps given in Claim \ref{24020607}, we derive 
\begin{equation}\label{24053001}
\begin{gathered}
\mathrm{tr}\,A_2=0,\;\; \mathrm{tr}\,A_1=-\mathrm{tr}\,A_4, \;\; (\mathrm{tr}\,A_1)^2=a-\frac{1}{2}, \;\; \mathrm{tr}\, (A_2A_3)=a-\frac{1}{2}
\end{gathered}
\end{equation}
where $a:=\det A_1$. Now we further change basis of $M$ and work with 
\begin{equation*}
\left(\begin{array}{cc}
A_1 & I \\
A_2A_3 & A_2A_4A_2^{-1}
\end{array}\right).
\end{equation*}
Since $\det A_2A_3=(1-a)^2$ and $\det A_1=\det A_4$, combined with \eqref{24053001}, we get that the primary matrix $P_M$ associated with $M$ is
\begin{equation*}
\left(\begin{array}{cc}
\frac{t\pm \sqrt{-\frac{1}{2}-3a}}{2} & 1\\
\frac{a-\frac{1}{2}\pm \sqrt{(a-\frac{1}{2})^2-4(1-a)^2}}{2} & \frac{-t\pm \sqrt{-\frac{1}{2}-3a}}{2} 
\end{array}\right) 
\end{equation*}
where $t:=\mathrm{tr}\,A_1$.

\begin{enumerate}
\item If $\det P_M$ is a root of unity, it is $\pm 1,\pm \sqrt{-1}$ or $\frac{\pm 1\pm \sqrt{-3}}{2}$, and one can check $a=\frac{5}{6}$ is the only possible value satisfying the criterion. However, in this case, it contradicts the fact that $\mathrm{tr}\,A_1=\pm\sqrt{a-\frac{1}{2}}\in \mathbb{Q}$. 

\item Suppose $\det P_M$ is not a root of unity. We only consider the following case
\begin{equation*}
\begin{gathered}
P_M=
\left(\begin{array}{cc}
\frac{t+ \sqrt{-\frac{1}{2}-3a}}{2} & 1\\
\frac{a-\frac{1}{2}+ \sqrt{(a-\frac{1}{2})^2-4(1-a)^2}}{2} & \frac{-t+ \sqrt{-\frac{1}{2}-3a}}{2}
\end{array}\right), 
\end{gathered}
\end{equation*}
as the remaining cases can be dealt with in an analogous manner. If two eigenvalues of $P_M$ are the same, then $\mathrm{Disc}\,P_M=0$. However, one can check there is no $a\in \mathbb{Q}$ satisfying this. If $\mathrm{tr}\,P_M\iota=0$, then $t=0$, so $a=\frac{1}{2}$ and $P_M=
\left(\begin{array}{cc}
\frac{\sqrt{-2}}{2} & 1\\
\frac{\sqrt{-1}}{2} & \frac{\sqrt{-2}}{2}
\end{array}\right)$. But in this case, it contradicts the fact that $P_M\in \mathrm{GL}_2(\mathbb{Q}(\sqrt{D}))$ for some $D<0$. Similarly, one confirms there is no $a$ such that $\det P_M\neq 0$, thus by Proposition \ref{24082301}, it follows that two eigenvalues of $P_M$ are all roots of unity, contradicting the assumption that $\det P_M$ is not a root of unity. 
 \end{enumerate}
\item Next we assume $A_2\parallel I$. Since $\det A_2=2(1-a)$, $A_2$ is either $\sqrt{2(1-a)}I$ or $-\sqrt{2(1-a)}I$. Without loss of generality, we assume the first one. Following the same steps given in the proof of Claim \ref{24021501}, one gets $A_3=-\sqrt{\frac{1-a}{2}}I$ and $A_4=aA_1^{-1}$. Since $A_1^2+A_2A_3=-A_3A_2-A_4^2-I$,\footnote{This follows from the 
expressions given in \eqref{25112103}.}  it follows that $A_1^2+(a-1)I=(1-a)I-a^2A_1^{-2}-I$. As $A_1+aA_1^{-1}=\sqrt{\frac{1}{2(1-a)}}I$, we conclude $a=\frac{1}{2}$ and $M$ is given as $\left(\begin{array}{cc}
A_1 & I \\
-\frac{I}{2} & \frac{A_1^{-1}}{2}
\end{array}\right)$ where $A_1^2-A_1+\frac{I}{2}=0$. However, in this case, $M^2-M+I=0$, contradicting the fact $\mathfrak{m}_M(x)=x^4+x^2+1$. 
\end{enumerate}

In conclusion, there is no $M$ of Type III such that $M^2$ is of the second form given in \eqref{24120403}. 

\item Similar to the previous case, by reproducing its proof almost verbatim, one can show that there is no $M$ of Type III such that $M^2$ is of the third form given in \eqref{24120403}. We omit the detailed proof here. 

\end{enumerate}
\end{enumerate}
\end{proof}

\subsection{$x^4-x^2+1$}

Finally, we consider $M$ satisfying $\mathfrak{m}_M(x)=x^4-x^2+1$.  

\begin{theorem}\label{230705251}
Let $M$ be a matrix in $\mathrm{Aut}\,\mathcal{X}$ (where $\mathcal{X}\in \mathfrak{Ghol}$) whose minimal polynomial is $x^4-x^2+1$. Then $M$ is not of Type I, and the following statements holds. 
\begin{enumerate}
\item If $M$ is of Type II, then $M$ is of the form 
\begin{equation}\label{25010803}
A_2\widetilde{\oplus}A_3 
\end{equation}
for some $A_2, A_3\in \mathrm{GL}_2(\mathbb{Q})$ satisfying $\mathfrak{m}_{A_2A_3}(x)=x^2-x+1$.

\item If $M$ is of Type III, then $M$ is either    
\begin{equation}\label{24110202}
\left(\begin{array}{cc}
A_1 & A_2\\
\frac{A_2^{-1}A_1^2}{3} & -A_2^{-1}A_1^{-1}A_2
\end{array}\right), \quad 
\left(\begin{array}{cc}
A_1 & A_2\\
\frac{3}{4}A_2^{-1} & A_2^{-1}A_1A_2
\end{array}\right) 
\quad  \text{or}\quad 
\left(\begin{array}{cc}
A_1 & A_2\\
\frac{A_2^{-1}}{2} & -\frac{A_2^{-1}A_1^{-1}A_2}{2}
\end{array}\right)
\end{equation}
where $\mathfrak{m}_{A_1}(x)$ is $x^2\pm \frac{3}{2}x+\frac{3}{4}$, $x^2+\frac{1}{4}$ or $x^2\pm x+\frac{1}{2}$ respectively. In particular, $M^2$ is of Type I (resp. Type III) and $\tau_i\in \mathbb{Q}(\sqrt{-3})$ (resp. $\mathbb{Q}(\sqrt{-1})$) in the first case (resp. second and third cases) for $i=1,2$.
\end{enumerate}
\end{theorem}

\begin{proof}
If $M$ is of Type I, given as $A_1\oplus A_4$, then both $\mathfrak{m}_{A_1}(x)$ and $\mathfrak{m}_{A_4}(x)$ must divide $x^4-x^2+1$. However, as $x^4-x^2+1$ is irreducible (over $\mathbb{Q}$) and $\deg \mathfrak{m}_{A_i}(x)\leq 2$ for $i=1,4$, it follows that such $\mathfrak{m}_{A_i}(x)$ do not exist. If $M$ is of Type II, then it is clearly of the form given in \eqref{25010803}. 

If $M$ is of Type III, then note that, since $\mathfrak{m}_{M^2}(x)=x^2-x+1$, $M^2$ is not of block anti-diagonal matrix by Theorem \ref{24040801}. If $M^2$ is of block diagonal matrix, then, following almost verbatim the proof of that case in Theorem \ref{23070525}, one finds $M$ must be the first type as given in \eqref{24110202}.

If $M^2$ is neither of block diagonal nor anti-diagonal matrix, then $M^2$ is one of the following: 
\begin{equation}\label{25012101}
\left(\begin{array}{cc}
\frac{I}{2} & \tilde{A_2}\\
-\frac{3}{4}\tilde{A_2}^{-1} & \frac{I}{2}
\end{array}\right), \;\; 
\left(\begin{array}{cc}
\tilde{A_1} & \tilde{A_2}\\
-\frac{\tilde{A_2}^{-1}}{2} & \tilde{A_2}^{-1}(I-\tilde{A_1})\tilde{A_2}
\end{array}
\right)\;\; \text{or}\;\;
\left(\begin{array}{cc}
\tilde{A_1} & \tilde{A_2}\\
-\frac{\tilde{A_2}^{-1}}{4} & \tilde{A_2}^{-1}(I-\tilde{A_1})\tilde{A_2}
\end{array}
\right)
\end{equation}
with $\det \tilde{A_2}=\frac{3}{4}, \mathfrak{m}_{\tilde{A_1}}(x)=x^2-x+\frac{1}{2}$ or $\mathfrak{m}_{\tilde{A_1}}(x)=x^2-x+\frac{3}{4}$ respectively by Theorem \ref{24040801}.
\begin{enumerate}
\item We assume the first case. By changing variables if necessary, we normalize $M^2$ and denote $M$ as 
\begin{equation*}
\left(\begin{array}{cc}
\frac{I}{2} & I\\
-\frac{3}{4}I & \frac{I}{2}
\end{array}\right)\quad \text{and}\quad \left(\begin{array}{cc}
A_1 & A_2\\
A_3 & A_4
\end{array}\right)
\end{equation*}
respectively. Then, closely following the proof of Claim \ref{24082401}, one verifies $\det A_1=\frac{1}{4}$ and $A_3=\frac{3}{4}A_2^{-1}$, thus $M$ is of the following form
\begin{equation*}
\left(\begin{array}{cc}
A_1 & A_2\\
\frac{3}{4}A_2^{-1} & A_2^{-1}A_1A_2
\end{array}\right).
\end{equation*}
Consequently, before normalization, $M$ is given by the second type described in \eqref{24110202}. 

\item Now we consider the second case in \eqref{25012101}. Similar to the above cases, by changing variables if necessary, we assume $M^2$ and $M$ are given as
\begin{equation*}
\left(\begin{array}{cc}
\tilde{A_1} & I\\
-\frac{I}{2} & I-\tilde{A_1}
\end{array}
\right)\quad \text{and}\quad \left(\begin{array}{cc}
A_1 & A_2\\
A_3 & A_4
\end{array}
\right)
\end{equation*}
respectively. Then 
\begin{equation}\label{240521011}
A_1^2+A_2A_3=\tilde{A_1}, \;\; A_1A_2+A_2A_4=I, \;\; A_3A_1+A_4A_3=-\frac{I}{2}, \;\; A_3A_2+A_4^2=-\tilde{A_1}+I
\end{equation}
and
\begin{equation}\label{240521021}
\det A_1=\det A_4, \;\; \det A_1+\frac{\det A_2}{2}=\det A_1+2\det A_3=1, \;\; A_4=-\frac{\det A_1}{1-\det A_1}A_3A_1^{-1}A_2. 
\end{equation}

\begin{enumerate}
\item First if $A_2\notparallel I$, using analogous steps from Claim \ref{24082401}, we derive the following from \eqref{240521011}-\eqref{240521021}:
\begin{equation}\label{240530011}
\begin{gathered}
\mathrm{tr}\,A_2=0, \quad  \det A_2=-\frac{1}{\mathrm{tr}\,A_2^{-1}A_1}=2(1-a), \quad \mathrm{tr}\,A_1=-\mathrm{tr}\,A_4, \\
(\mathrm{tr}\,A_1)^2=3a-\frac{1}{2} \quad\text{and}\quad \mathrm{tr}\, (A_2A_3)=-a+\frac{3}{2} 
\end{gathered}
\end{equation}
where $a:=\det A_1$. 

We further change basis of $M$ and normalize it as 
\begin{equation}\label{25010805}
\left(\begin{array}{cc}
A_1 & I \\
A_2A_3 & A_2A_4A_2^{-1}
\end{array}\right). 
\end{equation}
Since $\det A_2A_3=(1-a)^2$ and $\det A_1=\det A_4$ by \eqref{240521021}, the primary matrix $P_M$ associated with \eqref{25010805} is
\begin{equation*}
\left(\begin{array}{cc}
\frac{t\pm \sqrt{-\frac{1}{2}-a}}{2} & 1\\
\frac{-a+\frac{3}{2}\pm \sqrt{(a-\frac{3}{2})^2-4(1-a)^2}}{2} & \frac{-t\pm \sqrt{-\frac{1}{2}-a}}{2}
\end{array}\right) 
\end{equation*}
where $t=\mathrm{tr}\,A_1$. 

\begin{enumerate}
\item Suppose $\det P_M$ is not a root of unity. Similar to previous theorems, we only test the following case, 
\begin{equation*}
P_M=\left(\begin{array}{cc}
\frac{t+ \sqrt{-\frac{1}{2}-a}}{2} & 1\\
\frac{-a+\frac{3}{2}+ \sqrt{(a-\frac{3}{2})^2-4(1-a)^2}}{2} & \frac{-t+ \sqrt{-\frac{1}{2}-a}}{2}
\end{array}\right), 
\end{equation*}
as the remaining cases are treated analogously. First, note that there is no $a\in \mathbb{Q}$ satisfying $\mathrm{Disc}\,P_M=0$ or $\det P_M=0$. If $\mathrm{tr}\,(P_M\iota)=0$, then $t=0$, so $a=\frac{1}{6}$. However, it contradicts the fact that $P_M\in \mathrm{GL}_2(\mathbb{Q}(\sqrt{D}))$ for some $D<0$. In conclusion, the eigenvalues of $P_M$ are roots of unity by Proposition \ref{24082301}, contradicting the assumption that $\det P_M$ is not a root of unity.

\item If $\det P_M$ is a root of unity, it is one of $\pm 1, \pm \sqrt{-1}$ or $\frac{\pm 1\pm \sqrt{-3}}{2}$, and one gets $a=\frac{1}{2}$ is the only possible value satisfying the criterion. We claim 
\begin{claim}\label{24070810}
If $a=\frac{1}{2}$, then we have either  
\begin{equation*}
\mathfrak{m}_{A_1}(x)=x^2+x+\frac{1}{2}, \; A_2=-2A_1-I, \;  A_3=A_1+\frac{I}{2}, \; A_4=-\frac{A_1^{-1}}{2}
\end{equation*}
or 
\begin{equation}\label{24110201}
\mathfrak{m}_{A_1}(x)=x^2-x+\frac{1}{2}, \; A_2=-2A_1+I, \;  A_3=A_1-\frac{I}{2}, \; A_4=-\frac{A_1^{-1}}{2}.
\end{equation}
\end{claim}
\begin{proof}[Proof of Claim \ref{24070810}]
Since $a=\frac{1}{2}$, by \eqref{240530011}, we get $\mathfrak{m}_{A_2}(x)=x^2+1$, and $\mathfrak{m}_{A_1}(x)$ is either $x^2+x+\frac{1}{2}$ or $x^2-x+\frac{1}{2}$. Without loss of generality, let us assume $\mathfrak{m}_{A_1}(x)=x^2+x+\frac{1}{2}$. By \eqref{240521011}-\eqref{240521021}, 
\begin{equation}\label{25010809}
\begin{aligned}
A_1A_2+A_2A_4&=A_1A_2-A_2A_3A_1^{-1}A_2=I
\Longrightarrow 
A_3=(I+A_2^{-1}A_1)A_1
\end{aligned}
\end{equation}
and 
\begin{equation}\label{25010810}
A_4=-A_3A_1^{-1}A_2=A_2^{-1}(I+A_1A_2^{-1}). 
\end{equation}
Since $\mathrm{tr}\,(A_2^{-1}A_1)=-1$ from \eqref{240530011} and $\det A_2^{-1}A_1=\frac{1}{2}$, combined with \eqref{25010809}-\eqref{25010810}, we get 
\begin{equation*}
A_3=-\frac{A_1^{-1}A_2A_1}{2} \;\; \text{and}\;\;  A_4=-\frac{A_1^{-1}}{2}.
\end{equation*}

To show $A_2=-2A_1-I$, if we set $A_i:=\left(\begin{array}{cc}
a_i & b_i\\
c_i & d_i
\end{array}\right)$ ($i=1,2$), then  
\begin{equation}\label{24070801}
\begin{gathered}
a_1+d_1=-1, \; a_1(-1-a_1)-b_1c_1=\frac{1}{2}\;\; \text{and}\;\; a_2+d_2=0, \; -a_2^2-b_2c_2=1  
\end{gathered}
\end{equation}
from $\mathfrak{m}_{A_1}(x)=x^2+x+\frac{1}{2}$ and $\mathfrak{m}_{A_2}(x)=x^2+1$ respectively, and  
\begin{equation}\label{24070802}
a_1a_2+b_1c_2+b_2c_1+d_1d_2=1 
\end{equation}
from $\mathrm{tr}\,(A_2A_1)=1$. Since $\tau=\frac{-2a_1-1\pm \sqrt{-1}}{2b_1}=\frac{-a_2\pm \sqrt{-1}}{b_2}$, $2b_1$ is either $b_2$ or $-b_2$. If $2b_1=b_2$, then $a_2=2a_1+1$, and $c_2=2c_1$ by the last two equalities in \eqref{24070801}. Combined with \eqref{24070802}, it follows that 
\begin{equation*}
\begin{gathered}
4a_1^2+4a_1+4b_1c_1+1 =1. 
\end{gathered}
\end{equation*}
However, the above equality contradicts the second one in \eqref{24070801}. 
 
On the other hand, if $b_2=-2b_1$, then $a_2=-2a_1-1, c_2=-2c_1$ and $d_2=-2d_1-1$ by \eqref{24070801}, inducing $A_2=-2A_1-I$. 

Similarly, if $\mathfrak{m}_{A_1}(x)=x^2-x+\frac{1}{2}$, then $A_2=-2A_1+I$ and \eqref{24110201}.  
\end{proof}

By the claim, the original $M$, before normalization, takes either one of the following forms:  
\begin{equation*}
\left(\begin{array}{cc}
A_1 & (-2A_1-I)\tilde{A_2}\\
\tilde{A_2}^{-1}(A_1+\frac{I}{2}) & -\frac{\tilde{A_2}^{-1}A_1^{-1}\tilde{A_2}}{2}
\end{array}\right)\; \text{or}\; 
\left(\begin{array}{cc}
A_1 & (-2A_1+I)\tilde{A_2}\\
\tilde{A_2}^{-1}(-A_1+\frac{I}{2}) & -\frac{\tilde{A_2}^{-1}A_1^{-1}\tilde{A_2}}{2}
\end{array}\right)
\end{equation*}
with $\mathfrak{m}_{A_1}(x)$ being either $x^2+x+\frac{1}{2}$ or $x^2-x+\frac{1}{2}$ respectively. One can check each form above consistent with the third form given in \eqref{24110202}. 
\end{enumerate}

\item Next we assume $A_2\parallel I$. Since $\det A_2=2(1-a)$, $A_2$ is either $\sqrt{2(1-a)}I$ or $-\sqrt{2(1-a)}I$. Without loss of generality, we assume the first one. Following the same steps given in the proof of Claim \ref{24021501}, one can prove $A_3=-\sqrt{\frac{1-a}{2}}I$ and $A_4=aA_1^{-1}$. Since $A_1^2+2A_2A_3=-A_3A_2-A_4^2+I$ by \eqref{240521011}, it follows that $A_1^2-2(1-a)I-I+a^2A_1^{-1}=0$. As $A_1+aA_1^{-1}=\frac{1}{\sqrt{2(1-a)}}I$, we conclude $a=\frac{5}{6}$. However, in this case, it contradicts the fact that $\mathrm{tr}\,A_1=\frac{1}{\sqrt{2(1-a)}}\in \mathbb{Q}$. 
\end{enumerate}

\item The proof of the final case is analogous. Following the outline of the proof above, one can verify that there is no $M$ of Type III such that $M^2$ corresponds to the third form given in \eqref{25012101}. We omit the detailed proof here. 
\end{enumerate}
\end{proof}

\newpage
\section{Structure of $\mathrm{Aut}\,\mathcal{X}$}\label{25012005}

By amalgamating all the results from Sections \ref{24081804}-\ref{24081805}, this section explores the structure of the largest subgroup $\mathcal{G}$ of $\mathrm{Aut}\,\mathcal{X}$, whose elements are of Types I, II, and III. Using the theorems proven in this section, we will complete the proof of the main theorem in the next section. 

We first consider the case where $\mathcal{G}$ contains only elements of Types I and then address the remaining cases subsequently.  

\subsection{Type I}

The following theorem summarizes the results in Sections \ref{24081804}-\ref{24081805} concerning elements of Type I in $\mathrm{Aut}\,\mathcal{X}$. 

\begin{theorem}\label{24090201}
Let $\mathcal{X}\in\mathfrak{Ghol}$ and $M$ be an element in $\mathrm{Aut}\,\mathcal{X}$ of Type I, which is neither $\pm \iota$ nor $\pm I$. Then cusp shapes of $\mathcal{X}$ is either contained in $\mathbb{Q}(\sqrt{-3})$ or $\mathbb{Q}(\sqrt{-1})$. More specifically, the following statements hold. 
\begin{enumerate}

\item If the cusp shapes of $\mathcal{X}$ are contained in $\mathbb{Q}(\sqrt{-3})$. 
\begin{enumerate}
\item If $\mathfrak{m}_{M}(x)=x^2+x+1$ (resp. $x^2-x+1$), then $M$ is one of 
\begin{equation*}
A_1 \oplus A_4,\;  A_1\oplus (-I)
\;\; \text{or}\;\;
(-I) \oplus A_4\quad (\text{resp. } 
A_1 \oplus A_4,\;\;  
A_1\oplus I\;\; \text{or}\;\; I \oplus A_4)
\end{equation*}
for some $A_i\in \mathrm{GL}_2(\mathbb{Q})$ ($i=1,4$) satisfying $\mathfrak{m}_{A_i}(x)=x^2+x+1$ (resp. $x^2-x+1$). 

\item If $\mathfrak{m}_{M}(x)=x^3-1$ (resp. $x^3+1$), then $M$ is either 
\begin{equation*}
A_1\oplus I\;\; \text{or}\;\; 
I \oplus A_4\quad 
(resp. \;\;
A_1 \oplus(-I)
\;\; \text{or}\;\;
(-I) \oplus A_4)
\end{equation*}
for some $A_i\in \mathrm{GL}_2(\mathbb{Q})$ ($i=1,4$) satisfying $\mathfrak{m}_{A_i}(x)=x^2+x+1$ (resp. $x^2-x+1$). 

\item If $\mathfrak{m}_{M}(x)=x^3-2x^2+2x-1$ (resp. $x^3+2x^2+2x+1$)  then $M$ is either 
\begin{equation*}
A_1\oplus I
\;\; \text{or}\;\;  
I \oplus A_4\quad (resp. \;\;
A_1 \oplus (-I)
\;\; \text{or}\;\;
(-I) \oplus A_4)
\end{equation*}
for some $A_i\in \mathrm{GL}_2(\mathbb{Q})$ ($i=1,4$) satisfying $\mathfrak{m}_{A_i}(x)=x^2-x+1$ (resp. $x^2+x+1$).

\item If $\mathfrak{m}_{M}(x)=x^4+x^2+1$, then $M$ is either 
\begin{equation*}
A_1 \oplus (-A_4)\;\; \text{or}\;\;
(-A_1) \oplus A_4
\end{equation*}
for some $A_i\in \mathrm{GL}_2(\mathbb{Q})$ ($i=1,4$) satisfying $\mathfrak{m}_{A_i}(x)=x^2+x+1$. 

\end{enumerate}
\item If the cusp shapes of $\mathcal{X}$ are contained in $\mathbb{Q}(\sqrt{-1})$. 
\begin{enumerate}
\item If $\mathfrak{m}_{M}(x)=x^2+1$, then $M$ is of the form 
\begin{equation*}
A_1 \oplus A_4
\end{equation*}
for some $A_i\in \mathrm{GL}_2(\mathbb{Q})$ ($i=1,4$) satisfying $\mathfrak{m}_{A_i}(x)=x^2+1$. 

\item If $\mathfrak{m}_{M}(x)=x^3-x^2+x-1$ (resp. $x^3+x^2+x+1$), then $M$ is either 
\begin{equation*}
A_1\oplus I\;\; \text{or}\;\; 
I \oplus A_4
\quad (resp. \;\;
A_1\oplus (-I)\;\; \text{or}\;\; 
(-I)\oplus A_4)
\end{equation*}
for some $A_i\in \mathrm{GL}_2(\mathbb{Q})$ satisfying $\mathfrak{m}_{A_i}(x)=x^2+1$ ($i=1,4$). 
\end{enumerate}
\end{enumerate}
\end{theorem}

Using the theorem, we analyze the structure of a subgroup of $\mathrm{Aut}\mathcal{X}$ consisting of elements of Type I, as described below.

\begin{theorem}\label{24090302}
Let $\mathcal{X}\in\mathfrak{Ghol}$ and $\mathcal{G}$ be the largest subgroup of $\mathrm{Aut}\,\mathcal{X}$ generated by elements of Type I. If the cusp shape of $\mathcal{X}$ is contained in $\mathbb{Q}(\sqrt{-3})$ (resp. $\mathbb{Q}(\sqrt{-1}$), then $|\mathcal{G}|\leq 36$ (resp. $16$). In particular, the equality holds when $\mathcal{G}$ is generated by 
\begin{equation}\label{25032701}
A_1\oplus I\;\; \text{and}\;\; 
I \oplus A_4
\end{equation}
where $\mathfrak{m}_{A_i}(x)=x^2-x+1$ (resp. $x^2+1$) for $i=1,4$.  
\end{theorem}


\begin{proof}
Suppose the cusp shapes of $\mathcal{X}$ are contained in $\mathbb{Q}(\sqrt{-3})$. If $\mathcal{G}$ is generated by those given in \eqref{25032701}, then an element in $\mathcal{G}$ is of the form 
\begin{equation}\label{25121705}
A_1^i\oplus A_4^j
\end{equation}
for some $0\leq i,j\leq 5$, and so $|\mathcal{G}|=36$. Further, by Proposition \ref{24020501}(1) and the rigidity established in Lemma \ref{24090301}, every Type I element of $\textnormal{Aut}\,\mathcal{X}$ is of the form given in \eqref{25121705}. That is, $\mathcal{G}$ is the largest possible subgroup of $\textnormal{Aut}\,\mathcal{X}$ satisfying the condition, which confirms the result.

Similarly, one can treat the other case, where the cusp shapes of $\mathcal{X}$ is contained in $\mathbb{Q}(\sqrt{-1})$
\end{proof}

\subsection{Type II}\label{25021201}
The following theorem synthesizes the results from Sections \ref{24081804}-\ref{24081805} regarding elements of Type II in $\mathrm{Aut}\,\mathcal{X}$. 

\begin{theorem}\label{25012901}
For $\mathcal{X}\in\mathfrak{Ghol}$, if $M$ is an element in $\mathrm{Aut}\,\mathcal{X}$ of Type II, then it falls into one of the following categories. 
\begin{enumerate}
\item If $\mathfrak{m}_{M}(x)=x^2-1$ (resp. $x^2+1$), then $M$ is of the form 
\begin{equation*}
A_2\widetilde{\oplus}A_2^{-1} \;\; (\text{resp. } A_2\widetilde{\oplus}
(-A_2^{-1}))
\end{equation*}
for some $A_2\in \mathrm{GL}_2(\mathbb{Q})$.

\item If $\mathfrak{m}_{M}(x)=x^4+1$ (resp. $x^4+x^2+1$ or $x^4-x^2+1$), then $M$ is of the form 
\begin{equation*}
A_2\widetilde{\oplus}A_3 
\end{equation*}
for some $A_2, A_3\in \mathrm{GL}_2(\mathbb{Q})$ satisfying $\mathfrak{m}_{A_2A_3}(x)=x^2+1$ (resp. $x^2+x+1$ or $x^2-x+1$). In particular, in this case, the cusp shape of $\mathcal{X}$ is contained in $\mathbb{Q}(\sqrt{-1})$ (resp. $\mathbb{Q}(\sqrt{-3})$).
\end{enumerate}
\end{theorem}

In the following, we generalize Theorem \ref{24090302} to accommodate elements of Type II in $\mathcal{G}$. 

\begin{theorem}\label{24072401}
Let $\mathcal{X}\in\mathfrak{Ghol}$ whose cusp shapes are belong to the same quadratic field. Let $\mathcal{G}$ be the largest subgroup of $\mathrm{Aut}\,\mathcal{X}$ whose elements are either of Type I or II. 
\begin{enumerate}
\item If the cusp shapes of $\mathcal{X}$ are contained in $\mathbb{Q}(\sqrt{-3})$ (resp. $\mathbb{Q}(\sqrt{-1})$), then $|\mathcal{G}|\leq 72$ (resp. $32$). In particular, the equality holds when $\mathcal{G}$ is generated by 
\begin{equation}\label{24093002}
A_1\oplus I,\;\; I \oplus A_4\;\;\text{and} \;\;A_2\widetilde{\oplus}
A_2^{-1} 
\end{equation}
where $A_1, A_2, A_4\in \mathrm{GL}_2(\mathbb{Q})$ satisfying $\mathfrak{m}_{A_i}(x)=x^2-x+1$ (resp. $x^2+1$) for $i=1,4$. 

\item Otherwise, if the cusp shapes neither belong to $\mathbb{Q}(\sqrt{-3})$ nor $\mathbb{Q}(\sqrt{-1})$, then $|\mathcal{G}|\leq 8$. In particular, the equality holds when $\mathcal{G}$ is generated by
\begin{equation*}
\pm \iota\;\;\text{and} \;\;A_2\widetilde{\oplus}A_2^{-1} 
\end{equation*}
for some $A_2\in \mathrm{GL}_2(\mathbb{Q})$.
\end{enumerate}
\end{theorem}
\begin{proof}
We assume two cusp shapes of $\mathcal{X}$ are contained in $\mathbb{Q}(\sqrt{-3})$ and consider only this case, as the proof for the other cases are almost identical.  

Following a similar scheme presented in the proof of Theorem \ref{24090302}, suppose $\mathcal{G}$ is generated by the elements listed in \eqref{24093002}. We claim every element of $\mathcal{G}$ is one of the following forms:
\begin{equation}\label{24093015}
A_1^i \oplus A_4^j
\quad \text{or}\quad A_1^iA_2\widetilde{\oplus}A_2^{-1}A_1^j
\end{equation}
for some $0\leq i, j\leq 5$. Since
\begin{equation*}
(A_2\widetilde{\oplus}A_2^{-1})
(A_1 \oplus A_4)
(A_2\widetilde{\oplus}A_2^{-1})
=A_2A_4A_2^{-1} \oplus A_2^{-1}A_1A_2
\end{equation*}
is an element in $\mathcal{G}$ such that $\mathfrak{m}_{A_2A_4A_2^{-1}}(x)=\mathfrak{m}_{A_2^{-1}A_1A_2}(x)=x^2-x+1$, by Lemma \ref{24090301}, $A_2A_4A_2^{-1}$ is either $A_1$ or $A_1^{-1}$. If $A_2A_4A_2^{-1}=A_1$, then 
\begin{equation*}
(A_2\widetilde{\oplus}A_2^{-1})
(A_1\oplus A_4)
=A_2A_4\widetilde{\oplus}A_2^{-1}A_1
=A_1A_2\widetilde{\oplus}
A_2^{-1}A_1,
\end{equation*}
and the claim is readily verified from this observation. Analogously, one can prove the claim for the case where $A_2A_4A_2^{-1}=A_1^{-1}$. 

By Theorem \ref{24090302}, $\mathcal{G}$ contains every element of Type I. If $M\in \mathrm{Aut}\,\mathcal{X}$ is an element of Type II, distinct from $A_2\widetilde{\oplus}A_2^{-1}$, then, as $M(A_2\widetilde{\oplus}A_2^{-1})$ is of Type I, it is contained in $\mathcal{G}$ and so $M$ is as well. Consequently, $\mathcal{G}$ contains all possible elements of Types I and II in $\mathrm{Aut}\,\mathcal{X}$. The fact that $|\mathcal{G}|=72$ follows from the presentations given in \eqref{24093015}.   
\end{proof}

\subsection{Type III}\label{25021202}
In this subsection, we consider a subgroup of $\mathrm{Aut}\,\mathcal{X}$ that contains an element of Type III. As in the previous sections, we start by gathering the results of Sections \ref{24081804}-\ref{24081805} concerning all possible elements of Type III. 

\begin{theorem}\label{24092702}
Let $\mathcal{X}\in\mathfrak{Ghol}$ and $M$ be an element in $\mathrm{Aut}\,\mathcal{X}$ of Type III. 
\begin{enumerate}
\item First we assume $A_1\notparallel I$. Then cusp shapes of $\mathcal{X}$ is contained in $\mathbb{Q}(\sqrt{-3})$, $\mathbb{Q}(\sqrt{-2})$ or $\mathbb{Q}(\sqrt{-1})$. More specifically, the following cases may happen. 
\begin{enumerate}

\item Suppose the cusp shapes of $\mathcal{X}$ are contained in $\mathbb{Q}(\sqrt{-3})$. 
\begin{enumerate}
\item If $\mathfrak{m}_{M}(x)=x^2+1$ (resp. $\mathfrak{m}_{M}(x)=x^4+x^2+1$ and $M^2$ is of Type III), then $M$ is of the form 
\begin{equation}\label{24102103}
\left(\begin{array}{cc}
A_1 & A_2\\
-\frac{A_2^{-1}}{4} & -A_2^{-1}A_1A_2
\end{array}\right)\quad  (resp. \left(\begin{array}{cc}
A_1 & A_2\\
\frac{A_2^{-1}}{4} & A_2^{-1}A_1A_2
\end{array}
\right))
\end{equation}
for some $A_1, A_2\in \mathrm{GL}_2(\mathbb{Q})$ satisfying $\mathfrak{m}_{A_1}(x)=x^2+\frac{3}{4}$.

\item If $\mathfrak{m}_{M}(x)=x^3\pm1$ (resp. $\mathfrak{m}_{M}(x)=x^4+x^2+1$ and $M^2$ is of Type I), then $M$ is of the form 
\begin{equation*}
\left(\begin{array}{cc}
A_1 & A_2\\
-3A_2^{-1}A_1^2 & A_2^{-1}A_1A_2
\end{array}\right)\quad  (resp. \left(\begin{array}{cc}
A_1 & A_2\\
3A_2^{-1}A_1^2 & -A_2^{-1}A_1A_2
\end{array}\right))
\end{equation*}
for some $A_1, A_2\in \mathrm{GL}_2(\mathbb{Q})$ satisfying $\mathfrak{m}_{A_1}(x)=x^2\pm \frac{1}{2}x+\frac{1}{4}$.\footnote{Here and throughout the rest of the statements of the theorem, the $\pm$ signs are dependent and correspond to each other.}

\item If $\mathfrak{m}_{M}(x)=x^3\pm 2x^2+2x\pm 1$ (resp. $\mathfrak{m}_{M}(x)=x^4-x^2+1$ and $M^2$ is of Type I), then $M$ is of the form 
\begin{equation*}
\left(\begin{array}{cc}
A_1 & A_2\\
-\frac{1}{3}A_2^{-1}A_1^2 & A_2^{-1}A_1A_2
\end{array}\right)\quad  (resp. \left(\begin{array}{cc}
A_1 & A_2\\
\frac{1}{3}A_2^{-1}A_1^2 & -A_2^{-1}A_1A_2
\end{array}\right))
\end{equation*}
for some $A_1, A_2\in \mathrm{GL}_2(\mathbb{Q})$ satisfying $\mathfrak{m}_{A_1}(x)=x^2\pm \frac{3}{2}x+\frac{3}{4}$.
\end{enumerate}

\item Suppose the cusp shapes of $\mathcal{X}$ are contained in $\mathbb{Q}(\sqrt{-2})$. 
\begin{enumerate}
\item If $\mathfrak{m}_{M}(x)=x^2+1$ (resp. $\mathfrak{m}_{M}(x)=x^4+1$ and $M^2$ is of Type II), then $M$ is of the form 
\begin{equation}\label{24010205}
\left(\begin{array}{cc}
A_1 & A_2\\
-\frac{A_2^{-1}}{2} & -A_2^{-1}A_1A_2
\end{array}\right)\quad (resp. \left(\begin{array}{cc}
A_1 & A_2\\
\frac{A_2^{-1}}{2} & A_2^{-1}A_1A_2
\end{array}\right))
\end{equation}
for some $A_1, A_2\in \mathrm{GL}_2(\mathbb{Q})$ satisfying $\mathfrak{m}_{A_1}(x)=x^2+\frac{1}{2}$.

\item  If $\mathfrak{m}_{M}(x)=x^2\pm x+1$ (resp. $\mathfrak{m}_{M}(x)=x^4+1$ and $M^2$ is of Type III), then $M$ is of the form 
\begin{equation}\label{24010206}
\left(\begin{array}{cc}
A_1 & A_2\\
-\frac{A_2^{-1}}{4} & \frac{3A_2^{-1}A_1^{-1}A_2}{4}
\end{array}\right)\quad (resp. \left(\begin{array}{cc}
A_1 & A_2\\
\frac{A_2^{-1}}{4} & -\frac{3A_2^{-1}A_1^{-1}A_2}{4}
\end{array}\right))
\end{equation}
for some $A_1, A_2\in \mathrm{GL}_2(\mathbb{Q})$ satisfying $\mathfrak{m}_{A_1}(x)=x^2\pm x+\frac{3}{4}$.
\end{enumerate}

\item Suppose the cusp shapes of $\mathcal{X}$ are contained in $\mathbb{Q}(\sqrt{-1})$. 
\begin{enumerate}

\item If $\mathfrak{m}_{M}(x)=x^2+1$ (resp. $\mathfrak{m}_{M}(x)=x^4-x^2+1$ and $M^2$ is of Type III with $\mathfrak{m}_{\iota M^2}=x^2-1$), then $M$ is of the form 
\begin{equation}\label{24010201}
\left(\begin{array}{cc}
A_1 & A_2\\
-\frac{3}{4}A_2^{-1} & -A_2^{-1}A_1A_2
\end{array}\right)\quad 
(resp. \left(\begin{array}{cc}
A_1 & A_2\\
\frac{3}{4}A_2^{-1} & A_2^{-1}A_1A_2
\end{array}\right))
\end{equation}
for some $A_1, A_2\in \mathrm{GL}_2(\mathbb{Q})$ satisfying $\mathfrak{m}_{A_1}(x)=x^2+\frac{1}{4}$.

\item If $\mathfrak{m}_{M}(x)=x^3\pm x^2+x\pm 1$ (resp. $\mathfrak{m}_{M}(x)=x^4+1$ and $M^2$ is of Type I), then $M$ is of the form 
\begin{equation}\label{24121601}
\left(\begin{array}{cc}
A_1 & A_2\\
-A_2^{-1}A_1^2 & A_2^{-1}A_1A_2
\end{array}\right)\quad (resp. \left(\begin{array}{cc}
A_1 & A_2\\
A_2^{-1}A_1^2 & -A_2^{-1}A_1A_2
\end{array}\right))
\end{equation}
for some $A_1, A_2\in \mathrm{GL}_2(\mathbb{Q})$ satisfying $\mathfrak{m}_{A_1}(x)=x^2\pm x+\frac{1}{2}$.

\item If $\mathfrak{m}_{M}(x)=x^2\pm x+1$ (resp. $\mathfrak{m}_{M}(x)=x^4-x^2+1$ and $M^2$ is of Type III with $\mathfrak{m}_{\iota M^2}=x^4+1$), then $M$ is of the form 
\begin{equation}\label{24121602}
\left(\begin{array}{cc}
A_1 & A_2\\
-\frac{A_2^{-1}}{2} & \frac{A_2^{-1}A_1^{-1}A_2}{2}
\end{array}\right)\quad  (resp. \left(\begin{array}{cc}
A_1 & A_2\\
\frac{A_2^{-1}}{2} & -\frac{A_2^{-1}A_1^{-1}A_2}{2}
\end{array}\right))
\end{equation}
for some $A_1, A_2\in \mathrm{GL}_2(\mathbb{Q})$ satisfying $\mathfrak{m}_{A_1}(x)=x^2\pm x+\frac{1}{2}$.
\end{enumerate}
\end{enumerate}

\item If $A_1\parallel I$, then $\mathfrak{m}_{M}(x)$ is $x^2-1$, $x^2+ x+1$ or  $x^2-x+1$, and $M$ is of the form 
\begin{equation}\label{24010301}
\left(\begin{array}{cc}
\pm\frac{I}{2} & A_2\\
\frac{3}{4}A_2^{-1} & \mp\frac{I}{2}
\end{array}\right),\quad 
\left(\begin{array}{cc}
-\frac{I}{2} & A_2\\
-\frac{3}{4}A_2^{-1} & -\frac{I}{2}
\end{array}\right)\quad \text{or}\quad 
\left(\begin{array}{cc}
\frac{I}{2} & -A_2\\
\frac{3}{4}A_2^{-1} & \frac{I}{2}
\end{array}\right)
\end{equation}
for some $A_2\in \mathrm{GL}_2(\mathbb{Q})$ satisfying $\det A_2=\frac{3}{4}$ respectively. 
\end{enumerate}  
\end{theorem}

In \eqref{24102103}, if the minimal polynomial $M$ is $x^2+1$, then the minimal polynomial of $\iota M$ is $x^4+x^2+1$. In the same way, one verifies that the two matrices appearing in each statement of the theorem are always paired in this manner. 

Similar to the previous subsection, we analyze the problem by splitting it into several cases, depending on the quadratic field to which the cusp shapes of $\mathcal{X}$ belong. 

First, if the cusp shapes of $\mathcal{X}$ are contained in $\mathbb{Q}(\sqrt{-3})$, then we have the following theorem:

\begin{theorem}\label{24092901}
Let $\mathcal{X}$ be an element in $\mathfrak{Ghol}$ whose cusp shapes are contained in $\mathbb{Q}(\sqrt{-3})$. Let $\mathcal{G}$ be the largest subgroup of $\mathrm{Aut}\,\mathcal{X}$ generated by elements of Types I-III. If $\mathcal{G}$ contains an element $M$ of Type III such that $\mathfrak{m}_M(x)=x^3\pm 2x^2+2x\pm 1$, then $\mathcal{G}=\langle M, \pm \iota \rangle $. In particular, $\mathcal{G}$ contains every element of Type III exhibited in Theorem \ref{24092702}, and 
\begin{equation}\label{25020205}
\mathcal{G}=\{\pm \iota^{\epsilon_1}\eta^{\alpha_1}M^{\alpha_2}\iota^{\epsilon_2}\;|\; 0\leq \epsilon_1, \epsilon_2, \alpha_2\leq 1 \;\; \text{and}\;\;0\leq \alpha_1\leq 5\}
\end{equation}
where $\eta:=(\iota M^5)(\iota M^2)$. 
\end{theorem}

The precise meaning of the second claim in the theorem is as follows: if $N$ is any element of Type III whose minimal polynomial is one of those possible forms given in Theorem \ref{24092702}(1)(a), then $N\in \mathcal{G}$. However, not every element of Type III possesses this property. For instance, if $N$ is of Type III with $\mathfrak{m}_N(x)=x^2+1$, then there is no element in $\langle N, \pm\iota \rangle $ whose minimal polynomial is $x^3-2x^2+2x-1$. In some sense, $M$ in the theorem can be seen as a \textit{primitive} element that, together with $\pm \iota$, generates all the others.

For $\mathcal{X}$ with cusp shapes in $\mathbb{Q}(\sqrt{-2})$, we observe a similar pattern:

\begin{theorem}\label{25020403}
Let $\mathcal{X}$ be an element in $\mathfrak{Ghol}$ whose cusp shapes are contained in $\mathbb{Q}(\sqrt{-2})$, and $\mathcal{G}$ be the largest subgroup of $\mathrm{Aut}\,\mathcal{X}$ generated by elements of Types I-III. If $\mathcal{G}$ contains an element $M$ of Type III such that $\mathfrak{m}_M(x)=x^2\pm x+1$ and $\mathfrak{m}_{\iota M}(x)=x^4+1$, then $\mathcal{G}=\langle M, \pm \iota \rangle $. In particular, $\mathcal{G}$ contains every element of Type III exhibited in Theorem \ref{24092702}, and 
\begin{equation}\label{25020501}
\mathcal{G}=\{\pm \iota^{\epsilon_1}\eta^{\alpha_1}M^{\alpha_2}\iota^{\epsilon_2}\;|\; 0\leq \epsilon_1, \epsilon_2, \alpha_1\leq 1 \;\; \text{and}\;\;0\leq \alpha_2\leq 3\}
\end{equation}
where $\eta:=(\iota M^2)^2$. 
\end{theorem}

On the other hand, if the cusp shapes of $\mathcal{X}$ belong to $\mathbb{Q}(\sqrt{-1})$, then the structure exhibits slightly different behavior. In this case, there are two elements in $\mathrm{Aut }\mathcal{X}$ that are not compatible with each other, that is, they do not belong to the same subgroup of $\mathrm{Aut }\mathcal{X}$. The exact statement of the theorem is as follows:
 
\begin{theorem}\label{24111305}
Let $\mathcal{X}$ be an element in $\mathfrak{Ghol}$ whose cusp shapes are contained in $\mathbb{Q}(\sqrt{-1})$. Let $\mathcal{G}$ be the largest subgroup of $\mathrm{Aut}\,\mathcal{X}$ generated by elements of Types I-III. 

\begin{enumerate}
\item If $\mathcal{G}$ contains an element $M$ of Type III such that $\mathfrak{m}_{M}(x)$ is $x^2+1$, then $\mathcal{G}=\langle M, \pm \iota \rangle $. Further,    
\begin{equation}\label{24122603}
\begin{aligned}
\mathcal{G}=\{\pm \iota^{\epsilon_1}(\iota M)^{\alpha}\iota^{\epsilon_2}\;|\; 0\leq \epsilon_1, \epsilon_2\leq 1 \;\; \text{and}\;\;0\leq \alpha\leq 5\}.
\end{aligned}
\end{equation}

\item Suppose $\mathcal{G}$ contains Type III elements $M$ and $N$ with $\mathfrak{m}_{M}(x)=x^3\pm x^2+x\pm 1, \mathfrak{m}_{N}(x)=x^4-x^2+1$ and $\mathfrak{m}_{\iota N}(x)=x^2\pm x+1$. Then $\mathcal{G}=\langle M, N, \pm \iota \rangle $. In particular, $\mathcal{G}$ contains an abelian group $\mathcal{H}$ of $16$, whose elements are all of Type I, and 
\begin{equation}\label{24111302}
\mathcal{G}=\{\eta_1M^{\alpha}\eta_2\;|\; \eta_1, \eta_2\in \mathcal{H}, \;0\leq \alpha\leq 2 \}. 
\end{equation}
\end{enumerate}
\end{theorem}

Lastly, if the cusp shapes of $\mathcal{X}$ are not contained in any of the quadratic fields discussed above, then we have 
\begin{corollary}\label{25030302}
Let $\mathcal{X}$ be an element in $\mathfrak{Ghol}$ whose cusp shapes belong to the same quadratic field but not to $\mathbb{Q}(\sqrt{-3})$, $\mathbb{Q}(\sqrt{-2})$ or $\mathbb{Q}(\sqrt{-1})$. Let $\mathcal{G}$ be the largest subgroup of $\mathrm{Aut}\,\mathcal{X}$ generated by elements of Types I-III. If $\mathcal{G}$ contains an element $M$ of Type III, then either $\mathfrak{m}_M(x)$ or $\mathfrak{m}_{\iota M}(x)$ is $x^2-1$, and 
\begin{equation}\label{25030301}
\mathcal{G}=\{\pm\iota^{\epsilon_1}M^{\alpha}\iota^{\epsilon_2}\;|\; 0\leq \epsilon_1, \epsilon_2, \alpha\leq 1 \}.
\end{equation} 
\end{corollary}

\subsubsection{$\sqrt{-3}$}\label{24010501}

In this sub-subsection, we prove Theorem \ref{24092901}. To this end, one first needs to address the compatibility between two elements. Specifically, we consider the following questions: given two Type~III elements of $\mathcal{G}$ with the same minimal polynomial, what can be said about their relationship? More generally, how many Type~III elements with a fixed minimal polynomial can occur in $\mathcal{G}$?

In the proposition below, we establish the rigidity of a Type III element and, along the way, address these questions. Clearly, 
$$\mathfrak{m}_M(x)=\mathfrak{m}_{M^{-1}}(x)=\mathfrak{m}_{\iota M\iota}(x)=\mathfrak{m}_{\iota M^{-1}\iota}(x)
$$ 
for any Type III element $M$. In the following, we show any two elements of Type III with the same minimal polynomial must always arise in this way. 

\begin{proposition}\label{24080401}
Let $\mathcal{X}$ be an element in $\mathfrak{Ghol}$ whose cusp shapes are contained in $\mathbb{Q}(\sqrt{-3})$. Let $M$ and $N$ be two elements of Type III in $\mathrm{Aut}\,\mathcal{X}$ satisfying $\mathfrak{m}_{M}(x)=\mathfrak{m}_{N}(x)$. Then $N\in \langle M, \pm\iota \rangle $. More precisely, $N$ is one of $M, M^{-1}, \iota M\iota$ or $\iota M^{-1}\iota$.  
\end{proposition}

The idea of the proof essentially originates from Proposition \ref{24082301}. By Lemma \ref{24090301}, replacing $N$ with $N^{-1}$ if necessary, one may assume that two $(2\times 2)$-principal submatrices of $M$ and $N$ are the same. After normalizing $M$ and $N$ simultaneously, we find that the product of their two primary submatrices depends only on a single variable, whose exact value is determined by Proposition \ref{24082301}. The claim then follows from this computation. Once the claim is proved for one case, the remaining cases are handled either by an analogous manner or by a simple computational observation.

\begin{proof}[Proof of Proposition \ref{24080401}]
Note that, if $\mathfrak{m}_{M}(x)=\mathfrak{m}_{N}(x)$, then $\mathfrak{m}_{M\iota}(x)=\mathfrak{m}_{N\iota}(x)$, $\mathfrak{m}_{\iota M}(x)=\mathfrak{m}_{\iota N}(x)$, and $\mathfrak{m}_{-M}(x)=\mathfrak{m}_{-N}(x)$. Thus it is enough to verify the claim only for those whose minimal polynomials are
\begin{equation*}
x^2+1, \quad x^2-1, \quad x^3-1, \quad x^3-2x^2+2x-1.
\end{equation*}
The proofs of the first two cases are very similar, so we present only one of them.
\begin{enumerate}
\item Suppose $\mathfrak{m}_M(x)=\mathfrak{m}_N(x)=x^2+1$, and let $M$ and $N$ be given as  
\begin{equation*}
\left(\begin{array}{cc}
A_1 & A_2\\
-\frac{A_2^{-1}}{4} & -A_2^{-1}A_1A_2
\end{array}\right)\quad  \text{and}\quad  \left(\begin{array}{cc}
B_1 & B_2\\
-\frac{B_2^{-1}}{4} & -B_2^{-1}B_1B_2
\end{array}\right)
\end{equation*}
with $\mathfrak{m}_{A_1}=\mathfrak{m}_{B_1}(x)=x^2+\frac{3}{4}$ (by Theorem \ref{24092702}). By Lemma \ref{24090301}, $B_1$ is either $A_1$ or $-A_1$. If $B_1=-A_1$, by replacing $N$ with $-N$, we continue to assume $B_1$ in $N$ is $A_1$. To simplify the problem, we normalize $M$ and $N$ by
\begin{equation*}
\left(\begin{array}{cc}
A_1 & I\\
-\frac{I}{4} & -A_1
\end{array}\right)\quad  \text{and}\quad  \left(\begin{array}{cc}
A_1 & C_2\\
-\frac{C_2^{-1}}{4} & -C_2^{-1}A_1C_2
\end{array}\right)
\end{equation*}
respectively where $C_2:=B_2A_2^{-1}$ (and so $\det C_2=1$). We also assume $A_2\neq \pm B_2$, which implies $C_2\neq \pm I $. Then the primary matrices associated to $M$ and $N$ are 
\begin{equation}\label{24111501}
\left(\begin{array}{cc}
\frac{\sqrt{-3}}{2} & 1\\
-\frac{1}{4} & -\frac{\sqrt{-3}}{2}
\end{array}\right)\quad \text{and}\quad \left(\begin{array}{cc}
\frac{\sqrt{-3}}{2} & \lambda\\
-\frac{\lambda^{-1}}{4} & -\frac{\sqrt{-3}}{2}
\end{array}\right) 
\end{equation}
respectively where $\lambda:=\frac{\mathrm{tr}\,C_2\pm \sqrt{\mathrm{Disc}\,C_2}}{2}$, so the product of the two matrices in \eqref{24111501} is 
\begin{equation*}
P:=\left(\begin{array}{cc}
-\frac{3}{4}-\frac{\lambda^{-1}}{4} & \frac{\sqrt{-3}}{2}(-1+\lambda)\\
\frac{\sqrt{-3}}{8}(-1+\lambda^{-1}) & -\frac{3}{4}-\frac{\lambda}{4}
\end{array}\right).
\end{equation*}
Note that the eigenvalues $\lambda_i$ ($i=1,2$) of $P\iota$ are the solutions of $x^2-\frac{1}{4}(\lambda-\lambda^{-1})x-1=0$. Since $\lambda\in \mathbb{Q}(\sqrt{-3})\setminus \mathbb{Q}$, one can check two nonzero eigenvalues $\lambda_i$ ($i=1,2$) are distinct each other and $\mathrm{tr}\,P\neq 0$. Consequently, $\lambda_i$ ($i=1,2$) are roots of unity by Proposition \ref{24082301}, and thus $\deg \lambda_i\leq 2$ by Lemma \ref{24072601}. However, this induces $\lambda=\pm 1$, contradicting the assumption $C_2\notparallel I$. 

Similarly, when $\mathfrak{m}_M(x)=\mathfrak{m}_N(x)=x^2-1$, one shows $N$ is either $M$ or $\iota M\iota$ by the same procedure.

\item Now we assume $\mathfrak{m}_M(x)=\mathfrak{m}_N(x)=x^3-1$, and $M$ and $N$ are given as
\begin{equation*}
\left(\begin{array}{cc}
A_1 & A_2\\
-3A_2^{-1}A_1^2 & A_2^{-1}A_1A_2
\end{array}\right)\quad  \text{and}\quad  
\left(\begin{array}{cc}
B_1 & B_2\\
-3B_2^{-1}B_1^2 & B_2^{-1}B_1B_2
\end{array}\right)
\end{equation*}
with $\mathfrak{m}_{A_1}=\mathfrak{m}_{B_1}(x)=x^2-\frac{1}{2}x+\frac{1}{4}$ (by Theorem \ref{24092702}).  By Lemma \ref{24090301}, $B_1$ is either $A_1$ or $\frac{A_1^{-1}}{4}$. If $B_1=\frac{A_1^{-1}}{4}$, then by replacing $N$ with $N^{-1}$, we may assume without loss of generality that $B_1$ in $N$ is $A_1$. Note that $N$ is a conjugate of $M$ via
\begin{equation*}
N=(I \oplus B_2^{-1}A_2)M
(I \oplus B_2^{-1}A_2)^{-1}. 
\end{equation*}
Since 
\begin{equation*}
M\iota M^2=
\left(\begin{array}{cc}
\frac{I}{2} & 4A_1^2A_2\\
\frac{3}{2}A_2^{-1}A_1 & -\frac{I}{2}
\end{array}\right), 
\end{equation*}
one checks $\mathfrak{m}_{M\iota M^2}(x)=x^2-1$ (from Theorem \ref{24092702}) and, as $\mathfrak{m}_{N\iota N^2}=x^2-1$, concludes $N\iota N^2$ is either $M\iota M^2$ or $\iota M\iota M^2\iota $ by the previous claim. Consequently, $B_2^{-1}A_2=I$ or $-I$, implying $N=M$ or $\iota M\iota$ respectively.  

Similarly, when $\mathfrak{m}_M(x)=\mathfrak{m}_N(x)=x^3-2x^2+2x-1$, we have $\mathfrak{m}_{M^2}(x)=\mathfrak{m}_{N^2}(x)=x^3-1$ and attain the result by the same argument. 
\end{enumerate}
\end{proof}

Using the proposition, we establish Theorem \ref{24092901}. The proof of the theorem is two-fold: first, we prove the equality $\mathcal{G}=\langle M, \pm \iota \rangle $ and then confirm \eqref{25020205}. In the first claim, it will be shown that, for any minimal polynomial appearing in Theorem \ref{24092702}(1), $M$ and $\iota$ generate an element with that minimal polynomial. The proposition above then ensures that $\langle M, \pm \iota \rangle$ is the largest possible, and thus coincides with $\mathcal{G}$. The proof of the second claim will be computational in nature.

\begin{proof}[Proof of Theorem \ref{24092901}]
\begin{enumerate}
\item First, we show $\mathcal{G}=\langle M, \pm \iota \rangle $; to do so, it is enough to verify, for any $N\in\mathcal{G}$ of Type I, II or III, we have $N\in \langle M, \pm \iota \rangle $. 

Note that if $M$, with $\mathfrak{m}_M(x)=x^3\pm 2x^2+2x\pm 1$, is as in \eqref{24102103}, then $M^2$, $M\iota M^2$ and $M\iota M^5$ are given respectively by
\begin{equation*}
\begin{gathered}
\left(\begin{array}{cc}
\frac{2}{3}A_1^2 & 2A_1A_2\\
-\frac{2}{3}A_2^{-1}A_1^3 & \frac{2}{3}A_2^{-1}A_1^2A_2
\end{array}\right),\;
\left(\begin{array}{cc}
\frac{4}{3}A_1^3 & \frac{4}{3}A_1^3A_2\\
-\frac{4}{9}A_2^{-1}A_1^4 & \frac{4}{3}A_2^{-1}A_1^3A_2
\end{array}\right)\;\text{and}\;\left(\begin{array}{cc}
\frac{I}{2} & -\frac{3}{2}A_1^{-1}A_2\\
-A_2^{-1}\frac{A_1}{2} & -\frac{I}{2}
\end{array}\right)
\end{gathered}
\end{equation*}
where $\mathfrak{m}_{M^2}=x^3-1$, $\mathfrak{m}_{M\iota M^2}=x^2+1$ and $\mathfrak{m}_{M\iota M^5}=x^2-1$ by Theorem \ref{24092702}. Combining this with Proposition \ref{24080401}, we conclude $\langle M, \pm \iota \rangle $ contains every Type III element in $\mathcal{G}$, which also establishes the second claim of the theorem.

If $N\in\mathcal{G}$ is of Type I or II, thanks to Lemma \ref{25020503}, $NM$ is of Type III and so contained in $\langle M, \iota \rangle $, as noted above. Hence, $N\in \langle M, \iota \rangle $ follows. 

\item Now we verify the equality in \eqref{25020205}. First, note that $M^5$ and $\iota M^5\iota M^2$ are
$$\left(\begin{array}{cc}
\frac{3}{4}A_1^{-1} & \frac{16}{9}A_1^4A_2\\
\frac{A_2^{-1}}{4} & \frac{3}{4}A_2^{-1}A_1^{-1}A_2
\end{array}\right)\,\,\text{and}\,\,\left(\begin{array}{cc}
0 & 2A_2\\
-\frac{2}{3}A_2^{-1}A_1^2 & 0
\end{array}\right)$$
respectively. We denote $\iota M^5\iota M^2$ by $\eta$, and claim \eqref{25020205}. Since
\begin{equation*}
M\eta=\left(\begin{array}{cc}
-\frac{2}{3}A_1^2 & 2A_1A_2\\
-\frac{2}{3}A_2^{-1}A_1^3 & -\frac{2}{3}A_2^{-1}A_1^2A_2
\end{array}\right)=\eta M 
\end{equation*}
and
\begin{equation*}
M\iota M=(\frac{4}{3}A_1^2) \oplus (-\frac{4}{3}A_2^{-1}A_1^2A_2)=-\iota\eta^2, 
\end{equation*} 
combining these with the relations $\eta \iota=-\iota \eta$ and $\eta^6=-I$, it follows that any arbitrary element in $\mathcal{G}$, which is of the form
\begin{equation*}
\pm \iota^{\epsilon_1}M^{\alpha_1}\iota \cdots \iota M^{\alpha_n}\iota^{\epsilon_2}\;\; \text{where }0\leq \epsilon_1, \epsilon_2\leq 1\text{ and }0\leq \alpha_i\leq 5, 
\end{equation*}
is reduced to an element of the form
\begin{equation}\label{25010401}
\pm \iota^{\epsilon_1}\eta^{\alpha_1}M^{\alpha_2}\iota^{\epsilon_2}\;\;\text{where }0\leq \epsilon_1, \epsilon_2\leq 1, \;\text{and}\; 0\leq \alpha_1, \alpha_2\leq 5.
\end{equation}
Moreover, as $\eta M=-\iota M^2\iota$ and  
\begin{equation*}
\eta^3=(-\frac{8}{3}A_1^2A_2)\widetilde{\oplus}
(\frac{8}{9}A_2^{-1}A_1^4)=-M^3
\end{equation*}
\eqref{25010401} is further reduced to an element of the form given in \eqref{25020205}. This completes the proof of Theorem \ref{24092901}. 
\end{enumerate}
\end{proof}

\subsubsection{$\sqrt{-2}$}

Now we turn to the case where $\tau_1,\tau_2\in\mathbb{Q}(\sqrt{-2})$. The general scheme and strategy of the proof are very similar to those of Theorem \ref{24092901}. The following proposition is seen as a counterpart to Proposition \ref{24080401}, but its statement differs slightly.  

\begin{proposition}\label{250204031}
Let $\mathcal{X}$ be an element in $\mathfrak{Ghol}$ whose cusp shapes are contained in $\mathbb{Q}(\sqrt{-2})$. Let $M$ and $N$ be two elements of Type III in $\mathrm{Aut}\,\mathcal{X}$ with the same minimal polynomial. 
\begin{enumerate}
\item If $\mathfrak{m}_M(x)=x^4+1$, then $N\in \langle M, \pm\iota \rangle $. 
\item Suppose $\mathfrak{m}_M(x)=x^2-1$ and that $M$ and $N$ are given as
\begin{equation}\label{25031501}
\left(\begin{array}{cc}
\frac{I}{2} & A_2\\
\frac{3}{4}A_2^{-1} & -\frac{I}{2}
\end{array}\right)\quad \text{and}\quad 
\left(\begin{array}{cc}
\frac{I}{2} & B_2\\
\frac{3}{4}B_2^{-1} & -\frac{I}{2}
\end{array}\right)
\end{equation}
respectively for some $A_2, B_2\in \mathrm{GL}_2(\mathbb{Q})$. Then either $A_2=\pm B_2$ or $\mathfrak{m}_{B_2A_2^{-1}}(x)=x^2\pm \frac{2}{3}x+1$. In particular, in the latter case, $B_2A_2^{-1}$ is uniquely determined up to a scalar multiple of $\pm 1$. 
\end{enumerate}
\end{proposition}

\begin{proof}[Proof of Proposition \ref{250204031}]
The proof of the first claim is analogous to that of Proposition \ref{24080401}, so we only prove the second one. 

By Theorem \ref{24092702}, $M$ and $N$ are given as 
\begin{equation*}
\left(\begin{array}{cc}
\frac{I}{2} & A_2\\
\frac{3}{4}A_2^{-1} & -\frac{I}{2}
\end{array}\right)\quad  \text{and}\quad  \left(\begin{array}{cc}
\frac{I}{2} & B_2\\
\frac{3}{4}B_2^{-1} & -\frac{I}{2}
\end{array}\right)
\end{equation*}
respectively for some $A_2, B_2\in \mathrm{GL}_2(\mathbb{Q})$ satisfying $\det A_2=\det B_2=\frac{3}{4}$. To simplify the problem, if we normalize $M$ and $N$ as
\begin{equation*}
\left(\begin{array}{cc}
\frac{I}{2} & I\\
\frac{3I}{4} & -\frac{I}{2}
\end{array}\right)\quad  \text{and}\quad  \left(\begin{array}{cc}
\frac{I}{2} & C_2\\
\frac{3}{4}C_2^{-1} & -\frac{I}{2}
\end{array}\right)
\end{equation*}
respectively where $C_2:=B_2A_2^{-1}$, then the primary matrices $P_M$ and $P_N$ associated to $M$ and $N$ are 
\begin{equation}\label{25020404}
\left(\begin{array}{cc}
\frac{1}{2} & 1\\
\frac{3}{4} & -\frac{1}{2}
\end{array}\right)\quad \text{and}\quad \left(\begin{array}{cc}
\frac{1}{2} & \lambda\\
\frac{3}{4\lambda} & -\frac{1}{2}
\end{array}\right) 
\end{equation}
respectively where $\lambda:=\frac{\mathrm{tr}\,C_2\pm  \sqrt{\mathrm{Disc}\,C_2}}{2}(\in \mathbb{Q}(\sqrt{-2})\setminus \mathbb{Q})$. So their product $P:=P_MP_N$ is 
\begin{equation*}
\left(\begin{array}{cc}
\frac{1}{4}+\frac{3}{4\lambda} & \frac{\lambda}{2}-\frac{1}{2}\\
\frac{3}{8}-\frac{3}{8\lambda}  & \frac{3}{4}\lambda+\frac{1}{4}
\end{array}\right). 
\end{equation*}
If $\mathrm{tr}\,P\iota=0$, then $\lambda=\pm 1$ (i.e., $A_2=\pm B_2$). Otherwise, by Proposition \ref{24082301}, two nonzero eigenvalues of $P$, which are the solutions of 
\begin{equation}\label{25031505}
x^2-\big(\frac{1}{2}+\frac{3}{4\lambda}+\frac{3}{4}\lambda\big)x+1, 
\end{equation} 
are roots of unity. Since their degrees are at most $2$ by Lemma \ref{24072601}, an elementary computation shows that $\lambda=\frac{\pm 1\pm 2\sqrt{-2}}{3}$, that is, $\mathrm{tr}\,C_2=\pm \frac{2}{3}$. By Lemma \ref{24090301}, $C_2$ is uniquely determined up to its inverse and scalar multiplication by $\pm 1$. However, if both 
\begin{equation}\label{25112101}
\left(\begin{array}{cc}
\frac{I}{2} & C_2\\
\frac{3}{4}C_2^{-1} & -\frac{I}{2}
\end{array}\right)\quad  \text{and}\quad  \left(\begin{array}{cc}
\frac{I}{2} & \frac{1}{C_2}\\
\frac{3}{4}C_2 & -\frac{I}{2}
\end{array}\right)
\end{equation}
belong to $\mathrm{Aut}\,\mathcal{X}$ simultaneously, then we get a contradiction by applying Proposition \ref{24082301} to the product of the primary matrices associated to \eqref{25112101}. We omit the details because analogous procedures have been repeatedly shown. In conclusion, $C_2$ is unique up to a scalar multiple of $\pm 1$.  
\end{proof}

\noindent \textbf{Remark. }If the cusp shapes of $\mathcal{X}(\in \mathfrak{Ghol})$ are contained in a quadratic field other than $\mathbb{Q}(\sqrt{-2})$, then for $M, N\in \mathrm{Aut}\,\mathcal{X}$ given as in \eqref{25031501}, we always have $A_2=\pm B_2$ (so $N\in \langle M,\pm\iota \rangle $).\footnote{The reasoning follows directly from the proof of Proposition \ref{250204031}(2). More precisely, it is immediate from the fact that the solutions of the equation in \eqref{25031505} are roots of unity of degree at most $2$ if and only if $\lambda\in \mathbb{Q}(\sqrt{-2})$.} In this sense, the conclusion in Proposition \ref{250204031}(2) can be seen as the only exceptional case. \\

Now we prove Theorem \ref{25020403}. 
\begin{proof}[Proof of Theorem \ref{25020403}]
The proof of Theorem \ref{25020403} is similar to that of Theorem \ref{24092901}. Hence, we highlight the key computational steps instead of repeating the full details.

\begin{enumerate}
\item To claim $\mathcal{G}=\langle M, \pm \iota \rangle $, it is enough to show that any Type III element of $\mathcal{G}$ is contained in $\langle M, \pm \iota \rangle $. Note that $M^2, M^3$ and $(\iota M^2)^2$ are  
\begin{equation*}
\left(\begin{array}{cc}
A_1-\frac{I}{2} & (2A_1-I)A_2\\
A_2^{-1}(\frac{A_1}{2}-\frac{I}{4}) & A_2^{-1}(-A_1+\frac{I}{2})A_2
\end{array}\right),\; 
\left(\begin{array}{cc}
A_1-I & -A_2\\
-\frac{A_2^{-1}}{4} & A_2^{-1}A_1A_2
\end{array}\right)\;\text{and}\;
\left(\begin{array}{cc}
0 & -2A_2\\
\frac{A_2^{-1}}{2} & 0
\end{array}\right)
\end{equation*}
respectively. Thus, for $\eta:=(\iota M^2)^2$, $\eta M, \eta M^2$ and $\eta M^3$ are 
\begin{equation*}
\left(\begin{array}{cc}
-\frac{I}{2} & (-2A_1+2I)A_2\\
\frac{A_2^{-1}A_1}{2} & \frac{I}{2}
\end{array}\right),\; 
\left(\begin{array}{cc}
-A_1+\frac{I}{2} & (2A_1-I)A_2\\
A_2^{-1}(\frac{A_1}{2}-\frac{I}{4}) & A_2^{-1}(A_1-\frac{I}{2})A_2
\end{array}\right)\;\text{and}\;
\left(\begin{array}{cc}
\frac{I}{2} & -2A_1A_2\\
A_2^{-1}(\frac{A_1}{2}-\frac{I}{2}) & -\frac{I}{2}\\
\end{array}\right)
\end{equation*}
with the minimal polynomials $x^2-1, x^2+1$ and $x^2-1$ respectively, by Theorem \ref{24092702}. By Proposition \ref{250204031}, it follows that $\langle M, \pm \iota \rangle $ contains every Type III element of $\mathcal{G}$.  

\item Next we prove the equality in \eqref{25020501}. For convenience, let us temporarily denote the set on the right side of \eqref{25020501} by $\mathcal{\tilde{G}}$. Since every element of $\mathcal{G}$ is of the form 
\begin{equation*}
\pm \iota^{\epsilon_1}M^{\alpha_1}\iota\cdots \iota M^{\alpha_n}\iota^{\epsilon_2}
\end{equation*}
where $0\leq \epsilon_1, \epsilon_2\leq 1$ and $0\leq \alpha_i\leq 3$, to prove the equality, it is enough to show that $(\iota M^{\alpha_1})(\iota M^{\alpha_2})\in \mathcal{\tilde{G}}$ for every $1\leq \alpha_1, \alpha_2\leq 3$. For instance, since
\begin{equation*}
\iota M\iota M=\left(\begin{array}{cc}
A_1-I & A_2\\
-\frac{A_2^{-1}}{4} & -A_2^{-1}A_1A_2
\end{array}\right)=M^3\iota
\end{equation*}
and
\begin{equation*}
\iota M\iota M^{2}=\left(\begin{array}{cc}
-\frac{I}{2} & (2A_1-2I)A_2\\
-\frac{A_2^{-1}A_1}{2} & \frac{I}{2}
\end{array}\right)=\iota M^2\iota, 
\end{equation*}
both $\iota M\iota M$ and $\iota M\iota M^{2}$ are contained in $\mathcal{\tilde{G}}$. Similarly, one confirms the claim for the remaining cases. 
\end{enumerate}
\end{proof}

\begin{example}\label{25111910}
\normalfont{By Theorem \ref{25020403}, if $\mathcal{M}$ and $\mathcal{X}$ are given by $v2788$ and its analytic holonomy set, then $\mathcal{G}=\langle M, \iota \rangle$ where $M$ is
$$\left(\begin{array}{cccc}
\frac{1}{2} & \frac{1}{2} & -\frac{1}{2}  & 0\\
-1 & \frac{1}{2} & 0 & -\frac{1}{2}\\
\frac{1}{2} & 0 & \frac{1}{2} & -\frac{1}{2}\\
0 & \frac{1}{2} & 1 & \frac{1}{2}
\end{array}\right),$$ given in \eqref{25110405}. In fact, one can verify the following matrix
$$\left(\begin{array}{cccc}
0 & \frac{1}{2} & 0  & \frac{1}{2}\\
-1 & 0 & -1 & 0\\
0 & \frac{1}{2} & 0 & -\frac{1}{2}\\
-1 & 0 & 1 & 0
\end{array}\right),$$
given in \eqref{25112505}, is equal to $(M\iota)^2$ and so lies in $\mathcal{G}$. Using the structure of $\mathcal{G}$ described in Theorem \ref{25020403}, we can now completely classify every $H$ fitting into the setting of Example \ref{25101509}.}
\end{example}

\subsubsection{$\sqrt{-1}$: I}\label{24010501}

In this sub-subsection, we prove Theorem \ref{24111305}(1). The strategy of the proof closely parallels that of Theorems \ref{24092901}-\ref{25020403}. As in those cases, we first need the following proposition, whose proof is almost identical to that of Proposition \ref{24080401}, and will thus be omitted here.

\begin{proposition}\label{24010202}
Let $\mathcal{X}$ be an element in $\mathfrak{Ghol}$ whose cusp shapes are contained in $\mathbb{Q}(\sqrt{-1})$. Let $M$ be an element in $\mathrm{Aut}\,\mathcal{X}$ of Type III, satisfying $\mathfrak{m}_{M}(x)=x^2+1$. If $N$ is another element of Type III in $\mathrm{Aut}\,\mathcal{X}$ such that $\mathfrak{m}_{N}(x)$ is either $x^2+1$ or $x^2-1$, then $N\in \langle M, \pm \iota \rangle $. 
\end{proposition}

Using the proposition, we establish Theorem \ref{24111305}(1). 
\begin{proof}[Proof of Theorem \ref{24111305}(1)]
Similar to the previous cases, the proof consists of two parts: first, we show $\mathcal{G} = \langle  M, \pm \iota  \rangle $, and then establish \eqref{24122603}.

\begin{enumerate}
\item As seen in the proof of Theorem \ref{24092901}, thanks to Lemma \ref{25020503}, it is enough to verify $N\in \langle M, \pm\iota \rangle $ for those $N\in\mathcal{G}$ of Type III, in particular for those that satisfying $\mathfrak{m}_N(x)$ is one of the following:
\begin{equation}\label{25021403}
x^2+1, \quad x^2-1,  \quad x^2+x+1, \quad x^3+x^2+x+1. 
\end{equation}

For the first two cases, the result follows directly from Proposition \ref{24010202}. For $\mathfrak{m}_N(x) = x^2 + x + 1$, if we normalize $N$ and $\iota 
(\iota M)^2$ as  
\begin{equation*}
\left(\begin{array}{cc}
A_1 & I\\
-\frac{I}{2} & \frac{A_1^{-1}}{2}
\end{array}\right)\quad \text{and}\quad  
\left(\begin{array}{cc}
\frac{I}{2} & B_2\\
\frac{3}{4}B_2^{-1} & -\frac{I}{2}
\end{array}\right) 
\end{equation*}
respectively, then the primary matrices associated to $N$ and $\iota (\iota M)^2$ are
\begin{equation}\label{25021301}
\left(\begin{array}{cc}
\frac{-1+\sqrt{-1}}{2} & 1\\
-\frac{1}{2} & \frac{-1-\sqrt{-1}}{2}
\end{array}\right)\quad \text{and}\quad  
\left(\begin{array}{cc}
\frac{1}{2} & \lambda\\
\frac{3}{4\lambda} & -\frac{1}{2}
\end{array}\right)
\end{equation}
respectively where $\lambda:=\frac{\mathrm{tr}\,B_2\pm \sqrt{\mathrm{Disc}\,B_2}}{2}$. If we let $P$ be the product of two matrices in \eqref{25021301}, which is 
\begin{equation*}
\left(\begin{array}{cc}
\frac{-1+\sqrt{-1}}{4}+\frac{3}{4\lambda} & \frac{-1+\sqrt{-1}}{2}\lambda-\frac{1}{2}\\
-\frac{1}{4}+\frac{3(-1-\sqrt{-1})}{8\lambda} & \frac{-\lambda}{2}+\frac{1+\sqrt{-1}}{4}
\end{array}\right),    
\end{equation*}
then  
\begin{equation}\label{24122401}
\mathfrak{\chi}_P(x)=x^2-\big(\frac{\sqrt{-1}}{2}+\frac{3}{4\lambda} -\frac{\lambda}{2}\big)x-1.
\end{equation}
By appealing to Proposition \ref{24082301}, we get that the solutions of $\mathfrak{\chi}_P(x)=0$ must be roots of unity. Applying Lemma \ref{24072601} and elementary observation, however, one checks that there is no $\lambda\in \mathbb{Q}(\sqrt{-1})$ satisfying the required condition, i.e., $N\notin \mathcal{G}$. 

The last case $\mathfrak{m}_N(x)=x^3+x^2+x+1$ can be handled similarly. 

\item Next, we establish the equality in \eqref{24122603}. First, note that $M^2=-I$ and thus every element of $\mathcal{G}$ is of the following form: 
\begin{equation}\label{24122101}
\pm \iota^{\epsilon_1}M\iota \cdots M\iota^{\epsilon_2} 
\end{equation}
where $\epsilon_i$ is either $0$ or $1$ ($i=1,2$). Since $\mathfrak{m}_{\iota M}(x)=x^4-x^2+1$ (implying $(\iota M)^6=-I$), one can further reduce \eqref{24122101} to 
\begin{equation*}
\pm \iota^{\epsilon_1}(\iota M)^{\alpha}\iota^{\epsilon_2} 
\end{equation*}
where $0\leq \alpha\leq 5$ and each $\epsilon_i$ ($i=1,2$) is either $0$ or $1$. Consequently, the equality in \eqref{24122603} follows. 
\end{enumerate}
\end{proof}

\subsubsection{$\sqrt{-1}$: II}\label{24010502}
In this sub-subsection, we prove Theorem \ref{24111305}(2). As noted earlier, the statements of Theorem \ref{24111305} are somewhat different from those of Theorem \ref{24092901} or \ref{25020403}. In the previous theorems, there exists a primitive element such that every other element is contained in the group generated by this primitive element as well as $\pm\iota$. However, for $M$ and $N$ in the second claim of Theorem \ref{24111305}, neither, when combined with $\pm \iota$, generates the other. Nonetheless, we show below that both $M$ and $N$ are compatible in the sense that they can belong to the same group.

\begin{proposition}\label{24070301}
Let $\mathcal{X}$ be an element in $\mathfrak{Ghol}$ whose cusp shapes are contained in $\mathbb{Q}(\sqrt{-1})$. Let 
\begin{equation}\label{24092701}
M:=\left(\begin{array}{cc}
A_1 & A_2\\
-A_2^{-1}A_1^2 & A_2^{-1}A_1A_2
\end{array}\right)\quad\text{and}\quad N:=\left(\begin{array}{cc}
A_1 & A_2\\
-\frac{A_2^{-1}}{2} & \frac{A_2^{-1}A_1^{-1}A_2}{2}
\end{array}\right)
\end{equation}
be two elements in $\mathrm{Aut}\,\mathcal{X}$, both of Type III, such that $\mathfrak{m}_{A_1}(x)=x^2+x+\frac{1}{2}$. 
\begin{enumerate}
\item Then $\langle M, N, \pm \iota \rangle $ contains an abelian subgroup of order $16$ generated by $NM^{-1}$ and $(\iota M)^3$.

\item Suppose $L$ is an element in $\mathrm{Aut}\,\mathcal{X}$ of Type III such that $\mathfrak{m}_{L}(x)=x^3\pm x^2+ x\pm1$ or $x^4+1$ (resp. $x^2\pm x+1$ or $x^4-x^2+1$) with a form matching one of those given in \eqref{24121601} (resp. \eqref{24121602}). Then
\begin{equation}\label{25121401}
\langle L, N, \pm \iota \rangle =\langle M, N, \pm \iota \rangle \quad \text{(resp. }\langle M, L, \pm \iota \rangle =\langle M, N, \pm \iota \rangle). 
\end{equation}
\end{enumerate}
\end{proposition}

The first statement of the proposition is rather straightforward, and can be easily attained by directly computing $NM^{-1}$ and $(\iota M)^3$. Therefore, the key part of the proof lies in the second statement. To establish this, we proceed with a step-by-step approach, including a preliminary lemma. 

Suppose $M$ and $L$ are any two elements in $\mathrm{Aut}\,\mathcal{X}$ with the same minimal polynomial, for instance, $x^3+x^2+x+1$, and are thus represented as   
\begin{equation}\label{24112001}
\left(\begin{array}{cc}
A_1 & A_2\\
-A_2^{-1}A_1^2 & A_2^{-1}A_1A_2
\end{array}\right)\quad \text{and}\quad 
\left(\begin{array}{cc}
B_1 & B_2\\
-B_2^{-1}B_1^2 & B_2^{-1}B_1B_2
\end{array}\right)
\end{equation}
respectively, where $\mathfrak{m}_{A_1}(x)=\mathfrak{m}_{B_1}(x)=x^2+x+\frac{1}{2}$. Thanks to Lemma \ref{24090301}, $B_1$ is either $A_1$ or $ \frac{A_1^{-1}}{2}$. In the latter case, by replacing $L$ with $L^{-1}$, we still assume $B_1$ in \eqref{24112001} is $A_1$. 

In the following lemma, we establish a strong relationship between $A_2$ and $B_2$ in \eqref{24112001}. 
\begin{lemma}\label{24070311}
Let $M_N$ and $L$ be two elements of Type III in $\mathrm{Aut}\,\mathcal{X}$ such that $M_N$ is either 
\begin{equation}\label{24111801}
\left(\begin{array}{cc}
A_1 & A_2\\
-A_2^{-1}A_1^2 & A_2^{-1}A_1A_2
\end{array}\right)\quad \text{or}\quad 
\left(\begin{array}{cc}
A_1 & A_2\\
-\frac{A_2^{-1}}{2} & \frac{A_2^{-1}A_1^{-1}A_2}{2}
\end{array}\right), 
\end{equation}
and $L$ is
\begin{equation*}
\left(\begin{array}{cc}
A_1 & B_2\\
-B_2^{-1}A_1^2 & B_2^{-1}A_1B_2
\end{array}\right)
\end{equation*}
where $\mathfrak{m}_{A_1}(x)=x^2+x+\frac{1}{2}$. Then, for each possible form of $M_N$, we have either $A_2= \pm B_2$ or $\mathrm{tr}\,(A_2^{-1}B_2)=0$. 
\end{lemma}
\begin{proof}
We consider only the first form of $M_N$ among the two given in \eqref{24111801}, as the proof for the second case is analogous. 

To simplify the problem, we normalize $M_N$ and $L$ by 
\begin{equation}\label{24092706}
\left(\begin{array}{cc}
A_1 & I\\
-A_1^2 & A_1
\end{array}\right)\quad \text{and}\quad  \left(\begin{array}{cc}
A_1 & C_2\\
-C_2^{-1}A_1^2 & C_2^{-1}A_1C_2
\end{array}\right) 
\end{equation}
respectively where $C_2:=B_2A_2^{-1}$. Since $A_2\neq \pm B_2$ by the assumption, equivalently, $C_2\neq \pm I$. Without loss of generality, if we further assume the $(1,2)$-entries of $A_1$ is positive, then the primary matrices associated with $M_N$ and $L$ are 
\begin{equation}\label{24092801}
\left(\begin{array}{cc}
\frac{-1+\sqrt{-1}}{2} & 1\\
\frac{\sqrt{-1}}{2} & \frac{-1+\sqrt{-1}}{2}
\end{array}\right)\quad \text{and}\quad 
\left(\begin{array}{cc}
\frac{-1+\sqrt{-1}}{2} & \lambda\\
\lambda^{-1}\frac{\sqrt{-1}}{2} & \frac{-1+\sqrt{-1}}{2}
\end{array}\right)
\end{equation} 
respectively where $\lambda:=\frac{\mathrm{tr}\,C_2\pm \sqrt{\mathrm{Disc}\,C_2}}{2}(\neq \pm 1)$. Then the product of two matrices in \eqref{24092801} is   
\begin{equation*}
P:=\left(\begin{array}{cc}
\frac{\sqrt{-1}}{2}(-1+\lambda^{-1}) & \frac{-1+\sqrt{-1}}{2}(1+\lambda)\\
\frac{-1-\sqrt{-1}}{4}(1+\lambda^{-1}) & \frac{\sqrt{-1}}{2}(-1+\lambda)
\end{array}\right) 
\end{equation*}
and the eigenvalues $\lambda_i$ ($i=1,2$) of $P\iota$ are the solutions of $x^2+\frac{\sqrt{-1}}{2}(\lambda-\lambda^{-1})x+1=0$. Since there is no $\lambda\in \mathbb{Q}(\sqrt{-1})\setminus \mathbb{Q}$ such that the eigenvalues of $P\iota$ are distinct and $\mathrm{tr}\,P=0$, all $\lambda_i$ are roots of unity by Proposition \ref{24082301}, and satisfy $\deg \lambda_i\leq 2$ by Lemma \ref{24072601}. As a result, $\lambda=\pm\sqrt{-1}$, inducing $\mathrm{tr}\,C_2=0$. 
\end{proof}

Using the lemma, we proceed to prove Proposition \ref{24070301}. 
\begin{proof}[Proof of Proposition \ref{24070301}]
\begin{enumerate}
\item Note that 
\begin{equation}\label{25102001}
\begin{gathered}
NM^{-1}=I \oplus -2A_2^{-1}A_1^2A_2, \;\;
(\iota M)^3=2A_1^2 \oplus 2A_2^{-1}A_1^2A_2, 
\end{gathered}
\end{equation}
and thus the group generated by $NM^{-1}$ and $(\iota M)^3$ is an abelian group of order $16$.

\item We assume $\mathfrak{m}_L(x)=x^3+x^2+x+1$ and consider only this case, as the remaining cases are treated in the same manner. Without loss of generality, let $M$ and $L$ be normalized as in \eqref{24092706}. 
\begin{enumerate}
\item If $C_2$ is either $I$ or $-I$, then $L$ is either $M$ or $\iota M\iota$ respectively. In this case, the claim, \eqref{25121401}, is straightforward. 

\item Otherwise, if $C_2\neq \pm I$, since $\det C_2=1$, $\mathfrak{m}_{C_2}(x)=x^2+1$ by Lemma \ref{24070311}, and $LM^{-1}$ is 
\begin{equation*}
\left(\begin{array}{cc}
A_1 & C_2\\
-C_2^{-1}A_1^2 & C_2^{-1}A_1C_2
\end{array}\right)
\left(\begin{array}{cc}
\frac{A_1^{-1}}{2} & -\frac{A_1^{-2}}{2}\\
\frac{I}{2} & \frac{A_1^{-1}}{2}
\end{array}\right)
=\left(\begin{array}{cc}
\frac{I}{2}+\frac{C_2}{2} & -\frac{A_1^{-1}}{2}+\frac{C_2A_1^{-1}}{2}\\
-\frac{C_2^{-1}A_1}{2}+\frac{C_2^{-1}A_1C_2}{2} & \frac{C_2^{-1}}{2}+\frac{C_2^{-1}A_1C_2A_1^{-1}}{2}
\end{array}\right). 
\end{equation*}
Note that $LM^{-1}$ is a matrix of Type III. As $\mathfrak{m}_{\frac{I}{2}+\frac{C_2}{2}}(x)=x^2-x+\frac{1}{2}$, $\frac{I}{2}+\frac{C_2}{2}$ is either $-A_1$ or $-\frac{A_1^{-1}}{2}$, or, equivalently, $C_1$ is either $2A_1^2$ or $-2A_1^2$. If $C_1=2A_1^2$, then 
\begin{equation}\label{24111401}
\begin{aligned}
L=(I \oplus -2A_1^2)M(I\oplus 2A_1^2)&\Longrightarrow L\in \langle M, N, \pm \iota \rangle,\\
N^{-1}L=I \oplus 2A_1^2=\iota NM^{-1}&\Longrightarrow M\in \langle L, N, \pm \iota \rangle.
\end{aligned}
\end{equation}
Consequently, the equality in \eqref{25121401} follows.

Likewise, the same conclusion holds for $C_2=-2A_1^2$. 
\end{enumerate}
\end{enumerate}
\end{proof}

Now we complete the proof of Theorem \ref{24111305}. 

\begin{proof}[Proof of Theorem \ref{24111305}(2)]
Without loss of generality, by Proposition \ref{24070301}, we assume both $M$ and $N$ are as in \eqref{24092701}. 
\begin{enumerate}
\item We claim $\mathcal{G}=\langle M, N, \pm \iota\rangle $. As in the proofs of previous theorems, it suffices to verify that if a Type III element $L$, whose minimal polynomial is one of those given in \eqref{25021403}, lies in $\mathcal{G}$, then $L\in \langle M, N,  \pm \iota\rangle $. If $L\in \mathcal{G}$ with $\mathfrak{m}_{L}(x)=x^3+ x^2+x+1$ or $x^2+ x+1$, then $L\in \langle M, N, \pm \iota \rangle $ follows directly from Proposition \ref{24070301}. If $L\in \mathcal{G}$ with $\mathfrak{m}_L(x)=x^2-1$, then, as shown in the proof of Theorem \ref{24111305}(1), both $N$ and $L$ do not belong to the same subgroup of $\mathrm{Aut}\,\mathcal{X}$, which is a contradiction. Otherwise, if $L\in \mathcal{G}$ with $\mathfrak{m}_{L}(x)=x^2+1$, then $\mathcal{G}$ contains $L\iota L$ whose minimal polynomial is $x^2-1$. So it falls into the previous case, which leads to a contradiction.  

\item By Proposition \ref{24070301}(1), $\mathcal{G}$ contains an abelian subgroup of order $16$, so we only verify the equality in \eqref{24111302}. If we denote the righthand side in \eqref{24111302} by $\mathcal{\tilde{G}}$, then clearly $N\in \mathcal{\tilde{G}}$ from $NM^{-1}\in \mathcal{H}$. Thus to complete the proof, it is enough to show that $\mathcal{\tilde{G}}$ is closed under the group action, which reduces to proving  
\begin{equation*}
M\eta M, \;\; M^2\eta M \;\;\text{and}\;\;M\eta M^2\in \mathcal{\tilde{G}}\;\;\text{for any }\eta \in \mathcal{H}.
\end{equation*}
Further, instead of considering all $\eta\in \mathcal{H}$, we confirm the above claim only for $\eta\in\{I, \iota, \eta_1, \iota\eta_1\}$ where $\eta_1=I\oplus 2A_2^{-1}A_1^2A_2$, as the remaining elements of $\mathcal{H}$ are obtained by multiplying one of these by either $\pm I$ or $\pm\eta_2$ where $\eta_2=2A_1^2 \oplus 2A_2^{-1}A_1^2A_2$. 
\begin{enumerate}
\item First, 
$$M^3=\left(\begin{array}{cc}
\frac{A_1^{-1}}{2} & 2A_1^2A_2\\
\frac{A_2^{-1}}{2}  & -\frac{A_2^{-1}A_1^{-1}A_2}{2}
\end{array}\right)=-\iota \eta_2 M\iota\Longrightarrow M^3\in \mathcal{\tilde{G}}.$$ 

\item Second,  
$$M\eta_1M=-\eta_1 M\eta_1,\;\;
M\iota\eta_1 M=\iota\eta_1\eta_2M\eta_1 \;\;\text{and}\;\;M\iota M=\eta_2\iota,$$
implying $M\eta M\in \mathcal{\tilde{G}}$ for any $\eta\in \mathcal{H}$. 

\item Third, since 
$$M\eta_1M^2=-\iota M\eta_1,\;\;M\iota\eta_1M^2=-\iota M\iota\eta_1\;\; \text{and}\;\; M\iota M^2=\eta_2M$$ from 
\begin{equation}\label{25102006}
M^2=2A_1A_2\widetilde{\oplus}\frac{1}{2}A_2^{-1}A_1^{-1}, 
\end{equation}
we get $M\eta M^2\in \mathcal{\tilde{G}}$ for any $\eta\in \mathcal{H}$. 

\item Lastly,  
$$M^2\eta_1M
=\eta_1M\iota,\;\;M^2\iota\eta_1M=-\eta_1\iota M\iota\;\;\text{and}\;\; M^2\iota M=M\eta_2\iota,$$ 
implying $M^2\eta M\in \mathcal{\tilde{G}}$ for any $\eta\in \mathcal{H}$. 
\end{enumerate}
\end{enumerate}
\end{proof}

Finally, we prove Corollary \ref{25030302}: 
\begin{proof}[Proof of Corollary \ref{25030302}]
If $M$ is an element of Type III in $\mathrm{Aut}\,\mathcal{X}$, then either $\mathfrak{m}_M(x)$ or $\mathfrak{m}_{\iota M}(x)$ is $x^2-1$ by Theorem \ref{24092702}. If $N$ is another element of $\mathrm{Aut}\,\mathcal{X}$ such that $\mathfrak{m}_N(x)=x^2-1$, then $N\in \langle M, \pm\iota \rangle $ (see the remark after the proof of Proposition \ref{250204031}), hence concluding $\mathcal{G}=\langle M, \pm\iota \rangle $. Now the claim is established, once we show the set on the right side of \eqref{25030301} is closed under the group action. The verification of this fact is tedious and so we skip it here. 
\end{proof}

\newpage
\section{Proof of Theorem \ref{25101401}}\label{25012003}
In this section, by amalgamating all the results given in the previous section, we prove Theorem \ref{25101401}. 

We divide the problem into two cases depending on whether the two cusp shapes $\tau_1,\tau_2$ are contained in $\mathbb{Q}(\sqrt{-1})$ or not. The first case is relatively easier and follows immediately from Theorems \ref{24092901}-\ref{25020403} and Corollary \ref{25030302}. On the other hand, the latter case is more delicate and requires a detailed analysis of Theorem \ref{24111305}. 

Before proceeding with the proof, we recall the notation $A(p/q)$ where $A\in \mathrm{GL}_2(\mathbb{Q})$, introduced in Section \ref{intro} (after Theorem \ref{22050203}), and note the following equality: 
\begin{equation}\label{25031001}
A(p/q)=(aA)(p/q) \quad \text{for any  }a\in \mathbb{Q}\backslash\{0\}.
\end{equation}
 
\subsection{$\tau_1, \tau_2\notin\mathbb{Q}(\sqrt{-1})$}

The following is a detailed version of Theorem \ref{25101401} when $\tau_1, \tau_2\notin\mathbb{Q}(\sqrt{-1})$.

\begin{theorem}\label{25101601}
Let $\mathcal{M}$ be a two-cusped hyperbolic $3$-manifold whose cusp shapes belong to the same quadratic field but not to $\mathbb{Q}(\sqrt{-1})$, and let $\mathcal{M}_{p'_1/q'_1, p'_2/q'_2}$ and $\mathcal{M}_{p_1/q_1, p_2/q_2}$ be two Dehn fillings of $\mathcal{M}$ having the same pseudo complex volume with sufficiently large $|p'_i|+|q'_i|$ and $|p_i|+|q_i|$ ($i=1,2$). Suppose the complex lengths of the core geodesics of $\mathcal{M}_{p_1/q_1,p_2/q_2}$ are linearly independent over $\mathbb{Q}$. Then the following statements hold.
\begin{enumerate} 
\item If $\tau_1, \tau_2\in\mathbb{Q}(\sqrt{-3})$, then there exist unique elements $\sigma$ of order $3$ and $\rho$ in $\mathrm{GL}_2(\mathbb{Q})$ such that $(p'_1/q'_1, p'_2/q'_2)$ is either 
\begin{equation}\label{25031003}
\big(\sigma^i(p_1/q_1), \rho\circ\sigma^j\circ\rho^{-1}(p_2/q_2)\big)\quad \text{or}\quad  \big(\sigma^i\circ \rho^{-1}(p_2/q_2),\rho\circ\sigma^j (p_1/q_1)\big) 
\end{equation}
for some $0\leq i, j\leq 2$.

\item If $\tau_1,\tau_2\in\mathbb{Q}(\sqrt{-2})$, then there exist unique elements $\sigma$ of order $2$ and $\rho$ in $\mathrm{GL}_2(\mathbb{Q})$ such that $(p'_1/q'_1, p'_2/q'_2)$ is either  
\begin{equation*}
\big(\sigma^i(p_1/q_1), \rho\circ\sigma^i\circ\rho^{-1}(p_2/q_2)\big)\quad \text{or}\quad \big(\rho^{-1}(p_2/q_2), \rho(p_1/q_1)\big)
\end{equation*}
for some $0\leq i\leq 1$.
\item Otherwise, there exists unique $\rho\in \mathrm{GL}_2(\mathbb{Q})$ such that $(p'_1/q'_1, p'_2/q'_2)$ is either
\begin{equation}\label{25031008}
(p_1/q_1, p_2/q_2)\quad \text{or}\quad (\rho^{-1}(p_2/q_2), \rho(p_1/q_1)).
\end{equation} 
\end{enumerate}
\end{theorem}

To prove the theorem, we apply Theorem \ref{23070517}, in particular the relations \eqref{22031407} therein, to Theorems \ref{24092901}-\ref{25020403} and Corollary \ref{25030302}. To be more precise, we first list all the group elements appearing in these results, examine their actions on the Dehn filling coefficients as described in \eqref{22031407}, and then quotient by the scaling given in \eqref{25031001}.

We only consider the first case of the theorem in detail, as the remaining cases can be handled analogously.

\begin{proof}[Proof of Theorem \ref{25101601}]
As usual, let $\mathcal{G}$ denote the largest subgroup of $\mathrm{Aut}\,\mathcal{X}$, generated by elements of Types I-III. 
\begin{enumerate}
\item Suppose the two cusp shapes are contained in $\mathbb{Q}(\sqrt{-3})$. 
\begin{enumerate}
\item If $\mathcal{G}$ contains only elements of Types I-II, by Theorem \ref{24072401}, we assume the group is generated by \eqref{24093002} with its elements listed in \eqref{24093015}. According to Theorem \ref{22031403}, it means either 
\begin{equation*}
p'_1/q'_1=(A_1^T)^i(p_1/q_1),\quad p'_2/q'_2=A_2^T(A_1^T)^j(A_2^T)^{-1}(p_2/q_2) 
\end{equation*}
or 
\begin{equation*}
p'_1/q'_1=(A_1^T)^i(A_2^T)^{-1}(p_2/q_2),\quad p'_2/q'_2=A_2^T(A_1^T)^j(p_1/q_1)
\end{equation*}
for some $0\leq i, j\leq 5$. Since $A_1^3=-I$, letting $\sigma:=(A_1^T)^2$ and $\rho:=(A_2^T)^2$, we get the desired conclusion described in \eqref{25031003} by \eqref{25031001}.

\item If $\mathcal{G}$ contains an element of Type III, we assume without loss of generality that $\mathcal{G}=\langle M, \iota \rangle$ as described in Theorem \ref{24092901}. Moreover, by \eqref{25020205} and the observation in \eqref{25031001}, it is enough to consider the following elements in $\mathcal{G}$: 
$$\eta^i=
\begin{cases}
\left(\begin{array}{cc}
0 & 2\big(-\frac{4}{3}\big)^{\frac{i-1}{2}}A_1^{i-1}A_2\\
\frac{1}{2}\big(-\frac{4}{3}\big)^{\frac{i+1}{2}}A_2^{-1}A_1^{i+1} & 0
\end{array}\right)\text{ for }i=1,3,5 \\
\left(\begin{array}{cc}
\big(-\frac{4}{3}\big)^{\frac{i}{2}} A_1^i& 0\\
0 & \big(-\frac{4}{3}\big)^{\frac{i}{2}}A_2^{-1}A_1^iA_2
\end{array}\right) \text{ for }i=0, 2,4
\end{cases}$$
and 
$$ \eta^iM= 
\begin{cases}
\left(\begin{array}{cc}
-\frac{2}{3}\big(-\frac{4}{3}\big)^{\frac{i-1}{2}}A_1^{i+1} & 2\big(-\frac{4}{3}\big)^{\frac{i-1}{2}}A_1^{i}A_2\\
\frac{1}{2}\big(-\frac{4}{3}\big)^{\frac{i+1}{2}}A_2^{-1}A_1^{i+2} & \frac{1}{2}\big(-\frac{4}{3}\big)^{\frac{i+1}{2}}A_2^{-1}A_1^{i+1}A_2 
\end{array}\right)\text{ for }i=1,3,5 \\
\left(\begin{array}{cc}
\big(-\frac{4}{3}\big)^{\frac{i}{2}}A_1^{i+1} & \big(-\frac{4}{3}\big)^{\frac{i}{2}}A_1^iA_2\\
-\frac{1}{3}\big(-\frac{4}{3}\big)^{\frac{i}{2}}A_2^{-1}A_1^{i+2} & \big(-\frac{4}{3}\big)^{\frac{i}{2}}A_2^{-1}A_1^{i+1}A_2
\end{array}\right) \text{ for }i=0,2,4.
\end{cases}$$ 
Consequently, applying Theorem \ref{22031403} (and \eqref{25031001}), we find $(p'_1/q'_1, p'_2/q'_2)$ equals either
\begin{equation*}
\big((A_1^{T})^{2i}(p_1/q_1), A_2^T(A_1^{T})^{2i}(A_2^T)^{-1}(p_2/q_2)\big)
\end{equation*}
for some $0\leq i\leq 2$, or
\begin{equation*}
\big((A_1^T)^{2i}(A_2^T)^{-1}(p_2/q_2),A_2^T(A_1^T)^{2i+2}(p_1/q_1)\big)
\end{equation*}
for some $0\leq i\leq 2$. By letting $\sigma:=(A_1^T)^2$ and $\rho:=A_2^T$, the conclusion follows.
\end{enumerate}

\item For two cups shapes contained in $\mathbb{Q}(\sqrt{-2})$, the proof is very similar to the preceding case by invoking Theorem \ref{25020403}, so it will be omitted here. 

\item Otherwise, when two cusp shapes are not contained in any of the quadratic fields listed above, the conclusion in \eqref{25031008} follows from Theorem \ref{24072401}(2) as well as Corollary \ref{25030302}.   
\end{enumerate}
\end{proof}

\subsection{$\tau_1, \tau_2\in\mathbb{Q}(\sqrt{-1})$}\label{25111901}
In this section, we consider the remaining case of Theorem \ref{25101401}, that is, when $\tau_1, \tau_2 \in \mathbb{Q}(\sqrt{-1})$. As in the proof of Theorem \ref{25101601}, the argument relies on Theorem \ref{24111305}, but it requires a more subtle analysis to get the optimal bound. 

We begin by restating the remaining part of Theorem \ref{25101401}, not covered by Theorem \ref{25101601}, as follows.
\begin{theorem}\label{25101610}
Let $\mathcal{M}$ be a two-cusped hyperbolic $3$-manifold whose cusp shapes lie in $\mathbb{Q}(\sqrt{-1})$, and let $\mathcal{M}_{p_1/q_1, p_2/q_2}$ be a Dehn filling of $\mathcal{M}$ with $|p_i|+|q_i|$ ($i=1,2$) sufficiently large. If the complex lengths of the core geodesics of $\mathcal{M}_{p_1/q_1,p_2/q_2}$ are linearly independent over $\mathbb{Q}$, then the number of Dehn fillings of $\mathcal{M}$ whose pseudo complex volumes are the same as that of $\mathcal{M}_{p_1/q_1,p_2/q_2}$ is at most $8$.
\end{theorem}

Recall the definition of $\mathcal{S}$, initially given in the statement of Theorem \ref{23070303}(2) and refined in Theorem \ref{20103003}. To simplify notation for the proof of Theorem \ref{25101610} as well as that of Theorem \ref{23070505} later in the next section, we further stratify $\mathcal{S}$ in Theorem \ref{20103003} as follows.
 
\begin{definition}\label{25121101}
\normalfont{Let $\mathcal{M}$ be a two-cusped hyperbolic $3$-manifold whose cusp shapes belong to the same quadratic field. Let $\mathcal{M}_{p'_1/q'_1, p'_2/q'_2}$ and $\mathcal{M}_{p_1/q_1, p_2/q_2}$ be two Dehn fillings of $\mathcal{M}$ with the same pseudo complex volume and with sufficiently large coefficients. 
\begin{enumerate}
\item Suppose the complex lengths of the core geodesics of $\mathcal{M}_{p_1/q_1, p_2/q_2}$ are \textit{linearly dependent} over $\mathbb{Q}$. We define $\mathcal{S}_1$ to be the {\it smallest subset} of $\mathcal{S}$ such that $(p'_1/q'_1, p'_2/q'_2)$ is either
\begin{equation*}
(\sigma(p_1/q_1), \rho(p_2/q_2))\quad \text{or}\quad (\sigma(p_2/q_2), \rho(p_1/q_1))
\end{equation*} 
for some $(\sigma, \rho)\in \mathcal{S}_1$. 

\item Similarly, $\mathcal{S}_2$ is defined under the assumption that the complex lengths of the core geodesics of $\mathcal{M}_{p_1/q_1, p_2/q_2}$ are \textit{linearly independent} over $\mathbb{Q}$.
\end{enumerate}}
\end{definition}

Obviously, $\mathcal{S}=\mathcal{S}_1\cup \mathcal{S}_2$, and the identity is contained in both $\mathcal{S}_1$ and $\mathcal{S}_2$. In the next section, we  clarify the structure of the intersection, $\mathcal{S}_1\cap \mathcal{S}_2$, in order to prove Theorem \ref{23070505}. 

Having Definition \ref{25121101}, the conclusion of Theorem \ref{25101610} is simply restated as  
\begin{equation*}
|\mathcal{S}_2|\leq 8.
\end{equation*} 
Before embarking on the proof, we briefly explain how the following bound 
\begin{equation}\label{25112401}
|\mathcal{S}_2|\leq 12
\end{equation}
is directly derived from Theorem \ref{24111305}. We consider only the second case where $\mathcal{G}$ is given as in \eqref{24111302}, since the first one is easy to deal with.

Recall the definition of $\mathcal{H}$ in Theorem \ref{24111305}, whose generators are specified in \eqref{25102001}. Setting $\sigma:=A_1^T$ and $\rho:=A_2^T$, the elements in $\mathcal{H}$ induce the following types of elements in $\mathcal{S}_2$:
\begin{equation}\label{25102003}
(\sigma^{2i}, \rho\circ\sigma^{2j}\circ \rho^{-1}),\quad i,j\in \mathbb{Z}.
\end{equation} 

For $M$ as in Theorem \ref{24111305}(2) (explicitly described in \eqref{24092701}), the elements of the form $\eta_1M\eta_2$ where $\eta_1, \eta_2\in \mathcal{H}$ induce the following types of elements in $\mathcal{S}_2$: 
\begin{equation}\label{25102004}
(\sigma^{2i-1}, \rho\circ\sigma^{2j-1}\circ \rho^{-1}),\quad i,j\in \mathbb{Z}.
\end{equation}

Finally, recalling $M^2$ given in \eqref{25102006}, the elements of the form $\eta_1 M^2 \eta_2$ with $\eta_1, \eta_2 \in \mathcal{H}$ induce following types of elements in $\mathcal{S}_2$:
\begin{equation}\label{25102005}
(\sigma^{2i-1}\circ\rho^{-1},  \rho\circ\sigma^{2j-1} ), \quad i, j \in \mathbb{Z}.
\end{equation}

As the order of $\sigma$ is $4$ (due to $4A_1^4=-I$), \eqref{25112401} is derived from \eqref{25102003}-\eqref{25102005}.

The following is a consolidated and reformulated version of Theorems \ref{23070517}-\ref{22031403}.

\begin{theorem}\label{25120103}
Let $\mathcal{M}, \mathcal{X}, \mathcal{M}_{p'_{1}/q'_{1}, p'_{2}/q'_{2}}$ and $\mathcal{M}_{p_{1}/q_{1}, p_{2}/q_{2}}$ be the same as in Theorem \ref{23070517}. Let $\mathfrak{p}'$ and $\mathfrak{p}$ be the Dehn filling points on $\mathcal{X}$ corresponding to $\mathcal{M}_{p'_{1}/q'_{1}, p'_{2}/q'_{2}}$ and $\mathcal{M}_{p_{1}/q_{1}, p_{2}/q_{2}}$ respectively. Then there exists a finite set $\mathcal{G}(\subset \mathrm{Aut}\,\mathcal{X})$ such that $\mathfrak{p}'=M(\mathfrak{p})$ for some $M\in \mathcal{G}$. Moreover, when $M$ is of the form \eqref{21072901}, it satisfies \eqref{22031407}-\eqref{25120102}. 
\end{theorem}

When both cusp shapes of $\mathcal{M}$ are contained in $\mathbb{Q}(\sqrt{-1})$ and $\mathcal{G}$ is as given in \eqref{24111302}, a careful examination of the relations \eqref{22031407}-\eqref{25120102} shows that there exist $M, N \in \mathcal{G}$ such that, for the Dehn filling point $\mathfrak{p}$ introduced above, $M(\mathfrak{p})$ and $N(\mathfrak{p})$ cannot simultaneously correspond to Dehn-filled manifolds of the same pseudo complex volume. In other words, if $M(\mathfrak{p})$ corresponds to a \textit{manifold}, then $N(\mathfrak{p})$ must necessarily correspond either to an \textit{orbifold} or to a manifold whose pseudo complex volume is different from that of $M(\mathfrak{p})$, and vice versa. Consequently, when we restrict our attention strictly to counting the number of Dehn filled {\it manifolds}, the upper bound of $|\mathcal{S}_2|$ in this case drops from $12$ to $8$.

In the following three lemmas, $\mathcal{M}$ and $\mathcal{X}$ are the same as in Theorem \ref{25101610}, and $\mathcal{G}, \mathcal{H}$ and $M$ are the same as in Theorem \ref{24111305}(2).

The first lemma is readily obtained from Theorem \ref{22031403}.
\begin{lemma}\label{25061001}
Let $\mathcal{M}_{p'_1/q'_1, p'_2/q'_2}$ and $\mathcal{M}_{p_1/q_1, p_2/q_2}$ be two Dehn fillings of $\mathcal{M}$ that have the same pseudo complex volume, and let $\mathfrak{p}'$ and $\mathfrak{p}$ be their associated Dehn filling points on $\mathcal{X}$, respectively. Suppose $\mathfrak{p'}=(\eta_1M\eta_2)(\mathfrak{p})$ for some $\eta_1, \eta_2\in\mathcal{H}$ and $k_j$ ($1\leq j\leq 4$) are rational numbers, as obtained from Theorem \ref{25120103}. Then $k_j\neq \frac{1}{2}$ for each $1\leq j\leq 4$. Further, if $k_1\neq 1, 0$ and $\frac{1}{k_1}\in \mathbb{Z}$, then $\frac{1}{k_4}\notin \mathbb{Z}$. 
\end{lemma}
\begin{proof}
Note that if we assume $\eta_1M\eta_2$ is of the form given in \eqref{21072901}, then $\det A_j=\frac{1}{2}$ ($1\leq j\leq 4$). By \eqref{25120102}, we have
\begin{equation}\label{25061001}
\frac{1}{2k_1}+\frac{1}{2k_3}=1\Longrightarrow k_1+k_3=2k_1k_3, 
\end{equation}
which implies $k_1, k_3\neq \frac{1}{2}$. Analogously, $k_2, k_4\neq \frac{1}{2}$.

If both $\frac{1}{k_1}$ and $\frac{1}{k_4}$ are in $\mathbb{Z}$, since $k_3+k_4=1$ and
\begin{equation*}
k_4=1-k_3=1-\frac{k_1}{2k_1-1}=\frac{k_1-1}{2k_1-1},   
\end{equation*}
it follows that $\frac{1}{k_1-1}\in \mathbb{Z}.$ However, one can check there is no $k_1\in \mathbb{Q}\backslash\{0, 1,\frac{1}{2}\}$ satisfying the conditions.  
\end{proof}

In the following, we analyze the relations \eqref{22031407}-\eqref{25120102} in more detail, and show that, for a given Dehn filling point, not all elements in $\mathcal{G}$ realize Dehn filled manifolds of the same pseudo complex volume simultaneously. 

\begin{lemma}\label{25061205}
Let $\mathfrak{p}$ be a Dehn filling point on $\mathcal{X}$ corresponding to $\mathcal{M}_{p_1/q_1, p_2/q_2}$. Suppose one of the following conditions holds: 
\begin{enumerate}
\item either \begin{equation}\label{25061206}
(\pm A_1^{-1}A_2)\widetilde{\oplus}(\pm \frac{1}{2}A_2^{-1}A_1^{-1}) (\mathfrak{p})\quad \text{or}\quad 
(\pm 2A_1A_2)\widetilde{\oplus}(\pm A_2^{-1}A_1)(\mathfrak{p}) 
\end{equation}
is a point on $\mathcal{X}$ corresponding to a Dehn filling of $\mathcal{M}$ whose pseudo complex volume is equal to that of $\mathcal{M}_{p_1/q_1, p_2/q_2}$;

\item both 
\begin{equation}\label{25061208}
(\pm 2A_1A_2)\widetilde{\oplus}(\pm \frac{1}{2}A_2^{-1}A_1^{-1})(\mathfrak{p})\quad \text{and}\quad 
(\pm A_1^{-1}A_2)\widetilde{\oplus}
(\pm A_2^{-1}A_1)(\mathfrak{p})
\end{equation}
are points on $\mathcal{X}$ corresponding to Dehn fillings of $\mathcal{M}$ whose pseudo complex volumes are equal to that of $\mathcal{M}_{p_1/q_1, p_2/q_2}$. 
\end{enumerate}
Then, for any $\eta_1, \eta_2\in \mathcal{H}$, none of the points $(\eta_1 M \eta_2) (\mathfrak{p})$ corresponds to a Dehn filling whose pseudo complex volume is equal to that of $\mathcal{M}_{p_1/q_1,p_2/q_2}$. 
\end{lemma}

\begin{proof}
First suppose 
$$(A_1^{-1}A_2)\widetilde{\oplus}(\frac{1}{2}A_2^{-1}A_1^{-1})(\mathfrak{p})$$
is a point corresponding to a Dehn filling, say $\mathcal{M}_{p'_1/q'_1,p'_2/q'_2}$, whose pseudo complex volume is the same as that of $\mathcal{M}_{p_1/q_1, p_2/q_2}$. Then 
\begin{equation}\label{25101103}
\left(\begin{array}{c}
p'_1 \\
q'_1
\end{array}\right)=(A_1)^T(A_2^{-1})^T\left(\begin{array}{c}
p_2 \\
q_2
\end{array}\right)
\end{equation}
by Theorem \ref{22031403}. 

Second, suppose $\mathcal{M}(\mathfrak{p})$, which is
\begin{equation}\label{25061201} 
\left(\begin{array}{cc}
A_1 & A_2\\
-A_2^{-1}A_1^2  & A_2^{-1}A_1A_2
\end{array}\right)(\mathfrak{p}), 
\end{equation}
is a point corresponding to a Dehn filling of $\mathcal{M}$ whose pseudo complex volume is the same as that of $\mathcal{M}_{p_1/q_1,p_2/q_2}$. Let $k_i$ be $(1\leq i\leq 4)$ numbers associated with $\mathfrak{p}$ and \eqref{25061201}, as given by \eqref{22031407} in Theorem \ref{25120103}. Then 
\begin{equation*}
k_1(A_1^{-1})^T\left(\begin{array}{c}
p_1 \\
q_1
\end{array}\right)=k_2(A_2^{-1})^T\left(\begin{array}{c}
p_2 \\
q_2
\end{array}\right)\Longrightarrow 
\left(\begin{array}{c}
p_1 \\
q_1
\end{array}\right)=\frac{k_2}{k_1}(A_1)^T(A_2^{-1})^T\left(\begin{array}{c}
p_2 \\
q_2
\end{array}\right)
\end{equation*}
by \eqref{22031407}. So, combining it with \eqref{25101103},  
$$
\left(\begin{array}{c}
p_1 \\
q_1
\end{array}\right)
=\frac{k_2}{k_1}\left(\begin{array}{c}
p'_1 \\
q'_1
\end{array}\right).$$
Consequently, as both $(p_1, q_1)$ and $(p'_1, q'_1)$ are coprime pairs, $\frac{k_2}{k_1}=\pm 1$. But it contradicts the fact that $k_1+k_2=1$ and $k_1, k_2\neq \frac{1}{2}$ (by Lemma \ref{25061001}).

The remaining cases are treated analogously.  
\end{proof}

The following lemma is of a similar flavor to the above.

\begin{lemma}\label{25061203}
Let $\mathfrak{p}$ be a Dehn filling point on $\mathcal{X}$ corresponding to a Dehn filling $\mathcal{M}_{p_1/q_1, p_2/q_2}$. Suppose one of the following conditions holds: 
\begin{enumerate}
\item $$(\pm 2A_1^2)\oplus(\pm 2A_2^{-1}A_1^2A_2)(\mathfrak{p})$$ corresponds to a Dehn filling of $\mathcal{M}$ whose pseudo complex volume is equal to that of $\mathcal{M}_{p_1/q_1, p_2/q_2}$; 
\item both 
$$(\pm 2A_1^2)\oplus (\pm I)(\mathfrak{p})\quad \text{and}\quad 
(\pm I) \oplus (\pm 2A_2^{-1}A_1^2A_2)(\mathfrak{p})$$
correspond to Dehn fillings of $\mathcal{M}$ whose pseudo complex volumes are equal to that of $\mathcal{M}_{p_1/q_1, p_2/q_2}$. 
\end{enumerate}
Then, for any $\eta_1, \eta_2\in \mathcal{H}$, none of the points $(\eta_1 M \eta_2) (\mathfrak{p})$ corresponds to a Dehn filling whose pseudo complex volume is equal to that of $\mathcal{M}_{p_1/q_1,p_2/q_2}$. 
\end{lemma}
\begin{proof}
Suppose, for instance, both $(2A_1^2)\oplus(2A_2^{-1}A_1^2A_2)(\mathfrak{p})$ and $M(\mathfrak{p})$ are points corresponding to some Dehn fillings of $\mathcal{M}$ (say $\mathcal{M}_{p'_1/q'_1, p'_2/q'_2}$ and $\mathcal{M}_{\overline{p}_1/\overline{q}_1, \overline{p}_2/\overline{q}_2}$ respectively) whose pseudo complex volumes are the same as that of $\mathcal{M}_{p_1/q_1, p_2/q_2}$. 

First, applying Theorem \ref{25120103} to $\mathcal{M}_{p'_1/q'_1, p'_2/q'_2}$ and $\mathcal{M}_{p_1/q_1, p_2/q_2}$,  
$$\left(\begin{array}{c}
p'_1\\
q'_1
\end{array}\right)=\frac{1}{2}(A_1^T)^{-2}\left(\begin{array}{c}
p_1\\
q_1
\end{array}\right)\in\mathbb{Z}^2.$$
Since $\frac{1}{2}(A_1^T)^{-2}=I-(A_1^T)^{-1}$, this further implies 
\begin{equation}\label{25112503}
(A_1^T)^{-1}\left(\begin{array}{c}
p_1\\
q_1
\end{array}\right)\in\mathbb{Z}^2.
\end{equation} 

Now, applying Theorem \ref{25120103} to $\mathcal{M}_{\overline{p}_1/\overline{q}_1, \overline{p}_2/\overline{q}_2}$ and $\mathcal{M}_{p_1/q_1, p_2/q_2}$, there exists $k_i\in \mathbb{Q}$ $(1\leq i\leq 4)$ such that 
$$\left(\begin{array}{c}
\overline{p}_1\\
\overline{q}_1
\end{array}\right)=
k_1(A_1^T)^{-1}\left(\begin{array}{c}
p_1\\
q_1
\end{array}\right)\in\mathbb{Z}^2.$$ 
Combining this with \eqref{25112503}, it follows that $\frac{1}{k_1}\in \mathbb{Z}.$ By a similar argument, one has $\frac{1}{k_4}\in \mathbb{Z}$. But this contradicts Lemma \ref{25061001}

The other cases can be treated similarly. 
\end{proof}

By assembling all the three lemmas, we now prove Theorem \ref{25101610}.  

\begin{proof}[Proof of Theorem \ref{25101610}]
If $\mathcal{G}$ contains only elements of Type I or II, the proof is analogous to the previous case in Theorem \ref{25101601} and can be verified that $|\mathcal{S}_2|\leq 8$ by Theorem \ref{24072401}. 

Suppose $\mathcal{G}$ contains an element of Type III and let $\mathfrak{p}$ be a Dehn filling point associate to $\mathcal{M}_{p_1/q_1, p_2/q_2}$. We consider only the two cases listed in Theorem \ref{24111305}, as the remaining cases are easier to handle.   
\begin{enumerate}
\item If $\mathcal{G}=\langle M, \pm \iota \rangle $ where $\mathfrak{m}_{M}(x)=x^2+1$ as described in Theorem \ref{24111305}(1), then, by computing $(\iota M)^i$ for $1\leq i\leq 5$, one checks $(p'_1/q'_1, p'_2/q'_2)$ is either $(p_1/q_1, p_2/q_2)$ or $\big(\sigma^{-1}(p_2/q_2), \sigma(p_1/q_1)\big)$ for some $\sigma$ of order $2$. 

\item Suppose $\mathcal{G}=\langle M, N, \pm \iota \rangle $ where $\mathfrak{m}_{M}(x)=x^3\pm x^2+x\pm 1$ and $\mathfrak{m}_{N}(x)=x^4-x^2+1$ as in Theorem \ref{24111305}(2). First, if $(\iota M)^3(\mathfrak{p})$ corresponds to a Dehn filling of $\mathcal{M}$ whose pseudo complex volume is the same as that of $\mathcal{M}_{p_1/q_1, p_2/q_2}$, by Lemma \ref{25061203}, none of the points $(\eta_1 M \eta_2)(\mathfrak{p})$ corresponds to manifolds of the same pseudo complex volume for any $\eta_1, \eta_2\in \mathcal{H}$. Hence, the Dehn fillings having the same pseudo complex volume as $\mathcal{M}_{p_1/q_1, p_2/q_2}$ arise only from the remaining cases, namely, \eqref{25102003} or \eqref{25102005}. Consequently, $|\mathcal{S}_2|\leq 8$. Similarly, if $\mathfrak{p}$ satisfies any one of the conditions given in Lemmas \ref{25061205} and \ref{25061203}, then we obtain $|\mathcal{S}_2|\leq 8$. Otherwise, if $\mathfrak{p}$ satisfies neither of them, then we already exclude at least $5$ elements from $\mathcal{S}_2$, and thus $|\mathcal{S}_2|\leq 7$.


\end{enumerate}
\end{proof}

\newpage
\section{Proof of Theorem \ref{23070505}}\label{25102801}

Finally, in this section, we prove our first main theorem, Theorem \ref{23070505}, by synthesizing Theorems \ref{25082501}, \ref{25101601}, and \ref{25101610} obtained so far. 

For reader's convenience, we restate Theorem \ref{23070505} as follows. 
\\
\\
\noindent\textbf{Theorem \ref{23070505}} \textit{Let $\mathcal{M}$ be a two-cusped hyperbolic $3$-manifold, whose cusp shapes belong to the same quadratic field. We further suppose one of the following is true: }
\begin{itemize}
\item \textit{the potential function of $\mathcal{M}$ has a non-trivial term of degree $4$; or}
\item \textit{two cusps of $\mathcal{M}$ are SGI to each other.}
\end{itemize}
\textit{If $\tau_1, \tau_2\in\mathbb{Q}(\sqrt{-3})$ (resp. $\mathbb{Q}(\sqrt{-2})$ or $\mathbb{Q}(\sqrt{-1})$), then $|\mathcal{S}|\leq 18$ (resp. $3$ or $8$). Otherwise, $|\mathcal{S}|\leq 2$.}\\

To prove the theorem, we need to accommodate both $\mathcal{S}_1$ and $\mathcal{S}_2$ simultaneously, and, in particular, analyze their intersection. 

First note that, by the statement and proof of Theorem \ref{250825011}, $\sigma_1$ and $\sigma_2$ are uniquely determined, both of order $2$ or $3$ (depending on whether cusp shapes lie in $\mathbb{Q}(\sqrt{-1})$ or $\mathbb{Q}(\sqrt{-3})$), in the following three cases:
\begin{itemize}
\item two cusps of $\mathcal{M}$ are SGI to each other;
\item two cusp shapes of $\mathcal{M}$ are contained in $\mathbb{Q}(\sqrt{-3})$;
\item two cusp shapes of $\mathcal{M}$ are contained in $\mathbb{Q}(\sqrt{-1})$ and the coefficient of the term of bi-degree $(2,2)$ of the potential function of $\mathcal{M}$ is $0$. 
\end{itemize}

Thus each element of $\mathcal{S}_1$ acts on $(p_1/q_1, p_2/q_2)$ as a Type I element in the above cases, and therefore the conclusion of Theorem \ref{23070505} is covered by Theorems \ref{25101601}-\ref{25101610} in these cases.

Moreover, if the two cusp shapes are contained in neither $\mathbb{Q}(\sqrt{-1})$ nor $\mathbb{Q}(\sqrt{-3})$, then $|\mathcal{S}_1| = 1$, meaning $\mathcal{S}_1$ consists solely of the identity. Consequently, $\mathcal{S} = \mathcal{S}_2$, and the claim in this case is also immediate from Theorem \ref{25101601}. 

As a result, to complete the proof of Theorem \ref{23070505}, it is enough to consider the following remaining one:
\begin{itemize}
\item the two cusp shapes of $\mathcal{M}$ lie in $\mathbb{Q}(\sqrt{-1})$, and the coefficient of the term of bi-degree $(2,2)$ in the potential function of $\mathcal{M}$ is nonzero.
\end{itemize}

Recall that $\mathcal{G}$ denotes the largest subgroup of $\mathrm{Aut}\,\mathcal{X}$ whose elements are of Types I, II, and III. The following lemma simplifies the structure of $\mathcal{G}$ under the above assumption. 
\begin{lemma}\label{25102610}
Let $\mathcal{M}$ be a two-cusped hyperbolic $3$-manifold whose cusp shapes are contained in $\mathbb{Q}(\sqrt{-1})$, $\mathcal{X}$ be its analytic holonomy variety, and $\mathcal{G}$ be as above. If the coefficient of the bi-degree $(2,2)$ term in the potential function of $\mathcal{M}$ is nonzero, then $\mathcal{G}$ does not contain all Type I elements.
\end{lemma}
\begin{proof}
On the contrary, using the same notation in Theorem \ref{24072401}, we assume there exist $A_1, A_4$ of order $4$ such that 
\begin{equation*}
A_1^i\oplus A_4^j\in \mathcal{G}
\end{equation*}  
for any $0\leq i,j\leq 3$. Let $c_{4,0}u_1^4+c_{2,2}u_1^2u_2^2+c_{0,4}u_2^4$ be the homogeneous degree $4$ terms in the potential function of $\mathcal{M}$. By Proposition \ref{24020501}, it then satisfies
\begin{equation}\label{25011102}
\left(\begin{array}{c}
\lambda_1^{3i}2c_{4,0}u_1^3+c_{2,2}\lambda_1^{i}\lambda_2^{2j}u_1u_2^2\\
\lambda_2^{3j}2c_{0,4}u_2^3+c_{2,2}\lambda_1^{2i}\lambda_2^{j}u_2u_1^2
\end{array}\right)
=\left(\begin{array}{c}
\overline{\lambda_1}^i(2c_{4,0}u_1^3+c_{2,2}u_1u_2^2)\\
\overline{\lambda_2}^j(2c_{0,4}u_2^3+c_{2,2}u_2u_1^2)
\end{array}\right)
\end{equation}
for any $0\leq i,j\leq 1$ where $\lambda_k=\pm \sqrt{-1}$ ($k=1,2$). It is straightforward to check the equality holds if and only if $c_{2,2}=0$. 
\end{proof}

In the above proof, if $c_{2,2}\neq 0$ and the equation in \eqref{25011102} holds, then $i$ and $j$ must be either both even or odd. Having this observation, we now complete the proof of Theorem \ref{23070505}. 

\begin{proof}[Proof of Theorem \ref{23070505}]
We consider only the remaining case, described before the lemma. 

\begin{enumerate}
\item First, suppose $\mathcal{G}$ only contains elements of Types I and II. Since $c_{2,2}\neq 0$, by Lemma \ref{25102610} and the discussion after it, $\mathcal{G}$ is generated by 
\begin{equation*}
A_1^i \oplus A_4^j\quad \text{and}\quad 
A_2\widetilde{\oplus} A_2^{-1}
\end{equation*} 
where $0\leq i,j\leq 3$, and $i,j$ are either both odd or both even. Hence $|\mathcal{S}_2|\leq 4$ and, combined with $|\mathcal{S}_1|\leq 2$ (by Theorem \ref{250825011}), it follows that $|\mathcal{S}|\leq 5$.

\item Now suppose $\mathcal{G}$ contains an element of Type III. First, in the case of Theorem \ref{24111305}(1), we have $|\mathcal{S}_2|\leq 2$ and thus $|\mathcal{S}|\leq 3$ (since $|\mathcal{S}_1|\leq 2$). Second, in the case of Theorem \ref{24111305}(2), recall that $\mathcal{G}$ contains $\mathcal{H}$, which in turn contains all possible elements of Types I. But this contradicts Lemma \ref{25102610}. 
\end{enumerate}
\end{proof}

\vspace{5 mm}
\noindent Department of Mathematics, POSTECH\\
77 Cheong-Am Ro, Pohang, South Korea\\
\\
\text{Email Address}: bogwang.jeon@postech.ac.kr

\begin{thebibliography}{99}
\bibitem{za} E. ~Bombieri, D. ~Masser, U. ~Zannier, \textit{Anomalous subvarieties-structure theorems and applications}, IMRN \textbf{19} (2007), 1-33.
\bibitem{calegari1} D. ~Calegari, \textit{A note on strong geometric isolation in 3-orbifolds}, Bull. of Aust. Math. Soc. \textbf{53} 2 (1996), 271-280.

\bibitem{calegari} D. ~Calegari, \textit{Napoleon in isolation}, Proc. of Amer. Math. Soc. \textbf{129} 10, 3109-3119. 

\bibitem{feit} W. ~Feit, \textit{Orders of finite linear groups}, Proc. of the First Jamaican Conference on Group Theory and its Applications (Kingston, 1996) (Kingston), Univ. West Indies, 1996, 9-11.

\bibitem{CHDF} B. ~Jeon, \textit{Classification of hyperbolic Dehn fillings I}, Proc. of Lond. Math. Soc. \textbf{130} 1 (2025), e70017.

\bibitem{jeon3} B. ~Jeon, \textit{On the number of Dehn fillings of a given volume}, Trans. of Amer. Math. Soc. \textbf{374} (2021), 3947-3969.  

\bibitem{JO} B. ~Jeon, S. ~Oh, \textit{Hyperbolic Dehn filling, volume, and transcendentality}, arxiv.org/abs/2308.11574.

\bibitem{CHDFIII} B. ~Jeon, S. ~Oh, \textit{Classification of hyperbolic Dehn fillings III: Examples}, In preparation.

\bibitem{min} H. ~Minkowski, \textit{Untersuchungen $\ddot{u}$ber quadratische Formen. Bestimmung der Anzahl verschiedener Formen, welche ein gegebenes Genus enth$\ddot{a}$lt}, Acta Mathematica \textbf{7} (1885), 201-258.
\bibitem{min2}H. ~Minkowski, \textit{Zur Theorie der positiven quadratische Formen}, J. reine angew. Math. \textbf{101} (1887), 196-202.
\bibitem{walter} W. ~Neumann, \emph{Hilbert's 3rd Problem and Invariants of 3-manifolds}, In:``The Epstein Birthday Schrift", Geom. Topol. Monogr. \textbf{1} (1998), 383-411.

\bibitem{rigidity} W. ~Neumann, A.~Reid, \emph{Rigidity of cusps in deformations of hyperbolic 3-orbifolds}, Math. Ann. (1993), 223-237.

\bibitem{nz} W. ~Neumann, D. ~Zagier, \textit{Volumes of hyperbolic three-manifolds}, Topology \textbf{24} 3 (1985), 307-332.
\bibitem{PP} C. ~Petronio, J. ~Porti, \textit{Negatively oriented ideal triangulations and a proof of Thurston's hyperbolic Dehn filling theorem}, Expo. Math. \textbf{18} 1 (2000), 1-35.

\bibitem{thu} W.~Thurston, \textit{The Geometry and Topology of 3-manifolds}, Princeton University Mimeographed Notes (1979).
\bibitem{thu1} W.~Thurston, \textit{Three-dimensional manifolds, Kleinian groups and hyperbolic geometry}, Bull. of Amer. Math. Soc. \textbf{6} 3 (1982), 357-379.

\bibitem{yoshida} T. ~Yoshida, The $\eta$-invariant of hyperbolic 3-manifolds, Invent. Math. \textbf{81} (1985), 473-514.
\end{thebibliography}
\end{document}